\documentclass{amsbook}
\usepackage{amssymb}
\usepackage{enumerate}

\usepackage{graphicx}
\usepackage{amsmath}
\usepackage{amsfonts}
\usepackage{amsthm}
\usepackage{color}

\newcommand{\mtext}{}

\newtheorem{theorem}{Theorem}[section]

\newtheorem{corollary}[theorem]{Corollary}

\newtheorem{definition}[theorem]{Definition}
\newtheorem{example}[theorem]{Example}

\newtheorem{lemma}[theorem]{Lemma}

\newtheorem{proposition}[theorem]{Proposition}

\newtheorem{remark}[theorem]{Remark}

\numberwithin{equation}{section}

\makeindex

\title[Non-compact manifolds with general metric]{Inverse scattering on non-compact manifolds with general metric\\
$\ $
\\
$\ $
\\
{\rm Dedicated for the memory of Yaroslav Kurylev}}
\date{\today}

\begin{document}
\baselineskip 13pt

\author{Hiroshi Isozaki}
\address{Hiroshi Isozaki, Graduate School of Pure and Applied Sciences \\
University of Tsukuba,
Tsukuba, 305-8571, Japan}
\author{Matti Lassas}
\address{Matti Lassas, Department of Mathematics and Statistics \\
University of Helsinki, Finland}

\maketitle

\newpage

\begin{center}
{\Large{\bf Introduction and summary}}
\end{center}

\bigskip
The problem we address in this paper is the spectral theory and inverse problem associated with Laplacians on some class of  non-compact Riemannian manifolds (or more general manifolds admitting conic singularities). By observing behaviors of solutions to the Helmholtz equation on the manifold, we inroduce an analogue of Heisenberg's scattering matrix in quantum mechanics. We then show that the knowledge of the scattering matrix determines the topolgy and the metric of the manifold. We begin with a brief overview of our results, leaving  precise statements in the text.
\label{IntroFigure}
\begin{figure}[hbtp]
\centering
\includegraphics[height=6cm]{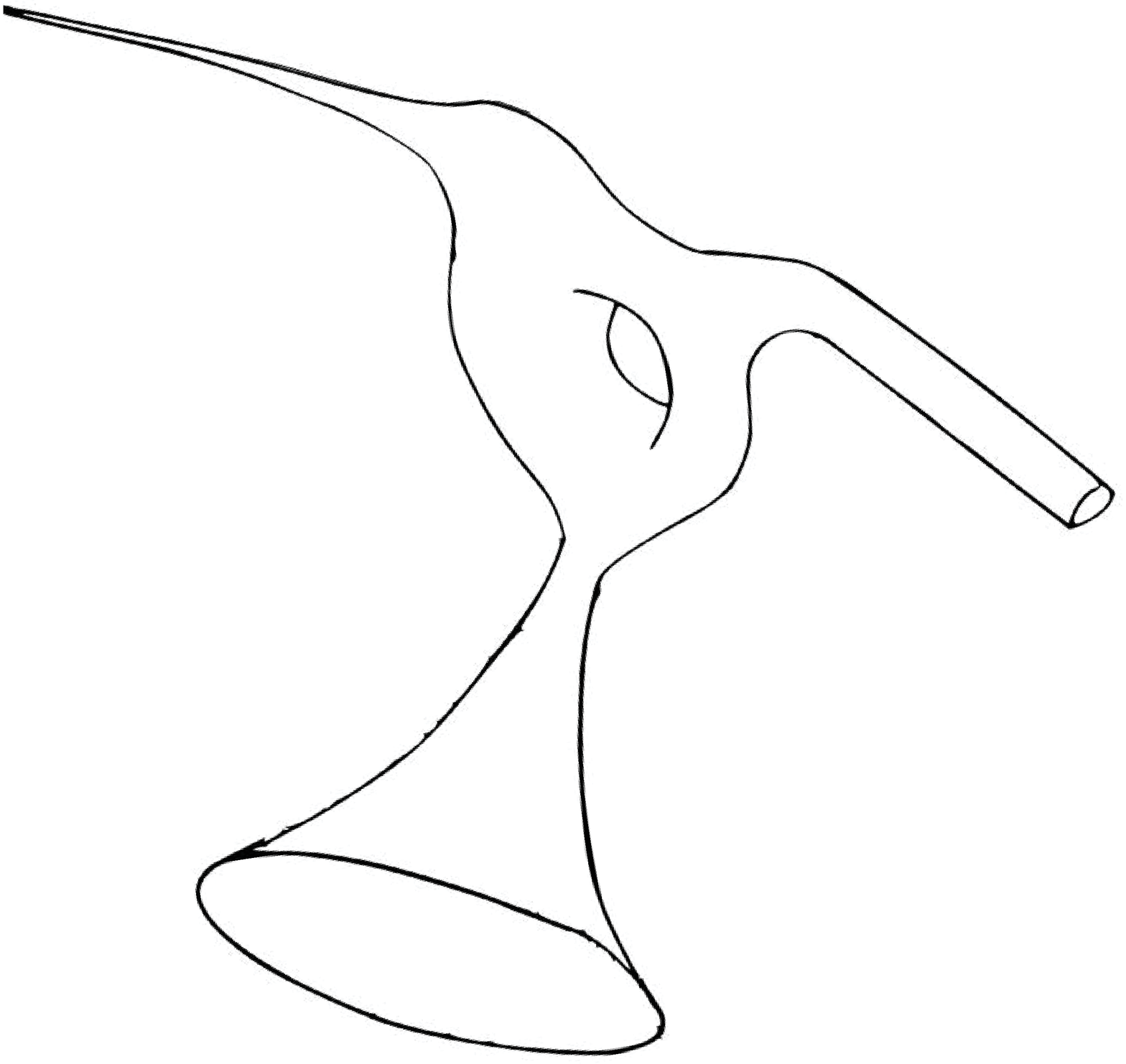}
\caption{Manifold $\mathcal M$}
\label{fig:math0}
\end{figure}

\subsection{Scattering of waves on non-compact manifolds}

We consider a connected, non-compact $n \, (\geq 2)$-dimensional Riemannian manifold $\mathcal M$ of the form
\begin{equation}
\mathcal M = \mathcal K  \cup \Big({\cup_{j=1}^{N}}\mathcal M_j\Big) \cup \Big({\cup_{j=N+1}^{N+N'}}\mathcal M_j\Big),
\label{S0manifoldmathcalM}
\end{equation}
where $\mathcal K$ is an open relatively compact subset and $\mathcal M_j$, $1 \leq j \leq N+N'$, is an open non-compact subset in 
 $\mathcal M$, henceforth called an {\it end}. 
As will be discussed later, we allow conical singularities for our manifolds.  For the sake of simplicity of explanation, however, we first  consider $C^{\infty}$-manifolds.

 We assume that each end $\mathcal M_j$ is diffeomorphic to 
$(0,\infty) \times M_j$, where $M_j$ is a compact $(n-1)$-dimensional Riemannian manifold endowed with metric $h_j(x,dx)$. Moreover, on each $\mathcal M_j$, the Riemannian metric of $\mathcal M$ is written in  the form 
\begin{equation}
ds^2 = (dr)^2 + \rho_j(r)h_j(r,x,dx),
\nonumber
\end{equation}
where $h_j(r,x,dx)$ is an $r$-dependent metric on $M_j$ satisfying\footnote{Until the end of subsection 0.6, we state only main parts of the assumptions. Precise assumptions are given in Subsection \ref{subsectionAssumptions}.}
\begin{equation}
h_j(r,x,dx) = h_j(x,dx) + O(r^{-\gamma_{j}}), \quad {\rm as} \quad r \to \infty,
\nonumber
\end{equation}
where $\gamma_{j} > 1$.  As for the behavior of $\rho_j(r)$, we assume that for $1 \leq j \leq N$, the section of $\mathcal M_j$ at $r$ is either exponentially or polynomially growing, 
\begin{equation}
\rho_j(r) = O(e^{c_{j}r}), \quad {\rm or} \quad \rho_j(r) = O(r^{2\beta_{j}}), \quad c_{j} > 0, \quad  \beta_{j}> 0,
\label{IntroGrowingEnd}
\end{equation}
and for $N+1 \leq j \leq N + N'$, the section of $\mathcal M_j$ at $r$ is either exponentially or polynomially decaying, 
\begin{equation}
\rho_j(r) = O(e^{-c_jr}), \quad {\rm or} \quad \rho_j(r) = O(r^{- 2\beta_j}), \quad c_j > 0, \quad  \beta_j > 0.
\label{IntroDecayingEnd}
\end{equation}
We represent the case (\ref{IntroGrowingEnd}) by 
\begin{equation}
\rho_j(r) \geq O(r^{2\beta_j})
\nonumber
\end{equation}
and, similary the case (\ref{IntroDecayingEnd}) by
\begin{equation}
\rho_j(r) \leq O(r^{-{2\beta_j}}).
\nonumber
\end{equation}
We put
\begin{equation}
E_j = \left(\frac{(n-1)c_j}{2}\right)^2, 
\nonumber
\end{equation}
when $\rho_j(r)$ is exponentially growing or decaying, and 
\begin{equation}
E_j = 0,
\nonumber
\end{equation}
when $\rho_j(r)$ grows or decays polynomially, but not exponentially. 
We put
$$
E = \mathop{{\min}_{1 \leq j \leq N+N'}}E_j.
$$
Let $H = - \Delta_{\mathcal M}$ be the Laplacian of $\mathcal M$. Then, the essential spectrum of $H$ is 
\begin{equation}
\sigma_e(H) = [E,\infty).
\nonumber
\end{equation}
Let
\begin{equation}
\mathcal E = \{E_1, \cdots, E_{N+N'}\}\cup \sigma_p(H),
\nonumber
\end{equation}
which is a discrete set with possible accumulation points in $\{E_1, \cdots, E_{N+N'}\}$.

We introduce a function space $\mathcal B^{\ast}(\mathcal M)$: $f \in \mathcal B^{\ast}(\mathcal M)$ if and only if $f \in L^2(\mathcal K)$ and 
\begin{equation}
\sup_{R>1}\frac{1}{R}\int_0^R\|f(r)\|^2_{L^2(M_j)}dr < \infty
\nonumber
\end{equation}
for  $1 \leq j \leq N + N'$. 
For two functions $f, g \in L^2_{loc}(\mathcal M)$, we denote $f \simeq g$ if they satisfy 
\begin{equation}
\lim_{R\to\infty}\frac{1}{R}\int_{0 < r < R}\|f(r) - g(r)\|^2_{L^2(M_j)}dr = 0
\nonumber
\end{equation}
on each $\mathcal M_j$, $1 \leq j \leq N + N'$.

 Our results are most transparent when $\rho_j(r) \geq r^{2\beta_j}$ with  $\beta_j > 1/2$ for $1 \leq j \leq N$. Take $\lambda \in \sigma_e(H)\setminus\mathcal E$ and let $\mathcal N(\lambda)$ be defined by
\begin{equation}
\mathcal N(\lambda) = \{u \in \mathcal B^{\ast}(\mathcal M)\, ; \, (- \Delta_{\mathcal M} - \lambda)u = 0\}. 
\nonumber
\end{equation}
Let 
\begin{equation}
{\bf h}_{\infty,j} = 
\left\{
\begin{split}
& L^2(M_j), \quad 1 \leq j \leq N, \\
& {\mathbb C}, \quad N+1 \leq j \leq N+N',
\end{split}
\right.
\nonumber
\end{equation}
and put
\begin{equation}
{\bf h}_{\infty}(\lambda) = \mathop{\oplus_{j=1}^{N+N'}}c_j(\lambda){\bf h}_{\infty,j},
\nonumber
\end{equation}
where $c_j(\lambda)$ is the characteristic function of $[E_{j},\infty)$.
 Note that for $a = (a_1,\cdots,a_{N+N'}) \in {\bf h}_{\lambda}$, 
 $a_j$ is an $L^2$-function on $M_j$ for $1 \leq j \leq N$, while $a_j \in \mathbb C$ for $N+1 \leq j \leq N + N'$, moreover $a_j = 0$ if $\lambda < E_j$.
  We put
\begin{equation}
\Phi_j(r,\lambda) = \int_{0}^r\phi_j(t,\lambda)\, dt,
\nonumber
\end{equation}
\begin{equation}
\phi_j(r,\lambda) = \sqrt{\lambda - \Big(\frac{(n-1)\rho'_j(r)}{2\rho_j(r)}\Big)^2}.
\nonumber
\end{equation}

{\it Theorem A}. 
 \label{PrefaceExpansiontheorem}
 {\it  For $1 \leq j \leq N$, assume $\rho_j(r) \geq O(r^{2\beta_j})$ with $\beta_j > 1/2$, and for $N+1 \leq j \leq N + N'$, assume $\rho_j(r) \leq O(r^{-2\beta_j})$ with  $\beta_j > 0$. 
Let $\lambda \in (E,\infty)\setminus \mathcal E$. Then, for any $a^{(in)} \in {\bf h}_{\infty}(\lambda)$, there exist unique $u \in \mathcal N(\lambda)$ and $a^{(out)} \in {\bf h}_{\infty}(\lambda)$ such that, $u$ behaves  as follows on each end $\mathcal M_j$. 

\smallskip
\noindent
(1) For $1 \leq j \leq N$, 
\begin{equation}
u \simeq \Big(\frac{\pi}{\sqrt{\lambda - E_{j}}}\Big)^{1/2}\rho_j(r)^{-(n-1)/2}\left(e^{- i\Phi_j(r,\lambda)}a^{(in)}_j(x) - e^{ i\Phi_j(r,\lambda)}a^{(out)}_j(x)\right).
\label{S0Expansionregularend}
\end{equation}

\noindent
(2) For $N+1 \leq j \leq N + N'$, 
\begin{equation}
u \simeq \Big(\frac{\pi}{\sqrt{\lambda - E_{j}}}\Big)^{1/2}\rho_j(r)^{-(n-1)/2}\Big(e^{- i\Phi_j(r,\lambda)}a^{(in)}_{j,0} - e^{ i\Phi_j(r,\lambda)}a^{(out)}_j,{0}\Big)e_{j,0},
\label{S0ExpansionCusp}
\end{equation}
where $e_{j,0}$ is the normalized eigenvector of the Laplace-Beltrami operator on $M_j$.

\noindent
(3) 
The operator
\begin{equation}
S(\lambda) : a^{(in)} \to a^{(out)}
\nonumber
\end{equation}
is unitary on ${\bf h}_{\infty}(\lambda)$.
}
 
 \bigskip
 Note that in Theorem A
 $a^{(in)}_j = a^{(out)}_j = 0$ if $\lambda < E_j$.

 The meaning of the above expansion is as follows. 
 For $u \in \mathcal N(\lambda)$, we put $v = e^{-i\sqrt{\lambda}t}u$ and
 \begin{equation}
 \begin{split}
 v^{(in)}_j &= \rho_j(r)^{-(n-1)/2}e^{i(-\Phi_j(r,\lambda)- \sqrt{\lambda}t)}a^{(in)},\\
  v^{(out)}_j &= \rho_j(r)^{-(n-1)/2}e^{i(\Phi_j(r,\lambda)-\sqrt{\lambda}t)}a^{(out)},
  \end{split}
 \nonumber
 \end{equation}
 for $1 \leq j \leq N$. We also put
 \begin{equation}
  \begin{split}
 v^{(in)}_j &= \rho_j(r)^{-(n-1)/2}e^{i(-\Phi_j(r,\lambda)- \sqrt{\lambda}t)}a_{j,0}e_{j,0},\\
  v^{(out)}_j &= \rho_j(r)^{-(n-1)/2}e^{i(\Phi_j(r,\lambda)-\sqrt{\lambda}t)}b_{j,0}e_{j,0},
  \end{split}
 \nonumber
 \end{equation}
  for $j = N+1,\cdots, N+N'$. Then, $v$ satisfies the wave equation $\partial_t^2 v + Hv = 0$ on $\mathcal M$, and 
  $v_j^{(out)}$, $v_j^{(in)}$ satisfies $\partial_t^2v_j + Hv_j =0$ asymptotically as $r \to \infty$ in $\mathcal M_j$. Let 
  $\psi_j^{(\pm)}(t) = \pm  \Phi_j(r,t) - \sqrt{\lambda}t$. Then, the phase $\psi_j^{(\mp)}(t)$ of the wave $v_j^{(in)}$, $v_j^{(out)}$ 
  is constant. The expansions in (\ref{S0Expansionregularend}) and (\ref{S0ExpansionCusp}) show that the wave front\footnote{In this case, it means a surface on which the phase is constant.} is diverging to infinity 
  as $t \to - \infty$ for $v^{(in)}$ and $t \to \infty$ for $v^{(out)}$. By this reason, in the expansion in Theorem A, the part with factor $- i\Phi_j(r,\lambda)$ is called {\it incoming} and that with factor $i\Phi_j(r,\lambda)$ is called {\it outgoing}. 
  
  Theorem A thus has the following interpretation. Omitting the time factor $e^{-i\sqrt{\lambda}t}$, 
  we send a wave $v^{(in)} = e^{-i\Phi^{(in)}}a^{(in)}$ in the remote past at infinity of $\mathcal M$, and observe the wave $v^{(out)} = e^{i\Phi^{(out)}}a^{(out)}$ coming back to infinity in the remote future. Then, the mapping, which is an analogue of Heisenberg's {\it S-matrix} in physics, 
 \begin{equation}
 S(\lambda) : a^{(in)} \to a^{(out)}
 \nonumber
 \end{equation}
 is unitary. 
 It is generally believed, and has been proved in various cases, that the S-matrix determines the whole physical system. Our aim is to prove this belief for the case of non-compact manifolds with (at most) exponetially growing or decaying ends.
 
 \subsection{Inverse scattering from regular end}

 For $j = 1, \cdots, N$, i.e. when $\mathcal M_j$ is growing to infinity, we call $\mathcal M_j$ {\it regular} end. 
 For $j = N+1, \cdots, N+N'$, i.e. when $\mathcal M_j$ is shrinking to a point, we call $\mathcal M_j$ {\it cusp}. 
 Letting
 \begin{equation}
 S_{jk}(\lambda) : a^{(in)}_j \to a_k^{(out)},
 \nonumber
 \end{equation}
 where $a^{(in)} = (a^{(in)}_1,\cdots,a^{(in)}_{N+N'})$, $a^{(out)} = (a_1^{(out)},\cdots,a^{(out)}_{N+N'})$, the S-matrix $S(\lambda)$ is an $(N+N')\times (N+N')$-matrix with operator entries $S_{jk}(\lambda)$. Our first main theorem asserts that the manifold $\mathcal M$ is determined by one entry $S_{jj}(\lambda)$ of $S(\lambda)$ associated with a regular end $\mathcal M_j$ for all energies $\lambda \in \sigma_e(H)\setminus\mathcal E$. 
 
 Assume that we are given two manifolds $\mathcal M^{(1)}$, $\mathcal M^{(2)}$ of the form (\ref{S0manifoldmathcalM}). Let 
 $S^{(i)}(\lambda) = \big(S^{(i)}_{jk}(\lambda)\big)$ be the S-matrix for $\mathcal M^{(i)}$ with size $(N^{(i)} + {N'}^{(i)})\times (N^{(i)} + {N'}^{(i)})$. Note that the number of ends of $\mathcal M^{(i)}$ is not assumed to be equal for $i = 1, 2 $ a-priori. 
 
 \medskip
{\it Theorem B.
\label{S0MainInverseTheorem}
Assume that $\rho_j(r) \geq O(r^{2\beta_jr})$ with $\beta_j > 1/2$ $(1 \leq j \leq N)$ on all regular ends and $\rho_j(r) \leq O(r^{-2\beta_j r})$ with $\beta_j > 0$ $(N+1 \leq j \leq N+N')$ on all cusps (Note that we are omitting the superscript $(i)$). Assume that 
a regular end $\mathcal M^{(1)}_1$ and a regular end $\mathcal M^{(2)}_1$ are isometric, and the associated $(1,1)$ components 
$S^{(1)}_{11}(\lambda)$ and $S^{(2)}_{11}(\lambda)$  coincide for all $\lambda \in \big(\sigma_e(H^{(1)})\setminus{\mathcal E}^{(1)}\big) \cap \big(\sigma_e(H^{(2)})\setminus{\mathcal E}^{(2)}\big)$.  Then, $\mathcal M^{(1)}$ and $\mathcal M^{(2)}$ are isometric.
}

\subsection{ Inverse scattering from slowly increasing end}

The spectral analysis becomes harder when the growth rate of the volume of the end becomes slower\footnote{It can be seen, for example, in the asymptotic expansion of solutions to the Helmholtz equation in (\ref{S0Expansionregularend}) and (\ref{S0Expansionregularend2}). The latter is more complicated than the former.}. Assume that on an regular end $\mathcal M_j$,  $\rho_j(r) \geq O(r^{2\beta_jr})$ with $0 < \beta_j \leq 1/2$. 
We fix such an end $\mathcal M_j$. Assume that $\mathcal M_j$ is diffeomorphic to $(0,\infty)\times M_j$, let $\Delta_{M_j}$ be the Laplac-Beltrami operator of $M_j$, $0 = \lambda_{0,j} < \lambda_{1,j} \leq \cdots  \to \infty$ the eigenvlues of $- \Delta_{M_j}$ and $P_{\ell,j}$ the eigenprojection associated with $\lambda_{\ell,j}$. Take a large constant $C$ which depends only on $\rho(r)$ and put for $\lambda, E > 0$
\begin{equation}
r_0(\lambda,E) = \left(\frac{2C(1 + E)}{\lambda}\right)^{1/\epsilon}, 
\nonumber
\end{equation}
$\epsilon > 0$ being a small constant. Take $\chi(r) \in C^{\infty}({\mathbb R})$ such that $\chi(r) = 0$ for $r < 1$, $\chi(r)=1$ for $r > 2$, and put
\begin{equation}
\chi_{\ell,j}(\lambda,r) = \chi\Big(\frac{r}{r_0(\lambda,\lambda_{\ell,j})}\Big).
\nonumber
\end{equation}
We also put
\begin{equation}
\varphi_j(\lambda,E,r) = \int_{r_0(\lambda,E)}^r\alpha_j(\lambda,E,s)ds,
\nonumber
\end{equation}
\begin{equation}
\alpha_j(\lambda,E,r) = \sqrt{\lambda - \Big(\frac{(n-1)\rho_j'}{2\rho_j}\Big)^2 - \frac{E}{\rho_j^2}}.
\nonumber
\end{equation}

\medskip
{\it Theorem C.
 If $0 < \beta_j \leq 1/3$, we assume that on $\mathcal M_j$, the metric is of the form
 $$
 ds^2 = (dr)^2 + \rho_j(r)h_j(x,dx).
 $$
  For the other regular ends, assume $\beta_j > 1/3$, and for the cusp ends, assume $\beta_j > 0$. 
Let $\lambda \in (E,\infty)\setminus \mathcal E$. 
Then, we have the same conclusion as in Theorem A 
except that 
on the end where $0 < \beta_j \leq 1/2$, we have the asymptotic expansion
\begin{equation}
\begin{split}
u \simeq & \Big(\frac{\pi}{\sqrt{\lambda - E_{j}}}\Big)^{1/2}\rho_j(r)^{-(n-1)/2}\\
& \times \sum_{\ell=0}^{\infty}\chi_{\ell,j}(\lambda,r)
\Big(e^{-i\varphi_j(\lambda,\lambda_{\ell,j},r)}a^{(in)}_{\ell, j} - e^{i\varphi_j(\lambda,\lambda_{\ell,j},r)}a^{(out)}_{\ell,j}\Big)
e_{\ell,j},
\end{split}
\label{S0Expansionregularend2}
\end{equation}
where $a^{(in)}_j = \sum_{\ell=0}^{\infty}a^{(in)}_{\ell,j}e_{\ell,j}(x)$, $a^{(out)}_j = \sum_{\ell=0}^{\infty}a^{(out)}_{\ell,j}e_{\ell,j}(x)$, and $e_{\ell,j}(x)$ is a normalized eigenvector of $- \Delta_{M_j}$ associated with the eigenvalue $\lambda_{\ell,j}$. 
}
 
 \medskip
 {\it Theorem D. 
Under the same assumption as in Theorem C, 
assume that 
a regular end $\mathcal M^{(1)}_1$ and a regular end $\mathcal M^{(2)}_1$ are isometric, and the associated $(1,1)$ components 
$S^{(1)}_{11}(\lambda)$ and $S^{(2)}_{11}(\lambda)$  coincide for all $\lambda \in \big(\sigma_e(H^{(1)})\setminus{\mathcal E}^{(1)}\big) \cap \big(\sigma_e(H^{(2)})\setminus{\mathcal E}^{(2)}\big)$.  Then, $\mathcal M^{(1)}$ and $\mathcal M^{(2)}$ are isometric.
}

\medskip
Comparing (\ref{PrefaceExpansiontheorem}) with (\ref{S0Expansionregularend2}), we see that the asymptotic behavior of solutions to the Helmholtz equation changes at the threshold $\beta_j = 1/2$. 

\subsection{Inverse scattering from cylindrical end}

An end $\mathcal M_j$ is said to be asymptotically cylindrical if the metric has the behavior
$$
ds^2 = (dr)^2 + h_j(r,x,dx),
$$
$$
h_j(r,x,dx) = h_j(x,dx) + O(r^{-\gamma_j}), \quad \gamma_j> 0.
$$
In this case, the expansion (\ref{S0Expansionregularend2}) is modified as follows:
\begin{equation}
\begin{split}
u \simeq & \Big(\frac{\pi}{\sqrt{\lambda }}\Big)^{1/2}\rho_j(r)^{-(n-1)/2}\\
& \times \sum_{\lambda_{\ell,j}<\lambda}
\Big(e^{-iy\sqrt{\lambda - \lambda_{\ell,j}}}a^{(in)}_{\ell,j} - e^{iy\sqrt{\lambda - \lambda_{\ell,j}}}a^{(out)}_{\ell,j} \Big)
e_{\ell,j}.
\end{split}
\nonumber
\end{equation}
Nemely, for a finite energy $\lambda$, we have only a finite number of scattering waves (channels), and the S-matrix becomes a matrix of finte size. Theorem D
holds also for this case, i.e. the manifold $\mathcal M$ is determined from the S-matrix associated with the cylindrical end (\cite{IKL09}).

\subsection{Inverse scattering from cusp}

It is known that the information of the S-matrix for the cusp does not determine the manifold (\cite{Ze}). 
To get more information,  we generalize the notion of  S-matrix by enlarging the solution spece of the Helmholtz equation $(-\Delta_{\mathcal M} - \lambda)u=0$ on the cusp. For this purpose, we assume that the end in question is a pure cusp. Namely, we assume that the metric is of the form 
$$
ds^2 = (dr)^2 + \rho_{N+N'}(r)h_{N+N'}(x,dx)
$$
on $\mathcal M_{N+N'}$. We put $N+N' = \kappa$ for the sake of simplicity. For the other cusp ends, we assume as before.

Let $0 = \lambda_{0,\kappa} \leq \lambda_{1,j\kappa} \leq \lambda_{2,\kappa} \leq \cdots$ be the eigenvalues of $- \Delta_{M_{\kappa}}$ with complete orthnormal system of eigenvectors $e_{\ell,\kappa}(x), \ell = 0,1,2,...$
We put
 \begin{equation}
 \Phi_{\kappa}(r,B) = \int_{r_0}^r \sqrt{\frac{B}{\rho_{\kappa}^2}- \lambda + \frac{(n^2-2n)}{4}\Big(\frac{\rho_{\kappa}'}{\rho_{\kappa}}\Big)^2 + \frac{(n-2)}{2}\Big(\frac{\rho_{\kappa}'}{\rho_{\kappa}}\Big)'}dr.
\nonumber
 \end{equation}
Then, there exist solutions $u_{\ell,\kappa,\pm}$ to the equation
\begin{equation}
- u'' - \frac{(n-1)\rho_{\kappa}'}{\rho_j}u' + \Big(\frac{\lambda_{\ell,\kappa}}{\rho_{\kappa}^2}-\lambda\Big)u = 0,
\nonumber
\end{equation}
which behave like\footnote{Here, $f \sim g$ means that $f/g \to 1$ as $r \to \infty$.}
\begin{equation}
u_{\ell,\kappa,\pm} \sim \rho_{\kappa}(r)^{-(n-2)/2}e^{\pm \Phi_{\kappa}(r,\lambda_{\ell,\kappa})}, \quad r \to \infty.
\nonumber
\end{equation}
Take any solution $u$ of the equation
\begin{equation}
(- \Delta_{\mathcal M} - \lambda)u = 0, \quad {\rm on} \quad \mathcal M_{\kappa}.
\nonumber
\end{equation}
Expanding it by $e_{\ell,\kappa}$, we have
\begin{equation}
(u(r,\cdot),e_{\ell,\kappa})_{L^2(M_j)} = a_{\ell,\kappa}u_{\ell,\kappa,+}(r) + 
b_{\ell,\kappa}u_{\ell,\kappa,-}(r).
\nonumber
\end{equation}
We introduce two spaces of sequences ${\bf A}_{\kappa,\pm}$ : 
\begin{equation}
{\bf A}_{\kappa,\pm} \ni \{c_{\ell,\pm}\}_{\ell=0}^{\infty} 
\Longleftrightarrow \sum_{\ell=0}^\infty|c_{\ell,\pm}|^2|u_{\ell,\kappa,\pm}(r)|^2 < \infty,\quad \forall r > 0.
\nonumber
\end{equation}
We take a 
 partition of unity $\{\chi_j\}$ on $\mathcal M$ such that $\chi_j(r) = 1$ on $\mathcal M_j\cap \{r > 2$, $\chi_j(r) = 0$ 
 on $\mathcal M_j\cap\{r < 1\}$ for $1 \leq j \leq N + N'$.
We define the generalized incoming  solution on the cusp end $\mathcal M_{\kappa}$ by
\begin{equation}
\Psi_{\kappa}^{(in)} = \chi_{\kappa}\sum_{\ell=0}^{\infty}a_{\ell,{\kappa}}u_{\ell,{\kappa},+}(r)e_{\ell,{\kappa}}(x),\quad \{a_{\ell,{\kappa}}\}_{\ell=0}^{\infty}\in {\bf A}_{{\kappa},+},
\label{Psijinpreface}
\end{equation}
which is growing as $r \to \infty$, and the generealized outgoing solution by
\begin{equation}
\Psi_{\kappa}^{(out)} = \chi_{\kappa}\sum_{\ell=0}^{\infty}b_{\ell,{\kappa}}u_{\ell,{\kappa},-}(r)e_{\ell,{\kappa}}(x),\quad \{b_{\ell,{\kappa}}\}_{\ell=0}^{\infty}\in {\bf A}_{{\kappa},-},
\label{Psijoutpreface}
\end{equation}
which is  decaying as $r \to \infty$. We also define the spaces of generalized scattering data by
\begin{equation}
{\bf h}_{\infty}^{(in)}(\lambda) = \Big({\mathop\oplus_{j=1}^N}c_j(\lambda)L^2(M_j) \Big)\oplus
 \Big({\mathop\oplus_{j=N+1}^{N+N'-1}}c_j(\lambda){\mathbb C}_j\Big)\oplus \Big(c_{\kappa}(\lambda){\bf A}_{{\kappa},+}\Big),
\nonumber
\end{equation}
\begin{equation}
{\bf h}_{\infty}^{(out)}(\lambda) = \Big({\mathop\oplus_{j=1}^N}c_j(\lambda)L^2(M_j) \Big)\oplus
 \Big({\mathop\oplus_{j=N+1}^{N+ N'-1}}c_j(\lambda){\mathbb C}_j\Big)\oplus \Big(c_{\kappa}(\lambda){\bf A}_{{\kappa},-}\Big),
\nonumber
\end{equation}
where $c_j(\lambda)$ is the characteristic function of the interval $(E_{j},\infty)$.

\medskip
{\it Theorem E. 
For any generalized incoming data $a^{(in)} \in {\bf h}_{\infty}^{(in)}(\lambda)$, there exist a  unique solution $u$ of the equation $(-\Delta_{\mathcal M} - \lambda)u=0$, and the  outgoing data $a^{(out)} \in {\bf h}_{\infty}^{(out)}(\lambda)$ such that
\begin{equation}
u- \Psi_{\kappa}^{(in)} \in \mathcal B^{\ast},
\nonumber
\end{equation}
\begin{equation}
u = \Psi_{\kappa}^{(in)} - \Psi_{\kappa}^{(out)}, \quad {\rm on} \quad \mathcal M_{\kappa},
\nonumber
\end{equation}
and on the ends $\mathcal M_j$, $1 \leq j \leq N+N'-1$, $u$ has the asymptotic form in Theorems A and C.
 Here, $\Psi_{\kappa}^{(in)}$ and  $\Psi_{\kappa}^{(out)}$ are written by (\ref{Psijinpreface}), (\ref{Psijoutpreface}) with $a_{\ell,{\kappa}}$, $b_{\ell,{\kappa}}$ replaced by the associated components of $a_{\kappa}^{(in)}$ and $a_{\kappa}^{(out)}$.
 }

\medskip
We call the mapping
\begin{equation}
\begin{split}
\mathcal S(\lambda) : {\bf h}_{\infty}^{(in)}(\lambda) \ni a^{(in)} \to a^{(out)} \in {\bf h}_{\infty}^{(out)}(\lambda)
\end{split}
\nonumber
\end{equation}
the {\it generalized scatteing matrix}. \index{generalized scattering matrix} Then, the inverse scattering theorem can be extended to the generalized S-matrix. Note that the $(\kappa,\kappa)$ component $\mathcal S_{\kappa\kappa}(\lambda)$ is an infinite matrix whose $(0,0)$ component is the usual S-matrix $S_{N+N',N+N'}(\lambda)$, which is a complex number of modulus 1. 

\medskip
{\it Theorem F.
Under the same assumption as in Theorem E, 
assume that 
the cusp ends $\mathcal M^{(1)}_{\kappa}$ and $\mathcal M^{(2)}_{\kappa}$ are isometric, and the associated $({\kappa},{\kappa})$ components 
of the generalized S-matrix  coincide for all $\lambda \in \big(\sigma_e(H^{(1)})\setminus{\mathcal E}^{(1)}\big) \cap \big(\sigma_e(H^{(2)})\setminus{\mathcal E}^{(2)}\big)$.  Then, $\mathcal M^{(1)}$ and $\mathcal M^{(2)}$ are isometric.
}

\medskip
As above, the number of ends of $\mathcal M^{(1)}$ and $\mathcal M^{(2)}$ are not assumed to be equal a-priori.

\

\subsection{Riemannian metric with continuous spectrum}

Properties of continuous spectrum of a manifold (i.e. that of the Laplacian) depend largely on its volume growth. 
Let $M$ be a compact manifold of dimension $n-1$.
We consider the Riemannian metric on $\mathcal M = (0,\infty)\times M$ of the form
\begin{equation}
ds^2 = (dr)^2 + g^{(M)}(r,x,dx),
\nonumber
\end{equation}
where\footnote{Using Einstein's summation convention, $a_{ij}b^{ij} = \sum_{i,j=1}^na_{ij}b^{ij}$.} $g^{(M)}(r,x,dx) = g^{(M)}_{\,\,\,\,\,\,\,\, ij}(r,x)dx^idx^j$ is a metric on $M$ depending smoothly on  $r > 0$. 
Identifying $g^{(M)}$ with an $(n-1)\times(n-1)$ matrix $\big( g^{(M)}_{\,\,\,\,\,\,\, ij}\big)$, we let 
\begin{equation}
g =g(r,x) = \det\big( g^{(M)}_{\,\,\,\,\,\,\, ij}\big). 
\nonumber
\end{equation}
For $\kappa \in {\mathbb R}$, let $S^{\kappa}$ be the set of  $C^{\infty}$ functions on $(0,\infty)\times M$ having the following property :
\begin{equation}
S^{\kappa} \ni f \Longleftrightarrow |\partial_r^{\ell}\partial_x^{\alpha}\, f(r,x)| \leq C_{\ell\alpha}(1 + r)^{\kappa - \ell }, \quad \forall \ell,  \alpha,
\label{IntoDefineSkappa}
\end{equation}
where $C_{\ell\alpha}$ is a constant. \index{$S^{\kappa}$}
This definition is naturally extended for tensor fields\footnote{We define the $S^{\kappa}$-norm for manifolds with conic singularties in Defintion \ref{S2DefineSingularSkappa}, which is actually used throughout this article. }. 

Now let us consider  $g = g(r,x)$ such that 
\begin{equation}
\frac{g'}{4g} - \frac{(n-1)c_0}{2} - c_1 r^{-\alpha} \in S^{-1- \epsilon}, 
\label{S1A-1}
\end{equation}
where  $  ' = \partial_r, $ and $c_0, c_1$, $\alpha$,  $\epsilon$ are real constants such that $ \alpha > 0, \epsilon > 0$. This constant $c_0$ is important, since it determines the infimum of the continuous spectrum of the Laplacian. 
Integrating the equation
$$
\frac{g'}{4g} - \frac{(n-1)c_0}{2} - c_1 r^{-\alpha} = O(r^{-1-\epsilon}),
$$
we obtain $g = \rho^{2(n-1)} O(1)$, where
\begin{equation}
\rho(r) = 
\left\{
\begin{split}
& \exp\big(c_0r + c_1'r^{1-\alpha}\big), \quad 0 < \alpha< 1,\\
& \exp(c_0r) \,r^{\beta},\quad \alpha=1, \\
& \exp(c_0r)\big(1 + O(r^{-\delta})\big), \quad
\alpha > 1, \quad  \delta = \min\{\alpha-1,\epsilon\},
\end{split}
\right.
\label{S1rho(r)}
\end{equation}
with
$c_1' = {\beta}/(1 - \alpha),  
\beta = {2c_1}/(n-1)$. 
We put  $ h(r,x,dx) = \rho(r)^{-2}g(r,x,dx)$, 
which is asymptotically equal to a metric independent of $r>0$, i.e.  $h(r,x,dx) \to h^{(M)}(x,dx)$ as $r \to \infty$. 
Thus, our typical example of the metric is written in terms of $\rho(r)$ in (\ref{S1rho(r)}) as
\begin{equation}
ds^2 = (dr)^2 + \rho(r)^2h(r,x,dx).
\nonumber
\end{equation}
Let $S(r) = \{r\}\times M$ be the section of $\mathcal M$ at $r$. Then  the volume of $S(r)$ is growing as $r \to \infty$ if either $c_0> 0$ or $ c_0=0,  \beta > 0$, 
and is shrinking to 0  if either $ c_0 < 0$, or $c_0=0,  \beta < 0$. 
Let us  call the metric of the form
\begin{equation}
ds^2 = (dr)^2 + \rho(r)^2h(r,x,dx),
\label{S1ds2=ds2+rho2d}
\end{equation}
a {\it perturbed warped product} metric, if it has the property
 $h(r,x,dx) \to h^{(M)}(x,dx)$ as $r \to \infty$. 
 In Chapter 1, \S \ref{Transformmetric}, we show that 
 the metric with cross term
\begin{equation}
ds^2= a(t,z)(dt)^2 + 2w(t)b_i(t,z)dtdz^i + w(t)^2c_{ij}(t,z)dz^idz^j,
\label{S1Geberalmetricform}
\end{equation}
\begin{equation}
w(t)^{-1} \in S^{-\kappa}, \ a(t,z)-1 \in S^{-\lambda}, \ b_i(t,z) \in S^{-\mu}, \ c_{ij}(t,z) - h_{ij}(z) \in S^{-\nu},
\nonumber
\end{equation}
is transformed into the form (\ref{S1ds2=ds2+rho2d}), where $z$ denote local coordinates on $M$. 

\subsection{Eaxmples of manifolds}

Simple examples can be constructed by taking $\mathcal M$ to be the surface of revolution with the metric induced from the Euclidean metric in ${\mathbb R}^{n+1}$ : $x_{n+1} = (x_1^2 + \cdots + x_n^2)^{1/(2\beta)}$. Negelecting the singularity at $x = 0$, the case $\beta=1$ corresponds to the conical surface, and the case $\beta=1/2$ to the parabola.

In (\ref{S1A-1}), we have restricted the growth order of $|\log g|$ at most linearly. This is a natural restriction, since outside this range, the Laplacian may not have continuous spectrum. To see it, let us consider the warped product metric
$$
ds^2 = (dr)^2 + \rho(r)^2h_M(x,dx).
$$
The  associated Laplacian is unitarily equivalent to 
(see (\ref{S7Deltagtranasformed}))
$$
- \partial_r^2 - \rho^{-2}\Lambda + \big(\frac{\rho'}{4\rho}\big)^2 + \big(\frac{\rho'}{4\rho}\big)',
$$
where $\Lambda$ is the Laplace-Beltrami operator for $h_M(x,dx)$. Letting $\{\lambda_{\ell}\}_{\ell=0}^{\infty}$ be the eigenvalues for $- \Lambda$, it is unitarily equivalent to 
\begin{equation}
{\mathop\oplus_{\ell=0}^{\infty}}\Big(- \partial_r^2 + \frac{\lambda_{\ell}}{\rho^2} + \big(\frac{\rho'}{4\rho}\big)^2 + \big(\frac{\rho'}{4\rho}\big)'\Big).
\nonumber
\end{equation}
If $\rho'/\rho = \pm r^{\epsilon}$, it has only the discrete spectrum. 

The metric of the warped product form includes many important examples. It is the hyperbolic metric when $\rho(r) = e^{\pm r}$, and the Euclidean metric when $\rho(r) = r^2$. The manifold with cylindrical end is the case $\rho(r)$ = constant.
Our assumption means that $\rho(r)$ is in between $e^{cr}$ and $e^{-cr}$ with $c > 0$. 
Therefore, as long as we start from the asymptotic expansion (\ref{S1A-1}), the class of the metric we employ in this paper seems to be optimal for the study of forward and inverse scattering on Riemannian manifolds.

\subsection{Conic singularities}
{\mtext Our another aim is to introduce a class of manifolds allowing cone-like singuralities. A simple example of manifolds with cone-like singularities is the sector
$$
S_{\alpha} = \{z \in \mathbb C\, ; \, 0 \leq {\rm arg}\,(z) \leq \alpha\},  \quad 0 < \alpha < 2\pi
$$
with two boundaries $\{z \in S_{\alpha}\, ; \, {\rm arg}\,(z) = 0\}\cup 
\{z \in S_{\alpha}\, ; \, {\rm arg}\,(z) = \alpha\}$ identified (see Figure  2 and  Example \ref{Example:CMGA1}). 
\begin{figure}[hbtp]
\setlength{\unitlength}{1.2mm}
\centering
\begin{picture}(0,30)
\put(-5,15) {\vector(1,0){10}}
\put(-45,15){\line(3,1){30}}
\put(-30,15){$\alpha$}
\put(-45,15){\line(3,-1){30}}
\bezier{800}(-35,18)(-32,15)(-35,12)
\put(15,15){\line(4,1){30}}
\put(15,15){\line(4,-1){30}}
\bezier{800}(45,23)(50,15)(45,7)
\bezier{800}(45,23)(40,15)(45,7)
\end{picture}
\label{fig:math1/2}
\caption{$S_{\alpha}$}
\end{figure}
One can induce the differential structure of $\mathbb R^2$ to $S_{\alpha}\setminus\{0\}$ to make it a $C^{\infty}$-manifold. One can also induce the Euclidean metric to $S_{\alpha}\setminus\{0\}$, which is not smooth at $z = 0$. 
If $2\pi/\alpha \in {\mathbb N}$, by the group action of rotation of angle $\alpha$, $S_{\alpha}$ becomes a $C^{\infty}$-manifold including $0$, although the metric is singular at $0$. This is a simple example of orbifold. 
Similarly, hyperbolic manifolds are orbifolds, and the singularity at the top of the cone does no harm for the spectral analysis (see e.g. \cite{ElGrMen}). Ths is no longer the case when $2\pi/\alpha \not\in {\mathbb Q}$. 
However, we can develope a theory of conic manifold with group action (CMGA) in order to allow this sort of singularities in the spectral analysis and inverse problems on more general class of manifolds.
}

Manifolds with singularities have been objects of long issue. For example, the regularity of solutions to the Dirichlet problem for Laplacians around corners, the behavior of solutions to the wave equation near cracks or thin sets are significant problems of classical physics. In differential geometry and
in global analysis, 
spaces having conic or more general singularities have been extensively studied,
 see e.g. \cite{MM87}. In particular, scattering problems have been studied from microlocal point of view e.g. in \cite{{Baskin}, Cheeger-Taylor1,Cheeger-Taylor2, {Kalka}, Mazzeo,Melrose-Vasy-Wunsch1,Melrose-Vasy-Wunsch2,Melrose-Wunsch, Galkowski17}.
 In this paper we use more classical functional analytical techniques and our focus is on the inverse scattering problems.

In  inverse spectral problems for elliptic
equations and in equivalent problems for heat and wave equations \cite{KKL04b}, singular spaces have appeared in the study of stability, see 
\cite{AKKLT,BKL2020,FIKLN,Liu-Oksanen}.
Indeed, for the stability problem one has to consider a class of spaces that 
is compact in suitable topology, for instance in the Gromov-Hausdorff or Lipschitz sense. 
In \cite{AKKLT}, a stability result 
for the inverse boundary spectral problem is studied in the class of manifolds  for which the curvature of the manifold, the diameter
and, in addition, the injectivity radius  are bounded. The Gromov-Hausdorff closure of these manifolds contain
manifolds which metric tensor is not smooth.
When the injectivity radius is not bounded from below,
the theory changes radically as  geometric collapses  can appear.
The study of collapsing manifolds has been
an important trend in  modern differential geometry. The celebrated papers of 
by Perelman  \cite{Per1} and Cheeger-Fukaya-Gromov  \cite{CFG}  are examples of the 
progress of collapse and
 the metric geometry of the non-smooth spaces.
 When the   manifolds converge  to a lower
dimensional object  we say that the geometry collapses.
The limit space can be very non-smooth  and not even a manifold, 
but generally an Alexandrov space \cite{BRP}. In the case that manifolds of dimension $n$ converge 
to a space whose Hausdorff dimension is $n-1$, the limiting objects are orbifolds.
The inverse interior spectral problems for compact orbifolds is studied in \cite{KLY09}.
A related inverse scattering problem for non-compact two-dimensional orbifolds is studied in \cite{IKL11}.
In this paper we extend this research for higher dimensional manifolds with more general singularities and study manifolds, which are non-compact and may 
have certain type of conic singularities. In particular, the class of  conic singularities studied in this paper contain
all  orbifold type singularities.
This enables us to extend the {\it boundary control method}\footnote{As will be explained in Subsection \ref{Subsection010Relatedworks},  this is the main tool for the inverse scattering procedure.}
on manifolds with conic singularities, since it is based on the combination of spectral theory for Laplacian, wave equation and geodesics on it.
We can then determine the manifold from the knowledge of S-matrix via the boundary control method. The details of the  metric  with conic singularities will be given in Definition \ref{Definitionconicchart} in \S \ref{Mfdandgroupaction}.


\subsection{Main subjects}

We assume that $\mathcal M$ is a regular conic manifold in the sense to be defined in Chapter 1, \S \ref{Mfdandgroupaction}. 
We study the spectral and scattering theory for the Laplace  operator $\Delta_{\mathcal M}$ on $\mathcal M$, in particular,  

\begin{itemize}
\item  Limiting absorption principle for the resolvent $R(z) = (- \Delta_{\mathcal M}-z)^{-1}$,

\item Spectral representation of $- \Delta_{\mathcal M}$,

\item Helmholtz equation for $- \Delta_{\mathcal M}$ and the S-matix,

\item Inverse scattering.

\end{itemize}
The first 3 subjects constitute the main  parts of the forward problem. Once we have solved them,  the inverse scattering problem can be solved in the same way as in \cite{IKL09}, \cite{IKL13}.

\medskip
Our main results are Theorem \ref{thm: main result for IP}  in Chapter 2, \S \ref{Uniquenessofinversescattering},  which asserts that the whole manifold is determined from the knowledge of the (generalized) S-matrix for all energies. 
To achieve this final goal, we prepare the following theorems, which are of independent importance from the view point of spectral theory:

\begin{enumerate}
\item Chapter 1 : Forward problem
\begin{itemize}
\item Theorem \ref{S2RellichTh1} (Rellich-Vekua theorem)
\item Theorem \ref{ModelLAPforbfh} (Limiting absorption principle for the resolvent)
\item Theorem \ref{GeneralizeFourierTh} (Resolvent asymptotics  and spectral representation)
\item Theorem \ref{S13SolspaceHelmholtzeq} (Characterization of the solution space of the Helmholtz equation)
\item Theorem \ref{Smatrixtheorem} (Asymptotic expansion of solutions to the Helmholtz equation and the S-matrix)
\end{itemize}
\item Chapter 2 : Inverse problem
\begin{itemize}
\item Theorem \ref{S15SmatrixdeterminesLW2} (The S-matrix and the generalized S-matrix determine the source-to-solution map)
\end{itemize}
\begin{itemize}
\item Theorem \ref{S15SmatrixdeterminesLW2cusp} (The generalized S-matrix determine the N-D map)
\end{itemize}
\end{enumerate}

\subsection{Related works}
\label{Subsection010Relatedworks}
Analysis of Laplacians on on non-compact manifolds and inverse spectral problems on compact or non-compact manifolds are our main concern. Since both of them are rather big subjects, we should better restrict ourselves to only recent developments. As for the latter subject, one should take notice of the evolution of multi-dimensional inverse problems of PDE which started in 1080s from the works of Sylvester-Uhlmann \cite{SylUhl87}, Nachman \cite{Nach88}, Khenkin-Novikov \cite{KheNov87}. There is a close connection between the inverse boundary value problem and the inverse scattering problem. For the survey of the recent development of inverse boundary value problems,  see \cite{U}, \cite{Uhl2}, \cite{MSalo} and the references therein. The other important stream is the boundary control method due to Belishev 
\cite{Be87} and its development to Riemannian manifold by Belishev-Kurylev \cite{BeKu92}. See \cite{Be97} for a survey of this method. 
See also an expository work \cite{KKL01}, and \cite{IsKu10} for the extension of BC method to inverse scattering problems. 
As for the spectral analysis on non-compact manifolds, a general framework based on microlocal calculus is provided by Melrose in \cite{Me95}, \cite{MM87}. Hyperbolic manifold is an interesting topic from analysis, geometry as well as number theory, and has been discussed in many articles
(\cite{Faddeev}, \cite{Bor07}, \cite{BJP05}, \cite{BP11}, \cite{Perry1}, \cite{Perry2},  \cite{Gu}, \cite{GM11}, \cite{GuZw95}, \cite{GuZw97}).
In particular, inverse scattering on hyperbolic manifolds by the boundary control method was studied by \cite{JSb1}, \cite{JSb2}, \cite{SaBa05}, \cite{GuSa}, \cite{GuSa2}. For the case including more slowly (polynomially) growing manifolds, see \cite{GST},  \cite{Kum13}, \cite{Kum13(1)}, and recent works \cite{ItoSki20},
\cite{ItoSkiI}, \cite{ItoSkiII}. These works are based on the functional analytic technique utilized in spectral theory.

We have been working on important typical cases of such metrics as cylindrical ends in \cite{IKL09},  asymptotically hyperbolic ends in  \cite{IKL13}, \cite{IKL13(2)}, \cite{IKL11},   the asymptotically Euclidean case in \cite{IsKu10}. In \cite{IKL11},  a relation between the generalized S-matrix and arithmetic surfaces was discussed. As was stated in \cite{IKL13(1)}, we  recognized that these typical metrics are embedded in a series of metrics, which can be understood in a unified scope\footnote{The main part of this paper, in particular for the case of smooth metric, has already been completed around 2014.}.  
Not only the behavior of the metric at infinity, but also the singularity of the manifold has been  our concern. 

In this paper, we  employ the method of integration by parts, which originates from the work of Eidus \cite{Eidus}. Since this technique is elementary, hence basic, it is flexible in applying to new problems. Moreover, although it is classical, its basic ideas are absorbed and transformed into new machieneries of this field. 
We do not pursue full generalities, although we recognize that many parts of our arguments can be performed under weaker assumptions.
	  Our intention is to try to simplify the whole theory
	  and to make the assumptions and statements as clear and understandable as possible.


\subsection{More precise assumptions of the metric}
\label{subsectionAssumptions}


Our assumptions on the metric change slightly from section to section according to the subjects. Therefore, let us summarize here the final assumptions on the metric. {\mtext After introducing CMGA in Definition \ref{Definitionconicchart} and $S^{\kappa}$ in Definition \ref{S2DefineSingularSkappa}, we assume that each end is of the form $(0,\infty)\times M_i$, where $M_i$ is CMGA of dimension $n-1$. We further assume as follows.}

\medskip
\noindent
 {\bf (A-1)}\index{(A-1)}  {\it The ends $\mathcal M_1, \cdots, \mathcal M_N$ are regular, and the ends $\mathcal M_{N+1},\cdots,\mathcal M_{N+N'}$ are cusp.}

\medskip
\noindent
{\bf (A-2)} \index{(A-2)} \label{AssumptionA2}
\  {\it  For each end $\mathcal M_j$, $j = 1, \cdots, N+N'$, there exist  constants $c_{0,j} \in {\mathbb R}$, $\alpha_{0,j} > 0, \gamma_{0,j} > 1$ and a metric $h_{M_j}(x,dx)$ on $M_j$ such that 
\begin{equation}
\frac{\rho'_j(r)}{\rho_j(r)} - c_{0,j} \in S^{-\alpha_{0,j}},
\label{rho'rho-c0}
\end{equation}
\begin{equation}
h_j(r,x,dx) - h_{M_j}(x,dx)  \in S^{-\gamma_{0,j}}.
\label{h-hMgamma}
\end{equation}
}
\medskip
\noindent
{\bf (A-3)}\index{(A-3)} \ 
\label{AssumptionA3} {\it  On  regular ends $\mathcal M_j$,   there exist constants $\beta_{0,j} > 0, r_{0,j}>0$ such that
\begin{equation}
\frac{\rho'_j(r)}{\rho_j(r)} \geq \frac{\beta_{0,j}}{r} \quad {\it  for} \quad r > r_{0,j}.
\label{S1growthCond}
\end{equation}
Moreover, the constants $\alpha_{0,j}$, $\gamma_{0,j}$ in (\ref{rho'rho-c0}), (\ref{h-hMgamma}) satisfy 
\begin{equation}
\left\{
\begin{split}
& \alpha_{0,j} > 0,  \ {\it for} \ {\it the} \ {\it regular} \ {\it ends} \ {\it with}\ \beta_{0,j} > \frac{1}{2},\ c_{0,j} \geq 0, \\
&\alpha_{0,j} >1, \ {\it for} \ {\it the} \ {\it regular} \ {\it ends} \ {\it with}\ 0 < \beta_{0,j} \leq \frac{1}{2},\ c_{0,j} \neq 0,  \\
&\alpha_{0,j} > \frac{1}{2},  \ {\it for} \ {\it the} \ {\it regular} \ {\it ends} \ {\it with}\ 0 < \beta_{0,j} \leq \frac{1}{2},\ c_{0,j} = 0. 
\end{split}
\right.
\label{Condalphabeta,1}
\end{equation}

On cusp ends $\mathcal M_j$, there exists a constant $\beta_{0,j} < 0$ such that
\begin{equation}
\rho_{j}(r) \leq C(1 + r)^{\beta_{0,j}}
\label{A-2CuspCondition}
\end{equation}
holds. Furthermore, for all cusp ends
\begin{equation}
\left\{
\begin{split}
& \alpha_{0,j} > 1, \quad {\it if}\quad c_{0,j} \neq 0, \\
& \alpha_{0,j} > 1/2,  \quad{\it if} \quad c_{0,j} = 0.
\end{split}
\right.
\label{A-1alpha11/2}
\end{equation}
}

\medskip
The results for the forward problem are proven under the assumptions (A-1), (A-2), (A-3).  The inverse scattering from the regular end with $\beta_j > 1/3$ is proven also under these assumptions.

To consider the inverse scattering from the regular end with $0 < \beta_j \leq 1/3$ and the cusp end, we impose the following 
assumtion (A-4). 

\medskip
\noindent
{\bf (A-4)}\label{(A-4)} {\it The metric on $\mathcal M_j$ is of the warped product form 
\begin{equation}
ds^2\Big|_{\mathcal M_j} = (dr)^2 + \rho_j(r)^2h_{M_j}(x,dx),
\label{S1Intro1/3warpedproductmetric}
\end{equation}
for either of the following cases (A-4-1) or (A-4-2):}

\medskip
\noindent
{\bf (A-4-1)} {\it There exists $1 \leq j \leq N$ such that for the regular end $\mathcal M_j$,  $\beta_j$ satifies
$$
0 < \beta_j \leq 1/3,
$$
and $\rho_j(r)$ satisfies (A-2) and (A-3).
}
\index{(A-4-1)}

\medskip
\noindent
{\bf (A-4-2)} {\it There exists $N + 1 \leq j \leq N + N'$ such that for the cusp end $\mathcal M_j$, $\rho_j(r)$ satisfies the condition (\ref{S7rhoasymp1}). 
}
\index{(A-4)}

In addition to the assumptions (A-1) $\sim$ (A-4), 
we also impose the assumptions (C-1) $\sim$ (C-4) (to appear in  Chapter 1, Subsection \ref{SubsectionConicmfdgroupaction}, 
 (D) (to appear in Chapter 1, Subsection \ref{SubsectionRegukarityofdomain}) and 
 (L) (to appear in Chapter 2, Section \ref{Uniquenessofinversescattering}).

\subsection{Thresholds}


The potential term of the Schr{\"o}dinger operator $- \Delta + V(x)$  in ${\mathbb R}^n$ with $V(x) = O(|x|^{-\alpha})$ as $|x| \to \infty$ is said to be short-range if $\alpha > 1$ and long-range if $0 < \alpha \leq 1$. The border line $\alpha = 1$ appears in the case of the Coulomb potential $V(x) = - C/|x|$, across which solutions to the Schr{\"o}dinger equation have different asymptotic behavior at infinity. In our case, this change of asymptotic behavior occurs across $\alpha = 1$ when $c_0 \neq 0$, and $\alpha = 1/2$ when $c_0=0$. 

We now explain the thresholds we will encounter in the spectral analysis of $\mathcal M$.

\smallskip
(i) In the first place, we have avoided the case in which $\rho(r)$ is  constant. This corresponds to the case of asymptotically cylindrical end, and we have studied it in \cite{IKL09}. Of course,  we  can include it in this paper, however, the statements of results are different in many places. Therefore, to simplify the description, it is better to deal with it separately.


\smallskip
(ii) For the regular ends, we have thresholds.

\begin{itemize}
\item  Transformation of the metric (\ref{S1Geberalmetricform}) into the perturbed warped product (\S \ref{Transformmetric}): $\beta = 1$, and $\beta = 1/2$. 
\end{itemize}

Note that the case $\beta = 1$ corresponds to the asymptotically Euclidean metric.

\begin{itemize}
\item Rellich-Vekua theorem (\S \ref{SectionRellichVekua}): $\beta = 1/3$. 
\end{itemize}

We have proved the Rellich-Vekua type theorem for the case $\beta > 1/3$, however it is open\footnote{As a matter of fact, we still do not know wherther it is a real border-line or not. } for $\beta \leq 1/3$. 



\begin{itemize}
\item  Resolvent asymptotics (\S \ref{ResolventAsymptotics1}, \S \ref{ResolventAsympto2}): 
$$
\alpha = 1\  (c_0 \neq 0), \quad \alpha = 1/2 \ (c_0=0), \quad \beta = 1/2, \quad \gamma = 1.
$$
\end{itemize}
The asymptotic behavior at infinity of the resolvent is strongly affected by all constants $c_0, \alpha_0, \beta_0, \gamma_0$. 
For $\beta > 1/2$, we allow the metric of long-range behavior, however, for  $0 < \beta \leq 1/2$, we deal with  only metrics of short-range behavior. Across the threshold $\alpha = 1$ for the case $c_0 \neq 0$, and $\alpha = 1/2$ for the case $c_0=0$, the asymptotic expansion at infinity of the resolvent changes its form. As for $\gamma$, we consider only the case $\gamma > 1$ in this article. For $\gamma \leq 1$, the behavior will be different from the one we give here. 
 
\subsection{Notations}

For Banach spaces $X$ and $Y$,
${\bf B}(X;Y)$ is the space of all bounded linear operators from $X$ to $Y$, and ${\bf B}(X) = {\bf B}(X;X)$. For a Hilbert space $\bf h$, $L^2(I;{\bf h};m(r)dr)$ is the $L^2$-space of ${\bf h}$-valued functions on an interval $I \subset {\mathbb R}$ with respect to the measure $m(r)dr$, whose inner product is defined by
\begin{equation}
(f,g) = \int_I (f(r),g(r))_{\bf h}m(r)dr.
\nonumber
\end{equation}
$H^s(\mathcal M)$ denotes the Sobolev space of order $s$ (with respect to $L^2$ derivatives) on a manifold $\mathcal M$. For a self-adjoint operator $A$, $\rho(A)$ denotes its resolvent set, and $\sigma(A)$, $\sigma_d(A)$, $\sigma_p(A)$, $\sigma_e(A)$, $\sigma_c(A)$ and $\sigma_{ac}(A)$ denote its spectrum, discrete spectrum, point spectrum (the set of eigenvalues), essential spectrum, continuous spectrum, and absolutely continuous spectrum, respectively (see e.g. \cite{Kato66}, \cite{Resi}).

For the reader's convenience, we give here a brief list of symbols used frequently in the text. 

\bigskip

\begin{center}
List of symbols
\end{center}
\medskip

\begin{center}
\begin{tabular}{llll|llll}
  $S^{\kappa}$ & (\ref{IntoDefineSkappa}) (\ref{DefineSkappanew}) & & & & $B(r), \Lambda(r)$& (\ref{S5B(r)andLambda(r)}) \\
 $D(k)$ &(\ref{DefDpmk})  & & & &  
   $v = D(k)u$ &  (\ref{Definevpm}) \\
  $w = \sqrt{B(r)}u$ & (\ref{Definevpm}) & & & &  
  $E_{0,i}$ & (\ref{DefE0i}) \\
 $E_{0,tot}$ &  (\ref{DefE0i}) & & & &  $\mathcal T$ & (\ref{E01cdotsE0N+N'})  \\  
 $\mathcal E$ & (\ref{S9ExceptionalSet})  & & & & $f \simeq g$ & (\ref{Definefsimeqg}) \\
  $f(r) \sim g(r)$&  (\ref{Definef(r)simg(r)}) & & & & 
   $f \asymp g$  & (\ref{S1AsypDefine}) 
   \end{tabular}
\end{center}

\bigskip

\subsection{Acknowledgement}
H. I. is supported by Grants-in-Aid for Scientific Research (S) 15H05740, (B)16H0394, Japan Society for the Promotion of Science. He expresses his gratitude to JSPS.

\bigskip
{\mtext We (H. I., M. L. and Yaroslav Kurylev) were working on the theme of this paper for more than 15 years, and know that there are still lots of room for improvements. However, H. I and M. L think that now it is the time to publish and dedicate this paper to Y. Kurylev, who passed away 
in 2019.}

 



\tableofcontents

\chapter{Spectral and scattering theory}

{\mtext We start with introducing manifolds with conical singularities (CMGA) and study the domain of its Laplacian. The manifold $\mathcal M$ is assumed to be CMGA and its ends $\mathcal M_i$ are of the form $(0,\infty)\times M_i$, where $M_i$ is CMGA of dimension $n-1$. Chapter 1 is devoted to the forward problem of scattering, and the main aim is to prove the limiting absorption principle of the resolvent and derive the spectral representation. Finally, the S-matrix is introduced.}

\section{Group action and conic manifolds}
\label{Mfdandgroupaction}
Our definition of conic manifold is rather involved. 
We begin with the notion of orbifolds, and then introduce conic manifolds admitting group action. 
In the following, $B^k(x_0,R)$ is an open ball in ${\mathbb R}^k$ of radius $R$ centered at $x_0 \in {\mathbb R}^k$. 


\subsection{Conic chart}
\begin{definition}
\label{Definitionofconicchart}
Given a topological space $U$, we call a triple $(\widetilde U, \Gamma, \pi)$ a conic chart of $U$ if it has the following properties:
\begin{itemize}
\item $\widetilde U$ is an open set in ${\mathbb R}^n$ having the form $
\widetilde U = \widetilde W \times \widetilde V$, where $\widetilde W = B^k(0,R_0)$, $\widetilde V = B^{n-k}(0,R_1)$ for some $R_0, R_1 > 0$ and $0 \leq k \leq n$ satisfying $k \neq n-1$.
\item $\Gamma$ is a finite group acting on $\widetilde U$, leaving  $\widetilde W$ invariant.
\item $\pi : \widetilde U \to U$ is defined to be $\pi = \widetilde{\Phi}\circ \widetilde \pi$, where $\widetilde{\pi} : \widetilde U \to \Gamma\backslash \widetilde U$ is the canonical projection to the set of orbits, and $\widetilde{\Phi} : \Gamma \backslash \widetilde U \to U$ is a homeomorphism.
\end{itemize}
\end{definition}

\begin{figure}[hbt]
\setlength{\unitlength}{0.7mm}
\begin{picture}(100,50)(0,-5)
\put(60,30){\makebox(30,15){$\widetilde U = \widetilde W \times \widetilde V$}}
\put(64,0){$\Gamma\backslash \widetilde U$}
\put(68,31){\vector(0,-1){24}}
\put(72,17){$\widetilde{\pi}$}
\put(27,0){$U$}
\put(60,1){\vector(-1,0){26}}
\put(49,-5){$\widetilde{\Phi}$}
\put(58,31){\vector(-1,-1){24}}
\put(45,25){$\pi$}
\end{picture}
\caption{Conical chart}
\label{fig.ConicalChart}
\end{figure}
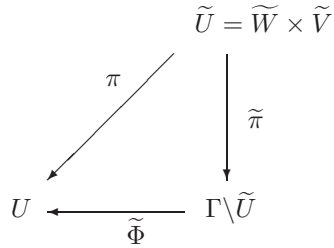


  Define $Y : \widetilde U \to \widetilde W$ and $Z : \widetilde U \to \widetilde V$ by
\begin{equation}
\widetilde U \ni x \to x = (y,z) = (Y(x),Z(x)).
\nonumber
\end{equation}
Then, by definition, $\Gamma$ leaves $\widetilde W$ invariant if 
\begin{equation}
Y(\gamma\cdot x) = Y(x), \quad \forall x \in \widetilde U, \quad \forall \gamma \in \Gamma.
\nonumber
\end{equation}
The isotropy group  for $x \in \widetilde U$ is defined by
\begin{equation}
\Gamma(x) = \{\gamma \in \Gamma\, ; \, \gamma\cdot x = x\}.
\nonumber
\end{equation}
For $k = 0$, we put
\begin{equation}
\widetilde U^{reg} = \{x \in \widetilde U\, ; \, \Gamma(x) = \{e\}\},
\nonumber
\end{equation}
where $\{e\}$ means the subgroup of $\Gamma$ consisting only of the unit.  For $0 < k < n$, we put
\begin{equation}
\widetilde U^{reg} = \Big(\widetilde W \times \big(B^{n-k}(0,R_1)\setminus\{0\}\big)\Big) \cap 
\big\{x \in \widetilde U \, ; \, \Gamma(x) = \{e\}\big\}.
\nonumber
\end{equation}
Finally, we put
\begin{equation}
\widetilde U^{sing} = {\widetilde U} \setminus \widetilde U^{reg},
\nonumber
\end{equation}
\begin{equation}
U^{reg} = \pi\big(\widetilde U^{reg}\big).
\nonumber
\end{equation}
Note that   $\widetilde U = \widetilde V$ for $k = 0$, in which case $\widetilde U^{sing}$ is either an empty set or $\{0\}$ in $\widetilde U$, and  $\widetilde U = \widetilde W$ for $k =n$, in which case $\widetilde U^{sing}$ is an empty set

Later, we introduce a Riemannian metric $\widetilde g$ on $\widetilde U$. Then a natural choice of  $\Gamma$ is a group of isometric transformations with respect to $\widetilde g$. Typical examples are a finite subgroup of $SO(2)$ when $\widetilde g$ is the Euclidean metric (Example \ref{Example:CMGA1}), and a finite subgroup of $SL(2,{\mathbb R})$ when $\widetilde g$ is the hyperbolic metric (Example \ref{S2ExampleModular}).

It is not necessary to define the differential structure of $U$, since the calculus on $U$ is done by lifting to  $\widetilde U$. 
Let $\mathcal F_{\widetilde U}$ and $\mathcal F_{U}$ be the set of all functions on $\widetilde U$ and $U$, and 
$\mathcal S_{\widetilde U}$ the set of $\Gamma$-invariant functions on $\widetilde U$, i.e. 
$$
v \in \mathcal S_{\widetilde U} \Longleftrightarrow v(\gamma\cdot x) = v(x), \quad \forall x  \in \widetilde U, \quad 
\forall \gamma \in \Gamma.
$$
For $v \in \mathcal F_{U}$, define
\begin{equation}
\pi^{\ast}v (x) = v(\pi(x)).
\label{Definepiastv}
\end{equation}
Then, since 
$\pi(\gamma\cdot x) = \pi(x)$ for all $\gamma \in \Gamma$ and $x \in \widetilde U$ by definition, 
we have $\pi^{\ast}v \in \mathcal S_{\widetilde U}$, and the map
\begin{equation}
\mathcal F_U \ni v \to \pi^{\ast}v  \in \mathcal S_{\widetilde U}
\label{Definepiastv}
\end{equation}
is a bijection. 
The operator $P$ defined by
\begin{equation}
 \mathcal F_{\widetilde U} \ni  u(x) \to Pu(x) = 
\frac{1}{\sharp \Gamma}\sum_{\gamma\in \Gamma}
u(\gamma\cdot x) \in \mathcal S_{\widetilde U}
\label{ConicprjectionP}
\end{equation}
is a surjection. 
Moreover, $P$
is an orthogonal projection on $L^2(\widetilde U)$, i.e.
$$
P^2 = P, \quad P^{\ast} = P.
$$
 The regularity of functions on $U$ is defined through $P$. Namely, 
$f$ on $U$ is in $C^m$ if and only if there exists $\widetilde f \in C^m(\widetilde U)$ such that $P\widetilde f = f$, which is equivalent to $\pi^{\ast}f \in C^m(\widetilde U)$. The same remark  applies to Sobolev spaces, which we explain in Subsection \ref{Subsection2.6Sobolevspace}.

\begin{example}
(An $n$-dimensional topological manifold). 
{\rm 
For an $n$-dimensional topological manifold $M$, any $p \in M$ has a coordinate neighborhood $(U, \varphi)$, where $U$ is an open neighborhood of $p$ in $M$ and $\varphi : U \to \widetilde U$ is a homeomorphsm with $\widetilde U$ an open set in ${\mathbb R}^n$. Then, 
$(\widetilde U, \{e\}, \varphi^{-1})$ is the conic chart, where $e$ is the unit in $SO(n)$.}
\end{example}

\begin{example} (Modular surface).  
\label{S2ExampleModular}
{\rm Let ${\bf H} = {\mathbb C}_+ = \{x + iy, ;\, x \in {\mathbb R}, \ y > 0\}$ be the upper-half plane. By taking a suitable discrete subgroup $\Gamma$ of $SL(2,{\mathbb R})$, we can define an action on ${\bf H}$ by
$$
\Gamma \times {\bf H} \ni \big(\gamma, z) \to \gamma\cdot z = \frac{az +b}{cz + d}, \quad 
\gamma = \left(
\begin{array}{cc}
a & b \\
c & d
\end {array}
\right),
$$
and obtain a 2-dimensional orbifold $\Gamma\backslash {\bf H}$ (see Definition \ref{Def:Orbifold}). A well-known example is the modular surface, which corresponds to $\Gamma = SL(2,{\bf Z})$. In this case, $\Gamma$ is generated by two elements $\gamma^{(T)}$ and $\gamma^{(I)}$, where
$$
\gamma^{(T)}\cdot z = z + 1, \quad \gamma^{(I)}\cdot z = - 1/z.
$$
Let $\mathcal M = SL(2, {\bf Z})\backslash{\bf H}$. Then, its fundametnal domain $\mathcal M^f$ is written as
$$
\mathcal M^f = \{z \in {\mathbb C}_+\, ; \, |z| \geq 1, \ |{\rm Re}\, z| \leq 1/2\}
$$
with boundary
$$
\partial{\mathcal M}^f = \partial{\mathcal M}_1^f \cup \partial{\mathcal M}_2^f,
$$ 
$$
\partial\mathcal M_1^f = L_- \cup L_+, \quad 
L_{\pm} = \Big\{\pm \frac{1}{2} + iy\, ; \, \frac{\sqrt{3}}{2} \leq y < \infty\Big\},
$$
$$
\partial{\mathcal M}_2^f = \Big\{e^{i\varphi}\, ; \, \frac{\pi}{3} \leq \varphi \leq \frac{2\pi}{3}\Big\}.
$$
To obtain $\mathcal M$ from $\mathcal M^f$, we glue $\partial M_1^f$ by the action of $\gamma^{(T)}$ : $z \to z + 1$, and $\partial M_2^f$ by the action of $\gamma^{I}$ : $e^{i\varphi} \to e ^{i(\pi - \varphi)}$. For $w \in {\bf H}$, let $\Gamma(w)$ be the isotropy group:
$$
\Gamma(w) = \{\gamma \in SL(2,{\bf Z})\, ; \, \gamma\cdot w = w\}.
$$
For $\gamma \in SL(2,{\bf Z})$, let $\langle\gamma\rangle$ be the cyclic group generated by $\gamma$. Then, for $w \in {\mathcal M}^f$, we have the following 4 cases:

\begin{itemize}
\item $\ \ \displaystyle{\Gamma(i) = \Big\langle
\left(
\begin{array}{cc}
0 & -1 \\
1 & 0
\end{array}
\right)\Big\rangle},$

\item $\ \ \displaystyle{\Gamma(e^{\pi i/3}) = \Big\langle
\left(
\begin{array}{cc}
0 & -1 \\
1 & -1
\end{array}
\right)\Big\rangle},$

\item $\ \ \displaystyle{\Gamma(e^{2\pi i/3}) = \Big\langle
\left(
\begin{array}{cc}
-1 & -1 \\
1 & 0
\end{array}
\right)\Big\rangle},$

\item $\ \ \displaystyle{\Gamma(w) = \big\langle
\pm I_2 \big\rangle}$, $\ \ w \neq i, e^{\pi i/3}, e^{2\pi i/3}$,
\end{itemize}

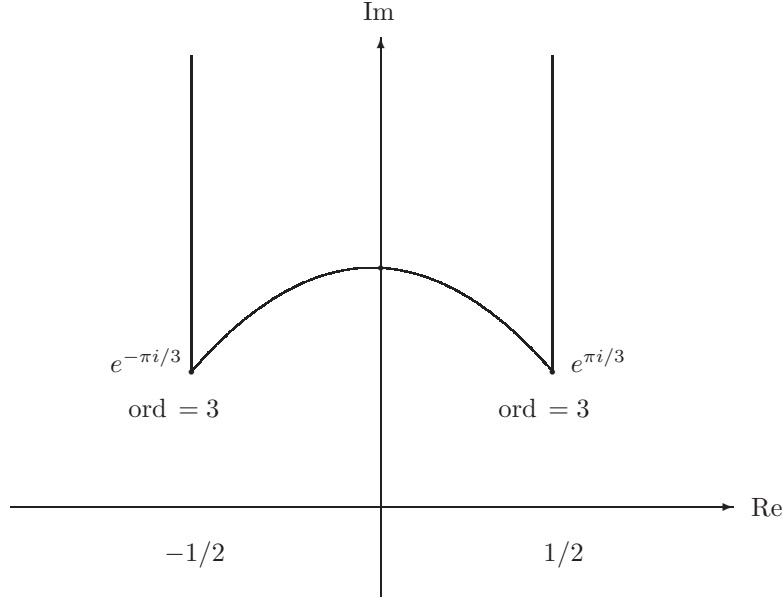
\begin{figure}[hbtp]
\setlength{\unitlength}{1.2mm}
\centering
\begin{picture}(0,70)
\put(-1,64){${\rm Im}$}
\put(1,0){\vector(0,1){62}}
\put(14,20){${\rm ord}\, = 3$}
\put(-27,20){${\rm ord}\, = 3$}
\put (20,25){\line(0,1){35}}
\put(19, 4){$1/2$}
\put (22, 25){$e^{\pi i/3}$}
\put(-23, 4){$-1/2$}
\put (-29, 25){$e^{-\pi i/3}$}
\put(-20,25){\line(0,1){35}}
\put(-40,10){\vector(1,0){80}}
\bezier{800}(-20,25)(0,48)(20,25)
\put(1,36.5){\circle*{0.6}}
\put(20,25){\circle*{0.6}}
\put(-20,25){\circle*{0.6}}
\put(42,9){${\rm Re}$}
\end{picture}

\label{fig:math2}
\caption{Fundamental domain for $PSL(2,{\bf Z})$}
\end{figure}

\noindent
where $I_2$ is the $2\times 2$ identity matrix. 
In the first case, $\Gamma(i)$ is order 2, in the 2nd and the 3rd cases, $\Gamma(e^{\pi i/3})$ and $\Gamma(e^{2\pi i/3})$ are order 3, and in the last case, the action $\gamma\cdot z$ is the identity. If $w \neq i, e^{\pi i/3}, e^{2\pi i/3}$, the manifold structure of $\mathcal M$ is easily defined by taking an open set $U = \widetilde U \subset {\mathbb C}_+$ such that $w \in U$, $i, e^{\pi i/3}, e^{2\pi i/3} \not\in U$ and $\varphi(z) = z$.


An element $\gamma \in SL(2,{\mathbb R})$ is said to be elliptic if it has only one fixed point in ${\mathbb C}_+$, which is equivalent to $|{\rm tr}\, \gamma| < 1/2$. Let $p$ be the fixed point of an elliptic $\gamma$. Then its isotropy group is cyclic. Let $n$ be the order of this cyclic group. Then, the generator $\gamma_0$ of this isotropy group satisfies
$$
\frac{w-p}{w - \overline{p}} = e^{2\pi i/n}\frac{z-p}{z - \overline{p}}, \quad 
w = \gamma_0\cdot z
$$
(see Lemma 2.4 of \cite{IKL11}). We can then take the local coordinates near $p$ by
$$
\zeta = \left(\frac{z - p}{z - \overline{p}}\right)^n, \quad \zeta(p) = 0.
$$
Then, we have
$$
z = \frac{p - \overline{p}\zeta^{1/n}}{1 - \zeta^n} = p + (p - \overline{p})\zeta^{1/n} + \cdots.
$$
(See Subsection 2.3 of  \cite{IKL11}). Letting $\zeta = \rho e^{i\theta}$,   the hyperbolic metric $\big((dx)^2 + (dy)^2\big)/y^2$ is rewritten as
$$
\frac{(dx)^2 + (dy)^2}{y^2} = (dr)^2 + \frac{1}{n^2}\sinh^2(r)(d\theta)^2.
$$
(See Subsection 2.4 of \cite{IKL11}). We have thus seen that around the elliptic fixed points, we can introduce a $C^{\infty}$-differentiable structure on $\mathcal M$ and a Riemannian metric except for fixed points. Note that:
\begin{itemize}
\item Although the covering space ${\mathbb C}_+$ has a global $C^{\infty}$-Riemannian metric $\dfrac{(dx)^2 + (dy)^2}{y^2}$, 
the induced metric on 
$\mathcal M$ is singular at elliptic fixed points.
 \end{itemize} 
}
\end{example}

\begin{example}
(3-dimensional orbifold). 
\label{Ex:3-dimorbifold}
{\rm An example of the 3-dimensional orbifold is given by the upper-half space model of the hyperbolic space ${\bf H}^3 = {\mathbb R}^3_+ = \{(x_1,x_2,y)\, y > 0\}$. Let us represent points in ${\bf H}^3$ by quarternions : $(x_1,x_2,y) \leftrightarrow x_1{\bf 1} + x_2{\bf i} + y{\bf k}$,
$$
{\bf 1} = \left(
\begin{array}{cc}
1 & 0 \\
0 & 1
\end{array}
\right), \quad 
{\bf i} = \left(
\begin{array}{cc}
i & 0 \\
0 & -i
\end{array}
\right), \quad 
{\bf j} = \left(
\begin{array}{cc}
0 & 1 \\
-1 & 0
\end{array}
\right), \quad 
{\bf k} = \left(
\begin{array}{cc}
0 & i \\
i & 0
\end{array}
\right).
$$
We identify $x_1{\bf 1} + x_2{\bf i}$ with $z = x_1 + ix_2 \in {\mathbb C}$. Then,
$$
x_1{\bf 1} + x_2{\bf i} + y{\bf j} = \left(
\begin{array}{cc}
z & y \\
- y & \overline{z}
\end{array}
\right) = : \zeta = z + yj.
$$
The action of $SL(2,{\mathbb C})$ is defined by
$$
SL(2,{\mathbb C})\times {\bf H}^3 \ni (\gamma,\zeta) \to \gamma\cdot\zeta := 
(a\zeta + b)(c\zeta + d)^{-1}, \quad 
\gamma = \left(
\begin{array}{cc} 
a & b \\
c & d
\end{array}
\right).
$$
Since the mapping $\gamma \to \gamma\cdot$ is 2 to 1, we consider 
$PSL(2,{\mathbb C}) = SL(2,{\mathbb C})/\{\pm I\}$. A counter part of the modular group is the Picard group
$$
\Gamma = PSL(2,{\bf Z}[i]) = \left\{
\left(
\begin{array}{cc}
a & b \\
c & d
\end{array}
\right)\, ; \, a, b, c, d \in {\bf Z}[i], \ ad - bc = 1
\right\},
$$
where ${\bf Z}[i] = {\bf Z} + i{\bf Z}$, the ring of Gaussian integers. As is well-known (see e.g. \cite{ElGrMen}, \cite{IKL13}), the fundamental domain of $\Gamma$ is 
$$
\mathcal M = \Gamma\backslash{\bf H}^3 = \Big\{z + yj\, ; \, |{\rm Re}\, z| \leq \frac{1}{2}, \ 0 \leq {\rm Im}\, z \leq \frac{1}{2}, \ |z|^2 + y^2 \geq 1\Big\}.
$$
The vertices of $\mathcal M$ are
$$
\infty, \quad - \frac{1}{2} + \frac{\sqrt{3}}{2}j, \quad \frac{1}{2} + \frac{\sqrt{3}}{2}j, \quad \frac{1}{2} + \frac{1}{2}i + \frac{\sqrt{2}}{2}j, \quad 
 -\frac{1}{2} + \frac{1}{2}i + \frac{\sqrt{2}}{2}j.
$$

\label{S2Example3dimorbifold}
\begin{figure}[hbtp]
\centering
\includegraphics[width=8cm]{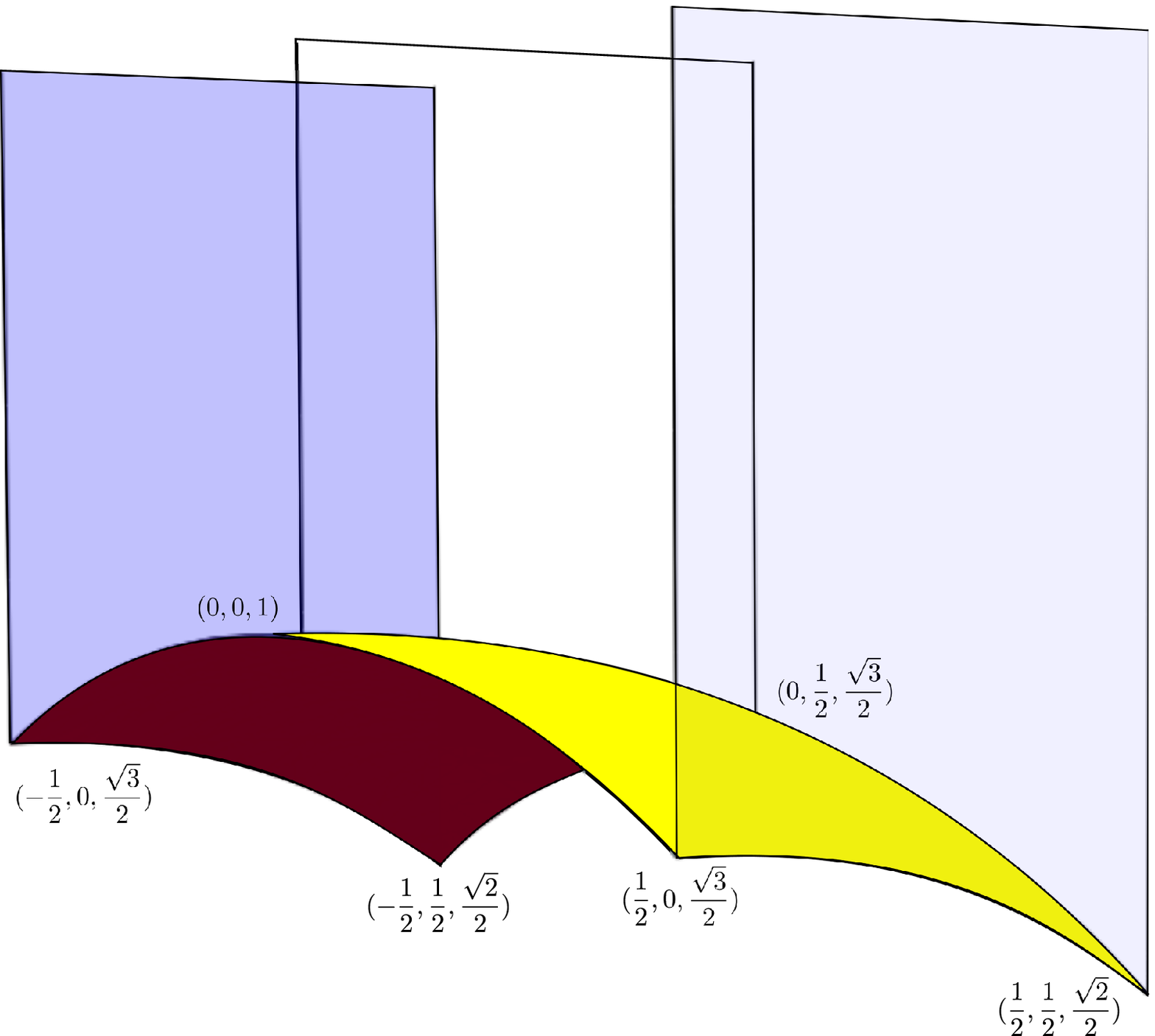}
\caption{Fundamental domain for $PSL(2,{\bf Z}[i])$}
\label{fig:math1}
\end{figure}

Letting $\mathcal M_{\pm} = \{z + yj \in \mathcal M\, ; \, 0 \leq \pm z \leq 1/2\}$, the boundary of the fundamental domain is split into 8 parts:
\begin{equation}
\begin{array}{lll}
S_1 = \mathcal M_- \cap \{{\rm Re}\, z =- 1/2\}, & \quad & S_2 = \mathcal M_+ \cap \{{\rm Re}\, z = 1/2\}, \\
S_3 = \mathcal M_- \cap \{{\rm Im}\, z = 0\}, & \quad & S_4 = \mathcal M_+ \cap \{{\rm Im}\, z = 0\},\\
S_5 = \mathcal M_- \cap \{{\rm Im}\, z = 1/2\}, & \quad & S_6 = \mathcal M_+ \cap \{{\rm Im}\, z = 1/2\},\\
S_7 = \mathcal M_- \cap \{|z|^2 + y^2 = 1\}, & \quad & S_8 = \mathcal M_+ \cap \{|z|^2 + y^2 = 1\}.
\end{array}
\nonumber
\end{equation}
Letting $L_{ij} = S_i\cap S_j$, the sets of sigular points are given by
\begin{equation}
\begin{array}{l}
 L_{13} = L_{24} = \{z = -1/2, y \geq \sqrt3/2\}, \ \ L_{34} = \{z = 0, y \geq 1\}\\
  L_{15} = L_{26} = \{z = -1/2 + i/2, y \geq 1/\sqrt2\}, \ \ L_{56} = \{z = i/2, y \geq \sqrt3/2\}\\
  L_{17} = L_{28} = \{|z|^2 + y^2 = 1, {\rm Re}\, z = -1/2, 0 \leq {\rm Im}\, z \leq 1/2\}, \\
  L_{78} =  \{|z|^2 + y^2 = 1,  {\rm Re}\, z = 0, 0 \leq {\rm Im}\, z \leq 1/2\}, \\
  L_{37} = L_{48} =  \{|z|^2 + y^2 = 1,  -1/2 \leq {\rm Re}\, z \leq 0,  {\rm Im}\, z =0\}, \\
  L_{57} = L_{68} =  \{|z|^2 + y^2 = 1,  -1/2 \leq {\rm Re}\, z \leq 0,  {\rm Im}\, z =1/2\}.
\end{array}
\nonumber
\end{equation}
We now let
\begin{equation}
\begin{array}{cccc}
\mathcal L_1 = L_{13}, & \mathcal L_2 = L_{15}, & \mathcal L_3 = L_{34}, &
\mathcal L_4 = L_{56}, \\
\mathcal L_5 = L_{17}, & \mathcal L_6 = L_{78}, & \mathcal L_7 = L_{37}, &
\mathcal L_8 = L_{57}.
 \end{array}
\nonumber
\end{equation}
The isotropy group for $\mathcal L_n$ are finite groups of rotations. In fact, 
$$
\Gamma(\mathcal L_i) = {\bf Z}_2, \quad i \neq 5, 8, \quad 
\Gamma(\mathcal L_i) = {\bf Z}_3, \quad i = 5, 8.
$$
The point is :
\begin{itemize}
\item The singular points form 1-dimensional curves. 
\item Some of the curves of singular points  go to infinity.
\item Some of the curves of singular points intersect.
\end{itemize}
 Since the isotropy groups for $\mathcal L_i (1 \leq i \leq 4)$ are rotation groups around them, for any $t > 1$, the horizontal slice $M_t = \mathcal M\cap\{y =t\}$ is a compact 2-dimensional orbifold with singular points $ - 1/2 + tj, -1/2 + i/2+ tj, yj, i/2 + tj$. Here, the covering is a disc in ${\mathbb R}^2$. For the details, see e.g. \cite{ElGrMen} and \cite{IKL13}.
}

\end{example}


\subsection{Orbifolds}
The notion of orbifold was first introduced by Satake \cite{Sa57}. 
We recall its definition for the sake of comparison. 

\begin{definition}
\label{Def:Orbifold}
An $n$-dimensional orbifold is a paracompact Hausdorff  space $X$ endowed with the set of local charts $\{U_i\, ; \, i \in I\}$ having the following properties:

\begin{itemize}
\item
 $X = \cup_{i\in I} U_i$. 
\item
 For each $U_i$, there exists a triple $(\widetilde U_i, \Gamma_i,\pi_i)$, where $\widetilde U_i \subset {\mathbb R}^n$ is an open set, $\Gamma_i$ is a finite group acting effectively (injectively) on $\widetilde U_i$, and $\pi_i : \widetilde U_i \to U_i$ is a continuous map inducing a homeomorphism $U_i \simeq \Gamma_i\backslash \widetilde U_i$.
 \item
  For any $U_i \subset U_j$, there is  a set $I_{ij} = \{(\phi,h)\}$, whose elements are called injections, where $\phi : \widetilde U_i \to \widetilde U_j$ is a smooth imbedding and $h : \Gamma_i \to \Gamma_j$ is an injective homeomorphism, such that $\phi$ is $h$-invariant and $\pi_i = \pi_j\circ \phi$. Moreover, $\Gamma_i\times \Gamma_j$ acts transitively on $I_{ij}$ by
$$
(g,g')\circ(\phi,h) = \big(g'\circ\phi\circ g^{-1}, {\rm Ad}\,(g')\circ h \circ {\rm Ad}\,(g^{-1})\big), 
$$
$$
 \forall g \in \Gamma_i, \ g' \in \Gamma_j, \ (\phi, h) \in I_{ij}.
$$
\item For each $p \in U_i\cap U_j$, there exists $U_k$ such that $p \in U_k \subset U_i\cap U_j$, and the injections are closed under composition for any $U_k \subset U_i \subset U_j$.
\end{itemize}
\end{definition}

Therefore, $(\widetilde U_i, \Gamma_i, \pi_i)$ is the conic chart in the sense of Definition \ref{Definitionofconicchart}. 


The examples \ref{S2ExampleModular} and \ref{Ex:3-dimorbifold} are 2  and 3 dimensional orbifolds arising in number theory.
For more detailed exposition of orbifolds, see e.g. 
\cite{Davis} and \cite{Thurs}. 

Let us repeat important points. Let $M$  be an orbifold with $C^\infty$-smooth coordinates.  Then, every point $x_0\in M$  has a neighborhood $U$ where $U$  is either
an open smooth manifold with smooth metric, or, $U$ has the structure
$U\equiv \Gamma\backslash \tilde U$,
where $(\widetilde U, \Gamma, \pi)$ is a conic chart of $U$. Theis means that (1)
$\widetilde U \subset {\mathbb R}^n$ is a local covering neighborhood having the product structure
$\widetilde U=\widetilde W\times \widetilde V$,
where $\widetilde W\subset B^{k}(0,R_1)\subset {\mathbb R}^{k} $ and 
$\widetilde V=B^{n-k}(0,R_0)\subset {\mathbb R}^{n-k}$,  and (2)  $\Gamma$  is a finite group acting on $\widetilde U$. Therefore, in $U$, there are local  coordinates 
$$
X :U\to \Gamma\backslash {\mathbb R}^n,\quad X (U)=\Gamma\backslash \widetilde U.
$$


\subsection{Conic manifolds with group action}
\label{SubsectionConicmfdgroupaction}
We now define the conic manifold which admits  group action,
{\mtext which is an extension of  orbifolds.}

\begin{definition}
\label{Definitionconicchart}
A {\it conic manifold admitting group action}, which is abbreviated to CMGA\index{CMGA}, is a topological space $M$ equipped with the following structure: 
{\it There exists a family of open covering $\{U_j\, ; \, j \in J\}$ of $M$ having the following properties (C-1) $\sim$ (C-4)}:

\smallskip
\noindent 
{\bf (C-1)}  {\it For any $j \in J$, $U_j$ has a conic chart $(\widetilde U_j, \Gamma_j, \pi_j)$, 
 $\widetilde U_j=\widetilde W_j\times \widetilde V_j$, where
  for some 
$0 \leq k \leq n,  \, k \neq n-1$,}
\begin{equation}
 \widetilde W_j = B^{k}(0,R_0)\subset {\mathbb R}^{k}, \quad  
\widetilde V_j=B^{n-k}(0,R_1)\subset {\mathbb R}^{n-k}. 
\label{S2ConstR0R1}
\end{equation}

\smallskip
\noindent
{\bf (C-2)}
{\it Define $Y_j : \widetilde U_j \to \widetilde W_j$, $Z_j : \widetilde U_j \to \widetilde V_j$ and $\widetilde U_j^{reg}$ by}
\begin{equation}
\widetilde U_j \ni x  \to x = (y,z) = (Y_j(x),Z_j(x)),
\nonumber
\end{equation}
\begin{equation}
\widetilde U_j^{reg}=\Big(\widetilde W_j\times \big(B_j^{n-k}(0,R_1)\setminus \{0\}\big)\Big)\cap\big\{x \in \widetilde U_j\, ; \, \Gamma_j(x) = \{e\}\big\}.
\label{DefinewidetildeUjreg}
\end{equation}
Then :

\smallskip
\noindent
(C-2-1)
 {\it The action of  $\gamma\in \Gamma_j$ keeps 
the $y$-coordinates invariant, i.e. }
\begin{equation}
Y_j(\gamma\cdot x)=Y_j(x), \quad \forall x\in \widetilde U_j, \quad \forall\gamma \in \Gamma_j.
\nonumber
\end{equation}

\smallskip
\noindent
(C-2-2) 
{\it There exists a $\Gamma_j$-invariant $C^\infty$-metric  $\widetilde g_j$ on $\widetilde U_j^{reg}$, i.e.}
\begin{equation}
{\gamma}_{\ast}\,\widetilde g_j=\widetilde g_j \quad {on} \quad \widetilde U_j^{reg}, \quad \forall \gamma \in \Gamma_j.
\nonumber
\end{equation}

\smallskip
\noindent
(C-2-3)
{\it In the spherical coordinates
 $Z_j(x)= s\omega = z$, $s = s(x) = |z|$, $\omega=\omega(x) = \dfrac{z}{|z|}$ such that
$ s \in (0,R_1)$ and $\omega\in S^{n-k-1}$ on $B^{n-k}(0,R_0)\setminus \{0\}$, $\widetilde g_j$
 has the form}
\begin{equation}
\begin{split}
 \widetilde g_j= & \sum_{p, q=1}^k a^{(j)}_{pq}(y,s,\omega)dy^pdy^q \\
&+ ds^2 +
s^2\sum_{\ell, m=1}^{n-k}  \,b^{(j)}_{\ell  m}(y,s,\omega)d\omega^{\ell} d\omega^{m}  + s\sum_{p=1}^k 
\sum_{\ell=1}^{n-k} \,h^{(j)}_{p\ell}(y,s,\omega)dy^p d\omega^{\ell}.
\end{split}
\label{Sec2Formofgj}
\end{equation}
The coefficients satisfy 
\begin{equation}
\left\{
\begin{split}
& a^{(j)}_{pq}(y,s,\omega)\to \widehat a^{(j)}_{pq}(y),\\ 
& b^{(j)}_{\ell m}(y,s,\omega)\to \widehat b^{(j)}_{\ell m}(y,\omega),\\ 
& h^{(j)}_{p\ell}(y,s,\omega)\to 0,
\end{split}
\right.
\label{Sec2abhcondition}
\end{equation}
{\it uniformly in $(y,\omega)$ as $ s \to 0$. Moreover, there exist  constants $C_1 \geq C_0> 0$ and a positive continuous function $T_j(y)$ such that}
\begin{equation}
 C_0\,  g_{S^{n-k-1}} \leq \sum_{\alpha, \beta = 1}^{n-k}\widehat b^{(j)}_{\ell m}(y,\omega)d\omega^{\ell} d\omega^{m} \leq   T_j(y)^2 g_{S^{n-k-1}}, 
\label{Sec2gcondition}
\end{equation}
{\it where  $ g_{S^{n-k-1}}$ is the standard metric of  ${S^{n-k-1}}$ and
$C_0 \leq T_j(y) \leq C_1$.}

\smallskip
\noindent
(C-2-4) For $e \neq \gamma \in \Gamma_j$, $e$ being the unit of $\Gamma_j$, 
\begin{equation}
{\rm cap}_2(\{x \in \widetilde U\, ; \, \gamma\cdot x = x\}) =0,
\label{Def2.1;diminvariantsetleqn-2}
\end{equation}
where ${\rm cap}_2(E)$ denotes the 2-capacity of a subset $E \subset {\mathbb R}^n$.

\smallskip 
\noindent
{\bf (C-3)} {\it If $U_j \cap U_k \neq \emptyset$, there exists $\ell \in J$ such that $U_{\ell} \subset U_j\cap U_k$. }

\smallskip
\noindent
{\bf (C-4)}  {\it If $U_{\ell} \subset U_k$, there exist an injective homomorphism $\mathcal I_{k\ell} : \Gamma_{\ell} \to \Gamma_k$, and a $C^{\infty}$ injective map $\widetilde I_{k\ell} : \widetilde U_{\ell} \to \widetilde U_k$.}
\end{definition}
 
 Recall that the 2-capacity of a subset $E \subset {\mathbb R}^n$ is defined by
 \begin{equation}
 {\rm cap}_2(E) = \inf\int_{{\mathbb R}^n}
\left(|u|^2 + |\nabla u|^2\right)dx,
 \nonumber
 \end{equation}
 where the infimum is taken over all $u \in H^2({\mathbb R}^n)$ such that $u \geq 1$ almost everywhere on a neighborhood of $E$.
 
To simplify the analysis below, we assume that the constants $R_0,R_1>0$ in (\ref{S2ConstR0R1}) are independent of $j$ and so are $C_0, C_1$ in (\ref{Sec2gcondition}).

We have another expression of (\ref{Sec2Formofgj}), which does not contain local coordinates on $S^{n-k-1}$: Omitting the subscript $j$, 
\begin{equation}
\label{radial coordinates 2}
\widetilde g=\sum_{p,q=1}^k \widetilde a_{pq}(y,z)dy^pdy^q +  \sum_{\ell,m=1}^{n-k} \widetilde b_{\ell m}(y,z)\,dz^{\ell}dz^{m} + 
\sum_{p=1}^k 
\sum_{\ell=1}^{n-k} \widetilde h_{p\ell}(y,z)dy^p dz^{\ell},
\end{equation}
where $z=(z^p)_{p=1}^{n-k}\in {\mathbb R}^{n-k}$,
$|z|=(\sum_{p=1}^{n-k}|z^p|^2)^{1/2}$, and
\begin{equation}
\left\{
\begin{split}
& \widetilde a_{pq}(y,z)=a_{pq}(y,|z|,\frac{z}{|z|}), \\
& \widetilde b_{\ell m}(y,z)= \frac{z^{\ell}z^{m}}{|z|^2} + \sum_{u, v=1}^{n-k} 
\big(\delta_{u \ell} - \frac{z^{u}z^{\ell}}{|z|^2}\big)b_{u v}(y,|z|,\frac {z}{|z|})
\big (\delta_{v m} - \frac{z^v z^{m}}{|z|^2}\big),\\
& \widetilde h_{p \ell}(y,z)=\sum_{u=1}^{n-k} h_{p u}(y,|z|,\frac{z}{|z|}) (\delta_{u\ell}-\frac{z^{u}z^{\ell}}{|z|^2}).
\end{split}
\right.
\nonumber
\end{equation}
By the conditions (\ref{Sec2abhcondition}), (\ref{Sec2gcondition}), there exist constants $R_0, C >0$ such that 
\begin{equation}
C|\xi|^2 \leq \widetilde g_x(\xi,\xi) \leq C^{-1}|\xi|^2, \quad |z| \leq R_0,
\label{S2metricgpositivecond}
\end{equation}
for any $\xi$, where $\widetilde g_x(\xi,\xi)$ is the metric $\widetilde g$  in (\ref{radial coordinates 2}), 
and we identify $\xi \in {\mathbb R}^n$ with the vector field
$\displaystyle{
\sum_{i=1}^k\xi^i\frac{\partial}{\partial y^i} + 
\sum_{i=k+1}^n\xi^i\frac{\partial}{\partial z^i}.
}$

Let us check the validity of (\ref{Def2.1;diminvariantsetleqn-2}) for  the case of linear action by $SO(n)$.

\begin{lemma}
\label{Lemma:dimLgammaleqn-2}
Let a finite group  $\Gamma \subset SO(n)$ act linearly on ${\mathbb R}^n$, and  for $e \neq \gamma \in \Gamma$, put $L_{\gamma} = \{x \in {\mathbb R}^n\, ; \, \gamma\cdot x = x\}$. Then, ${\rm dim}\, (L_{\gamma}) \leq n-2$.
\end{lemma}

\begin{proof}
Assume that ${\rm dim}\, (L_{\gamma}) = n-1$. Then, $1$ is an eigenvalue with multilicity at least $n-1$. Since $\gamma \in SO(n)$, the remaining eigenvalue must be $1$.  This is a contradiction. 
\end{proof}

\begin{lemma}
\label{LemmaSO(n)capapcity0}
If $\Gamma_j$ is a finite subgroup of $SO(n)$ acting linearly on $\widetilde U_j$, the condition (\ref{Def2.1;diminvariantsetleqn-2}) is satisfied.
\end{lemma}

\begin{proof}
Lemma \ref{Lemma:dimLgammaleqn-2} implies that $(n-2)$-dimensional Hausdorff measure of $\{x \in \widetilde U\, ; \, \gamma\cdot x = x\}$ is finite. Then, its 2-capacity is 0. See \cite{Heino07}, p. 18. 
\end{proof}

For Hausdorff measures and dimensions, see \cite{Mattila}.


\subsection{Orbifold and CMGA}
The examples \ref{S2ExampleModular} and \ref{Ex:3-dimorbifold} are both orbifold and CMGA. 
The difference between  CMGA and orbifold is the singularity of the metric. For  the case of orbifold, instead of  (\ref{DefinewidetildeUjreg}),
 the metric $\widetilde g$ on the covering space is assumed to be smooth in a neighborhood of $z=0$. 
 

\begin{example}
\label{Example:CMGA1}
A simple example of CMGA is constructed from a sector 
$$
U = \big\{se^{i\psi} \, ; \, s \geq  0, \ \psi \in [0,2\pi\kappa]\big\}, \quad 0 < \kappa < 1.
$$
We regard $U$ as a cone by identifying $\psi =\psi_0$ with $\psi = 0$.
Letting $\theta = \psi/\kappa$, and $\omega = 
(\cos\theta,\sin\theta)$, which varies over whole $S^1$, 
we equip $U$, which is homeomorphic to ${\mathbb R}^2$, with the standard polar coordinates $(s,\theta) \in [0,\infty)\times S^1$. Then, $x = s\cos\theta$, $y = s\sin\theta$ in the rectangular coordinates. We define a metric $g$ on $U$ by
\begin{equation}
g = (ds)^2 + \kappa^2s^2(d\theta)^2.
\label{ExampleofgonM0}
\end{equation}
 Then, $U$ is a $C^{\infty}$ manifold with the standard differential structrure of ${\mathbb R}^2$, but its metric is singular at the origin. 
In fact,  $g$ is written as
\begin{equation}
\begin{split}
g &= \frac{x^2 + \kappa^2y^2}{x^2 + y^2}(dx)^2 - 
\frac{2(\kappa^2-1)xy}{x^2 + y^2}dxdy + \frac{\kappa^2x^2 + y^2}{x^2 + y^2}(dy)^2 \\
& = (dx)^2 + (dy)^2 + (\kappa^2-1)(d\theta)^2.
\end{split}
\nonumber
\end{equation}

If $1/\kappa$ is a natural number, then $U$ is an orbifold. In fact,  put  $\widetilde U = {\mathbb R}^2$, $\Gamma$ = the group generated by the rotation of angle $2\pi/n$, where $n = 1/\kappa$, and $\pi : se^{i\psi} \to 
se^{in\psi}$. Then, $(\widetilde U, \Gamma, \pi)$ is a conic chart of $U$. Endow $\widetilde U$ with the metric $\widetilde g = (dx)^2 + (dy)^2$. Then, $g$ is induced from $\widetilde g$.

If $\kappa \not\in {\mathbb Q}$, $U$ is not an orbifold, since the group generated by a rotation of angle $2\pi\kappa$ is an infinite group. However, we can regard $U$ as a CMGA, where we take $\widetilde U = {\mathbb R}^2$ equipped with the metric (\ref{ExampleofgonM0}), $\Gamma = \{e\}$ (the unit group) and $\pi$ to be identity. 
\end{example}

\begin{example}
We consider a higher dimensional analogue of Example \ref{Example:CMGA1}. Let
\begin{equation}
\widetilde U = \widetilde W \times \widetilde V = {\mathbb R}^{k}\times {\mathbb R}^{n-k},
\nonumber
\end{equation}
and $\Gamma$  a finite subgroup of $SO(n)$ acting on $\widetilde U$. As in Definition \ref{Definitionconicchart}, assume that this action leaves $y$-coordinates invariant, i.e. 
\begin{equation}
Y(\gamma\cdot x) = Y(x), \quad \forall x \in \widetilde U, \quad \forall \gamma \in \Gamma.
\nonumber
\end{equation}
We then have
\begin{equation}
\Gamma\backslash \widetilde U = \widetilde W \times 
\big(\Gamma\backslash \widetilde V\big).
\nonumber
\end{equation}
We define a $\Gamma$-invariant metric $\widetilde g$ by
\begin{equation}
\widetilde g = (dy)^2 + (ds)^2 + a(y)s^2g_{S^{n-k-1}},
\nonumber
\end{equation}
where $a(y)$ is any positive $C^{\infty}$-function on $\widetilde W$. Then, $\widetilde g$ is a smooth metric on
\begin{equation}
U \setminus \big({\mathbb R}^k\times \{0\}\big) = 
{\mathbb R}^k \times \Big((0,\infty)\times \big(\Gamma\backslash S^{n-k-1}\big)\Big),
\nonumber
\end{equation}
where we deal with $\Gamma\backslash S^{n-k-1}$ as an orbifold. This space has the structure described in the left-hand side of Figure \ref{fig:math4}. As is seen in its right-hand side, this CMGA has  screen type singularities.
\end{example}

\begin{figure}[hbtp]
\setlength{\unitlength}{1.2mm}
\centering
\begin{picture}(5,65)
\put(-16,65){${\mathbb R}^k$}
\put(-15,0){\vector(0,1){62}}
\put(-15,55){\line(2,-1){20}}
\put(-15,55){\line(-2,-1){20}}
\bezier{800}(-35,45)(-15,36)(5,45)
\put(-42,40){$\Gamma\backslash S^{n-k-1}$}
\put(-15,32){\line(5,-2){20}}
\put(-15,32){\line(-5,-2){20}}
\bezier{800}(-35,24)(-15,15)(5,24)
\put(-15,8){\line(6,-1){20}}
\put(-15,8){\line(-6,-1){20}}
\bezier{800}(-35,5)(-15,0)(5,5)
\put(-15,55){\line(1,-1){13}}
\put(-15,55){\line(-1,-1){13}}
\put(-15,32){\line(1,-1){11}}
\put(-15,32){\line(-1,-1){11}}
\put(-15,8){\line(2,-1){10}}
\put(-15,8){\line(-2,-1){10}}

\put(29,65){${\mathbb R}^k$}
\put(30,0){\vector(0,1){62}}
\put(30,55){\line(1,-1){13}}
\put(30,8){\line(1,-1){13}}
\put(43,-5){\line(0,1){47}}
\put(30,55){\line(-1,-1){13}}
\put(30,8){\line(-1,-1){13}}
\put(17,-5){\line(0,1){47}}
\put(15,48){$U^{sing}$}
\end{picture}

\label{fig:math4}
\caption{CMGA}
\end{figure}

\subsection{Global distance on $M$}

Assume that a  CMGA, denoted by $M$,  is covered by locally finite coordinate neighborhoods $U_j$, $j \in J = \{1,2,\cdots\}$, (allowing the case $\sharp J < \infty$). We define
\begin{equation}
\begin{split}
& M^{reg}=\bigcup_{j=1}^\infty \pi_j\Big(
  \big(\widetilde W_j\times (B^{n-k}(0,R_1)\setminus \{0\})\big)\cap \big\{x \in \widetilde U_j \, ; \, 
  \Gamma_j(x) = \{e\}\big\}\Big), \\
& M^{sing}=M\setminus M^{reg}.
\end{split}
\nonumber
\end{equation}
 We also denote 
\begin{equation}
 M_j^{reg} = \pi_j\Big(
  \big(\widetilde W_j\times (B^{n-k}(0,R_0)\setminus \{0\})\big)\cap \big\{x \in \widetilde U_j\, ; \, \Gamma_j(x) = \{e\}\big\}\Big),
\label{DefineMjreg}
\end{equation}
 $$
 U^{reg}_j=
\pi_j( \widetilde U^{reg}_j).
$$

The metric tensors $g_j$ in (\ref{Sec2Formofgj}) define a smooth metric on the regular part  $M^{reg}$ of $M$ and we denote this metric by $g$.
These coordinates also determine a topology on $M$.
As $M$ is the topological closure of  $M^{reg}$, we define 
\begin{equation}
\label{extension of metric}
d_M(x,y)=\inf_{\gamma} \hbox{Length}_g(\gamma([0,1])\cap M^{reg}),
\end{equation}
where the infimum is taken over paths $\gamma:[0,1]\to M$  that
are  piecewise $C^1$-smooth on the lifted local coordinates,
$\gamma(0)=x$, $\gamma(1)=y$, and $\gamma([0,1)])\cap  M^{sing}$ is a finite set. Note that
$\gamma([0,1])\cap M^{reg}_j$ is rectifiable in all local coordinate charts $U^{reg}_j$. Also, 
when $\gamma((s_1,s_2))\subset U_j^{reg}\subset M^{reg}$,
 the above length  in (\ref{extension of metric}) is written using local coordinates as
$$
\hbox{Length}_g(\gamma((s_1,s_2)))=\int_{s_1}^{s_2}\big(g_{jk}(\gamma(s))\dot\gamma^j(s)\dot\gamma^k(s)\big)^{1/2}ds.
$$


\subsection{Integration on CMGA}
Integration over $M$ is actually done on its covering as explained below.
Let $\{U_j\}_{j \in J}$ be a covering of $M$, and $(\widetilde U_j, \Gamma_j, \pi_j)$ the conic chart of $U_j$. 
We can constrcut a partition of unity $\{\chi_j\}_{j\in J}$ on $M$ and  $\widetilde \chi_j \in C^{\infty}(\widetilde U_j)$ 
satisfying
${\rm supp}\, \chi_j \subset U_j$ and $\chi_j(\pi_j(x)) = \widetilde \chi_j(x)$. Then, assuming that 
integration over $M$ is defined, 
for any $v, w \in L^2(M)$, 
\begin{equation}
(v,w)_{L^2(M)} = \big(\sum_{j=1}^{\infty}\chi_jv,w)_{L^2(M)} = \sum_{j=1}^{\infty}(\chi_jv,w)_{L^2(M)}, 
\nonumber
\end{equation}
since $\sum_{j=1}^{\infty}\chi_j = 1$. This reduces the computation on each patch $U_j$. 
We omit the subscript $j$ and assume that  the supports of $v, w$ are contained in $U$.
As was noticed before, the functions on $U$ are ragarded as $\Gamma$-invariant functions on $\widetilde U$.  
Therefore, $L^2(U)$ should be regarded as the closed subspace of $L^2(\widetilde U)$ consisting of $\Gamma$-invariant functions. With in mind,  we {\it define} integration over $U$ by the one over $\widetilde U$, and put
\begin{equation}
(v,w)_{L^2(U)} = \frac{1}{\sharp\Gamma} (\pi^{\ast}v,\pi^{\ast}w)_{L^2(\widetilde U)}.
\nonumber
\end{equation}


\subsection{Sobolev spaces}
\label{Subsection2.6Sobolevspace}
We fix a coordinate patch $U$ in $M$.
Letting $\widetilde g = \widetilde g_{ab}dx^adx^b$ be the $\Gamma$-invariant Riemannian metric on $\widetilde U^{reg}$, we define a quadratic form $\widetilde q$ with form domain $D(\widetilde q) = C_0^{\infty}(\widetilde U)\cap \mathcal S_{\widetilde U}$ by
\begin{equation}
\begin{split}
\widetilde q(\widetilde v,\widetilde w)  = \frac{1}{\sharp\Gamma}(\widetilde v,\widetilde w)_{L^2(\widetilde U)} +  \frac{1}{\sharp\Gamma}\big(\widetilde g^{ab}
\partial_a\widetilde v,\partial_b\widetilde w\big)_{L^2(\widetilde U)},
\end{split}
\label{DefineqonU}
\end{equation}
where $\partial_a = \partial/\partial x^a$,  $\widetilde v, \widetilde w \in C_0^{\infty}(\widetilde U)\cap \mathcal S_{\widetilde U}$. 
Letting $v(\pi(x)) = \widetilde v(x)$, $w(\pi(x)) = \widetilde w(x)$, we put
\begin{equation}
q(v,w) = \widetilde q(\widetilde v,\widetilde w),
\nonumber
\end{equation}
and call it a quadratic form on $C_0^{\infty}(U)$. Note the abuse of notation, since $C_0^{\infty}(U)$ is not defined. 
For convenience, we write the right-hand side of (\ref{DefineqonU})  as
\begin{equation}
\int_Uv\overline{w}\sqrt{g}dx +  \int_{U} g^{ab}(x)
\big(\partial_a v(x)\big)\big(\overline{\partial_bw(x)}\big)\sqrt{g(x)}\, dx.
\nonumber
\end{equation}
Since the metric $\widetilde g$ is $\Gamma$-invariant, so is $\widetilde g^{ab}(\partial_a\widetilde v)(\partial_b\widetilde v)$ for $\Gamma$-invariant $\widetilde v$. Taking account of (\ref{DefineqonU}), the usual calculus for integral, e.g. integration by parts,  can be applied also to this formal expression.
Returning to the total manifold $M$, the form $q$ is rewritten formally as
\begin{equation}
\begin{split}
q(v,w)  &= \sum_j\int_{U_j}\chi_j v\overline{w}\sqrt{g_j}dx \\
 &+ \sum_{j} \int_{U_j}g_j^{ab}
\big(\partial_a(\chi_jv)\big)\overline{\partial_b w}\sqrt{g_j}\, dx
\end{split}
\nonumber
\end{equation}
for $v, w$ whose lifts $\widetilde v, \widetilde w$ are smooth and compactly supported. 
Actually, it is represented as
\begin{equation}
\begin{split}
\widetilde q(\widetilde v,\widetilde w) &= \sum_j\frac{1}{\sharp\Gamma_j}\int_{\widetilde U_j}\widetilde \chi_j(\pi_j^{\ast}v)(\overline{\pi_j^{\ast}w})\sqrt{\widetilde g_j}dx \\
 &+ \sum_{j}\frac{1}{\sharp\Gamma_j} \int_{\widetilde U_j} \widetilde g_j^{ab}
\big(\partial_a(\widetilde\chi_j\pi_j^{\ast}v)\big)\big(\overline{\partial_b \pi_j^{\ast}w}\big)\sqrt{\widetilde g_j}\, dx. 
\end{split}
\label{qdefinedonwidetildU}
\end{equation}
Since the  computation is reduced to that on $\widetilde U_j$, the usual theory of quadratic form is 
applied to $q$. Then, due to (\ref{Sec2abhcondition}), $q$ is closable. 
The Sobolev space $H^1(M)$ is defined to be the domain of $\overline{q}$ = the closure of $q$:
\begin{equation}
H^1(M) = D\big(\overline{q}\big).
\label{H1M=D(q)}
\end{equation}
This is more clearly stated in the following lemma.

\begin{lemma}
\label{LemmaH1M=H1Mreg}
$ \ H^1(M) = H^{1}(M^{reg})$.
\end{lemma}

\begin{proof}
Recall that $v \in H^1(M^{reg})$ means that $\pi^{\ast}v \in H^1(\widetilde U_j\setminus \widetilde U_j^{sing})$ for all $j$. Therefore $H^1(M) \subset H^1(M^{reg})$ is obvious. To show the converse inclusion, we use the following fact.

\begin{lemma}
Let $\Omega$ be a domain in ${\mathbb R}^n$, and assume that there exists a subset $S \subset \Omega$ whose 2-capacity is equal to 0. Then, there exists a bounded operator ${\rm Ext}\, : H^1(\Omega\setminus S) \to 
H^1(\Omega)$ such that ${\rm Ext}\, v = v$ on $\Omega\setminus S$ for any $v \in H^1(\Omega\setminus S)$.
\end{lemma}

See \cite{KilKinMar00} Theorem 4.6 (and its proof). 
Therefore, if we show that the 2-capacity of $\widetilde U_j^{sing}$ is 0, the inclusion $H^1(M) \supset H^1(M^{reg})$ follows. 
Now, omitting the subscript $j$, $\widetilde U^{sing}$ consists of two parts. 
$$
\widetilde U^{sing,1} = \widetilde W \times \{0\}, \quad \widetilde U^{sing,2} = \{x \in \widetilde U\, ; \, \Gamma(x) \neq \{e\}\}.
$$
If $k = n$, $\widetilde U$ has no singularities. If $ k \neq n$, then $k \leq n-2$ by the assumption (C-1). Therefore $\widetilde U^{sing,1}$ has a finite $(n-2)$-dimensional Hausdorff measure, hence its 2-capasity is 0, as was seen in the proof of Lemma \ref{LemmaSO(n)capapcity0}.
Noting that
\begin{equation}
\widetilde U^{sing,2} \subset {\mathop\cup_{e\neq \gamma\in\Gamma}}\{x \, ; \, \gamma\cdot x = x\},
\nonumber
\end{equation}
and using  (\ref{Def2.1;diminvariantsetleqn-2}), we see that the 2-capcity of
$\widetilde U^{sing,2}$ is 0.
\end{proof}

 Returning to the quadratic form (\ref{qdefinedonwidetildU}), let $L'$ be the associated Friedrichs extension, i.e.  the unique self-adjoint operator $L'$ satisfying $D(\sqrt{L'}) = H^1(M)$ and
\begin{equation}
\overline{q}(u,v) = (L'u,v), \quad \forall u \in D(L'), \quad \forall  v \in D(\overline{q}).
\nonumber
\end{equation}
We put $L = L'-1 = - \Delta_M$. Then, we have
\begin{equation}
H^1(M) = D(\sqrt{1 - \Delta_M}).
\nonumber
\end{equation}
We define the Sobolev space $H^s_{g}(M)$ by $H^s_g(M) = D(L')^{s/2},  s \in {\mathbb R}$. In other words, 
\begin{equation}
H^s_g(M) = D((1 - \Delta_M)^{s/2}), \quad s \in {\mathbb R}.
\nonumber
\end{equation}
We write $H^s_g(M)$ as $H^s(M)$ for the sake of simplicity. 
We need a uniformity for the bounds of the metric on each $\widetilde U_j$ to define $L'$.  However, we do not pursue this condition here, since  in our later applications the number of charts are finite.
 The projection $P$ defined 
on $L^2(\widetilde U)$
is naturally extended to $H^m(\widetilde U)$, $m = 0, 1, 2,\cdots$. 

In the case of orbifolds, the metric $\widetilde g$ on the covering space is assumed to be smooth even on $\{z = 0\}$. Therefore, there is no problem for the regularity, and we have
\begin{equation}
{\mathop\cap_{m=0}^{\infty}}H^m(M) \subset C^{\infty}(M).
\nonumber
\end{equation}
However, for the case of  CMGA, although
 \begin{equation}
{\mathop\cap_{m=0}^{\infty}}H^m(M) \subset C^{\infty}(M^{reg})
\nonumber
\end{equation}
holds, $u \in {\mathop\cap_{m=0}^{\infty}}H^m(M) $ may not be regular on $M^{sing}$, since $\widetilde g$ is not assumed to be smooth around $\{z =0\}$.


\section{Laplacian on  conic manifold}
\label{Laplacian on  conic manifold}


\subsection{Regularity of the domain of Laplacian}
\label{SubsectionRegukarityofdomain}
Let $L = - \Delta$ be the Laplacian on $M$ defined in the previous section. We study the regularity of the elements in $D(L)$ near $M^{sing}$.
Recall that $M$ is covered by $\{U_j\, ; \, j \in J\}$. 
By taking local coordinates $\theta \in {\mathbb R}^{n-k-1}$ on $S^{n-k-1}$, we introduce the following norm for functions defined on $\widetilde U_j$:
\begin{equation}
p_j^{(m)}(f) = \sup_{|y|<R_0,0 < s< R_1,\theta}\sum_{|\alpha| + \beta + |\gamma|\leq m}
s^{-|\gamma|}\big|\partial_y^{\alpha}\partial_s^{\beta}\partial_{\theta}^{\gamma}f(y,s,\theta)\big|.
\label{DefinepjmS2}
\end{equation}

We make the condition (\ref{Sec2abhcondition}) more precise. 

\medskip
\noindent
{\it {\bf (D)} 
Letting $c(y,s,\omega)$ be any of $a^{(j)}_{pq}(y,s,\omega)$, $b^{(j)}_{\ell m}(y,s,\omega)$ and $h^{(j)}_{p\ell}(y,s,\omega)$, we assume that 
\begin{equation}
\sup_{|y|<R_0, 0 < s < R_1,\theta}p_j^{(2)}(c) < \infty,
\nonumber
\end{equation}
and $c(y,0,\omega) = \lim_{s\to 0} c(y,s,\omega)$ satisfies 
$c(y,0,\omega) = c(y)$ if $c = a^{(j)}_{pq}$, $c(y,0,\omega) = 0$ if $c = h^{(j)}_{p\ell}$. Moreover, (\ref{Sec2gcondition}) is satisfied.}

\medskip
The collection of norms $\{p_j^{(m)}\}_{j\in J}$, together with a partition of unity $\{\chi_j\}_{j\in J}$ on $M$, endows a norm on $M$, which we denote by $|\cdot|_m$:
\begin{equation}
|f|_m = \sup_{j \in J} p_j^{(m)}(\chi_jf).
\label{||mnormonM}
\end{equation}

\medskip
Take a conic chart $(\widetilde U_j, \Gamma_j, \pi_j)$. Omittng the subscript $j$, we rewrite the metric (\ref{Sec2Formofgj}) in terms of local coordinates $\theta$ on $S^{n-k-1}$  as follows:
\begin{equation}
\Big(\widetilde g_{ij}\Big) = \left(
\begin{array}{ccc}
a & 0 & s\, ^t h\\
0 & 1 & 0 \\
sh & 0 & s^2b
\end{array}
\right), 
\label{S1gijandT}
\end{equation}
where $a = \big(a_{pq}\big), b  = \big(b_{\alpha\beta}\big), h = \big(h_{p\alpha}\big)$, with $ 1 \leq p, q \leq k$, $1 \leq \alpha, \beta \leq n-k-1$. Note that $a$ is the same as in (\ref{Sec2Formofgj}), and $b$ and $h$ are different, since we have used $d\omega = \sum_{\ell=1}^{n-k-1}\frac{\partial \omega}{\partial \theta^{\ell}}d\theta^{\ell}$. However,  $a, b, h$ still have the properties in (D).

We then have
\begin{equation}
T\Big(\widetilde g_{ij}\Big)T = 
\left(
\begin{array}{ccc}
a & 0 & ^t h\\
0 & 1 & 0 \\
h & 0 & b
\end{array}
\right), \quad 
 T = 
\left(
\begin{array}{ccc}
1 & 0 & 0\\
0 & 1 & 0 \\
0 & 0 & 1/s
\end{array}
\right).
\label{TgijtileT=abh}
\end{equation}
The inverse of this matrix is computed as follows:
\begin{equation}
\left(
\begin{array}{ccc}
a & 0 & ^t h\\
0 & 1 & 0 \\
h & 0 & b
\end{array}
\right)^{-1} = \left(
\begin{array}{ccc}
A & 0 & \, ^tH\\
0 & 1 & 0 \\
H & 0 & B
\end{array}
\right),
\nonumber
\end{equation}
\begin{equation}
A = \big(a - \, ^thb^{-1}h\big)^{-1}, \quad 
H = - b^{-1}hA, \quad B = \big(b - h\,a^{-1}\, ^th\big)^{-1}.
\label{S2AHBform}
\end{equation}

In fact, consider the equation
\begin{equation}
\left(
\begin{array}{ccc}
a & 0 & \, ^th \\
0 & 1 & 0 \\
h & 0 & b
\end{array}
\right)\left(
\begin{array}{ccc}
A & 0 & \, ^tH \\
0 & 1 & 0 \\
H & 0 & B
\end{array}
\right) = \left(
\begin{array}{ccc}
1 & 0 & 0 \\
0 & 1 & 0 \\
0 & 0 & 1
\end{array}
\right), 
\nonumber
\end{equation}
that is
\begin{equation}
\begin{split}
& aA + \, ^thH = 1,  \quad  a\, ^t H + \, ^thB = 0, \\
& hA + bH = 0, \quad  h\, ^tH + bB = 1.
\end{split}
\nonumber
\end{equation}
Since $\det b \neq 0$ by (\ref{Sec2gcondition}), we have $H = - b^{-1}hA$, which implies $(a - \, ^th b^{-1}h)A = 1$. Then, $A$ is written as 
in  (\ref{S2AHBform}). $B$ can be computed similarly.

Therefore, we have 
\begin{equation}
\Big(\widetilde g_{ij}\Big)^{-1} = \Big(\widetilde g^{ij}\Big) = 
\left(
\begin{array}{ccc}
A & 0 & ^tH/s\\
0 & 1 & 0 \\
H /s& 0 & B/s^2
\end{array}
\right).
\label{S3g-1ijformula}
\end{equation}
We put
\begin{equation}
g =\det\big(\widetilde g_{ij}\big).
\nonumber
\end{equation}
Then, by (\ref{S3g-1ijformula}),  the Laplacian is formally decomosed  as
\begin{equation}
L = L_A + L_S + L_{B} + L_H,
\label{L=LA+LB+LS+LH}
\end{equation}
\begin{equation}
 L_A = - \frac{1}{\sqrt{g}}\partial_p\big(\sqrt{g}a^{pq}\partial_q\big),  \quad 
L_{B} = - \frac{1}{\sqrt{g}}\frac{1}{s^2}\partial_{\alpha}\big(\sqrt{g}b^{\alpha\beta}\partial_{\beta}), 
\nonumber
\end{equation}
\begin{equation}
L_S =  - \frac{1}{\sqrt{g}} \partial_s\big(\sqrt{g}\partial_s\big), \quad
L_H = - \frac{1}{\sqrt{g}}\Big(\frac{1}{s}\partial_p\big(\sqrt{g}h^{p\alpha}\partial_{\alpha}) +  \frac{1}{s}\partial_{\alpha}\big(\sqrt{g}h^{\alpha p}\partial_{p})\Big),
\nonumber
\end{equation}
where $\partial_p = \partial/\partial y^p$, $\partial_s = \partial/\partial s$ and $\partial_{\alpha} = \partial/\partial\theta^{\alpha}$, $\theta^{\alpha}$ being the local coordinate on $S^{n-k-1}$.
 Actually, we decompose the  associated quadratic form as
\begin{equation}
q = q_A + q_S + q_{B} + q_H,
\label{q=qA+qS+etc}
\end{equation}
\begin{equation}
q_A(u,v) = \big(a^{pq}\partial_qu,\partial_pv\big), \quad 
q_{B}(u,v) = \big(b^{\alpha\beta}\frac{\partial_{\beta}u}{s},\frac{\partial_{\alpha}v}{s}\big),
\nonumber
\end{equation}
\begin{equation}
q_S(u,v) = \big(\partial_su,\partial_sv\big), \quad 
q_H(u,v) = \big(h^{p\alpha}\frac{\partial_{\alpha}u}{s},\partial_pv\big) + \big(h^{\alpha p}\partial_{p}u,\frac{\partial_{\alpha}v}{s}\big).
\nonumber
\end{equation}
Here,  $q$ is $q(u,v) - (u,v)_{L^2(M)}$ where $q(u,v)$ is defined by (\ref{qdefinedonwidetildU}), and we omit $\sharp\Gamma$ for the sake of simplicity. Recall that in the computation below  $M^{sing}$ and its neighborhood are understood to be lifted to $\widetilde U$.

To study the regularity of $u \in D(L)$ near $M^{sing}$, it is sufficient to consider on a small coordinate patch intersecting $M^{sing}$. We take local coordinates $y, s , \theta$ near $M^{sing}$, where $\theta$ is a local coordinate on $S^{n-k-1}$. Letting $\delta$ be a small positive constant, we can thus  identify the $\delta$-neighborhood of $M^{sing}$ with the following set:
\begin{equation}
M^{sing}_{\delta} = \{(y, s, \theta)  \, ; \, y \in \mathcal O^k, \ 0 \leq s < \delta, \ \theta \in \mathcal O^{n-k-1}\},
\nonumber
\end{equation}
where $\mathcal O^k$ and $\mathcal O^{n-k-1}$ are bounded open sets in ${\mathbb R}^k$ and ${\mathbb R}^{n-k-1}$. 
  Letting $c(y,s,\theta)$ be any of $a_{pq}$, $b_{\ell m}$ and $h_{p\ell}$, we have by the assumption (D) :
\begin{equation}
\mathop{\sup_{y,s,\theta}}|s|^{-\gamma}\big|\partial_y^{\alpha}\partial_s^{\beta}\partial_{\theta}^{\gamma}c(y,s,\theta)\big| < \infty
\nonumber
\end{equation}
for $|\alpha| + \beta + |\gamma| \leq 2$. Moreover,  $a_{pq}$, $b_{\ell m}$ and $h_{p\ell}$ are extended to $s \leq 0$ as $C^2$-functions.  Taylor expansion and (\ref{Sec2abhcondition}) yield
\begin{equation}
\begin{split}
& s ^{-|\gamma|}\Big|\partial_y^{\alpha}\partial_s^{\beta}\partial_{\theta}^{\gamma}\big( a_{pq}(y,s,\theta) - \widehat a_{pq}(y)\big)\Big| \leq Cs, \\
& s ^{-|\gamma|}\Big|\partial_y^{\alpha}\partial_s^{\beta}\partial_{\theta}^{\gamma}\big(b_{\ell m}(y,s,\theta) - \widehat b_{\ell m}(y,\theta)\big)\Big| \leq Cs, \\
& s ^{-|\gamma|}\Big|\partial_y^{\alpha}\partial_s^{\beta}\partial_{\theta}^{\gamma}\big(h_{p\ell}(y,s,\theta)\big)\Big| \leq Cs,
\end{split}
\label{S2Blowupcond}
\end{equation}
for $|\alpha| + \beta + |\gamma| \leq 1$.

Recall that the inner product of $L^2(\widetilde U)$, hence $L^2(M)$,  is defined by
\begin{equation}
(u,v) = \int u\overline{v} \sqrt{g}\, dydsd\theta.
\label{S3:Innerproduct(u,v)}
\end{equation}
By (\ref{TgijtileT=abh}), 
$$
g = s^{2(n-k-1)}\big(\det\, a\det\, b + o(1)\big), \quad s \to 0.
$$
Therefore, letting
\begin{equation}
g(y,s,\theta) = s^{2(n-k-1)}g_0(y,s,\theta),
\label{S3:g=s2(n-k-1)g0}
\end{equation}
we have
$$
C \leq g_0 \leq C^{-1}
$$
for a constant $C > 0$.  
Taking this into account, we introduce the following norm and semi-norm:
\begin{equation}
\|u\|_{a,b,c} = \left(\sum_{|\alpha|\leq a,   \beta \leq b,  |\gamma| \leq c}
\int_{M_{\delta}^{sing}} s^{-2|\gamma|}\big|\partial_y^{\alpha}
\partial_s^{\beta}\partial_{\theta}^{\gamma} u\big|^2 s^{n-k-1}dydsd\theta\right)^{1/2},
\label{S2abcnorm}
\end{equation}
\begin{equation}
|u|_{0,2,0} = \left(
\int_{M_{\delta}^{sing}} \Big|\partial_s^2 u + \frac{n-k-1}{s}
\, \partial_s u\Big|^2
s^{n-k-1} dydsd\theta\right)^{1/2}.
\label{S2020norm}
\end{equation}
We put
\begin{equation}
\|u\| = \|u\|_{0,0,0},
\nonumber
\end{equation}
\begin{equation}
\|u\|_{H^1} = \|u\|_{1,0,0} + \|u\|_{0,1,0} + \|u\|_{0,0,1}.
\nonumber
\end{equation}
Note that for $u$ supported in $M^{sing}_{\delta}$, $\|u\|_{0,0,0}$ is equivalent to $\|u\|_{L^2(M)}$. 
We also use the following inner product:
\begin{equation}
(u,v)_{0} = \int u{\overline v} \, s^{n-k-1}dydsd\theta.
\label{S3:Innerproduct(u,v)0}
\end{equation}

\begin{definition}
\label{S2Definewidetilde HMsingdelta}
Let $\widetilde H^2(M^{sing}_{\delta})$ be the set of functions $u$ such that 
$$
\|u\|_{\widetilde H^2(M^{sing}_{\delta})} : = |u|_{0,2,0} + \sum_{a + b + c \leq 2, b \leq 1}\|u\|_{a,b,c} < \infty.
$$
\end{definition}
 
The aim of this subsection is to show the following theorem.


\begin{theorem}
\label{S2DLregularitytheorem}
$D(L) = H^2_{loc}(M^{reg})\cap \widetilde H^2(M^{sing}_{\delta})$. 
\end{theorem}

The proof of this theorem is done by the standard argument for elliptic regularity, however, it requires careful computation.   
In order to make the singularity with respect to $s$ of the volume element of the quadratic form $q$ more explicit, we make the gauge transformation 
$v \to g_0^{1/4}v$, where $g_0$ is in (\ref{S3:g=s2(n-k-1)g0}). In view of (\ref{S3:Innerproduct(u,v)}) and (\ref{S3:Innerproduct(u,v)0}), we 
rewrite $q$ in (\ref{q=qA+qS+etc}) as
$$
q(g_0^{-1/4}u,g_0^{-1/4}v) = Q(u,v) = Q_A + Q_S + Q_{B} + Q_H,
$$
\begin{equation}
\begin{split}
 Q_A(u,v) &= \big(a^{pq}(\partial_q + \partial_q\log g_0^{-1/4})u,(\partial_p+ \partial_p\log g_0^{-1/4})v\big)_{0}, \\
 Q_{B}(u,v) & = \big(b^{\alpha\beta}\Big(\frac{\partial_{\beta}}{s} + \frac{\partial_{\beta}}{s}\log g_0^{-1/4}\Big)u
,\Big(\frac{\partial_{\alpha}}{s}+ \frac{\partial_{\alpha}}{s}\log g_0^{-1/4}\Big)v\big)_{0},\\
 Q_S(u,v) &= \big((\partial_s + \partial_s\log g_0^{-1/4})u,(\partial_s + \partial_s\log g_0^{-1/4})v\big)_{0}, \\
 Q_H(u,v) &= \big(h^{p\alpha}\Big(\frac{\partial_{\alpha}}{s} + \frac{\partial_{\alpha}}{s}\log g_0^{-1/4}\Big)u,(\partial_p + \partial_p\log g_0^{-1/4})v\big)_{0} \\
& + \big(h^{\alpha p}(\partial_{p} + \partial_p\log g_0^{-1/4})u,\Big(\frac{\partial_{\alpha}}{s} + \frac{\partial_{\alpha}}{s}\log g_0^{-1/4}\Big)v\big)_{0}.
\end{split}
\nonumber
\end{equation}
Note $g_0^{1/4}(\partial_pg_0^{-1/4}u) = \partial_pu + (\partial_p\log g_0^{-1/4})u,  etc.$. 
We put
\begin{equation}
L_{0} = g_0^{1/4}Lg_0^{-1/4}, \quad D(L_0) = g_0^{1/4}D(L).
\nonumber
\end{equation}
Then we have for $u \in D(L_{0})$ and $v \in g_0^{1/4}D(\sqrt{L})$
\begin{equation}
Q(u,v) = (L_{0}u,v)_{0}.
\nonumber
\end{equation}

Here, we recall  a well-known lemma on Freidrich's mollifier. Take $\rho(x) \in C_0^{\infty}({\mathbb R}^m)$ such that $\rho(x) = 0$ for $|x| > 1$ and $\int_{{\mathbb R}^m}\rho(x)dx=1$. For a sufficiently small $\delta > 0$, we put $\rho_{\delta}(x) = \delta^{-m}\rho(x/\delta)$ and let the operator $\rho_{\delta}\ast$ be defined by
$$
\rho_{\delta}\ast : u \to \rho_{\delta}\ast u(x) = \int_{{\mathbb R}^m}\rho_{\delta}(x-y)u(y)dy.
$$
 We also put
$$
\big[\rho_{\delta}\ast,a(x)\frac{\partial}{\partial x_j}\big]u = 
\rho_{\delta}\ast\big(a(x)\frac{\partial u}{\partial x_j}\big) - a(x)\frac{\partial}{\partial x_j}\big(\rho_{\delta}\ast u\big).
$$


\begin{lemma}
\label{LemmaMolifier}
Let $a(x) \in C^1({\mathbb R}^m)$ be such that 
$$
|a|_{\mathcal B^1} : = {\mathop{\sup}_{x\in{\mathbb R}^m}}\sum_{|\alpha| \leq 1}|\partial_x^{\alpha}a(x) | < \infty.
$$
 Then, there exists a constant $C > 0$ independent of $\delta > 0$ such that
\begin{equation}
\Big\|\big[\rho_{\delta}\ast,a(x)\frac{\partial}{\partial x_j}\big]u\Big\|_{L^2({\mathbb R}^m)} \leq C|a|_{\mathcal B^1}\|u\|_{L^2({\mathbb R}^m)},
\label{S2FriedricjsLemma1}
\end{equation}
and as $\delta \to 0$ 
\begin{equation}
\big[\rho_{\delta}\ast,a(x)\frac{\partial}{\partial x_j}\big]u \to 0, \quad {in} \quad L^2({\mathbb R}^m)
\label{S2FriedricjsLemma2}
\end{equation}
for any $u \in L^2({\mathbb R}^m)$. 
\end{lemma}

\begin{proof}
We have
$$
\big[\rho_{\delta}\ast,a(x)\frac{\partial}{\partial x_j}\big]u = \int_{{\mathbb R}^m}u(y)\frac{\partial}{\partial y_j}
\left\{
\rho_{\delta}(x-y)\big(a(x) - a(y)\big)\right\}dy.
$$
Noting that
$$
|a(x) - a(y)| \leq C|a|_{\mathcal B^1}|x-y|, \quad 
\int_{|x-y|<\delta}|x-y|\Big|\frac{\partial \rho_{\delta}}{\partial y_j}(x-y)\Big|dy < C,
$$
we obtain (\ref{S2FriedricjsLemma1}). To prove (\ref{S2FriedricjsLemma2}), we only have to consider the case in which $u \in C_0^{\infty}({\mathbb R}^m)$. Then, we have
\begin{equation}
\begin{split}
& a(x)\frac{\partial}{\partial x_j}\rho_{\delta}\ast u - 
\rho_{\delta}\ast\big(a(x)\frac{\partial}{\partial x_j}u\big)  = \int_{{\mathbb R}^n}\big(a(x) - a(y)\big)\rho_{\delta}(x-y)\frac{\partial u}{\partial y_j}(y)dy.
\end{split}
\nonumber
\end{equation}
Noting that
$$
\int|x-y|\rho_{\delta}(x-y)\Big|\frac{\partial u}{\partial y_j}(y)\Big|dy \leq C\delta\int\rho_{\delta}(x-y)\Big|\frac{\partial u}{\partial y_j}(y)\Big|dy,
$$
we obtain (\ref{S2FriedricjsLemma2}). 
\end{proof}

\medskip
It is sufficient to prove Theorem \ref{S2DLregularitytheorem} with $D(L)$ replaced by $D(L_{0})$. Let $u \in D(L_{0})$. 
Since $g_0^{1/4}u \in D(\sqrt{L})$, we have by (\ref{H1M=D(q)})
\begin{equation}
\sum_{|\alpha| + \beta + |\gamma| \leq 1}
\int s^{-2|\gamma}\big|\partial_y^{\alpha}
\partial_s^{\beta}\partial_{\theta}^{\gamma} u\big|^2 s^{n-k-1}dydsd\theta < \infty.
\label{DsqrtLnearMsing}
\end{equation}

We consider the 2nd order derivatives of $u \in D(L_{0})$. 
Take $\rho(y,\theta) \in C^{\infty}_0({\mathbb R}^{n-1})$ such that
 $\rho(y,\theta) = 0$ for $|y|^2 + |\theta|^2 > 1$ and 
$\int_{{\mathbb R}^{n-1}}\rho(y,\theta)dyd\theta = 1$.
We put $\rho_{\delta}(y,\theta)= \delta^{-(n-1)}\rho(y/\rho,\theta/\rho)$. 
Note that the operator $\rho_{\delta}\ast$ has the property in Lemma \ref{LemmaMolifier} with $L^2({\mathbb R}^m)$-norm replaced by $\|\cdot\| = \|\cdot\|_{0,0,0}$.

\begin{lemma}
\label{Lemma2.9}
If $u \in D(L_{0})$, then $\rho_{\delta}\ast u \in D(L_{0})$, and $L_0(\rho_{\delta}\ast u) \to f$ in $L^2(M)$ as $\delta \to 0$.
\end{lemma}

\begin{proof}
 If $u \in D(L_{0})$, letting $L_{0}u = f$, we have for $v \in g_0^{1/4}D(\sqrt{L})$, $Q(u,v) = (f,v)_{0}$. 
We replace $v$ by $\rho_{\delta}\ast v$. Then, since $\rho_{\delta}\ast$ is self-adjoint and commutes with $\partial_q$ and $\partial_{\theta}$, we have
\begin{equation}
\begin{split}
Q_A(u,\rho_{\delta}\ast v) & = Q_A(\rho_{\delta}\ast u,v) + (f_{1\delta},v)_0,\\
f_{1\delta}& = (\partial_p + \partial_p\log g_0^{-1/4})^{\ast}
[\rho_{\delta}\ast,a^{pq}(\partial_q + \partial_q\log g_0^{-1/4})]u\\
& + ([\partial_p + \partial_p\log g_0^{-1/4},\rho_{\delta}{\ast}]^{\ast}a^{pq}(\partial_q + \partial_q\log g_0^{-1/4})u, 
\end{split}
\nonumber
\end{equation}
\begin{equation}
\begin{split}
Q_B(u,\rho_{\delta}\ast v) &= Q_B(\rho_{\delta}\ast u,v) +(f_{2\delta},v)_0, \\
f_{2\delta} &= \big(\frac{\partial_{\alpha}}{s} + \frac{\partial_{\alpha}}{s}\log g_0^{-1/4}\big)^{\ast}\big[\rho_{\delta}\ast,b^{\alpha\beta}\big(\frac{\partial_{\beta}}{s} + \frac{\partial_{\beta}}{s}\log g_0^{-1/4}\big)\big]u \\
&+ \big[\frac{\partial_{\alpha}}{s} + \frac{\partial_{\alpha}}{s}\log g_0^{-1/4},\rho_{\delta}\ast\big]^{\ast}b^{\alpha\beta}\big(\frac{\partial_{\beta}}{s} + \frac{\partial_{\beta}}{s}\log g_0^{-1/4}\big)u,
\end{split}
\nonumber
\end{equation}
\begin{equation}
\begin{split}
Q_{H}(u,\rho_{\delta}\ast v) &= Q_H(\rho_{\delta}\ast u,v) + (f_{3\delta},v)_0,\\
f_{3\delta}& = (\partial_p + \partial_p\log g_0^{-1/4})^{\ast}
\big[\rho_{\delta}\ast,h^{p\alpha}\big(\frac{\partial_{\alpha}}{s} + \frac{\partial_{\alpha}}{s}\log g_0^{-1/4}\big)\big]u \\
& + 
\big[\partial_p + \partial_p\log g_0^{-1/4},\rho_{\delta}\ast]^{\ast}h^{p\alpha}\big(\frac{\partial_{\alpha}}{s} + \frac{\partial_{\alpha}}{s}\log g_0^{-1/4}\big)u \\
& + \big(\frac{\partial_{\alpha}}{s} + \frac{\partial_{\alpha}}{s}\log g_0^{-1/4}\big)^{\ast}[\rho_{\delta}\ast,h^{\alpha p}(\partial_p + \log g_0^{-1/4})]u\\
& + \big[\frac{\partial_{\alpha}}{s} + \frac{\partial_{\alpha}}{s}\log g_0^{-1/4},\rho_{\delta}\ast\big]^{\ast}h^{\alpha p}(\partial_p + \log g_0^{-1/4})u,
\end{split}
\nonumber
\end{equation}
\begin{equation}
\begin{split}
Q_S(u,\rho_{\delta}\ast v) &= Q_S(\rho_{\delta}\ast u,v) + (f_{4\delta},v)_0,\\
f_{4\delta} & = ([\partial_s\log g_0^{-1/4},\rho_{\delta}\ast]^{\ast}(\partial_s + \partial_s\log g_0^{-1/4})u.
\end{split}
\nonumber
\end{equation}
Summing up these 4 terms, we obtain
$$
Q(u,\rho_{\delta}\ast v) = Q(\rho_{\delta}\ast u,v) + (\sum_{i=1}^4f_{i\delta},v)_0.
$$
On the other hand, we have
$$
Q(u,\rho_{\delta}\ast v) = (L_0u,\rho_{\delta}\ast v)_0 = (f,\rho_{\delta}\ast v)_0 = (\rho_{\delta}\ast f,v)_0.
$$
Therefore, we have
$$
Q(\rho_{\delta}\ast u,v) = (\rho_{\delta}\ast f,v)_0 - (\sum_{i=1}^4f_{i\delta},v)_0.
$$
Since this holds for all $v \in g_0^{-1/4}D(\sqrt{L})$, we see that $\rho_{\delta}\ast u \in D(L_{0})$ and
\begin{equation}
L_{0} (\rho_{\delta}\ast u) = \rho_{\delta}\ast f - \sum_{i=1}^4f_{i\delta}.
\nonumber
\end{equation}
Since $u \in H^1(M)$, by Lemma \ref{LemmaMolifier}, 
$L_0(\rho_{\delta}\ast u) \to f$ in $\| \cdot \|_{0,0,0}$ norm 
as $\delta \to 0$. 
This proves the lemma. 
\end{proof}

As a formal differential operator, $L_0$ is rewritten as
\begin{equation}
L_0 = L_A^{(0)} + L^{(0)}_S + L^{(0)}_{B} + L^{(0)}_H + L^{(0)}_R,
\nonumber
\end{equation}
\begin{equation}
 L^{(0)}_A = - \partial_p\big(a^{pq}\partial_q\big),  \quad L^{(0)}_S =  - 
s^{-(n-k-1)} \partial_s\big(s^{n-k-1}\partial_s\big), 
\nonumber
\end{equation}
\begin{equation}
L^{(0)}_{B} = - \frac{\partial_{\alpha}}{s}\big(b^{\alpha\beta}\frac{\partial_{\beta}}{s}),  \quad 
L^{(0)}_H = - \Big(\partial_p\big(h^{p\alpha}\frac{\partial_{\alpha}}{s}) +  \frac{\partial_{\alpha}}{s}\big(h^{\alpha p}\partial_{p})\Big),
\nonumber
\end{equation}
where $L^{(0)}_R$ consists of first order and zeroth order terms of differential operators $\partial_p$, $\frac{\partial_{\alpha}}{s}$, $\partial_s$ with bounded coefficients. 
We also decompose the quadratic form $Q$ as
\begin{equation}
Q(u,v) = Q^{(0)}_A(u,v) + Q^{(0)}_S(u,v) + Q^{(0)}_B(u,v) + Q^{(0)}_H(u,v) + Q^{(0)}_R(u,v),
\nonumber
\end{equation}
\begin{equation}
\begin{split}
 Q^{(0)}_A(u,v) &= \big(a^{pq}\partial_qu,\partial_pv\big)_{0}, \\
 Q^{(0)}_S(u,v) &= \big(\partial_su,\partial_sv\big)_{0}, \\
 Q^{(0)}_{B}(u,v) & = \big(b^{\alpha\beta}\frac{\partial_{\beta}}{s}u
,\frac{\partial_{\alpha}}{s}v\big)_{0},\\
 Q^{(0)}_H(u,v) &= \big(h^{p\alpha}\frac{\partial_{\alpha}}{s}u,\partial_pv\big)_{0} + \big(h^{\alpha p}\partial_{p}u,\frac{\partial_{\alpha}}{s}v\big)_{0}, \\
Q^{(0)}_R(u,v) &= (L^{(0)}_Ru,v)_0.
\end{split}
\nonumber
\end{equation}

We are going to estimate ${\rm Re}\,Q_E^{(0)}(u,L_F^{(0)}u)$, $F \neq S$, from below.  In the following, $o(H^2)$ denotes a term which is estimated as follows: For any $\epsilon > 0$, there exists a constant $C_{\epsilon} > 0$ such that
\begin{equation}
\big|o(H^2)| \leq \epsilon \sum_{a+ b + c \leq 2, b \leq 1}\|u\|^2_{a,b,c} + 
C_{\epsilon}\|u\|_{H^1}^2.
\label{Lemma2.10Proof0}
\end{equation}
We put
\begin{equation}
m(\delta) = \max_{p, \alpha}\sup_{0<s<\delta, y, \theta}\big(|h^{p\alpha}(y,s,\theta)| + |h^{\alpha p}(y,s,\theta)|\big).
\nonumber
\end{equation}
By the assumption (D), $m(\delta) \to 0$ as $\delta \to 0$.

\begin{lemma}
\label{Lemma2.10}
There exists a constant $C > 0$ such that if $u$ is in $ D(L_0)$ and three times differentiable with respect to $y$ and $\theta$, we have
\begin{equation}
\|u\|_{2,0,0}^2 \leq C{\rm Re}\, Q^{(0)}_A(u,L_A^{(0)}u) + o(H^2),
\label{Q0ALAu}
\end{equation}
\begin{equation}
\|u\|_{1,0,1}^2 \leq C{\rm Re}\, Q_B^{(0)}(u,L^{(0)}_Au) + o(H^2),
\label{IntroLemmaQBLA}
\end{equation}
\begin{equation}
\|u\|_{1,1,0}^2 \leq C{\rm Re}\, Q_S^{(0)}(u,L_A^{(0)}u) + o(H^2),
\label{Lemma2.10QSoLAu}
\end{equation}
\begin{equation}
\big| Q_H^{(0)}(u,L_A^{(0)}u)\big| \leq Cm(\delta)\sum_{a+b+c\leq 2, b\leq 1}\|u\|_{a,b,c}^2,
\label{Lemma2.10QH0LAu}
\end{equation}
\begin{equation}
 Q_R^{(0)}(u,L_A^{(0)}u) = o(H^2).
\label{Lemma2.10QRLA}
\end{equation}
\end{lemma}

\begin{proof}
By integration by parts if necessary,  ${\rm Re}\, Q_R^{(0)}(u,L_A^{(0)}u) $ consists of the terms like e.g.,
\begin{equation}
(\cdots \partial_qu,\cdots \partial_p\frac{\partial_{\beta}}{s}u)_0, \quad 
(\cdots \partial_qu,\cdots \frac{\partial_{\beta}}{s}u)_0,
\nonumber
\end{equation}
which proves (\ref{Lemma2.10QRLA}). 
Similar integration by parts yields  (\ref{Lemma2.10QH0LAu}), (\ref{Lemma2.10QSoLAu}),   (\ref{Q0ALAu}).
We prove (\ref{IntroLemmaQBLA}). By integration by parts,
\begin{equation}
Q_B^{(0)}(u,L_A^{(0)}u) = (b^{\alpha\beta}\partial_p\frac{\partial_{\beta}}{s}u,a^{pq}\partial_q\frac{\partial_{\alpha}}{s}u)_0 + o(H^2).
\nonumber
\end{equation}
Let $C(\eta)$ be the cube in $(y,\theta)$ of size $\eta$ with center at the origin. 
We take a partition of unity $\{\chi_j\}$ such that $\sum_{j}\chi_j^2 = 1$ and each $\chi_j$ is a translation of fixed $\chi$ which has support inside $C(2\eta)$ and $\chi=1$ on $C(\eta)$. Let $C_j(\eta)$ be the associated translation of $C(\eta)$ and $x_{j}$  the center of $C_j(\eta)$. 
Inserting $1 = \sum_j\chi_j^2$, we have
\begin{equation}
Q_B^{(0)}(u,L_A^{(0)}u) = \sum_j\big(b^{\alpha\beta}\partial_p\frac{\partial_{\beta}}{s}(\chi_ju),a^{pq}\partial_q\frac{\partial_{\alpha}}{s}(\chi_ju)\big)_0 + 
o(H^2).
\nonumber
\end{equation}
We put 
\begin{equation}
a^{pq}_j = a^{pq}(x_j), \quad b^{\alpha\beta}_j = b^{\alpha\beta}(x_j). 
\nonumber
\end{equation}
Given $\epsilon > 0$, 
taking $\eta$ small enough independently of $j$, we have
\begin{equation}
\begin{split}
& {\rm Re}\, \big(b^{\alpha\beta}\partial_p\frac{\partial_{\beta}}{s}(\chi_ju),a^{pq}\partial_q\frac{\partial_{\alpha}}{s}(\chi_j u)\big)_0 \\
& \geq {\rm Re}\, \big(b_j^{\alpha\beta}\partial_p\frac{\partial_{\beta}}{s}(\chi_ju), a^{pq}_j\partial_q\frac{\partial_{\alpha}}{s}(\chi_j u)\big)_0 - 
\epsilon\|\chi_ju\|^2_{1,0,1}.
\end{split}
\label{Lemma2.10Proof1}
\end{equation}
Since $\big(a^{pq}_j\big), \big(b^{\alpha\beta}_j\big)$ are positive definte matrices corresponding to independent variables $y, \theta$, by passing to the Fourier transform with respect to $y$ and $\theta$, we have
\begin{equation}
{\rm Re}\,\big(b^{\alpha\beta}\partial_p\frac{\partial_{\beta}}{s}(\chi_ju), a^{pq}_j\partial_q\frac{\partial_{\alpha}}{s}(\chi_j u)\big)_0 \geq 
C_0\big((- \frac{\Delta_{\theta}}{s^2})\chi_j u,(- \Delta_y)\chi_j u\big)_0,
\label{Lemma2.10Proof2}
\end{equation}
where $\Delta_y$ and $\Delta_{\theta}$ are the Laplacians with respect to $y$ and $\theta$, and $C_0$ is a positive constant independent of $j$.
It is easy to show that
\begin{equation}
\|\chi_j\|^2_{1,0,1}  \leq C\Big( 
\big((- \frac{\Delta_{\theta}}{s^2})\chi_ju,(- \Delta_y)\chi_ju\big)_0
+ \|\chi_j u\|^2_{H^1}\Big).
\label{Lemma2.10Proof3}
\end{equation}
The inequalities (\ref{Lemma2.10Proof1}), (\ref{Lemma2.10Proof2}) and (\ref{Lemma2.10Proof3}) prove (\ref{IntroLemmaQBLA}). 
\end{proof}

Similarly, one can prove the following  two lemmas . 

\begin{lemma}
\label{Lemma2.11}
There exists a constant $C > 0$ such that if $u$ is in $ D(L_0)$ and three times differentiable with respect to $y$ and $\theta$, we have
\begin{equation}
\|u\|_{1,0,1}^2 \leq C{\rm Re}\, Q_A^{(0)}(u,L^{(0)}_Bu) + o(H^2),
\label{IntroLemmaQALB}
\end{equation}
\begin{equation}
\|u\|_{0,0,2}^2 \leq C{\rm Re}\, Q^{(0)}_B(u,L_B^{(0)}u) + o(H^2),
\label{Q0BLBu}
\end{equation}
\begin{equation}
\|u\|_{1,1,0}^2 \leq C{\rm Re}\, Q_S^{(0)}(u,L_B^{(0)}u) + o(H^2),
\label{Lemma2.10QSoLBu}
\end{equation}
\begin{equation}
\big| Q_H^{(0)}(u,L_B^{(0)})\big| \leq Cm(\delta)\sum_{a+b+c\leq 2, b\leq 1}\|u\|_{a,b,c}^2,
\label{Lemma2.10QH0LBu}
\end{equation}
\begin{equation}
Q_R^{(0)}(u,L_B^{(0)}u) = o(H^2).
\label{Lemma2.10QRLB}
\end{equation}
\end{lemma}

\begin{lemma}
There exists a constant $C > 0$ such that if $u$ is in $ D(L_0)$ and three times differentiable with respect to $y$ and $\theta$, we have for $E \neq H$
$$
\big|Q_E^{(0)}(u,L_H^{(0)}u)\big| \leq Cm(\delta)\sum_{a+b+c\leq 2, b\leq 1}\|u\|_{a,b,s}^2.
$$
\end{lemma}

The following estimates are our main purpose.

\begin{lemma}
\label{Lemma2.12}
There exists a constant $C > 0$ such that for any $u \in D(L_0)$
\begin{equation}
\sum_{a+b+c\leq 2, b \leq 1}\|u\|_{a,b,c} \leq 
C\big(\|L_0u\| + \|u\|_{H^1}\big),
\label{Lemma2.10estimate1}
\end{equation}
\begin{equation}
\|\big(\partial_s^2 + \frac{n-k-1}{s}\partial_s\big) u\| \leq C\big(\|L_0u\| + \|u\|_{H^1}\big).
\label{Lemma2.10estimate2}
\end{equation}
\end{lemma}

\begin{proof}
 By virue of Lemma \ref{Lemma2.9}, we only have to prove this lemma when $u$ is sufficiently smooth with respect to $y$ and $\theta$. Let $L_0u = f$. We put $v = L_A^{(0)}u$
in the quadratic form $Q(u,v) = (f,v)_0$. Then,
\begin{equation}
\begin{split}
&Q_A^{(0)}(u,L_A^{(0)}u) + Q_S^{(0)}(u,L_A^{(0)}u) + Q_B^{(0)}(u,L_A^{(0)}u) \\
& + Q_H^{(0)}(u,L_A^{(0)}u) + Q_R^{(0}(u,L_A^{(0)}u) = 
(f,L_A^{(0)}u)_0.
\end{split}
\label{Lemma3.8ProofInnerproduct}
\end{equation}
 By Lemma \ref{Lemma2.10}, we have
$$
- {\rm Re}\, Q_S^{(0)}(u,L_A^{(0)}u) \leq o(H^2), \quad - {\rm Re}\, Q_B^{(0)}(u,L_A^{(0)}u) \leq o(H^2), 
$$
$$
- {\rm Re}\, Q_H^{(0)}(u,L_A^{(0)}u) \leq Cm(\delta)\sum_{a+b+c\leq 2, b\leq 1}\|u\|_{a,b,c}^2, \quad - {\rm Re}\, Q_R^{(0)}(u,L_A^{(0)}u) \leq o(H^2).
$$
We take $\epsilon > 0$ small enough, and put
$$
o_{\epsilon}(H^2) = \epsilon \sum_{a+ b + c \leq 2, b \leq 1}\|u\|^2_{a,b,c} + 
C_{\epsilon}\|u\|_{H^1}^2.
$$
Noting that $Q_A^{(0)}(u,L_A^{(0)}u) = \|L_A^{(0)}u\|^2$, we have by (\ref{Lemma3.8ProofInnerproduct}), choosing $\delta$ small enough
$$
\|L_A^{(0)}u\|^2 \leq C\|f\|^2 + o_{\epsilon}(H^2).
$$
Using (\ref{Q0ALAu}), we then have
\begin{equation}
\|u\|_{2,0,0}^2 \leq C\|f\|^2 + o_{\epsilon}(H^2).
\label{u200fromabove}
\end{equation}
By the similar argumens, we have
$$
\|L_B^{(0)}u\|^2 \leq C\|f\|^2 + o_{\epsilon}(H^2),
$$
\begin{equation}
\|u\|_{0,0,2}^2 \leq C\|f\|^2 + o_{\epsilon}(H^2),
\label{u002fromabove}
\end{equation}
$$
\|L_H^{(0)}u\|^2 \leq C\|f\|^2 + o_{\epsilon}(H^2),
$$
\begin{equation}
\|u\|_{1,0,1}^2 \leq C\|f\|^2 + o_{\epsilon}(H^2),
\label{u101fromabove}
\end{equation}
$$
\|L_R^{(0)}u\|^2 \leq C\|f\|^2 + o_{\epsilon}(H^2).
$$
Since $L_0 =   L_A^{(0)} + L_S^{(0)} +  L_B^{(0)} +  L_H^{(0)} +  L_R^{(0)}$, we then have
$$
\|L_S^{(0)}u\|^2 \leq C\|f\|^2 + o_{\epsilon}(H^2).
$$
This yields
\begin{equation}
\|\big(\partial_s^2 + \frac{n-k-1}{s}\big)u\|^2 \leq C\|f\|^2 + o_{\epsilon}(H^2).
\label{partials2+n-k-1overspartials}
\end{equation}
Summing up (\ref{u200fromabove}), (\ref{u200fromabove}), (\ref{u200fromabove}) and (\ref{partials2+n-k-1overspartials}), we get finally
$$
\|u\|_{\widetilde H^2}^2 \leq C\|f\| +  \epsilon\|u\|_{\widetilde H^2}^2 + C_{\epsilon}\|u\|^2_{H^1}.
$$
This proves the lemma. 
 \end{proof}

 Theorem \ref{S2DLregularitytheorem} now follows from Lemma \ref{Lemma2.12}. In particular, for any $u \in D(L)$, 
\begin{equation}
Lu = L_Au + L_Su + L_Bu + L_Hu,
\nonumber
\end{equation}
where each term of the right-hand side makes sense in the sense of distribution and belongs to $L^2(M)$. 


\subsection{Limit metric and its perturbation} 
\label{LimitMetricPerturbation}
We return to our manifold $\mathcal M$  of the form (\ref{S0manifoldmathcalM}). We need to change the definition of $S^{\kappa}$ in (\ref{IntoDefineSkappa}) as follows. 
Let $\{\chi_j\}_{j\in J_M}$ be the set of partition of unity on $M$, and define 
$p_j^{(m)}(f)$ by (\ref{DefinepjmS2}). We put
\begin{equation}
\begin{split}
p_M(f) &= \sum_{j\in J_M}p_j^{(2)}(f) \\
&= \sum_{j \in J_M}\sup_{y,s,\theta}\sum_{|\alpha| + \beta + |\gamma|\leq 2}s^{-|\gamma|}
\big|\partial_y^{\alpha}\partial_s^{\beta}\partial_{\theta}^{\gamma}
\big(\chi_jf(y,s,\theta)\big)\big|.
\end{split}
\nonumber
\end{equation}

\begin{definition}
\label{S2DefineSingularSkappa}
We define $S^{\kappa}$ \index{$S^{\kappa}$} to be the set of $C^{\infty}((0,\infty);C^2(M^{reg})\big)$-functions $f$ on $\mathcal M$  satisfying
\begin{equation}
p_M\big(\partial_r^{\ell} f(r)\big) \leq C(1 + r)^{\kappa - \ell}, \quad  \forall \ell \geq 0.
\label{DefineSkappanew}
\end{equation}
\end{definition}

The assumptions on $\mathcal M$ are as follows.  

\begin{itemize}
\item $\mathcal K$ {\it is a relatively compact  $n$-dimensional regular conic manifold.}

\item  
{\it Each end  $\mathcal M_i$ $(i = 1,\cdots, N + N')$ is an $n$-dimensional CMGA with the following properties.  There exist an $(n-1)$-dimensional CMGA $M_i$ and a family of metrics $h_i(r,x,dx)\ (r > 0)$ on $M_i$ such that 
$\mathcal M_i$ is diffeomorphic to $(0,\infty)\times M_i$ and equipped with the metric 
\begin{equation}
ds_i^2 = (dr)^2 + \rho_i(r)^2h_i(r,x,dx),
\nonumber
\end{equation}
where  $\rho_i(r)$ satisfies (A-2), (A-3).
Moreover, there exists a metric $h_{M_i}(x,dx)$ as a CMGA on $M_i$ such that}
\begin{equation}
h_i(r,x,dx) - h_{M_i}(x,dx) \in S^{-\gamma_{0,i}}, \quad \gamma_{0,i} > 1.
\nonumber
\end{equation}
\end{itemize}

\begin{remark}
The above assumption is stronger than actually needed. In Definition \ref{S2DefineSingularSkappa}, we have only to assume that $f \in S^{\kappa}$ if and only if $C^{2}((0,\infty);C^2(M^{reg})\big)$-functions $f$ on $\mathcal M$  satisfying
\begin{equation}
p_M\big(\partial_r^{\ell} f(r)\big) \leq C(1 + r)^{\kappa - \ell}, \quad  0 \leq  \ell \leq 2.
\nonumber
\end{equation}
For the assumption for the metric $g_{ij}$, this is still stronger. However, in order not to make the assumption too complicated, we proceed under the above condition.
\end{remark}

Our main concern for the conic singularities in the previous subsection is the local structure of the singular set $M^{sing}$ and the regularity of the domain of the Laplacian. They are obviously invariant by the multiplication by $C^1$-functions to the metric. Then we can assume without loss of generality that for each end $\mathcal M_i$ 
\begin{equation}
\mathcal M_i^{sing} = (0,\infty) \times M_i^{sing},
\nonumber
\end{equation}
and $M_i^{sing}$ has the structure described in the previuous section.
We omit the subscript $i$ for the sake of simplicity. 
Around $\mathcal M^{sing}$, letting $x = (y,z)$ be the local coordinate on $M$, $s = |z|$, $\omega = z/s$ and $\theta$ the local coordinate on $S^{n-k-1}$, the Riemannian metric of $\mathcal M$ is rewritten as
\begin{equation}
\begin{split}
g   = & (dr)^2  + \rho(r)^2\Big(\sum_{p,q=1}^ka_{pq}(r,y,s,\theta)dy^pdy^q + ds^2\\
&+ s^2\sum_{\ell,m=1}^{n-1-k}b_{\ell m}(r,y,s,\theta)d\omega^{\ell}d\omega^m + s\sum_{p=1}^{k}\sum_{\ell=1}^{n-1-k}h_{p\ell}(r,y,s,\theta)dy^pd\omega^{\ell}\Big),
\end{split}
\label{S2g=s2(n-k-1)g1}
\end{equation}
where the coefficients $a_{pq}(r,y,s,\theta)$, $b_{\ell m}(r,y,s,\theta)$, $h_{p\ell}(r,y,s,\theta)$ satisfy the assumptions (C-1), (C-2) and (C-3) uniformly with respect to $r > 0$. Moreover, there exist $a^{\infty}_{pq}(y,s,\theta)$, $b^{\infty}_{\ell m}(y,s,\theta)$, $h^{\infty}_{p\ell}(y,s,\theta)$ such that 
\begin{equation}
p_M(\partial_r^m\alpha-\partial_r^m \alpha^{\infty})  \leq Cr^{-\gamma_0- m},
\nonumber
\end{equation}
where $\alpha$ is any of $a_{pq}(r,y,s,\theta)$, $b_{\ell m}(r,y,s,\theta)$, $h_{p\ell}(r,y,s,\theta)$, and 
$\alpha^{\infty}$ is any of $a^{\infty}_{pq}(y,s,\theta)$, $b^{\infty}_{\ell m}(y,s,\theta)$, $h^{\infty}_{p\ell}(y,s,\theta)$.
Using (\ref{S2g=s2(n-k-1)g1}), we also see that
$g$ satisfies
\begin{equation}
g = \rho(r)^{2(n-1)}s^{2(n-k-2)}g_1, \quad 
C < g_1 < C^{-1}
\label{gformathcalMestimate}
\end{equation}
for a constant $C > 0$. 


\subsection{Laplacians}
The Laplacian $- \Delta_{\mathcal M}$ on $\mathcal M$ has a self-adjoint realization in $L^2(\mathcal M)$ through the quadratic form 
$\sum_{i,j} (g^{ij}\partial_i u,\partial_j v)$, 
which is denoted by $H$. By  Lemma \ref{LemmaH1M=H1Mreg} and Theorem \ref{S2DLregularitytheorem}, it satisfies 
\begin{equation}
D(\sqrt{H}) = H^1(\mathcal M), 
\label{DsqrtHinSection3}
\end{equation}
\begin{equation}
D(H) = H^2_{loc}(\mathcal M^{reg})\cap {\widetilde H}^2(\mathcal M^{sing}_{\delta}).
\label{DHinsection3}
\end{equation}
We take $\chi_{reg}, \chi_{sing} \in C^{\infty}(\mathcal M)$ such that $\chi_{reg} + \chi_{sing} = 1$ on $\mathcal M$, $\chi_{reg} = 0$ on a small neighborhood of $\mathcal M^{sing}$. Choose $\delta > 0$ so that 
${\rm supp}\,\chi_{sing} \subset \mathcal M_{\delta}^{sing}$. 
We put
\begin{equation}
P_{a,b,c,d}(u) = 
{\Big(}\sum_{\substack{|\alpha| \leq a, \beta \leq b\\|\gamma|\leq c, \kappa \leq d}}s^{n-k-1-2|\gamma|}\rho(r)^{-2(|\alpha| + \beta + |\gamma|)}
\Big|\partial_{y}^{\alpha}\partial_s^{\beta}\partial_{\theta}^{\gamma}\partial_r^{\kappa}u\Big|^2{\Big)}^{1/2},
\nonumber
\end{equation}
and for $t \in {\mathbb R}$, define the $\widetilde H^{2,t}(\mathcal M^{sing}_{\delta})$ norm by
\begin{equation}
\begin{split}
\|u\|^2_{\widetilde H^{2,t}(\mathcal M^{sing}_{\delta})} =& 
\sum_{\substack{a+b+c+d \leq 2\\b\leq 1}}
 \int_{{\mathcal M}^{sing}_{\delta}}(1 + r)^{2t}P_{a,b,c,d}(u)^2\rho^{n-1}(r)drdydsd\theta \\
& + \int_{{\mathcal M}^{sing}_{\delta}}(1 + r)^{2t}
\Big|\partial_s^2u + \frac{n-k-1}{s}\partial_su\Big|^2s^{n-k-1}\rho^{n-1}(r)drdydsd\theta.
\end{split}
\nonumber
\end{equation}
We also put
\begin{equation}
\|u\|_{\widetilde H^{2}(\mathcal M^{sing}_{\delta})} = \|u\|_{\widetilde H^{2,0}(\mathcal M^{sing}_{\delta})},
\nonumber
\end{equation}
\begin{equation}
\|u\|_{\widetilde H^2(\mathcal M)} = \|\chi_{reg}u\|_{H^2(\mathcal M^{reg})} + 
\|\chi_{sing}u\|_{{\widetilde H}^2({\mathcal M}_{\delta}^{sing})}.
\label{S2DefineH2tildenorm}
\end{equation} 
By Lemma \ref{Lemma2.12}, the following elliptic regularity theorem holds.

\begin{theorem}
\label{S2ModifiedEllipticregularity}
There exists a constant $C > 0$ such that for any $u \in D(H)$
$$
\|u\|_{\widetilde H^2(\mathcal M)} \leq C(\|Hu\|_{L^2(\mathcal M)} + \|u\|_{L^2(\mathcal M)}).
$$
\end{theorem}

By our assumption, each end $\mathcal M_i$ is diffeomorphic to 
$(0,\infty)\times M_i$. Recalling that $M_i$ has two metrics $h_i(r,x,dx)$ and $h_{M_i}(x,dx)$, let $\Lambda_i(r)$ and $\Lambda_i$ be the associated Laplace-Beltrami operators on $M_i$, respectively.  By the same arguments as above, $-\rho(r)^{-2}\Lambda_i(r)$ and $- \rho(r)^{-2} \Lambda_i$ have self-adjoint realizations in $L^2(M_i)$ as Friedrichs extensions, which are denoted by $B_i(r)$ and $B_i$. Therefore,
\begin{equation}
- \Delta_{\mathcal M} = - \frac{\partial^2}{\partial r^2} - \frac{\partial_rg_i}{2g} \frac{\partial}{\partial r} + B_i(r) 
\quad {\it on} \quad \mathcal M_i,
\nonumber
\end{equation}
\begin{equation}
B_i(r) = \rho_i(r)^{-2}\Lambda_i(r) \quad {\it on} \quad M_i.
\nonumber
\end{equation}
They satisfy
\begin{equation}
D(\sqrt{B_i(r)}) = D(\sqrt{B_i}) = H^1(M_i).
\nonumber
\end{equation}
We impose the Dirichlet boundary condition at $r = 1$, and let $H_i$ and $H_{0,i}$  be the Laplacians on $\widetilde{\mathcal M}_i = (1,\infty)\times M_i$ associated with metrics $(dr)^2 + \rho_i(r)^2h_i(r,x,dx)$ and $(dr)^2 + \rho_i(r)^2h_i(x,dx)$, respectively. They are the Friedrichs extensions of 
the Laplace operators restricted to $C_0^{\infty}(\widetilde{\mathcal M}_i)$. Hence (\ref{DsqrtHinSection3}) and (\ref{DHinsection3}) hold also for 
$H_i$ and $H_{0,i}$.

We derive  resolvent equations for $H_i$ and  $H$ restricted to $\mathcal M_i$. Let 
$$
R(z) = (H - z)^{-1}, \quad R_{i}(z) = (H_i - z)^{-1}.
$$
Take $\chi_i \in C^{\infty}(\mathcal M)$ such that 
\begin{equation}
\chi_i = \left\{
\begin{split}
&1 \quad {\rm on}\quad  (2,\infty)\times M_i \subset \mathcal M_i, \\
& 0\quad {\rm on}  \quad\big((0,1)\times M_i\big) \cup\big(\mathcal M \setminus \mathcal M_i\big).
\end{split}
\right.
\nonumber
\end{equation}
Since the elements in the domains of $H$ and $H_i$ are twice differentiable with respect to $r$, and $\chi_i$ depends only on $r$,  we have 
\begin{equation}
\chi_i R_{i}(z)f \in D(H), \quad \forall f \in L^2(\widetilde{\mathcal  M}_i),
\label{S1chiRAiD(A)}
\end{equation}
\begin{equation}
\chi_iR_{i}(z)f = R(z)\chi_if + R(z)[H,\chi_i]R_{i}(z)f.
\label{S1chiiRAi=RA}
\end{equation}
Conversely, we have
\begin{equation}
\chi_i R(z)f \in D(H_i), \quad \forall f \in L^2(\mathcal M),
\label{S1chiiRAinDAI}
\end{equation}
\begin{equation}
\chi_iR(z)f = R_{i}(z)\chi_if + R_{i}(z)[H_i,\chi_i]R(z)f.
\label{S1chiiRA=RAi}
\end{equation}

We have thus finished preliminary consideration for the Laplacian on CMGA. We add here 
two more facts which are useful to study  the cusp end. 


\subsection{Functional calculus}
We introduce a formula due to Helffer-Sj{\"o}strand \cite{HelSjos} on the representation of functions of self-adjoint operators in terms of their resolvents. We use the following notation. 
$$
{\mathbb C} \ni z = x + iy, \quad 
\overline{\partial_z} = \frac{1}{2}\big(\partial_x + i\partial_y\big), \quad dzd\overline{z} = - 2i dxdy.
$$

\begin{lemma}
\label{C1S3AlmostAnalyticExt}
If $f(x) \in C^{\infty}({\mathbb R})$ satisfies for some $s \in {\mathbb R}$
\begin{equation}
|f^{(k)}(x)| \leq C_k (1 + |x|)^{s-k}, \quad \forall k \geq 0,
\label{S2Lemma2.15f(x)Condition}
\end{equation}
there exists $F(z) \in C^{\infty}({\mathbb C})$, called an almost analytic extension of $f(x)$, satisfying
\begin{equation}
\left\{
\begin{split}
& F(x) = f(x), \quad x \in {\mathbb R}, \\
& \big|F(z)\big| \leq C(1 + |z|)^s, \quad z \in {\mathbb C}, \\
& \big|\overline{\partial_z}F(z)\big| \leq C_N\big|{\rm Im}\, z\big|^N(1 + |z|)^{s - 1 - N}, \quad \forall N \geq 0, \quad z \in {\mathbb C}, \\
& {\rm supp}\, F(z) \subset \{ z \in {\mathbb C} \, ; \, \big|{\rm Im}\, z\big| \leq 2 + 2\big|{\rm Re}\, z\big|\}.
\end{split}
\right.
\nonumber
\end{equation}
One can take $F(z) \in C_0^{\infty}({\mathbb C})$, if $f(x) \in C_0^{\infty}({\mathbb R})$.
\end{lemma}

\begin{lemma}
\label{C1S3HelfferSjostrandFormula}
If $f(x)$ satisfies (\ref{S2Lemma2.15f(x)Condition}) for some $s < 0$,
\begin{equation}
f(A) = \frac{1}{2\pi i}
\int_{{\mathbb C}}\overline{\partial_z}F(z) (z - A)^{-1}dzd\overline{z}
\label{HelfferSjostrandFormula}
\end{equation}
holds  for any self-adjoint operator $A$, where $F$ is an almost analytic extension of $f$.
\end{lemma}

For the proof, see \cite{DeGe}, p. 392.

This formula is suitable to deal with the perturbation $f(B) - f(A)$. For two self-adjoint operators $A$ and $B$, assume that $D(A) = D(B)$ and $(A- B)(i + A)^{-1}$ and $(A-B)(i + B)^{-1}$ are bounded operators.
Then, we have for any $\varphi \in C_0^{\infty}({\mathbb R})$
\begin{equation}
\varphi(B) - \varphi(A) = \frac{1}{2\pi i}\int_{{\mathbb C}}
\overline{\partial_z}\widetilde{\varphi}(z)\,(z - B)^{-1}(B-A)(z - A)^{-1}dzd\overline{z},
\label{varphiB-varphiA}
\end{equation}
where $\widetilde{\varphi} \in C_0^{\infty}({\mathbb C})$ is an almost analytic extension of $\varphi$. This is formally obvious by the resolvent equation. We show the convegence of the integral of the right-hand side. In fact, we have 
$$
\|(B-A)(z - A)^{-1}\| \leq \|(B-A)(i + A)^{-1}\|\|(i + A)(z - A)^{-1}\| \leq C\|(i + A)(z - A)^{-1}\|.
$$
By using the spectral decomposition of $A$, we have
$$
\|(i + A)(z - A)^{-1}\| \leq \sup_{\lambda\in{\mathbb R}}
\left|\frac{i + \lambda}{z - \lambda}\right| \leq 
C|{\rm Im}\, z|^{-1}(1 + |z|).
$$
Using the estimate
$$
\left|\overline{\partial_z}\widetilde{\varphi}(z)\right| \leq 
C|{\rm Im}\, z|^2(1 + |z|)^{-4},
$$
we see that the right-hand side of (\ref{varphiB-varphiA}) is a bounded operator.

Returning to our manifold $\mathcal M$, we pick up one end 
$(0,\infty)\times M$ equipped with the metric $ds^2 = (dr)^2 + \rho(r)^2h(r,x,dx)$. Let $h_M(x,dx)$ be the limit metric on $M$ satisfying
$$
h(r,x,dx) - h_M(x,dx) \in S^{-\gamma}
$$
for $\gamma > 0$. Let $\Lambda(r)$ be the Laplace-Beltrami operator on $M$ associated with the metric $h(r,x,dx)$. 

\begin{lemma}
\label{varphiLambdaderivative}
For any $\varphi \in C_0^{\infty}({\mathbb R})$, $\varphi(\Lambda(r))$ is strongly differentiable on $L^2(M)$ with respect to $r > 0$ and satisfies
\begin{equation}
\Big\|\big(\frac{d}{dr}\big)^{\ell}\varphi(\Lambda(r))\Big\| \leq C_n(1 + r)^{-\ell-\gamma}, \quad \forall \ell \geq 1.
\label{S2VarphiLambdaderivativedecay}
\end{equation}
\begin{equation}
\Big\|\Big[ \varphi(\Lambda(r)),\frac{g'}{g}\Big]\Big\| \leq 
C (1 + r)^{-1 - \gamma}, 
\label{S2varphLambdacommutatordecay}
\end{equation}
\end{lemma}

\begin{proof}
By (\ref{S2DLregularitytheorem}), $D(\Lambda(r))$ is independent of $r>0$. In view of (\ref{L=LA+LB+LS+LH}), 
 for any $u \in D(\Lambda(r))$, each term in 
$\Lambda(r)u$ is differentiable with respect to $r$. Let us check it for the most delicate term $L_Su$.  Letting $h(r,x) = s^{2(n-k-2)}h_0(r,x)$, we have $C \leq h_0 \leq C^{-1}$ and
$$
L_S = - \sqrt{h_0}\Big(\partial_s^2 + \frac{n-k-2}{s}\partial_s\Big) - \frac{\partial_s\sqrt{h_0}}{\sqrt{h_0}}\partial_s.
$$
Therefore, $L_Su$ is differentiable with respect to $r$. This proves that $(z - \Lambda(r))^{-1}$ is strongly differentiable with respect to $r > 0$. 
 Using (\ref{varphiB-varphiA}) and arguing as above, we obtain (\ref{S2VarphiLambdaderivativedecay}). 

 By virtue of (\ref{HelfferSjostrandFormula}), we have
\begin{equation}
\Big[\varphi(\Lambda(r)),\frac{g'}{g}\Big] = 
\frac{1}{2\pi i}\int_{{\mathbb C}}\overline{\partial_z}\widetilde
\varphi(z)
(z - \Lambda(r))^{-1}\Big[\Lambda(r),\frac{g'}{g}\Big]
(z - \Lambda(r))^{-1}dzd\overline{z}.
\nonumber
\end{equation}
Since $g'/g \in S^{-\gamma}$,  $\big[\Lambda(r),g'/g\big]$ is a 1st order differential operator with coefficients decaying like $r^{-1-\gamma}$. This proves (\ref{S2varphLambdacommutatordecay}). 
\end{proof}


\subsection{1-dimensional equation}
We summarize here basic facts about the 1-dimensional Helmholtz equation, since they elucidate the role of radiation condition and are also utilized in the spectral analysis of cusp. 


\begin{lemma}
\label{Lemma2.181dimHelmholtz}
Consider the equation on ${\mathbb R}^1$:
\begin{equation}
- u'' - k^2 u = f \quad {\it for} \quad - \infty < r < \infty.
\nonumber
\end{equation}
Assume that $k > 0$, $u = 0$ for $r < 1$, and $f \in L^1(0,\infty)$. \\
\noindent
(1) If
$\mathop{\underline{\lim}_{r\to\infty}}\big|u'(r) - ik u(r)\big| = 0$, 
we have
\begin{equation}
u(r) =  \frac{i}{2k} \int_0^{\infty}e^{ik|r-s|}f(s)ds,
\label{S21dimGreenfunction}
\end{equation}
\begin{equation}
\big|u'(r) - iku(r)\big| \leq \int_r^{\infty}|f(s)|ds.
\label{S21dimRadcond}
\end{equation}
(2) If furthermore $\mathop{\underline{\lim}_{r\to\infty}}|u(r)| = 0$, we have
\begin{equation}
|u(r)| \leq \frac{1}{k}\int_r^{\infty}|f(s)|ds,
\label{Subsection2.8Inequality2}
\end{equation}
\begin{equation}
\int_0^{\infty}(1 + r)^{2(s-1)}|u|^2dr  \leq C_{k,s}\int_0^{\infty}(1 + r)^{2s}|f(r)|^2dr,
\label{Subsection2.8Inequality2.5}
\end{equation}
for any $s > 1/2$, where the constant $C_{k,s}$ is independent of $k$ when $k$ varies over a compact interval in $(0,\infty)$.
\end{lemma}
\begin{proof} 
Let $g_0(k,r,s) = \frac{i}{2k}e^{ik|r-s|}$. For $0 < r < a$, by integration by parts
$$
- \int_0^au''(s)g_0(k,r,s)ds = u(r) + 
\frac{e^{ik(a-r)}}{2ik}\left(u'(a) - ik u(a)\right) 
+ k ^2\int_0^au(s)g_0(k,r,s)ds.
$$
Using the equation, we have
\begin{equation}
u(r) + \frac{e^{ik(a-r)}}{2ik}\left(u'(a) - iku(a)\right) 
= \frac{i}{2k}\int_0^ae^{ik|r-s|}f(s)ds.
\label{S2.8Helmholtz1}
\end{equation}
Letting $a \to \infty$, we obtain (\ref{S21dimGreenfunction}). 
This and (\ref{S2.8Helmholtz1}) yield
\begin{equation}
u'(a) - ik u(a) = e^{ik(r-a)}\int_a^{\infty}e^{ik|r-s|}f(s)ds,
\nonumber
\end{equation}
which implies (\ref{S21dimRadcond}). 
If $\mathop{\underline{\lim}_{r\to\infty}}|u(r)| = 0$,
we have $\displaystyle{\int_{0}^{\infty}e^{-iks}f(s)ds = 0}$, hence
\begin{equation}
u(r) = \int_{r}^{\infty}\frac{\sin k(r-s)}{k}f(s)ds,
\label{Subsection2.8Inequality3}
\end{equation}
which yields (\ref{Subsection2.8Inequality2}). 

Let us recall  well-known Hardy's inequality: {\it Let $h(r) \in L^1((0,\infty);dr)$ and put 
$\displaystyle{
w(r) = \int_r^{\infty}h(t)dt.}$ 
Then for $s > 1/2$,}
\begin{equation}
\int_0^{\infty}r^{2(s-1)}|w(r)|^2dr \leq \frac{4}{(2s-1)^2}\int_0^{\infty}
r^{2s}|h(r)|^2dr.
\label{Lemma7.4HardyInequality}
\end{equation}
(See e.g. \cite{IsKu10}, p. 106). This and (\ref{Subsection2.8Inequality3}) imply (\ref{Subsection2.8Inequality2.5}).
\end{proof}


\section{Transformation of the metric}
\label{Transformmetric}
We pick up one end $\mathcal M_i$ of $\mathcal M$, omit the subscript $i$ and consider the  perturbed warped product  metric
\begin{equation}
ds^2 = (dr)^2 + \rho(r)^2h(r,x,dx)
\label{S2meric}
\end{equation}
on $(0,\infty) \times M$
having the property
\begin{equation}
h(r,x,dx) - h_{M}(x,dx) \in S^{-\gamma}.
\label{S12hij-hMij}
\end{equation}
The perturbation term $h(r,x,dx) - h_{M}(x,dx)$ is said to be {\it short-range} if $\gamma > 1$ in (\ref{S12hij-hMij}), and 
{\it  long-range} if $0 < \gamma \leq 1$. 
In the study of spectral properties of the associated Laplace operator, the former can be dealt with by the standard perturbation technique, however the latter requires involved analysis. In this section,  given a metric
\begin{equation}
ds^2 = a(t,z)(dt)^2 + 2w(t)b_i(t,z)dtdz^i +
 w(t)^{2}c_{ij}(t,z)dz^idz^j,
\label{S72ds2}
\end{equation}
where $z = (z^1,\cdots,z^{n-1})$ are local coordinates on $M$, 
we seek the condition under which we can transform this metric into the perturbed warped product form (\ref{S2meric}).
The decay order $\gamma$ in (\ref{S12hij-hMij}) will be affected by the growth order of the volume of the manifold.   For the sake of simplicity, we consider the case in which $\mathcal M$ has no conic singularities hence the metric is $C^{\infty}$.
Let $\{U_j\}_{j\in J}$ be a finite open covering of $M$, each $U_j$ being diffeomorphic to a bounded open subset $V_j \subset {\mathbb R}^{n-1}$. 
 We assume that there exists an $r_0 > 0$ such that for $r > r_0$, $\mathcal M$ is covered by $\{(r_0,\infty)\times U_j\}_{j \in J}$. We take one of $U_j$,  assume that $z \in V_j$ and omit the subscript $j$.

Assume that
\begin{equation}
w(t)^{-1} \in S^{-\kappa},
\label{Condw(t)}
\end{equation}
\begin{equation}
a(t,z) -1 \in S^{-\lambda},
\label{Condatz-1}
\end{equation}
\begin{equation}
b_i(t,z) \in S^{-\mu},
\label{condbitz}
\end{equation}
\begin{equation}
c_{ij}(t,z) - h_{ij}(z) \in S^{-\nu},
\label{S7cij-hij}
\end{equation}
where  $\kappa, \lambda, \mu, \nu $ are constants such that
\begin{equation}
\kappa > 1/2, \quad \lambda > 1, \quad \mu > 0,   \quad \nu > 0, \quad \kappa + \mu > 1,
\label{S72DecayCond}
\end{equation}
$h_{ij}(z)dz^idz^j$ is a $C^{\infty}$ metric on $M$, and $S^{\kappa}$ is defined as in (\ref{IntoDefineSkappa}). 
Letting $y = (t,z)$, we rewrite (\ref{S72ds2})  as
$$
ds^2 = g_{ij}dy^idy^j, 
\quad 
\big(g_{ij}\big) = 
\left(
\begin{array}{cc}
1 & 0\\
0 & w
\end{array}
\right)
\left(
\begin{array}{cc}
a & \ ^t b\\
b & c
\end{array}
\right)
\left(
\begin{array}{cc}
1 & 0\\
0 & w
\end{array}
\right),
$$
where $^t b = (b_1,\cdots,b_{n-1})$, 
$c = \left( c_{ij}\right)$. 
Therefore, its inverse is written as
$$
\left(g^{ij}\right) = 
\left(
\begin{array}{cc}
1 & 0\\
0 & w^{-1}
\end{array}
\right)
\left(
\begin{array}{cc}
{\widetilde a} & \ ^t{\widetilde b}\\
{\widetilde b} & {\widetilde c}
\end{array}
\right)
\left(
\begin{array}{cc}
1 & 0\\
0 & w^{-1}
\end{array}
\right),
$$
where
\begin{equation}
\widetilde a-1 \in S^{-\lambda}, \quad \widetilde b \in S^{-\mu}, \quad \widetilde c -\widetilde  h\in S^{-\nu},
\quad
\widetilde h = \big(\widetilde h^{ij}\big) = \big({\widetilde h}_{ij}\big)^{-1}.
\label{a-1bc-hdecay}
\end{equation}
Define the classical Hamiltonian $H$ by
\begin{equation}
H(t,z,\tau,\zeta) = \widetilde a\tau^2 + 2w^{-1}\widetilde b^i\tau\zeta_i + w^{-2}\widetilde c^{ij}\zeta_i\zeta_j.
\nonumber
\end{equation}
We rewrite it as
\begin{equation}
H(t,z,\tau,\zeta) = \tau^2 + K(t,z,\tau,\zeta), 
\nonumber
\end{equation}
\begin{equation}
K = (\tilde a-1)\tau^2 + 2w^{-1}\tilde b^i\tau\zeta_i + w^{-2}\tilde c^{ij}\zeta_i\zeta_j.
\nonumber
\end{equation}
By (\ref{a-1bc-hdecay}), we have as a function of $t, z$
\begin{equation}
K \in S^{-m}, \quad
m = \min\{\lambda, \kappa + \mu, 2\kappa\}.
\label{S7EstimateK}
\end{equation}
The assumption (\ref{S72DecayCond}) implies
\begin{equation}
m = 1 + \epsilon_0, \quad \epsilon_0 >0.
\label{m=1+epsilon0}
\end{equation}
Letting $\ \dot {} = \dfrac{d}{dr}$, we solve the Hamilton equation
\begin{equation}
\left\{
\begin{split}
& \dot t = \frac{\partial H}{\partial\tau}, \quad \dot z = \frac{\partial H}{\partial \zeta}, \\
& \dot \tau = - \frac{\partial H}{\partial t}, \quad \dot \zeta = - \frac{\partial H}{\partial z},
\end{split}
\right.
\label{S7HamiltonEq}
\end{equation}
with the condition at infinity
\begin{equation}
t(r) - r \to 0,  \quad z(r) \to x, \quad \tau(r) \to 1/2, \quad 
\zeta(r) \to 0,
\label{S7HamiltonEqBCInfty}
\end{equation}
where $x \in V$. 
 Note that, since $\frac{\partial H}{\partial \tau} \to 2\tau$ as $t \to \infty$, the condition $\tau(r) \to 1/2$ is compatible with $t(r) -r \to 0$. 
Then, $\tau$, $z$ and $\zeta$ should satisfy the integral equations
\begin{equation}
\begin{split}
2\tau - 1 &= 2 \int_r^{\infty}\frac{\partial H}{\partial t}dr', \\
z - x &= - \int_r^{\infty}\frac{\partial H}{\partial \zeta}dr', \\
\zeta &= \int_r^{\infty}\frac{\partial H}{\partial z}dr',
\end{split}
\label{IntegralEqfortauzzeta}
\end{equation}
and $t$ should satisfy the integro-differential equation
$$
\frac{dt}{dr} = 1 + 2\int_r^{\infty}\frac{\partial H}{\partial t}dr' + \frac{\partial K}{\partial \tau}.
$$
Letting
\begin{equation}
\begin{split}
X(r) &= X(r,x) = \big(t(r)-r,z(r)-x,\tau(r)-1/2,\zeta(r)\big), \\
\widetilde X(r) &= \widetilde X(r,x) = \big(t(r),z(r),\tau(r),\zeta(r)\big),\\
X_{\infty}(r)  &= X_{\infty}(r,x) = (r,x,1/2,0),
\end{split}
\nonumber
\end{equation}
(hence $\widetilde X(r) = X(r) + X_{\infty}$), 
 we put
\begin{equation}
\begin{split}
U_0(X(r)) &= 2 \int_r^{\infty}\frac{\partial H}{\partial t}(\widetilde X(r'))dr' + \frac{\partial K}{\partial\tau}(\widetilde X(r)),\\
U(X(r)) &= \int_r^{\infty}
\left(- U_0(\widetilde X(r')), - \frac{\partial H}{\partial \zeta}(\widetilde X(r')), 
\frac{\partial H}{\partial t}(\widetilde X(r')), \frac{\partial H}{\partial z}(\widetilde X(r'))\right)dr'.
\end{split}
\nonumber
\end{equation}
Then, $t$ should satisfy
\begin{equation}
t- r = - \int_{r}^{\infty}U_0(X(r'))dr'.
\label{Intgralequationfort}
\end{equation}
Note that
$$
\int_r^{\infty}U_0(X(r'))dr' = \int_r^{\infty}\left(2(r'-r)\frac{\partial K}{\partial t}(\widetilde X(r')) + \frac{\partial K}{\partial\tau}(\widetilde X(r'))\right)dr'.
$$
By (\ref{IntegralEqfortauzzeta}) and (\ref{Intgralequationfort}),  
the differential equation (\ref{S7HamiltonEq}) with the condition at infinity (\ref{S7HamiltonEqBCInfty}) is converted to the integral equation
\begin{equation}
X(r,x) = U(X_{\infty}(r,x) + X(r,x)).
\label{S7NonlinearEq}
\end{equation} 
Define the norm $\|\cdot\|$ by
$$
\|X\| = \sup_{(r,x) \in (r_0,\infty)\times V}|X(r,x)|.
$$
The conditions (\ref{S7EstimateK}) and (\ref{m=1+epsilon0}) imply that for  a sufficiently small $\epsilon, r_0^{-1}>0$, $U$ is a contraction mapping in the ball 
$$
B_{\epsilon,r_0} = \{X \in C((r_0,\infty)\times V)\, ; \, \|X\| \leq \epsilon\}.
$$
We then have:


\begin{lemma}\label{S7ExistenceX(r)}
There exists a unique solution $X(r,x)$ of the equation (\ref{S7NonlinearEq}). It satisfies
\begin{equation}
|\partial_r^{\ell}\partial_x^{\alpha} X(r,x)| \leq C_{\ell\alpha}r^{-\ell-\epsilon_0}, \quad \forall \ell, \ \alpha.
\label{S7X(r)estimate}
\end{equation}
Moreover, the differential of the map $(t,z) \to (r,x)$ is $I + O(r^{-\epsilon_0})$. 
\end{lemma}


\begin{lemma}\label{S7closedformlemma} As a 2-form with respect to $r, x$, we have
$$
 d\tau\wedge dt + \sum_{i=1}^{n-1}d\zeta_i\wedge dz^i=0. 
$$
\end{lemma}

\begin{proof}
Let $y=(t,z)$, $\eta=(\tau,\zeta)$ and $\theta = (r,x)$. Then
$$
d\tau\wedge dt + \sum_{i=1}^{n-1}d\zeta_i\wedge dz^i = \sum_{j<k}\left[\eta,y\right]_{jk}d\theta^j\wedge d\theta^k,
$$
$$
\left[\eta,y\right]_{jk}= \frac{\partial \eta}{\partial \theta^j}\cdot\frac{\partial y}{\partial \theta^k} - 
\frac{\partial\eta}{\partial \theta^k}\cdot\frac{\partial y}{\partial \theta^j}. 
$$
Noting that
$$
\frac{\partial}{\partial r}\left(\frac{\partial \eta}{\partial \theta^j}\cdot\frac{\partial y}{\partial \theta^k}\right) = - \frac{\partial^2H}{\partial y^i\partial y^m}\frac{\partial y^m}{\partial \theta^j}\frac{\partial y^i}{\partial\theta^k} + 
\frac{\partial^2H}{\partial \eta_i\partial \eta_m}\frac{\partial \eta_i}{\partial \theta^k}\frac{\partial\eta_m}{\partial\theta^j}
$$
is symmetric with respect to $j$ and $k$, we have
$
\frac{\partial}{\partial r}\left[\eta,y\right]_{jk} =0.
$
Lemma \ref{S7ExistenceX(r)} implies $\left[\eta,y\right]_{jk} \to 0$ as $r\to\infty$. Hence $\left[\eta,y\right]_{jk} =0$, which proves Lemma \ref{S7closedformlemma}. 
\end{proof}

\medskip
By Lemma \ref{S7ExistenceX(r)}, the map $(r,x) \to (t,z)$ is a global diffeomorphism on $(r_0,\infty)\times M$. 
We invert it to get 
$r=r(t,z)$, $x=x(t,z)$, $\tau = \tau(t,z)$, $\zeta = \zeta(t,z)$.
 Lemma \ref{S7closedformlemma} implies
$$
\frac{\partial \zeta_j}{\partial z^k}=\frac{\partial\zeta_k}{\partial z^j}, \quad 
\frac{\partial \zeta_j}{\partial t}=\frac{\partial\tau}{\partial z^j}, \quad
1 \leq j, k \leq n-1.
$$
Recalling that
\begin{equation}
\tau(r) - \frac{1}{2}=\int_r^{\infty}\frac{\partial H}{\partial t}(\widetilde X(r'))dr' = O(r^{-1-\epsilon_0}),
\label{S7taur-1/2estimate}
\end{equation}
we define
$$
\varphi(t,z) = \frac{t}{2}- \int_t^{\infty}\Big(\tau(t',z)-\frac{1}{2}\Big)dt'.
$$


\begin{lemma}\label{S7psivarphilemma} There exists $t_0>0$ such that for $t > t_0$, \\
(1)  $\ \partial_t\varphi(t,z) = \tau(t,z)$, \\
(2)  $\ \partial _z\varphi(t,z) = \zeta(t,z)$, \\
(3) $\ H(t,z,\partial_t\varphi(t,z),\partial_z\varphi(t,z)) = 1/4$, \\
(4) $\ |\partial_t^{\ell}\partial_z^{\alpha}\big(\varphi(t,z) - t/2\big)| \leq C_{\ell\alpha}t^{-\epsilon_0 - \ell}, \quad 
\forall \ell, \alpha, $ \\
(5) $\ \varphi(t,z) = r(t,z)/2$.
\end{lemma}

\begin{proof}
The assertion (1) is obvious, and (2) follows from
\begin{equation}
\begin{split}
\frac{\partial \varphi}{\partial z^j} =- \int_t^{\infty}\frac{\partial\tau}{\partial z^j}dt'=- \int_t^{\infty}\frac{\partial\zeta_j}{\partial t'}dt' = \zeta_j(t,z).
\end{split}
\nonumber
\end{equation}
Since the energy is conserved, $H(t,z,\tau,\zeta)$ is constant along the orbit, which turns out to be $1/4$ by letting $r\to \infty$. The assertion (4) follows from (\ref{S7taur-1/2estimate}) and (\ref{S7EstimateK}). Using (1), (2) and (3), we have
\begin{equation}
\begin{split}
\frac{\partial\varphi}{\partial r} & = \frac{\partial \varphi}{\partial t}\frac{\partial t}{\partial r} +  \frac{\partial\varphi}{\partial z}\cdot\frac{\partial z}{\partial r} = \tau\frac{\partial t}{\partial r} + \zeta\frac{\partial z}{\partial r} \\
&= \tau\frac{\partial H}{\partial \tau} + \zeta\cdot\frac{\partial H}{\partial \zeta} = 2H = 1/2.
\end{split}
\nonumber
\end{equation}
Therefore, $\varphi(t,z) - r/2$ is independent of $r$. On the other hand $\varphi(t,z) - r/2 \to 0$ as $r\to\infty$. This proves $\varphi=r/2$. 
\end{proof}

\medskip
The diffeomorphism $(r,x) \to (t,z)$ induces  $r$-dependent local coordinates
$z = z(r,x)$ on $M$. As $r \to \infty$, they converge to local coordinates $z(x)$ on $M$.


\begin{theorem}\label{Metrictransformed}
Assume (\ref{Condw(t)}) $\sim$ (\ref{S72DecayCond}), and $h_{ij}$ be as in (\ref{S7cij-hij}). Then,
in the coordinate system $(r,x)$, the Riemannian metric (\ref{S72ds2}) is written as
\begin{equation}
ds^2 = (dr)^2 + w(r)\overline{h}(r,x,dx)
\label{Stds2rewritten}
\end{equation}
where $\overline h(r,x,dx)$ is a Riemmanian metric on $M$ and satisfies
\begin{equation}
\overline h_{ij}(r,x) - h_{ij}(z(x)) \in S^{-\min\{\nu,\epsilon_0\}}.
\end{equation}

\end{theorem}

\begin{proof}
We put $y = (t,z)$, $\overline{y} = (r,x)$. Then, the Hamiltonain is written as
$$
H = g^{ij}(y)\eta_i\eta_j = \overline{g}^{ij}(\overline{y})\overline{\eta}_i\overline{\eta}_j,
$$
where $\eta = (\tau,\zeta)$. 
Using Lemma \ref{S7psivarphilemma}, we have
$$
\overline{g}^{00} = g^{ij}\frac{\partial{\overline y}^0}{\partial y^i}\frac{\partial{\overline y}^0}{\partial y^j} = g^{ij}\frac{\partial r}{\partial y^i}\frac{\partial r}{\partial y^j} = 4g^{ij}\frac{\partial \varphi}{\partial y^i}\frac{\partial \varphi}{\partial y^j} = 4H=1,
$$
$$
\overline{g}^{0k} = g^{ij}\frac{\partial{\overline y}^0}{\partial y^i}\frac{\partial{\overline y}^k}{\partial y^j} = g^{ij}\frac{\partial r}{\partial y^i}\frac{\partial x^k}{\partial y^j} = 0, \quad 1 \leq k \leq n-1,
$$
Here, in the 2nd line, we have used
$$
 0 = \frac{\partial x^k}{\partial r} = \frac{\partial x^k}{\partial y^j}\frac{\partial y^j}{\partial r} = \frac{\partial x^k}{\partial y^j}
 g^{ij}\eta_i= \frac{1}{2}\frac{\partial x^k}{\partial y^j}
 g^{ij}\frac{\partial r}{\partial y^i}.
 $$
Therefore, the Riemannain metric has the form
$$
ds^2 = (dr)^2 + \sum_{1\leq i, j \leq n-1}\overline{g}_{ij}dx^idx^j.
$$
We observe the matrix $\big(\overline{g}_{ij} \big)_{1 \leq i, j \leq n-1}$ with
$$
\overline{g}_{ij} = g_{00}\frac{\partial t}{\partial{\overline y}^i}\frac{\partial t}{\partial {\overline y}^j} + 2g_{0k}\frac{\partial t}{\partial{\overline y}^i}\frac{\partial z^k}{\partial {\overline y}^j} + 
g_{k\ell}\frac{\partial z^k}{\partial{\overline y}^i}\frac{\partial z^{\ell}}{\partial {\overline y}^j}.
$$
Since $g_{00} = O(1)$, $g_{0k} = O(w)$ and $g_{k\ell} = w^{2}c_{k\ell}$, we have
\begin{equation}
\overline{g}_{ij} = w^{2}\left(c_{ij} + O(r^{-\epsilon_0})\right),
\nonumber
\end{equation}
from which we can derive the lemma. 
\end{proof}

\medskip
 The above theorem says that the metric of the form (\ref{S72ds2}) can be transformed to (\ref{Stds2rewritten}) if
 $$
 w(t) \sim \exp\Big(c_0t + \frac{\beta}{1-\alpha}t^{1-\alpha}\Big), \quad {\rm or}\quad t^{\beta} \quad {\rm with} \quad \beta>1/2.
 $$
Let us consider the case $w(t) = t^{\beta}$ for large $t$. 
Then, (\ref{Stds2rewritten}) is a short-range perturbation of $(dr)^2 + w(r)^2h_M$ when $\min\{\nu,\epsilon_0\} > 1$, and long-range perturbation when $\min\{\nu,\epsilon_0\} \leq 1$. Since $\kappa = \beta$, (\ref{Stds2rewritten}) is a short-range perturbation of $(dr)^2 + w(r)^2h_M$ only when 
 $\beta>1$. 
  For $1/2 < \beta \leq 1$, it is a long-range perturbation. The border-line case appears when the metric is Euclidean.
  This fact is pointed out by Bouclet \cite{Bouclet}. He also mentions that there appears a conformal factor in front of $h_{ij}(z(x))$ in Theorem \ref{Metrictransformed}. This is because he solves the Hamilton equation as an initial value problem with data on a surface $\{t = t_0\}$. If, as has been done above, we solve it as the Cauchy problem from infinity, the conformal factor does not appear.


\section{Rellich-Vekua type theorem}\label{SectionRellichVekua}


\subsection{Volume growth condition}
The aim of this section is to derive a decay rate of solutions to the Helmholtz equation near infinity of an end $(0,\infty)\times M$, which is crucial to study the spectral theory  for the Laplacian, in particular for the discreteness or the non-existence of eigenvalues embedded in the continuous spectrum. This property is closely related to the volume growth of the manifold at infinity. We consider only the growing metric satisfying the following assumption:

\medskip
\noindent
{\it (VG)
 There exist a non-negative constant $c_0$ and positive constants $\alpha_0, \beta_0, \gamma_0, r_0$ such that}
\begin{equation}
\left\{
\begin{split}
& \frac{\rho'}{\rho}-c_0 \in S^{-\alpha_0},  \\
& \frac{\rho'}{\rho} \geq \frac{\beta_0}{r}, \quad r > r_0,\\
& h^{ij}(r,x) - h^{ij}_M(x) \in S^{-\gamma_0}.
\end{split}
\right.
\label{S7CondC}
\end{equation}

\medskip
Note that (A-2) and (A-3) imply (VG), and that (VG) yields
\begin{equation}
\rho(r) \geq \rho(r_0)\Big(\frac{r}{r_0}\Big)^{\beta_0}, \quad r > r_0.
\nonumber
\end{equation}
 Letting $h = \det h(r,x,dx)$, we have $g = \rho^{2(n-1)}h$, and the assumption on $h^{ij}$ implies $h' = O(r^{-1-\gamma_0})$. Therefore, we can reformulate (\ref{S7CondC}) in terms of $g$:

\medskip
\noindent
{\it (VG)' There exist a non-negative constant $c_0$ and positive constants $\alpha_0, \beta_0, r_0$  such that}
\begin{equation}
\left\{
\begin{split}
& \frac{g'}{4g}-\frac{(n-1)c_0}{2} \in S^{-\alpha_0}, \\
& \frac{g'}{4g} \geq \frac{(n-1)\beta_0}{2r}, \quad r>r_0,\\
& h^{ij}(r,x) - h^{ij}_M(x) \in S^{-\gamma_0}.
\end{split}
\right.
\label{S7CondC'}
\end{equation}

\medskip
Let $S(r) = \{r\} \times M$ and $dS(r)$ be the surface element of $S(r)$ induced from the metric $ds^2$.  The volume element $dV$ of the end $(0,\infty)\times M$ is then written as
$$
dV = drdS(r).
$$
Let $\|\cdot\|$ and $(\;,\;)$ be the norm and the inner product of $L^2(M)$. 
 We say that $u$ is locally in $D(-\Delta_{\mathcal M})$ if $\chi(r)u \in D(-\Delta_{\mathcal M})$ for any $\chi(r) \in C_0^{\infty}\big((0,\infty)\big)$. Then the following theorem holds.


\begin{theorem}\label{S2Rapiddecay} 
Assume  (VG) with 
\begin{equation}
  \alpha_0 > 0, \quad \beta_0 > 0,  \quad \gamma_0 > 0.
\label{Theorem4.1Cond}
\end{equation}
Let
$E_0 =\big(
(n-1)c_0/2\big)^2$. 
Suppose $u$ is locally in $D(-\Delta_{\mathcal M})$ and satisfies 
$$
(- \Delta_{\mathcal M} - E)u = 0, \quad {for } \quad r > R
$$
for some constants $E > E_0$ and $R > 0$. 
If $u$ satisfies
$$
\liminf_{r\to\infty}\, r^{\gamma}\int_{S(r)}\Big(|u'|^2 +\big|u\big|^2\Big)dS(r)= 0
$$
for a constant $\gamma >0$, then 
\begin{equation}
\int_{r>R}r^m\left(\|u'\|^2 + (Bu,u) + \|u\|^2\right)dV < \infty
\end{equation}
for any $m > 0$, where the operator $B = B(r)$ is defined in (\ref{S7Laplacian1}).
\end{theorem}

As will be shown later, this theorem yields the discreteness of embedded eigenvalues in the continuous spectrum for metrics having growth order $\beta > 0$. Much more significant is the following Rellich-Vekua type theorem which proves the non-existence of embedded eigenvalues and also plays an important role in the inverse problem, however with the trade-off of losing slowly growing metrics of  order  $0 < \beta \leq 1/3$.


\begin{theorem}\label{S2RellichTh1} 
Assume  (VG) with 
\begin{equation}
   \alpha_0 > 0, \quad \beta_0 > 1/3, \quad \gamma_0 > 0. \\
\label{Theorem4.2Cond}
\end{equation}
Let
$E_0 =\big(
(n-1)c_0/2\big)^2$. 
Suppose $u$ is locally in $D(-\Delta_{\mathcal M})$ and satisfies 
$$
(- \Delta_{\mathcal M} - E)u = 0, \quad {for } \quad r > R
$$
for some constants $E > E_0$ and $R > 0$. 
If $u$ satisfies
$$
\liminf_{r\to\infty}\, \int_{S(r)}\Big(|u'|^2 +\big|u\big|^2\Big)dS(r)= 0,
$$
 there exists a constant $R_1 > 0$ such that $u = 0$ for $r > R_1$.
\end{theorem}

 For the other results on this type of theorem on non-compact manifolds,  see Kumura \cite{Kum13}, \cite{Kum13(1)}, Ito-Skibsted \cite{ItoSkiI} and the references therein.
 Note that these works do not deal with the manifold with conic singularities. 
 Below, we derive this theorem from growth properties of solutions to an abstract ordinary differential equation with operator-valued coefficients. 
 

\subsection{Abstract differential equations}
\label{Subsec4.2}

 Let $X$ be a Hilbert space and consider the following differential equation for an $X$-valued function $u(t)$: 
\begin{equation}
- u^{\prime\prime}(t) + B(t)u(t) + V(t)u(t) - \lambda u(t) = 0, \quad t > 0,
\label{eq:Chap2Sect3Diffeq}
\end{equation}
$\lambda > 0$ being a constant. We assume the following conditions (B) and (V):

\medskip
\noindent
  {\it (B) \ For each $t>0$, $B(t)$ is a non-negative self-adjoint operator on $X$  having the following properties:
\begin{itemize}
\item the domain $D = D(B(t)))$ of $B(t)$ is independent of $t > 0$, 
\item the domain $D_1 = D(\sqrt{B(t)})$ of $\sqrt{B(t)}$ is independent of $t >0$,
\item for any $x \in D$, $(B(t)x,x)$ is differentiable with respect to $t > 0$, 
\item there exists a closed form $b_t(\cdot,\cdot)$ on $D_1$ such that
\begin{equation}
\frac{d}{dt}(B(t)x,x) = b_t(x,x), \quad \forall x \in D, \quad \forall t > 0,
\nonumber
\end{equation}
\item there exist positive constants $\delta, C, \epsilon_0$ such that
\begin{equation}
tb_t(x,x) \leq - \delta (B(t)x,x) + Ct^{-\epsilon_0}(x,x), \quad 
\forall x \in D, \quad \forall t > 0.
\label{tbt(x,x)leq-delatB(t)+Ct-psilon0}
\end{equation}
\end{itemize}
We denote this quadratic form $b_t(x,y)$  by
$\big(\dfrac{d}{dt}B(t)x,y\big)$ or $\big(B'(t)x,y\big)$ and write (\ref{tbt(x,x)leq-delatB(t)+Ct-psilon0}) as}
\begin{equation}
t\frac{dB(t)}{dt} \leq - \delta B(t) + Ct^{-\epsilon_0}, \quad \forall t > 0.
\label{S8C3tdBtdt}
\end{equation}
{\it (V)\  For each $t > 0$,  $V(t)$ is bounded self-adjoint on $X$, 
$V(t) \in C^1\big((0,\infty);{\bf B}(X)\big)$,  and satisfies
\begin{equation}
\|V(t)\| + t\big\|\frac{dV(t)}{dt}\big\|\leq C\, t^{- \epsilon_0},
\quad \forall t \geq 1,
\label{S7C-4normVt}
\end{equation}
for a constants $C> 0$, where $\|\cdot\|$ is the norm in ${\bf B}(X)$.}

\medskip
By the solution $u(t)$ of (\ref{eq:Chap2Sect3Diffeq}), we mean that $u(t) \in C^1((0,\infty);X)$, $u''(t) \in X$ exists almost everywhere on $(0,\infty)$, $u''(t) \in L^1_{loc}((0,\infty);X)$ and 
$$
u'(t) = u'(t_0) + \int_{t_0}^tu''(s)ds
$$
holds for all $t >t_0 >  0$. Moreover, $u(t) \in D$ for any $t > 0$, and the equation (\ref{eq:Chap2Sect3Diffeq}) is satisfied on $(0,\infty)$. 


\begin{theorem}\label{Growthproperty}
Assume (B),  (V)  and for some $\gamma > 0$
$$
\liminf_{t\to\infty}\, t^{\gamma}(\|u'(t)\| + \|u(t)\|) = 0.
$$
Then 
$$
\int_0^{\infty}t^m\left(\|u'\|^2 + (Bu,u) + \|u\|^2\right)dt < \infty
$$
for all $m > 0$, where $B = B(t)$ is defined in (\ref{S7Laplacian1}).
\end{theorem}
 
The Rellich type theorem is proven by restricting the range of $\delta$.


\begin{theorem}\label{AbstractRllich1}
Assume (B) with $\delta > 2/3$,  (V)  and
$$
\liminf_{t\to\infty}(\|u'(t)\| + \|u(t)\|) = 0.
$$
Then there exists $t_2 > 0$ such that $u(t) = 0$,  for $t > t_2$.
\end{theorem}

\medskip
 Let us derive Theorems \ref{S2Rapiddecay} and \ref{S2RellichTh1} from Theorems 
\ref{Growthproperty} and \ref{AbstractRllich1}. Let  $g = \rho^{2(n-1)} h$, where
 \begin{equation}
 h = h(r,x) = \det\big(h_{ij}(r,x)\big), \quad h_M = h_M(x) = \det\big(h_{ij,M}(x)\big).
 \label{S7Definehrx}
 \end{equation}
 Then, the Laplacian is rewritten as
\begin{equation}
- \Delta_{\mathcal M} =- \partial_r^2 - \frac{\partial_r g}{2g}\partial_r + B(r), \quad B(r) = -
\frac{1}{\sqrt{g}}\partial_{x_i}\Big(\rho(r)^{-2}\sqrt{g}h^{ij}(r,x)\partial_{x_j}\Big).
\label{S7Laplacian1}
\end{equation} 
Letting $' = \partial_r$ and noting that
$$
g^{1/4}\Big(\partial_r^2 + \frac{g'}{2g}\partial_r\Big)g^{-1/4} = 
\partial_r^2 - \Big(\frac{g'}{4g}\Big)^2 -\Big(\frac{g'}{4g}\Big)', 
$$
we transform $\Delta_{\mathcal M}$ as
\begin{equation}
\begin{split}
- \Big(\frac{g}{h_M}\Big)^{1/4}\Delta_{\mathcal M}\Big(\frac{g}{h_M}\Big)^{-1/4}= -\partial_r^2 + B_0(r)+ C(r),
\end{split}
\label{S7Deltagtranasformed}
\end{equation}
\begin{equation}
B_0(r) = \Big(\frac{g}{h_M}\Big)^{1/4} B(r)\Big(\frac{g}{h_M}\Big)^{-1/4}, 
\label{S7B(r)}
\end{equation}
\begin{equation}
C(r) =  \Big(\frac{g'}{4g}\Big)^2 + \Big(\frac{g'}{4g}\Big)'.
\label{S7V(r)}
\end{equation}
Let $X = L^2(M)$ equipped with the inner product
$$
(u,v)_{X} = \int_Mu(x)\overline{v(x)}\sqrt{h_M(x)}dx.
$$
Then, the transformation
$$
L^2((0,\infty)\times M) \ni u \to v = \Big(\frac{g}{h_M}\Big)^{1/4}u \in L^2((0,\infty);X;dr)
$$
is unitary, 
$B_0(r)$ is self-adjoint on $X$ with domain characterized in Theorem \ref{S2DLregularitytheorem}. We put 
$$
V(t) = C(t) - E_0, \quad 
\lambda = E - E_0.
$$
Then, by virtue of Theorem \ref{S2DLregularitytheorem}, the solution of the equation $- \Delta_{\mathcal M}u = Eu$ is transformed to an $X$-valued solution of the ordinary differential equation 
$$
(- \partial_t^2 + B_0(t) + V(t)\big)v=\lambda v, \quad v = \big(g/h_M\big)^{1/4}u.
$$

We show that this differential equation satisfies the assumptions (B), (V). Note that by (\ref{S7CondC'})
\begin{equation}
\frac{g'}{4g} - \sqrt{E_0} \in S^{-\alpha_0},
\label{S4fracg'4g-sqrtE0=S-epsilon0}
\end{equation}
which implies
$$
\|V(t)\| \leq Ct^{-\alpha_0}, \quad \|\partial_tV(t)\| \leq Ct^{-1-\alpha_0}.
$$
By Theorem \ref{S2DLregularitytheorem}, $D(B_0(t))$ is independent of $t > 0$, and by (\ref{H1M=D(q)}), $D(\sqrt{B_0(t)}) = H^1(M)$. Passing to the quadratic form, 
we see that $(B_0(t)x,x)$ is differentiable with respect to $t > 0$. 


\begin{lemma}\label{LemmatB(t)leqct-epsilon}
Let $0 < \delta < 2\beta_0$. Then, there exist positive  constants $C, t_0$ such that 
$$
t\frac{d B_0(t)}{dt} + \delta B_0(t) \leq Ct^{-\gamma_0}, \quad t > t_0.
$$
\end{lemma}

\begin{proof}
Rewrite the left-hand side as
$$
t\frac{dB_0(t)}{dt} + \delta B_0(t) = t\frac{dB_0(t)}{dt}+ 2\beta_0B_0(t) - (2\beta_0 -\delta)B_0(t).
$$
Since $B_0(r) = \big(h/h_M\big)^{1/4}B(r)\big(h/h_M\big)^{-1/4}$, we have
for $u \in D(B_0(t))$ 
$$
(B_0(t)u,u) =  \sum_{i,j}\big(\frac{h^{ij}}{\rho(t)^2}\partial_{x_j}u,\partial_{x_i}u\big) + 
\sum_{|\alpha| \leq 1}\big(b_{1,\alpha}\partial_x^{\alpha}u,\partial_x^{\alpha}u\big),
$$
where $b_{1,\alpha}$ behaves like
\begin{equation}
b_{1,\alpha} = O(t^{-1-\gamma_0}\rho(t)^{-2}).
\nonumber
\end{equation}
Therefore, 
$$
(\frac{dB_0(t)}{dt}u,u) =- \sum_{i,j}\big(c^{ij}(t,x)\partial_{x_j}u,\partial_{x_i}u\big) + \cdots, \quad c^{ij}(t,x) = \partial_t\Big(\frac{h^{ij}}{\rho(t)^2}\Big).
$$
By (\ref{S7CondC}), 
$$
c^{ij}(t,x) = -2\Big(\frac{h^{ij}}{\rho(t)^2}\Big)\frac{\rho'(t)}{\rho(t)} + 
O(t^{-1-\gamma_0}\rho(t)^{-2}).
$$
We then have, again using (\ref{S7CondC}),
$$
\big(c^{ij}\big) \leq - \frac{2\beta_0}{t}\Big(\frac{h^{ij}}{\rho(t)^2}\Big) + 
O(t^{-1-\gamma_0}\rho(t)^{-2}).
$$
Therefore, we have in the sense of quadratic form 
$$
t\frac{dB_0(t)}{dt}+ 2\beta_0B_0(t) \leq B_1(t),
$$
where $B_1(t)$ is a 2nd order differential operator on $X$ whose coefficients decay like $t^{-\gamma_0}\rho(t)^{-2}$. By the positive definiteness of the matrix $\big(h^{ij}\big)$, we have
$$
B_1(t) \leq (2\beta_0 - \delta)B_0(t)  + Ct^{-\gamma_0}
$$
for a constant $C > 0$, which proves the lemma. 
\end{proof}

Therefore, (B) and (V) are justified. 


\subsection{Proof of Theorem \ref{Growthproperty}} 

The following method based on integration by parts is essentially due to Eidus \cite{Eidus}.
Let $(\cdot,\cdot),  \|\cdot\|$ be the inner product and the norm of $X$, respectively.
For $0 < a < b < \infty$, we put
\begin{eqnarray*}
J_{(a,b)}(u,v) &=& \int_a^b(u(t),v(t))\,dt, \\
S_t(u,v)&=& (u(t),v(t)), \\
S(u,v)\Big|_a^b &= &S_b(u,v) - S_a(u,v).
\end{eqnarray*}
Note the formula
$$
S(u,v)\Big|_a^b = J_{(a,b)}(u',v) + J_{(a,b)}(u,v').
$$

Take a real-valued $C^{\infty}$-function $d(t)$ on $(0,\infty)$ and put 
$v(t) = e^{d(t)}u(t)$, where $u$ is a solution to (\ref{eq:Chap2Sect3Diffeq}). Then, it satisfies
\begin{equation}
- v''(t) + 2d'(t)v'(t) + (B(t) + q(t))v(t)=0, 
\label{S11Equationforv}
\end{equation}
$$
q(t) = V(t) + d''(t) - d'(t)^2 - \lambda.
$$
In the following arguments, we compute under the additional asumption that $(v'(t),v(t))$, etc.  are real-valued, for the sake of simplicity. 
This is not essential at all. In fact, we have only to take the real part in the follwing formulas. Or, in the practial applications, say $X=L^2({\bf h})$, we have only to take the real part of $v$. 

We start from the following identities.


\begin{lemma}\label{S9IntIdentities}
 Let $v(t) = e^{d(t)}u(t)$, and $\psi = \psi(t)$  a positive $C^{\infty}$-function on $(0,\infty)$. Then we have  for $0<a<b$
\begin{equation}
\begin{split}
&  S\big(\psi v',v'\big)\Big|_a^b - S\big(\psi (B+q)v,v\big)\Big|_a^b \\
&	   =   J_{(a,b)}\big((\psi' + 4d'\psi)v',v'\big) - J_{(a,b)}\big((\psi (B+q))'v,v\big), 
\label{S11Identity1} 
\end{split}
\end{equation}
\begin{equation}
\begin{split}
& 
 S\big(\psi v',v\big)\Big|_a^b -S\big(\psi d'v,v\big)\Big|_a^b - 
\frac{1}{2}S(\psi'v, v)\Big|_a^b\\   
	&  = J_{(a,b)}\big(\psi v',v'\big)+ J_{(a,b)}\big(\psi (B+ q)v,v\big) 
	 \\
&  	\ \ \ \ \ \ \ - J_{(a,b)}\big(v,(d''\psi + d'\psi'+ \frac{1}{2}\psi'')v\big).
\label{S11Identity2}
\end{split}
\end{equation}
\end{lemma}

\begin{proof}
Take the  inner product of  (\ref{S11Equationforv})  and $\psi v'$.
By integration by parts, we have 
\begin{equation}
\begin{split}
S(v',\psi v')\Big|_a^b & =  2 J_{(a,b)}(v'',\psi v') + J_{(a,b)}(v',\psi'v'), \\
S((B+q)v,\psi v)\Big|_a^b & = 2 J_{(a,b)}((B+q)v,\psi v') + J_{(a,b)}((\psi (B+q))'v,v),
\end{split}
\nonumber
\end{equation}
which
yield (\ref{S11Identity1}).  Next note that
\begin{equation}
J_{(a,b)}(v'',\psi v) = S(v',\psi v)\Big|_a^b - \frac{1}{2}S(v,\psi 'v)\Big|_a^b + \frac{1}{2}J_{(a,b)}(v,\psi'' v) - J_{(a,b)}(v',\psi v').
\nonumber
\end{equation}
Using this and taking the inner product of (\ref{S11Equationforv}) with $\psi v$, we obtain 
\begin{equation}
\begin{split}
& S(v',\psi v)\Big|_a^b - \frac{1}{2}S(v,\psi'v)\Big|_a^b + 
\frac{1}{2}J_{(a,b)}(v,\psi''v) - J_{(a,b)}(v',\psi v')\\
& = 2J_{(a,b)}(d'v',\psi v) + J_{(a,b)}(\psi(B + q)v,v).
\end{split}
\nonumber
\end{equation}
By integration by parts, 
$$
2J_{(a,b)}(d'v',\psi v) = S(v,d'\psi v)\Big|_a^b - J_{(a,b)}(v, (d''\psi + d'\psi')v).
$$
These three formulas imply
(\ref{S11Identity2}). 
\end{proof}

\medskip
Our main task is to increase the decay order of $\|u(t)\|$ step by step.


\begin{lemma} \label{S9Bvvdecayestimate}
Assume that for some $\alpha < \delta$
\begin{equation}
\liminf_{t\to\infty}\, t^{\alpha}\big(\|u'(t)\|^2 + \|u(t)\|^2\big)=0.
\label{DecayLmmaassump}
\end{equation}
 Then we have
\begin{equation}
\liminf_{t\to\infty}\, t^{\alpha}(B(t)u(t),u(t)) = 0.
\nonumber
\end{equation}
\end{lemma}

\begin{proof}
Take $\psi = t^{\alpha}$, $d=0$ in (\ref{S11Identity1}). Noting that $q = V - \lambda$, we  have
\begin{equation}
\begin{split}
& S(t^{\alpha}u',u')\Big|_a^b - S(t^{\alpha}(B+V)u,u)\Big|_a^b 
+ \lambda S(t^{\alpha}u,u)\Big|_a^b \\
=& \alpha J_{(a,b)}(t^{\alpha-1}u',u') - J_{(a,b)}((t^{\alpha}(B+V))'u,u) + \lambda\alpha J_{(a,b)}(t^{\alpha-1}u,u).
\end{split}
\label{Lemma4.7Eq1}
\end{equation}
Note that by (\ref{S8C3tdBtdt}), 
$$
-(t^{\alpha} B)' = -t^{\alpha-1}(tB' + \delta B) - (\alpha-\delta)t^{\alpha -1}B \geq - Ct^{\alpha-1-\epsilon_0}  + (\delta-\alpha)t^{\alpha-1}B.
$$
By the assumption (\ref{S8C3tdBtdt}), the right-hand side of (\ref{Lemma4.7Eq1}) is estimated from below by
\begin{equation}
\begin{split}
\alpha J_{(a,b)}(t^{\alpha-1}u',u') - CJ_{(a,b)}(t^{\alpha-1-\epsilon_0}u,u) \\
+ (\delta - \alpha)J_{(a,b)}(t^{\alpha-1}Bu,u) + \lambda\alpha J_{(a,b)}(t^{\alpha-1}u,u).
\end{split}
\nonumber
\end{equation}
Taking $a$ large enough (independently of $u$), this is estimated from below by
\begin{equation}
\begin{split}
\kappa_0\Big(J_{(a,b)}(t^{\alpha-1}u',u') + J_{(a,b)}(t^{\alpha-1}Bu,u) + 
J_{(a,b)}(t^{\alpha-1}u,u)\Big),
\end{split}
\nonumber
\end{equation}
where
\begin{equation}
\kappa_0 = \min\{\alpha,\delta-\alpha,\lambda\alpha/2\}.
\nonumber
\end{equation}
We have by (\ref{Lemma4.7Eq1})
\begin{equation}
\begin{split}
& S_b(t^{\alpha}u',u') - S_b(t^{\alpha}Vu,u) + S_a(t^{\alpha}(B + V)u,u) + 
\lambda S_b(t^{\alpha}u,u) \\
& \geq   S_b(t^{\alpha}Bu,u) + S_a(t^{\alpha}u',u') + \lambda S_a(t^{\alpha}u,u) \\
& + \kappa_0\Big(J_{(a,b)}(t^{\alpha-1}u',u') + J_{(a,b)}(t^{\alpha-1}Bu,u) + 
J_{(a,b)}(t^{\alpha-1}u,u)\Big).
\end{split}
\label{Lemma4.7ProofInequality2}
\end{equation}
Noting that
$$
CS_b(t^{\alpha - \epsilon_0}u,u) \geq - S_b(t^{\alpha}Vu,u), \quad
S_b(t^{\alpha}Bu,u) \geq 0,
$$
 we have by using the assumption of the lemma and letting $b \to \infty$ in (\ref{Lemma4.7ProofInequality2}) along a suitable sequence
\begin{equation}
\begin{split}
& S_a(t^{\alpha}(B + V)u,u)   \geq    S_a(t^{\alpha}u',u') + \lambda S_a(t^{\alpha}u,u) \\
 & \hskip5mm + \kappa_0\Big(J_{(a,\infty)}(t^{\alpha-1}u',u') + J_{(a,\infty)}(t^{\alpha-1}Bu,u) + 
J_{(a,\infty)}(t^{\alpha-1}u,u)\Big).
\end{split}
\label{Lemma4.7Ineq1}
\end{equation}
This implies $J_{(a,\infty)}(t^{\alpha-1}Bu,u) < \infty$, hence $\liminf_{t\to\infty}S_t(t^{\alpha}Bu,u)=0$. In fact, if $g(t) \geq 0$, 
$\displaystyle{\int_a^{\infty}}g(t)dt < \infty$ implies that $\liminf_{t\to\infty}tg(t)=0$. \end{proof}

The following lemma \ref{S11Bv,vestmatedformbelow} follows easily from (\ref{Lemma4.7Ineq1}).


\begin{lemma}\label{S11Bv,vestmatedformbelow}
	Assume $\alpha<\delta$ and (\ref{DecayLmmaassump}). Then, there exist constants $C, a_0>0$ such that
\begin{equation}
\begin{split}
 t^{\alpha}((B + V)u,u)\Big|_{t=a} \geq &\  t^{\alpha}\Big(\|u'\|^2 + 
\lambda \|u\|^2\Big)\Big|_{t=a} \\
& + C\int_a^{\infty}t^{\alpha-1}\Big(\|u'\|^2 + (Bu,u) + \|u\|^2\Big)dt
\end{split}
\label{S11Formula7}
\end{equation} 
for any $a>a_0$.
\end{lemma}


\begin{lemma}\label{S11polynomweighinequality}
Assume $\alpha<\delta$ and (\ref{DecayLmmaassump}). 
Then, there exist constants $C, r_0>0$ such that for $r > r_0$
\begin{equation}
\begin{split}
& C\int_r^{\infty}(t-r)t^{\alpha-1}\Big(\|u'\|^2 + (Bu,u) + \|u\|^2\Big)dt \\
& \leq t^{\alpha}\Big(\|u'\|^2 + \|u\|^2\Big)\Big|_{t=r}
+ \int_r^{\infty}t^{\alpha-2}\|u\|^2dt.
\end{split}
\label{Lemma4.7Inequality}
\end{equation}
\end{lemma}

\begin{proof}
We put
\begin{equation}
f(t) = t^{\alpha-1}\Big(\|u'\|^2 + (Bu,u) + \|u\|^2\Big).
\label{Lemma4.7Definef(t)}
\end{equation}
Integrating (\ref{S11Formula7}) with respect to $a$ over $(r,b)$, we have
\begin{equation}
\begin{split}
& C\int_r^{b}(t-r)f(t)dt + C\int_b^{\infty}(b-r)f(t)dt \\
& \leq 
\int_r^{b}\big(- (t^{\alpha}u',u') + (t^{\alpha}(B+V)u,u) - \lambda (t^{\alpha}u,u)\big)dt.
\end{split}
\label{Lemma4.8Inequality1}
\end{equation} 
Taking $\psi=t^{\alpha}$ and $d=0$ in (\ref{S11Identity2}), we see that the right-hand side is equal to
\begin{equation}
\begin{split}
& S(t^{\alpha}u',u)\Big|_r^b - \frac{\alpha}{2}S( t^{\alpha-1}u,u)\Big|_r^b \\
& - 2J_{(r,b)}(t^{\alpha}u',u')  + \frac{\alpha(\alpha-1)}{2}J_{(r,b)} (t^{\alpha-2}u,u).
 \end{split}
\nonumber
\end{equation} 
We then have
\begin{equation}
\begin{split}
& C\int_r^b(t-r)f(t)dt + 2J_{(r,b)}(t^{\alpha}u',u') \\
& \leq (t^{\alpha}u',u)\Big|_{t=b} - (t^{\alpha}u',u)\Big|_{t=r} 
+ \frac{\alpha}{2}(t^{\alpha-1}u,u)\Big|_{t=r}  + \frac{\alpha(\alpha-1)}{2}J_{(r,b)}(t^{\alpha-2}u,u).
\end{split}
\nonumber
\end{equation}
Therefore, letting $b\to\infty$ along a suitable sequence, we obtain
\begin{equation}
\begin{split}
& C\left(\int_r^{\infty}(t-r)f(t)dt + \int_r^{\infty}t^{\alpha}\|u'\|^2dt\right) \\
& \leq t^{\alpha}\big(\|u'\|^2 + \|u\|^2\big)\Big|_{t=r} + 
\int_r^{\infty}t^{\alpha - 2}\|u\|^2dt.
\end{split}
 \nonumber
\end{equation}
This proves the lemma.
\end{proof}

 
\begin{lemma}\label{S11polynomweightlemma}
Let $u$ be as in Theorem \ref{Growthproperty}. Then
$$
\int_1^{\infty}t^m\Big(\|u'\|^2 + (Bu,u) + \|u\|^2\Big)dt < \infty, \quad \forall m > 0.
$$
\end{lemma}

\begin{proof}
We take $0<\alpha<\min\{2\gamma,\delta\}$, and prove this lemma in the form
\begin{equation}
\int_1^{\infty}t^{\alpha-1+m}\Big(\|u'\|^2 + (Bu,u) + \|u\|^2\Big)dt < \infty, \quad \forall m \geq 0.
\label{S11Lemma11.6Variant}
\end{equation}
By the assumption of the theorem, (\ref{DecayLmmaassump}) holds for $0<\alpha <2\gamma$. In the proof of Lemma \ref{S9Bvvdecayestimate}, i.e. in (\ref{Lemma4.7Ineq1}), we have already proven (\ref{S11Lemma11.6Variant})  for $m=0$. 
Let $f(t)$ be as in (\ref{Lemma4.7Definef(t)}) and put
$$
g(t) = t^{\alpha-1}\big(\|u'\|^2 + \|u\|^2).
$$ 
We show that for $m \geq 1$
\begin{equation}
\begin{split}
& C\int_r^{\infty}\frac{(t-r)^{m+1}}{(m+1)!}f(t)dt \\
& \leq 
\int_r^{\infty}\frac{(t-r)^{m-1}}{(m-1)!}tg(t)dt  + \int_r^{\infty}\frac{(t-r)^m}{m!}g(t)dt
\end{split}
\label{Lemma4.8Inequality}
\end{equation}
and the right-hand side is finite. 
By (\ref{Lemma4.7Inequality}), we have
$$
C\int_r^{\infty}(t-r)f(t)dt \leq rg(r) + \int_r^{\infty}g(t)dt,
$$
where the right-hand side is finite. Integrating this inequality, we obtain 
(\ref{Lemma4.8Inequality}) for $m=1$, where the right-hand side is finite. 
 Repeating this procedure, we prove (\ref{Lemma4.8Inequality}), hence (\ref{S11Lemma11.6Variant}) for all $m \geq 1$. 
\end{proof}

The proof of Theorem \ref{Growthproperty} is now completed. 

If we assume $\delta > 2$, which holds when $\rho(r)$ grows like $O(r^{\beta})$ with $\beta > 1$, i.e. faster than the Euclidean metric, one can show that $u=0$ near infinity. 
In fact, let $u$ be as in Theorem \ref{Growthproperty},  and put $d = m\log t$, $v = e^{d}u=t^mu$. Take $0<\alpha<\min\{2\gamma,\delta\}$. Since $u$ decays rapidly at infinity, one can take $\gamma > 2$. Then,  (\ref{DecayLmmaassump}) is satisfied.
In (\ref{S11Identity1}), we take $\psi = t^{\alpha}$, $d = m\log t$, $a = r$ and let $b \to \infty$. Since $\delta > 2$, one can take $\alpha > 2$. Then, we have
\begin{equation}
\begin{split}
& - S_r(t^{\alpha}v',v') + S_r(t^{\alpha}(B+q)v,v) \\
=&\  J_{(r,\infty)}\big((\alpha + 4m)t^{\alpha-1}v',v'\big) 
 - J_{(r,\infty)}((t^{\alpha}(B+q))'v,v).
\end{split}
\end{equation}
Since $q = V-\lambda  - (m+m^2)/t^2$, we have, by the assumption (\ref{tbt(x,x)leq-delatB(t)+Ct-psilon0})
$$
- (t^{\alpha}(B+q))' \geq - Ct^{\alpha-1-\epsilon_0} + (\delta - \alpha)t^{\alpha - 1}B + (\alpha -2)(m^2 + m)t^{-3+\alpha} + \lambda\alpha t^{\alpha-1}
$$
for a constant $\kappa_0 > 0$.
Then, we have
\begin{equation}
\begin{split}
& - S_r(t^{\alpha}v',v') + S_r(t^{\alpha}(B+q)v,v)) \\
& \geq
\kappa_0(J_{(r,\infty)}(t^{\alpha-1}v',v') + 
J_{(r,\infty)}(t^{\alpha-1}Bv,v) + 
J_{(r,\infty)}(t^{\alpha-1}v,v)).
\end{split}
\nonumber
\end{equation}
 In particular, there exists $t_0>0$ such that
$$
- S_t(v',v') + S_t((B+q)v,v)\geq 0, \quad \forall t\geq t_0.
$$
Using $v = t^mu$, $q = V -\lambda - (m+m^2)/t^2$, we then have
\begin{equation}
\begin{split}
-(2m^2 + m)S_t\Big(\frac{u}{t},\frac{u}{t}\big) - 2mS_t\big(\frac{u}{t},u'\big) - S_t(u',u') \\
+ S_t(Bu,u) + S_t((V-E)u,u) \geq 0, \quad \forall t \geq t_0.
\end{split}
\nonumber
\end{equation}
Since this holds for any $m > 0$, we see that $u(t)=0$ for $t\geq t_0$. 

\medskip
This observation shows that the rapid growth of the volume facilitates the spectral analysis. 
To deal with  slowly growing metrics, we need more elaborate consideration.


\subsection{Proof of Theorem
\ref{AbstractRllich1}}
 We use the method in Saito \cite{Saito}, which orginates from the work of Kato \cite{Kato59}. Although the proof here is apprently different from the one in the previous section, they are actually closely related.
 In the following, $\|\cdot\|_X$ is simply written as $\|\cdot\|$.
We show that if ${\rm supp}\,u(t)$ is unbounded,
\begin{equation}
\liminf_{t\to\infty}\,(\|u'(t)\|^2 + \|u(t)\|^2) > 0.
\label{S6GrowthPositive}
\end{equation}
To prove this, we consider the following two cases. 
 We put
\begin{equation}
(Ku)(t) = \|u'(t)\|^2 + \lambda\|u(t)\|^2 - (B(t)u(t),u(t)) - (V(t)u(t),u(t)).
\nonumber
\end{equation}
\medskip
\noindent
{\it Case 1}.  There exists a sequence $t_1 < t_2 < \cdots \to \infty$ such that
\begin{equation}
 (Ku)(t_n) > 0, \ n = 1, 2, \cdots.
\end{equation}


\begin{lemma}\label{GrowthLemma1}
There exist constants $C_1, T_1 > 0$ such that
$$
\frac{d}{dt}(Ku)(t) \geq - {C_1}{(1 + t)^{-1-\epsilon}}(Ku)(t), \quad 
\forall t > T_1.
$$
\end{lemma}

\begin{proof}
By choosing $\epsilon$ small enough, we can assume that 
\begin{equation}
 \|V'(t)\| \leq C(1 + t)^{-1-2\epsilon}.
 \label{eq:Chap2Sect3Vprimet}
\end{equation}
By the equation (\ref{eq:Chap2Sect3Diffeq})
\begin{eqnarray*}
\frac{d}{dt}(Ku)(t) & = & 2{\rm Re}\,\Big[(u'',u') + \lambda(u,u') - (Bu,u') 
- (Vu,u')\Big]   - ((B' + V')u,u) \\
& =&   - ((B' + V')u,u). 
\end{eqnarray*}
By (\ref{eq:Chap2Sect3Vprimet}), there exists $t_0 = t_0(\epsilon) > 0$ such that for $t > t_0$ 
$$
|(V'(t)u,u)| \leq \frac{\epsilon}{2}(1 + t)^{-1-\epsilon}\|u\|^2.
\nonumber
$$
By Lemma \ref{LemmatB(t)leqct-epsilon},
$$
- (B'u,u) \geq \frac{\delta}{t}(Bu,u) - \frac{C}{t^{1+\epsilon_0}}\|u\|^2.
$$
By virtue of the above estimates, there is $C_{\epsilon} > 0$ such that 
for $t > t_0$
\begin{eqnarray*}
\frac{d}{dt}(Ku)(t) 
&\geq& - Ct^{-1-\epsilon}(\|u'\|^2 + \|u\|\|u'\| + \frac{\epsilon}{2}\|u\|^2)
 + \frac{\delta}{t}(Bu,u) \\
& \geq& - C_{\epsilon}t^{-1-\epsilon}\|u'\|^2  - C\epsilon t^{-1-\epsilon}\|u\|^2
 + \frac{\delta}{t}(Bu,u). 
\end{eqnarray*}
We rewrite the right-hand side as
\begin{eqnarray*}
& &-C_{\epsilon}t^{-1-\epsilon}(\|u'\|^2 + \lambda\|u\|^2) + (C_{\epsilon}\lambda- C\epsilon)t^{-1-\epsilon}\|u\|^2 + \frac{\delta}{t}(Bu,u) \\
&=& - C_{\epsilon}t^{-1-\epsilon}(Ku)(t) \\
& &+ (C_{\epsilon}\lambda- C\epsilon)t^{-1-\epsilon}\|u\|^2 - 
C_{\epsilon}t^{-1-\epsilon}(Vu,u)  + \frac{\delta}{t}(Bu,u).
\end{eqnarray*}
Choose $C_{\epsilon}$ large enough so that $C_{\epsilon}\lambda - C\epsilon \geq \frac{1}{2}C_{\epsilon}\lambda$. Using (\ref{S7C-4normVt}), choose $t_0 = t_0(\epsilon,C_{\epsilon})$ such that, for $t > t_0$, $\frac{\lambda}{2}\|u\|^2 - (Vu,u) \geq 0$. Thus, the 3rd line is non-negative for  $t > t_0$. Hence the lemma is proved. \end{proof}

\medskip
Let us prove (\ref{S6GrowthPositive}) for the Case 1.
Let $T_1$ be as in  Lemma \ref{GrowthLemma1}. Then for some $T > T_1$, $(Ku)(T) > 0$. We show that  $(Ku)(t) \geq 0, \ \forall t > T$.
In fact Lemma \ref{GrowthLemma1} implies
\begin{equation}
\frac{d}{dt}\left\{\exp\left(C_1\int_T^t(1 + s)^{-1-\epsilon}ds\right)
(Ku)(t)\right\} \geq 0, \quad \forall t > T.
\nonumber
\end{equation}
Hence,
\begin{equation}
(Ku)(t) \geq \exp\left(- C_1\int_T^t(1 + s)^{-1-\epsilon}ds\right)
(Ku)(T), \quad \forall t > T.
\nonumber
\end{equation}
This then implies that, for $t >t(\lambda)$,
\begin{equation}
\begin{split}
\|u'(t)\|^2 + \lambda\|u(t)\|^2 &= Ku(t) + (B(t)u(t),u(t)) + (V(t)u(t),u(t)) \\
&\geq \exp\left(- C_1\int_T^t(1 + s)^{-1-\epsilon}ds\right)
(Ku)(T)   - C t^{-\epsilon}\|u(t)\|^2.
\end{split}
\nonumber
\end{equation}
Therefore, we arrive at
\begin{equation}
\liminf_{t\to\infty}\,(\|u'(t)\|^2 + \|u(t)\|^2) \geq 
\frac{1}{2}\exp\left(- C_1\int_T^{\infty}(1 + s)^{-1-\epsilon}ds\right)
(Ku)(T) > 0.
\nonumber
\end{equation}
Note that in  this case, we do not use the assumption $\delta > 2/3$.

\medskip
Next let us consider Case 2:

\medskip
\noindent
{\it Case 2}. There exists $T_1 > 0$ such that $(Ku)(t) \leq 0$ for all $t > T_1$.

\medskip
To deal with this case, take $\beta, \gamma, m,  d(t)$ such that
\begin{equation}
m > 0, \quad \frac{1}{3} < \gamma < 1, \quad  2\gamma < \beta < \delta,  \quad d(t) = \frac{m}{1 - \gamma}
t^{1-\gamma},
\label{S&mbetagamma}
\end{equation}
and put
\begin{equation}
(Nu)(t) = t^{\beta}\left[K(e^{d(t)}u) + \frac{m^2 - \log t}{t^{2\gamma}}\|e^{d(t)}u\|^2\right].
\nonumber
\end{equation}


\begin{lemma} \label{Nupositive}
If ${\rm supp}\,u(t)$ is unbounded, there exist constants
 $m_1 \geq 1$, $T_2 \geq T_1$ such that
\begin{equation}
(Nu)(t) \geq 0, \quad \forall t \geq T_2, \quad \forall m \geq m_1.
\nonumber
\end{equation}
\end{lemma}

\begin{proof}
Letting $w(t) = e^{d(t)}u(t)$, we have by  a direct computation, 
\begin{equation}
\begin{split}
w' & = e^du' + mt^{-\gamma}w, \\
w'' & = e^du'' + mt^{-\gamma}e^du' + mt^{-\gamma}w' - \gamma mt^{-\gamma - 1}w 
\\ 
&= Bw + Vw - \lambda w + 2mt^{-\gamma}w'   -  (\gamma m t^{-\gamma -1} +m^2t^{-2\gamma})w.
\end{split}
\nonumber
\end{equation}
Hence,
\begin{equation}
\begin{split}
& \frac{d}{dt}(Kw)\\
 & = 2{\rm Re}\,(w'' + \lambda w - Vw - Bw,w') - 
(\big(B' + V'\big)w,w) \\
&=  4mt^{-\gamma}\|w'\|^2 
 - 2(\gamma mt^{-\gamma - 1} + m^2t^{-2\gamma}){\rm Re}\,(w,w')  - (\big(B' + V'\big)w,w) . 
\end{split}
\label{eq:Chpa2Sect3ddtKw}
\end{equation}
We then have
\begin{equation}
\begin{split}
  \frac{d}{dt}(Nu) 
 & = \beta t^{\beta -1}Kw + t^{\beta}\frac{d}{dt}(Kw) + \frac{(\beta - 2\gamma)(m^2 - \log t) - 1}{t^{2\gamma - \beta+1}}\|w\|^2 \\
& \ \ \ + \frac{2(m^2-\log t)}{t^{2\gamma-\beta}}{\rm Re}\,(w',w)
\end{split}
\nonumber
\end{equation}
Using (\ref{eq:Chpa2Sect3ddtKw}), we have
\begin{equation}
\begin{split}
& t^{1-\beta}\frac{d}{dt}(Nu) \\
&  = (4mt^{1-\gamma} + \beta)\|w'\|^2 
+\big (\beta \lambda + \frac{(\beta-2\gamma)(m^2-\log t)-1}{t^{2\gamma}}\big)\|w\|^2 \\
&   - 2(\gamma m t^{-\gamma} + t^{1-2\gamma}\log t)\,{\rm Re}\,(w',w) 
 -  ((\beta V + tV')w,w) \\
&  - ((\beta B + tB')w,w) \\
& =: I_1 + I_2 + I_3.
\end{split}
\label{eq:Chap2Sect3ddtNt}
\end{equation}
For large $t > 0$, $I_1$ is estimated from below as
\begin{equation}
 I_1 \geq (4mt^{1-\gamma} + \beta)\|w'\|^2 + \big(\frac{\beta \lambda }{2} + 
 (\beta - 2\gamma)t^{-2\gamma}m^2\big)\|w\|^2.
 \nonumber
\end{equation}
By (\ref{S7C-4normVt}), 
$I_2$ is estimated from below as
\begin{eqnarray*}
 I_2 &\geq& - 2(\gamma m t^{-\gamma} + t^{1-2\gamma}\log t)\|w\|\|w'\| - Ct^{-\epsilon}\|w\|^2 \\
 &\geq& - \epsilon m^2 t^{-2\gamma}\|w\|^2 - \frac{\gamma^2}{\epsilon}\|w'\|^2 \\
 & & - 2t^{1-2\gamma}\log t\|w\|\|w'\| - Ct^{-\epsilon}\|w\|^2,
 \nonumber
\end{eqnarray*}
where $\epsilon > 0$ is chosen arbitrarily  small. Note that the constant $C$ is  independent of $m$.

Since $\beta < \delta$, 
$I_3$ is estimated from below as
\begin{eqnarray*}
I_3 \geq - Ct^{-\epsilon}\|w\|^2.
\end{eqnarray*}
Choosing $\beta - 2\gamma \geq \epsilon$ and putting the above estimates together, we have
\begin{equation}
t^{1-\beta} \frac{d}{dt}(Nu) \geq 3mt^{1-\gamma}\|w'\|^2 + \frac{\beta \lambda}{3}\|w\|^2 - 
 2t^{1-2\gamma}\log t\,\|w\|\|w'\|.
 \nonumber
 \end{equation}
Finally, we use the inequality
\begin{equation}
  t^{1-2\gamma}\log t\|w\|\|w'\| \leq \epsilon t^{1-\gamma}\|w'\|^2 + 
 C_{\epsilon}t^{1-3\gamma}(\log t)^2\|w\|^2
 \nonumber 
\end{equation}
and $1 - 3\gamma < 0$. Then there is $t_0 > 0$ independent of $m$ such that
\begin{equation}
  \frac{d}{dt}(Nu)(t) \geq t^{\beta-1}\Big(2mt^{1-\gamma}\|w'\|^2 + \frac{\beta \lambda}{4}\|w\|^2\Big) \geq 0 
 \label{eq:Chap2SEct3ddtNutgeq3m}
\end{equation}
for $t > t_0$. 

On the other hand, $Nu(t)$ can be rewritten as
\begin{equation}
 \begin{split}
  (Nu)(t) & = t^{\beta}e^{2d}\big[\|mt^{-\gamma}u + u'\|^2 + \lambda\|u\|^2 \\
  & \ \ \ - (Bu,u) - (Vu,u) + t^{-2\gamma}(m^2 - \log t)\big]\|u\|^2 \\
  &= t^{\beta}e^{2d}\big[2t^{-2\gamma}\|u\|^2m^2 + 2t^{-\gamma}{\rm Re}\,(u,u')m \\
  & + (Ku - t^{-2\gamma}\|u\|^2\log t)\big].
 \end{split}
 \label{eq:Chap2Sect3Nutte2d}
\end{equation}
By the assumption of the lemma, ${\rm supp}\,u(t)$ is unbounded. Therefore, there is $T_2 > t_0$ such that $\|u(T_2)\| > 0$.
By choosing $m_1$ large enough, we then have 
\begin{equation}
(Nu)(T_2) > 0, \quad \forall m > m_1.
\label{eq:Chap2Sect3NuT2}
\end{equation}
The inequalities (\ref{eq:Chap2SEct3ddtNutgeq3m}) and (\ref{eq:Chap2Sect3NuT2}) prove the lemma. \qed

\bigskip

Since $Ku\leq 0$, Lemma \ref{Nupositive} and (\ref{eq:Chap2Sect3Nutte2d}) show that, for large $t$,
\begin{equation}
2t^{-2\gamma}\|u(t)\|^2m^2 + 2t^{-\gamma}{\rm Re}\,(u(t),u'(t))m
- t^{-2\gamma}\|u(t)\|^2\log t \geq 0,
\nonumber
\end{equation}
which together with
\begin{equation}
 \frac{d}{dt}\|u(t)\|^2 = 2{\rm Re}\,(u(t),u'(t)),
 \nonumber
\end{equation}
yields, for large $t > 0$, 
\begin{equation}
\frac{d}{dt}\|u(t)\|^2 \geq t^{-\gamma}\left(\frac{1}{m}\log t - 
2m\right)\|u(t)\|^2 \geq 0.
\label{eq:Chap2Sect3ddtnormut}
\end{equation}
Since the support of $u(t)$ is unbounded, we can choose $T$ large enough so that $\|u(T)\| > 0$.
In view of (\ref{eq:Chap2Sect3ddtnormut}), we then have
\begin{equation}
 \|u(t)\| \geq \|u(T)\| > 0, \quad \forall t > T,
 \nonumber
\end{equation}
which proves (\ref{S6GrowthPositive}). 
\end{proof}

 
\section{Integral identities}
\label{SectionIntegralidentities}
The next aim is to derive  resolvent estimates for $- \Delta_{\mathcal M}$ on one end $(0,\infty)\times M$, which we denoted  as $\mathcal M = (0,\infty)\times M$ for the sake of simplicity. Our method is based on the two integral identities to be proved in Lemmas \ref{IntIdentity1} and \ref{IntIdentity2} below.
The basic assumption on the metric
$ds^2 = (dr)^2 + \rho(r)^2h(r,x,dx)$ in this section 
is the existence of 
the constants $\alpha_0, \gamma_0 > 0$ such that 
\begin{equation}
 \frac{\rho'(r)}{\rho(r)} - c_0  \in S^{-\alpha_0}, 
\label{rho'/rhosimbeta/r}
\end{equation}
\begin{equation}
h^{ij}(r,x) - h_{M}^{ij}(x) \in S^{-\gamma_0},
\label{rh0prime-rho,hinS-gamma}
\end{equation}
where $S^{\kappa}$ is defined in Definition \ref{S2DefineSingularSkappa}.
Let $g = \rho^{2(n-1)}h$, $h = h(r,x) = \det (h_{ij}(r,x))$.


\subsection{Preliminaries} 
We regard $L^2(\mathcal M)$ as the $L^2$-space of $L^2(M)$-valued functions over $(0,\infty)$. We did it already in Subsection \ref{Subsec4.2}, however, we employ here a slightly different formulation in order to take account of the growth order of $\rho(r)$ more explicitly.
We define the $L^2$-space over $M$  by
\begin{equation}
{\bf h}(r) = L^2(M;\sqrt{h(r,x)}dx)
\label{Definebhf(r)}
\end{equation}
with the $r$-dependent inner product and norm 
$$
(\varphi,\psi)_{{\bf h}(r)} = \int_M\varphi(x)\overline{\psi(x)}\sqrt{h(r,x)}dx, \quad \|\varphi\|_{{\bf h}(r)} = \sqrt{(\varphi,\varphi)_{{\bf h}(r)}}.
$$ 
Here, we identify the local coordinate $x$ with a point in $M$.
Note that the space ${\bf h}(r)$ is independent of $r$ as a set. 
The inner product of the Hilbert space  $L^2\big((0,\infty)\times M;\sqrt{g}\,drdx\big)$ is rewritten as
\begin{equation}
(u,v) = \int_0^{\infty} (u(r),v(r))_{{\bf h}(r)}\, \rho^{n-1}(r)dr.
\label{S31diminnerproduct}
\end{equation}
For $0 < a < b < \infty$, we put
\begin{eqnarray}
M_r(u,v)&=& (u(r),v(r))_{{\bf h}(r)}\,\rho^{n-1}(r) = \int_Mu(r,x)\overline{v(r,x)}\sqrt{g(r,x)}dx, 
\nonumber\\
M(u,v)\Big|_a^b &= &M_b(u,v) - M_a(u,v), \nonumber
\\
I_{(a,b)}(u,v) &=& \int_a^b(u(r),v(r))_{{\bf h}(r)}\,\rho^{n-1}(r)dr 
= \int_a^bM_r(u,v)dr,\nonumber
 \\
I(u,v) &= &I_{(0,\infty)}(u,v). 
\nonumber
\end{eqnarray}
   By integration by parts, we have
\begin{equation}
I_{(a,b)}(u',v) = M(u,v)\Big|_a^b - I_{(a,b)}(u,v')
- I_{(a,b)}(\frac{g'}{2g}u,v),
\label{S3partialrast=-partial+}
\end{equation}
which implies that the formal adjoint of $\partial_r$ is
\begin{equation}
\partial_r^{\ast} = - \partial_r - \frac{g'}{2g},
\label{S5partialrast}
\end{equation}
and, for  a real-valued $C^1$-function $\varphi$,
\begin{equation}
{\rm Re}\,I_{(a,b)}(v',\varphi v) = 
M\big(\frac{\varphi}{2} v,v\big)\Big|_a^b - 
I_{(a,b)}\big((\frac{\varphi'}{2} + \frac{g'\varphi}{4g}) v,v\big).
\label{S3Intbyparts}
\end{equation}

Recall that the Laplacian on $\mathcal M$ is written as
\begin{equation}
- \Delta_{\mathcal M} = - \partial_r^2 - \frac{g'}{2g}\partial_r + B(r), \quad ' = \partial_r,
\nonumber
\end{equation}
\begin{equation}
B(r)  = - \rho(r)^{-2}\Lambda(r), \quad \Lambda(r) = \frac{1}{\sqrt{h}}\partial_{x_i}\Big(\sqrt{h}h^{ij}\partial_{x_j}\Big).
\label{S5B(r)andLambda(r)}
\end{equation}
\index{$B(r)$}
For any $r>0$, $\Lambda(r)$ is self-adjoint on ${\bf h}(r)$ with domain described in Theorem \ref{S2DLregularitytheorem}, which is independent of $r$.  We rewrite $\Delta_{\mathcal M}$ into a form which is more convenient for our computation.
Put for $z \in {\mathbb C}$
\begin{eqnarray}
E_0 &=& \big((n-1)c_0/2\big)^2,
\nonumber \\
k &=& \sqrt{z - E_0}, 
\nonumber \\
Q &=& \Big(\frac{g'}{4g}\Big)^2+ \Big(\frac{g'}{4g}\Big)' - E_0.
\nonumber
\end{eqnarray}
Here and in the sequel, we take the branch of $\sqrt{\cdot}$ on ${\mathbb C}\setminus[0,\infty)$ in such a way that ${\rm Im}\,\sqrt{\cdot} \geq 0$, i.e. $\sqrt{z} = \sqrt{r}e^{i\theta/2}$ if $z = re^{i\theta}$, $0 \leq \theta < 2\pi$. 
Then,  $- \Delta_{\mathcal M} - z$ is rewritten as 
\begin{equation}
\begin{split}
- \Delta_{\mathcal M} - z & = 
- \Big(\partial_r + \frac{g'}{4g}\Big)^2 + B(r) +  Q -k^2.
\end{split}
\label{S5DeltaM-zrewrittena1}
\end{equation}
The assumption  (\ref{rh0prime-rho,hinS-gamma}) implies, 
\begin{equation}
Q \in S^{-\alpha_0}.
\label{S5QinS-gamma}
\end{equation}
Take a complex-valued function $\psi$ and introduce a differential operator
\begin{equation}
D_r = \partial_r + \frac{g'}{4g} - i\psi.
\nonumber
\end{equation}
Then 
\begin{equation}
- \Delta_{\mathcal M} - z = - D_r^2 - 2i\psi D_r + B(r) - i\psi' + \psi^2 + Q-k^2.
\label{S2DeltagwrittenbyDr}
\end{equation}

We construct an approximate solution of
$$
- i\psi' + \psi^2 + Q -k^2 = 0
$$
by putting $\psi_{-1}=0$ and
\begin{equation}
\psi_m(r,x,k) = 
 \chi(r/R_m)\sqrt{ k^2 - Q + i\partial_r\, \psi_{m-1}}, \quad m \geq 0,
\label{S5Definepsimrxk}
\end{equation}
where $\chi \in C^{\infty}({\mathbb R})$, $\chi(t) = 1$ $(t > 2)$, $\chi(t) = 0$ $(t < 1)$.


\begin{lemma}\label{psimlemma}
For $k \neq 0$ and $m \geq 0$, by choosing $R_m$ large enough, $\psi_m$ is $C^{\infty}$ with respect to $r$ and has the following properties. \\
\noindent
(1) $\ \psi_m \in S^0$, $\partial_r\psi_m \in S^{-1-\alpha_0}$.\\
\noindent
(2) $ -i\partial_r\psi_m+ (\psi_m)^2 + Q - k^2 \in 
S^{-m-1- \alpha_0}. $ \\
\noindent
(3) $\ \psi_m - \sqrt{k^2 - Q} \in S^{-1-\alpha_0}$.
\end{lemma}

\begin{proof}
The assertion (1) is proven by induction on $m$. For large $r$ we have
$$
- i\partial_r\psi_m + (\psi_m)^2 + Q - k^2 = -i(\partial_r \psi_m - \partial_r\psi_{m-1}), \quad m \geq 0,
$$
$$
i(\psi_m - \psi_{m-1})= - \frac{\partial_r\psi_{m-1} - \partial_r\psi_{m-2}}{\psi_{m}+ \psi_{m-1}}, \quad m \geq 1.
$$
One can then derive (2) from (1). The assertion (3) follows from (\ref{S5Definepsimrxk}) and (1).
\end{proof}

\medskip
We take $m$ large enough, and put
\begin{equation}
\psi(r,x,k) = \psi_m(r,x,k),
\end{equation}
\begin{equation}
V = - i\psi' + \psi^2 + Q - k^2,
\label{DefVpm}
\end{equation}
\begin{equation}
 D(k)  = \partial_r + \frac{g'}{4g} - i\psi(r,x,k).
 \label{DefDpmk}
\end{equation}
\index{$D(k)$}
For a solution $u$ of the equation
\begin{equation}
(- \Delta_{\mathcal M} - z)u = (- \partial_r^2 - \frac{g'}{2g}\partial_r + B(r) - z)u = f,
\label{S2Equation1}
\end{equation}
 we put 
\begin{equation}
v = D(k)u,  \quad w = \sqrt{B(r)}u =  \rho(r)^{-1}\sqrt{-\Lambda(r)}\, u, 
\label{Definevpm}
\end{equation}
\index{$v = D(k)u$} \index{$w = \sqrt{B(r)}u$} 
where $B(r)$ and $\Lambda(r)$ are defined by (\ref{S5B(r)andLambda(r)}).


\subsection{The 1st identity}


\begin{lemma}
\label{IntIdentity1} 
Let $u$ be a solution to (\ref{S2Equation1}), and $v$, $w$ as in (\ref{Definevpm}). Then, for a non-negative $C^1$-function $\varphi(r)$, we have
\begin{equation}
\begin{split}
& I_{(a,b)}\big(\varphi'v,v\big) - I_{(a,b)}\big( \varphi' w,w\big) + M\big(\varphi w,w\big)\Big|_a^b - M\big({\varphi}v,v\big)\Big|_a^b\\
 =& - 2I_{(a,b)}\big(({\rm Im}\,\psi)\varphi v,v\big) - 2I_{(a,b)}\big(({\rm Im}\,\psi)\varphi w,w\big) \\
& + 2{\rm Re}\,I_{(a,b)}\big(\varphi f,v\big) - 2{\rm Re}\, I_{(a,b)}\big(\varphi Vu,v\big)  + I_{(a,b)}\big(\varphi B'u, u\big).
\end{split}
\label{EqualityinLemmaIntIdentity1} 
\end{equation}
\end{lemma}

\begin{proof}
By (\ref{S2DeltagwrittenbyDr}), we have
$$
- D(k)^2u - 2i\psi D(k)u + B(r) u = f - Vu,
$$
which implies
\begin{equation}
- \big(\partial_r + \frac{g'}{4g}\big)v
+ Bu
 =\ f + i\psi v - Vu.
\label{S3Eq1}
\end{equation}
Taking the inner product with $\varphi v$, 
\begin{equation}
\begin{split}
& I_{(a,b)}\big(-(\partial_r + \frac{g'}{4g})v + Bu,\varphi v\big) \\
& = 
 I_{(a,b)}(\varphi f, v) + iI_{(a,b)}(\psi\varphi v, v)  - I_{(a,b)}(\varphi Vu, v).
\end{split}
\label{S3Eq3}
\end{equation}
On the other hand,   by (\ref{S3Intbyparts}),
\begin{equation}
2 {\rm Re}\, I_{(a,b)}\big(-(\partial_r + \frac{g'}{4g})v, \varphi v\big) \\
= I_{(a,b)}\big(\varphi'v,v) - M\big({\varphi}v,v\big)\Big|_a^b.
\label{S3Eq4}
\end{equation}
Using
\begin{equation}
\begin{split}
 {\rm Re}\,I_{(a,b)}(Bu,\varphi u')
 = &\ M\big(\frac{\varphi}{2}w,w \big)\Big|_a^b - 
I_{(a,b)}(\frac{\varphi'}{2}w,w) \\
& - {\rm Re}\,I_{(a,b)}\big(\frac{g'}{4g}\varphi Bu, u\big) - I_{(a,b)}\big(\frac{\varphi}{2} B' u,u\big),
\end{split}
\end{equation}
we have
\begin{equation}
\begin{split}
{\rm Re}\,I_{(a,b)}(Bu,\varphi v)
= & 
M\big(\frac{\varphi}{2}w,w\big)\Big|_a^b - 
I_{(a,b)}\big(\frac{\varphi'}{2} w,w\big)  \\
& - I_{(a,b)}(\frac{\varphi}{2} B'u,u)  + \,I_{(a,b)}\big(({\rm Im}\,\psi)\varphi w,w\big).
\end{split}
\label{S3Eq5}
\end{equation}
By (\ref{S3Eq4}) and (\ref{S3Eq5}),  the real part of  the left-hand side of (\ref{S3Eq3}) is equal to
\begin{equation}
\begin{split}
& M\big(\frac{\varphi}{2}w,w\big)\Big|_a^b - 
M\big(\frac{\varphi}{2}v,v\big)\Big|_a^b 
 +  I_{(a,b)}\big(\frac{\varphi'}{2}v,v\big) \\
& - I_{(a,b)}\big(\frac{\varphi'}{2}
w,w\big) 
 -  I_{(a,b)}\big(\frac{\varphi}{2} B'u, u\big) 
 + I_{(a,b)}\big(({\rm Im}\, \psi)\varphi w, w\big).
\end{split}
\label{S3Eq6}
\end{equation}
In view of  (\ref{S3Eq3}) and (\ref{S3Eq6}), we have
\begin{equation}
\begin{split}
& M\big(\frac{\varphi}{2}w,w\big)\Big|_a^b - 
M\big(\frac{\varphi}{2}v,v\big)\Big|_a^b  + I_{(a,b)}\big(\frac{\varphi'}{2}v,v\big) \\
& 
- I_{(a,b)}\big(\frac{\varphi'}{2}
w,w\big)  -I_{(a,b)}\big(\frac{\varphi}{2} B'u, u\big) + I_{(a,b)}\big(({\rm Im}\,\psi)\varphi w,w\big)\\
= &\, {\rm Re}\, I_{(a,b)}\big(\varphi f,v\big) 
- I_{(a,b)}\big(({\rm Im}\,\psi)\varphi v,v\big) - {\rm Re}\, I_{(a,b)}\big(\varphi Vu,v\big).
\end{split}
\nonumber
\end{equation}
This proves the lemma. 
\end{proof}


\subsection{The 2nd identity}

 
\begin{lemma}\label{IntIdentity2}
The solution $u$ of  (\ref{S2Equation1}) satisfies
\begin{equation}
\begin{split}
{\rm Im}\,M(v,u)\Big|_a^b = &\ - 2\, I_{(a,b)}(({\rm Re}\,\psi\, {\rm Im}\,\psi)u,u) - M(({\rm Re}\,\psi) u,u)\Big|_a^b
+ I_{(a,b)}(u,({\rm Re}\,\psi')u)  \\
&  - {\rm Im}\,I_{(a,b)}(f,u) 
+ {\rm Im}\, I_{(a,b)}(Vu,u).
\end{split}
\nonumber
\end{equation}
\end{lemma}

\begin{proof}
Note that ${\rm Im}\,M_r(v,u) = {\rm Im}\,M_r(u',u) - M_r(({\rm Re}\,\psi) u,u)$, and 
$$
\partial_rM_r(v,u) = M_r(v',u) + M_r(v,u') + M_r(v,\frac{g'}{2g}u).
$$
Integrating over $(a,b)$, we have
\begin{equation}
M(v,u)\Big|_a^b = I_{(a,b)}(v',u) + I_{(a,b)}(v,u') + I_{(a,b)}(\frac{g'}{2g}v,u).
\label{Lemma5.3M(v,u)first}
\end{equation}
Using the equation (\ref{S3Eq1}), we have
\begin{equation}
\begin{split}
{\rm Im}\,I_{(a,b)}(v',u)  + {\rm Im}\,I_{(a,b)}(\frac{g'}{2g}v,u)
&= {\rm Im}\, I_{(a,b)}(\frac{g'}{4g}v,u) - {\rm Im}\,I_{(a,b)}(f,u) \\
& - {\rm Im}\, I_{(a,b)}(i\psi v,u) + {\rm Im}\, I_{(a,b)}(Vu,u).
\end{split}
\nonumber
\end{equation}
Hence, by (\ref{Lemma5.3M(v,u)first}), 
\begin{equation}
\begin{split}
{\rm Im}\, M(v,u)\Big|_a^b =&\  {\rm Im}\, I_{(a,b)}\big(v,u' + \frac{g'}{4g}u + i\overline{\psi}u\big)\\
& - {\rm Im}\, I_{(a,b)}(f,u) + {\rm Im}\, I_{(a,b)}(Vu,u) \\
=&\ {\rm Im}\, I_{(a,b)}(v,2i({\rm Re}\,\psi) u) - {\rm Im}\, I_{(a,b)}(f,u) + {\rm Im}\, I_{(a,b)}(Vu,u).
\end{split}
\label{Lemma5.3Eq1}
\end{equation}
We have, by integration by parts, 
\begin{equation}
\begin{split}
{\rm Im}\, I_{(a,b)}(v,2i({\rm Re}\,\psi)u) 
= & \ -2 {\rm Re}\, I_{(a,b)}(u',({\rm Re}\,\psi)u) - 2 I_{(a,b)}\big(\frac{g'}{4g}u,({\rm Re} \, \psi) u\big) \\
& - 2\, I_{(a,b)}\big(({\rm Re}\,\psi\, {\rm Im}\, \psi)u,u\big) \\
= & - M(u,({\rm Re}\,\psi)u)\Big|_a^b + I_{(a,b)}\big(u,({\rm Re}\,\psi')u\big) \\
& - 2\, I_{(a,b)}\big(({\rm Re}\,\psi\, {\rm Im}\, \psi)u,u\big).
\end{split}
\label{Lemma5.3Eq2}
\end{equation}
The lemma then follows from (\ref{Lemma5.3Eq1}) and (\ref{Lemma5.3Eq2}). 
 \end{proof}


\section{A priori-estimates  on each end}\label{SectionAprioiestimatesends}

We begin to estimate the resolvent of the Laplacian on each end. We assume (\ref{rho'rho-c0}) and (\ref{h-hMgamma}) for
\begin{equation}
\alpha_{0} > 0, \quad \gamma_{0} > 0,
\label{S6alpha0j>0gamma0j>0}
\end{equation}
omitting the subscript $j$. 
 We pick up one end $(0,\infty)\times M$, and consider the solution $u$ of $(- \Delta_{\mathcal M} - z)u = f$ with support in $(0,\infty)\times M$. For the sake of simplicity, we denote $(0,\infty)\times M$ as $\mathcal M$.
We use the same notations as in \S 5, and define the following  Besov type function spaces introduced by Agmon-H{\"o}rmander \cite{AgHo76} in the case of Euclidean space. 


\subsection{Function spaces}\label{spaceBBast}
For an inteval $I \subset (0,\infty)$, we put
$$
L^2(I) = L^2(I;{\bf h}(r);\rho^{n-1}(r)dr).
$$
Let $\mathcal B$ be the set of functions $f$ satisfying
\begin{equation}
\|f\|_{\mathcal B} = \sum_{j=0}^{\infty}2^{j/2}\|f\|_{L^2(I_j)}< \infty,
\label{S6Bnorm}
\end{equation}
where $I_0 = (0,1]$, $I_j = (2^{j-1},2^j]$, $j \geq 1$.
The dual space of $\mathcal B$ is identified with the set of functions $v(r)$ satisfying
$$
\|v\|^{\ast} = \sup_{j\geq0}\,2^{-j/2}\|v\|_{L^2(I_j)} < \infty.
$$
Letting
\begin{equation}
\|v\|_{\mathcal B^{\ast}} = \Big(\sup_{R>1}\frac{1}{R}\int_{0}^R\|v(r)\|_{{\bf h}(r)}^2\,\rho^{n-1}(r)dr\Big)^{1/2},
\label{S6Bastnorm}
\end{equation}
one can show the existence of a constant $C > 0$ such that
$$
C^{-1}\|v\|^{\ast} \leq \|v\|_{\mathcal B^{\ast}} \leq C\|v\|^{\ast}.
$$
Therefore, we employ  $\|\cdot\|_{\mathcal B^{\ast}}$ as the norm of the dual space of $\mathcal B$. Note that the coupling of $f \in \mathcal B$ and $v \in \mathcal B^{\ast}$ is given by
$$
(f,v) = \int_0^{\infty}(f(r),v(r))_{{\bf h}(r)}\,\rho^{n-1}(r)dr.
$$
We introduce a closed subspace $\mathcal B^{\ast}_0$  of $\mathcal B^{\ast}$ as follows :
\begin{equation}
\mathcal B^{\ast}_0 \ni v \Longleftrightarrow \lim_{R\to\infty}\frac{1}{R}\int_0^R\|v(r)\|_{{\bf h}(r)}^2\,\rho^{n-1}(r)dr = 0.
\label{Bast01dim}
\end{equation}
For $s \in {\mathbb R}$, let $L^{2,s}$ be the set of functions $v(r)$ satisfying
\begin{equation}
\|v\|_s = \left(\int_0^{\infty}\|v(r)\|_{{\bf h}(r)}^2\,(1 + r)^{2s}\rho^{n-1}(r)dr\right)^{1/2} < \infty.
\label{S6L2s}
\end{equation}
For $s = 0$, this norm is denoted by $\|\cdot\|$. 
For $s >1/2$, the following inclusion relations hold :
\begin{equation}
L^{2,s} \subset \mathcal B \subset L^{2,1/2} \subset L^2 \subset L^{2,-1/2} \subset \mathcal B^{\ast} 
\subset L^{2,-s}.
\label{S6InclusionRelation}
\end{equation}
For example, that $L^{2,-1/2} \subset \mathcal B^{\ast} 
\subset L^{2,-s}$ is proven as follows:
\begin{equation}
\begin{split}
\|f\|_{-s}^2 &\leq C\sum_{j=0}^{\infty}2^{-2sj}\|f\|_{L^2(I_j)}^2 \leq
C\Big(\sum_{j=0}^{\infty}2^{-j(2s-1)}\Big)\sup_{j\geq0}2^{-j}\|f\|^2_{L^2(I_j)} \\
& \leq 
 C\sum_{j=0}^{\infty}2^{-j}\|f\|^2_{L^2(I_j)} \leq C\|f\|^2_{-1/2}.
\end{split}
\nonumber
\end{equation}

\begin{lemma}\label{LemmaL2-1/2isBast0}
(1) \ $\ L^{2,-1/2} \subset \mathcal B^{\ast}_0$. \\
\noindent
(2) If $v \in \mathcal B^{\ast}_0$, we have
$\mathop{\underline{\lim}_{r\to\infty}}M_r(v,v) = 0.$\\
\noindent
(3) If $u \in \mathcal B^{\ast}$ and $v \in \mathcal B^{\ast}_0$, we have $\displaystyle{\underline{\lim}_{r\to\infty}|M_r(u,v)| = 0}$.
\end{lemma}

\begin{proof}
 We prove (1). For $v \in L^{2,-1/2}$ and  $\epsilon > 0$, there exists $R_{\epsilon} > 1$ such that
$$
\int_{R_{\epsilon}}^{\infty}\frac{1}{r}\|v(r)\|^2_{{\bf h}(r)}\rho^{n-1}(r)dr < \epsilon.
$$
Then, if $R > R_{\epsilon}$,
$$
\frac{1}{R}\int_{R_{\epsilon}}^R\|v(r)\|^2_{{\bf h}(r)}\rho^{n-1}(r)dr \leq 
\int_{R_{\epsilon}}^R\frac{1}{r}\|v(r)\|^2_{{\bf h}(r)}\rho^{n-1}(r)dr < \epsilon.
$$
Therefore, for $R > R_{\epsilon}$,
$$
\frac{1}{R}\int_{0}^R\|v(r)\|^2_{{\bf h}(r)}\rho^{n-1}(r)dr \leq 
\frac{1}{R}\int_{0}^{R_{\epsilon}}\|v(r)\|^2_{{\bf h}(r)}\rho^{n-1}(r)dr + \epsilon,
$$
which implies $\displaystyle{\overline{\lim}_{R\to\infty}}\frac{1}{R}\int_{0}^R\|v(r)\|^2_{{\bf h}(r)}\rho^{n-1}(r)dr \leq \epsilon$. 

If $\mathop{\underline{\lim}_{r\to\infty}}M_r(v,v) > 0$, there exist constants $C, r_0 > 0$ such that $M_r(v,v) > C$ for $r > r_0$, hence $v \not\in \mathcal B_0^{\ast}$.  This proves (2).

It follows from the inequality
\begin{equation}
\begin{split}
& \frac{1}{R}\int_0^R|(u(r),v(r))_{{\bf h}(r)}|\rho^{n-1}(r)dr \\
\leq & \Big(\frac{1}{R}\int_0^R\|u(r)\|_{{\bf h}(r)}^2\rho^{n-1}(r)dr \Big)^{1/2}
 \Big(\frac{1}{R}\int_0^R\|v(r)\|_{{\bf h}(r)}^2\rho^{n-1}(r)dr \Big)^{1/2}.
\end{split}
\nonumber
\end{equation}
that $\frac{1}{R}\int_0^R |(u(r),v(r))_{{\bf h}(r)}| \rho^{n-1}(r)dr \to 0$, which implies (3). 
\end{proof}

Finally, for a non-negative integer $m$,  the Sobolev space of order $m$ 
on $\mathcal M = (0,\infty)\times M$ is denoted by $H^m(\mathcal M)$. For an interval $I \subset (0,\infty)$, $H^m(I\times M)$ is defined similarly. For $s \in {\mathbb R}$, $H^{m,s}(\mathcal M)$ is the set of functions $v$ such that $(1 + r^2)^{s/2}v \in H^m(\mathcal M)$ equipped with norm
$$
\|v\|_{H^{m,s}} = \|(1 + r^2)^{s/2}v\|_{H^m(\mathcal M)}.
$$
 Similarly, for $s \in {\mathbb R}$, we introduce $\widetilde H^{2,s}(\mathcal M)$-norm by
\begin{equation}
\begin{split}
\|u\|^2_{\widetilde H^{2,s}(\mathcal M)} = 
\sum_{m\leq 2}\|( 1 + r)^{s}\partial_r^mu\|^2_{L^2(\mathcal M)}  + 
\int_0^{\infty}(1 + r)^{2s}\|u\|^2_{{\widetilde{\mathcal H}}^2(M)}\rho^{n-1}(r)dr,
\end{split}
\nonumber
\end{equation}
where $\|\cdot\|_{\widetilde H^2}$-norm is defined by (\ref{S2DefineH2tildenorm}).
For an interval $I \subset (0,\infty)$, $\widetilde{H}^{2,s}(I\times M)$ is defined similarly.


\subsection{A-priori estimates}
In this subsection, we derive some a-priori estimates for the solution $u$ to 
\begin{equation}
(- \Delta_{\mathcal M} - z)u = f, \quad {\rm on} \quad (0,\infty)\times M
\label{S6Equation}
\end{equation}
satisfying $u=0$ for $r < 1$. 


\begin{lemma}\label{S4Apriori1}
(1) For any  $a,\epsilon >0$, there exists a constant $C_{a,\epsilon}>0$  such that  
\begin{equation}
\|u\|_{\widetilde H^2((0,a)\times M)} \leq C_{a,\epsilon}\big(\|u\|_{L^2((0,a+\epsilon)\times M)} + \|f\|_{L^2((0,a+\epsilon)\times M)}\big).
\nonumber
\end{equation}
(2) For any $s, s' \in {\mathbb R}$ satisfying $s' \leq s$, there exists a constant $C_{s,s'} > 0$ such that if $u, f \in L^{2,s}$ and $u' \in L^{2,s'}$,
\begin{equation}
\|u\|_{\widetilde H^{2,s}}  + \|B(r)u\|_s 
 \leq C_{s,s'}\big(\|u\|_s + \|f\|_s\big).
\nonumber
\end{equation}
In the above inequalities, the constants $C_{a,\epsilon}$ and $C_{s,s'}$ are independent of $z$ if $z$ varies over a compact set in ${\mathbb C}$.
\end{lemma}

\begin{proof}
Take $\chi \in C^{\infty}({\mathbb R})$ such that $\chi(r) = 1$ for $r < a$ and $\chi(r) = 0$ for $r > a + \epsilon$. Applying Theorem \ref{S2ModifiedEllipticregularity} to $\chi u$, we obtain (1). To prove (2),  we have only to consider the case in which $s' = s-1$. 
Let $I_R = (0,Ra)$, $J_R = (0,R(a+\epsilon))$ and $K_R = (Ra,R(a+\epsilon))$. 
Applying  Theorem \ref{S2ModifiedEllipticregularity} to $(1 + r)^{s}\chi (r/R)u$, we have
\begin{equation}
\|u\|_{\widetilde H^{2,s}(I_R\times M)} \leq C\big(\|Hu\|_{L^{2,s}(J_R\times M)} + \|u\|_{L^{2,s}(J_R\times M)} + \|u'\|_{L^{2,s-1}(J_R\times M)}\big).
\label{S6Lemma6.Inequality}
\end{equation}
Here, we note $ \|u'\|^2_{L^{2,s-1}(J_R\times M)} = \|u'\|^2_{L^{2,s-1}(I_R\times M)} + \|u'\|^2_{L^{2,s-1}(K_R\times M)}$ and
\begin{equation}
\begin{split}
\int_0^R\int_M(1 + r)^{2(s-1)}u(r)\overline{u''(r)}\sqrt{g}drdx &= 
\int_M(1 + R)^{2(s-1)}u(R)\overline{u'(R)}\sqrt{g}dx \\
& - \int_0^R\int_M\left((1+r)^{2(s-1)}u\sqrt{g}\right)'
\overline{u'(r)}drdx.
\end{split}
\nonumber
\end{equation}
This implies
\begin{equation}
\begin{split}
\|u'\|^2_{L^{2,s-1}(I_R\times M)} \leq \ & 2 (1 + R)^{2(s-1)}|M_R(u,u')| \\
& + \epsilon\|u''\|^2_{L^{2,s-1}(I_R\times M)} + C_{\epsilon}\|u\|^2_{L^{2,s-1}(I_R\times M)}. 
\end{split}
\nonumber
\end{equation}
Letting $R \to \infty$ in (\ref{S6Lemma6.Inequality}) along a suitable sequence, we then obtain
$$
\|u\|_{\widetilde H^{2,s}(\mathcal M)} \leq C\big(\|Hu\|_{L^{2,s}(\mathcal M)} + \|u\|_{L^{2,s}(\mathcal M)}\big).
$$
Using the equation (\ref{S2Equation1}), we also have
$$
\|B(r)u\|_{L^{2,s}(\mathcal M)} \leq C\big(\|Hu\|_{L^{2,s}(\mathcal M)} + \|u\|_{L^{2,s}(\mathcal M)}\big).
$$
These two inequalities imply (2). 
\end{proof}

We define the $\|\cdot\|_{{\mathcal B}^{\ast}_m}$ norm by
\begin{equation}
\|u\|_{{\mathcal B}^{\ast}_m} = \sum_{|\alpha|\leq m}
\|\partial^{\alpha}u\|_{\mathcal B^{\ast}}.
\nonumber
\end{equation}

\begin{lemma}\label{Lemma6.2Bastapriori}
If $u, u', f \in \mathcal B^{\ast}$, we have
$$
\|u\|_{\mathcal B^{\ast}_1} \leq C\big(\|f\|_{\mathcal B^{\ast}} + \|u\|_{\mathcal B^{\ast}}\big),
$$
where the constant $C$ is independent of $z$ when $z$ varies over a compact set in ${\mathbb C}$. 
\end{lemma}

\begin{proof}
Take $\chi \in C^{\infty}({\mathbb R})$ such that $\chi(t) = 1$ for $t < 1$, $\chi(t) = 0$ for $t > 2$, and put $u_R = \frac{1}{\sqrt R}\chi_R(r)u$, where $\chi_R(r) = \chi(r/R)$. Then, we have
\begin{equation}
(- \Delta_{\mathcal M} - z)u_R = \frac{1}{\sqrt R}\chi_Rf 
- \frac{1}{\sqrt R}\big(2\chi_R'u' + \frac{g'}{2g}\chi'_Ru + \chi_R''u\big).
\nonumber
\end{equation}
Taking the inner product with $u_R$ and integrating by parts, we have
$$
\frac{1}{R}\|(\chi_Ru)'\|^2 + \frac{1}{R}\|\sqrt{B}\chi_Ru\|^2 \leq C\big(\|f\|^2 _{\mathcal B^{\ast}} + \|u\|^2 _{\mathcal B^{\ast}} + \frac{1}{R^2}|(c_R(r)u',u)|\big),
$$
where $c_R(r)$ is the characteristic function of $(0,2R)$. For any $\epsilon >0$, there exists $R_{\epsilon} > 0$ such that
$$
\frac{1}{R^2}|(c_R(r)u',u)| \leq \epsilon \big( \|u'\|^2_{\mathcal B^{\ast}}+ \|u\|^2_{\mathcal B^{\ast}}\big), \quad {\rm for} \quad R > R_{\epsilon}.
$$
For $1 < R < R_{\epsilon}$, applying Lemma \ref{S4Apriori1} (1), we have
$$
\frac{1}{R^2}|(c_R(r)u',u)| \leq C_{\epsilon}\big(\|f\|^2_{\mathcal B^{\ast}} + \|u\|^2_{\mathcal B^{\ast}}\big).
$$
We then have
$$
\frac{1}{R}\|(\chi_Ru)'\|^2 + \frac{1}{R}\|\sqrt{B}\chi_Ru\|^2 \leq \epsilon \|u'\|^2_{\mathcal B^{\ast}} + C_{\epsilon}\big(\|f\|^2 _{\mathcal B^{\ast}} + \|u\|^2 _{\mathcal B^{\ast}}\big),
$$
which yields the lemma. 
\end{proof}

\medskip
For energies below $E_0 = \big((n-1)c_0/2\big)^2$, the following inequality holds.


\begin{lemma}\label{S7Lowenerguestimate}
Take any compact interval $[c_1,c_2] \subset (-\infty,E_0)$. Let $s \geq 0$, $s' = s -  \min\{\alpha_0,1\}$. Then, there exists a constant $C_{s}>0$ such that for a solution $u$ to (\ref{S6Equation}), if ${\rm Re}\, z \in [c_1,c_2]$, and  $u, u' \in L^{2,s'}$, $f \in L^{2,s}$,
$$
\|u\|_{H^{1,s}}  \leq C_{s}(\|f\|_s + \|u\|_{s'}).
$$
\end{lemma}

\begin{proof}
 By (\ref{S5partialrast}) and (\ref{S5DeltaM-zrewrittena1}), we have
$$
- \Delta_{\mathcal M} - z = \big(\partial_r + \frac{g'}{4g}\big)^{\ast}\big(\partial_r + \frac{g'}{4g}\big) + B(r) + Q + E_0 -z.
$$
Hence for any $h \in D(H)$,
$$
{\rm Re}\,\big((- \Delta_{\mathcal M} - z)h,h\big) 
\geq (E_0-c_2)\|h\|^2 + (Qh,h).
$$
Take $\chi \in C^{\infty}({\mathbb R})$ such that $\chi(t) = 1$ for $t < 1$ and $\chi(t) = 0$ for $t > 2$. Then, letting $u_R = \chi(r/R)u$, and $h = (1 + r)^su_R$, we have by (\ref{S5QinS-gamma})
$$
\|u_R\|_s^2 \leq \epsilon \|u_R\|_s^2 + C_{\epsilon}(\|f\|_s^2 + \|u\|_{s'}^2 + \|u'\|_{s'}^2).
$$ 
Letting $R \to \infty$, and using Lemma \ref{S4Apriori1}, we obtain the lemma.
\end{proof}


\subsection{Regular ends} 
We begin with the case of regular ends.  
In addition to (\ref{S6alpha0j>0gamma0j>0}), we assume (\ref{S1growthCond}), which implies
\begin{equation}
\rho(r)^{-1} \leq C(1 + r)^{-\beta_0}, \quad \beta_0 > 0.
\label{S6.3Assumption}
\end{equation}
Moreover, there exists a constant $r_0>0$ such that 
\begin{equation}
\rho'(r) > 0 \quad {\rm for} \quad r > r_0.
\label{rho'/rhaopm}
\end{equation}
We fix an interval $[c_1, c_2] \subset (E_0,\infty)$, and  put
\begin{equation}
J_{\pm} = \big\{z \in {\mathbb C}\, ; \, c_1 \leq {\rm Re}\,z \leq c_2, \   0 < \pm {\rm Im}\,z < 1\big\}.
\label{S4DefineJpm}
\end{equation}
We derive  a-priori estimates for  a solution $u$ to (\ref{S2Equation1}) satisfying 
\begin{equation}
f \in \mathcal B, \quad u, u' \in \mathcal B^{\ast}, \quad z \in J_{\pm}.
\nonumber
\end{equation}
 Using (\ref{DefDpmk}), we put as in the previous section
$$
v = D(k)u, \quad 
w = \sqrt{B(r)}u,
$$
\index{$v = D(k)u$} \index{$w = \sqrt{B(r)}u$}
 and assume that 
\begin{equation}
u = 0 \quad {\rm for} \quad r < 1,
\label{Condu=0r<1}
\end{equation} 
\begin{equation}
v \in \mathcal B^{\ast}_0.
\label{CondMrvv=01}
\end{equation}
In the following, $C$'s denote  constants  independent of $z \in J_{\pm}$.


\begin{lemma} \label{Lemma6.5IabIneq1}
Let $0 < \delta < 2\beta_0$. Then, for any non-negative $C^1$-function $\varphi(r)$, we have
\begin{equation}
\begin{split}
&I_{(a,b)}(\varphi'v,v) + I_{(a,b)}\big((\delta \varphi r^{-1} - \varphi')w,w\big)  + M(\varphi w,w)\Big|_a^b - M(\varphi v,v)\Big|_a^b \\
&  \leq  \, C\Big(I_{(a,b)}(\varphi r^{-1-\alpha_0}v,v) + I_{(a,b)}(\varphi r^{-1-\alpha_0}w,w)   + |I_{(a,b)}(\varphi f,v)| \\
& \ \ \ \ \ \ \ \ \ \ \ + |I_{(a,b)}(\varphi Vu,v)| + 
I_{(a,b)}(\varphi r^{-1-\gamma_0}u,u)\Big).
\end{split}
\nonumber
\end{equation}
\end{lemma}

\begin{proof}
Adding $I_{(a,b)}\big(\frac{\delta\varphi}{r}Bu,u\big)$ to the both sides of 
(\ref{EqualityinLemmaIntIdentity1}) 
\begin{equation}
\begin{split}
& I_{(a,b)}\big(\varphi'v,v\big) + I_{(a,b)}\big(\big(\frac{\delta\varphi}{r} -\varphi'\big) w,w\big) + M\big(\varphi w,w\big)\Big|_a^b - M\big({\varphi}v,v\big)\Big|_a^b\\
 =& - 2I_{(a,b)}\big(({\rm Im}\,\psi)\varphi v,v\big) - 2I_{(a,b)}\big(({\rm Im}\,\psi)\varphi w,w\big) \\
& + 2{\rm Re}\,I_{(a,b)}\big(\varphi f,v\big) - 2{\rm Re}\, I_{(a,b)}\big(\varphi Vu,v\big)  + I_{(a,b)}\big(\varphi \big(B' + \frac{\delta}{r}B\big)u, u\big).
\end{split}
\label{EqualityinLemmaIntIdentity12} 
\end{equation}
The proof of Lemma \ref{LemmatB(t)leqct-epsilon} works also for $B(r)$. This and 
 Lemma \ref{psimlemma} (3)  imply 
$$
- {\rm Im}\,\psi \leq Cr^{-1-\alpha_0}, \quad B' + \frac{\delta}{r}B \leq Cr^{-1-\gamma_0}
$$
for  $0 < \delta < 2\beta_0$. The lemma then follows from  (\ref{EqualityinLemmaIntIdentity12}). 
\end{proof}


\begin{lemma}\label{Lemma6.3w-1/2leqfu}
Let $\epsilon_0 = \min(\alpha_0,\gamma_0)$ and $\frac{1}{2} <s \leq \frac{1}{2} (1 + \epsilon_0)$. Then, we have
\begin{equation}
\|v\|_{\mathcal B^{\ast}} + \|w\|_{-1/2} \leq C\big(\|f\|_{\mathcal B} +  \|u\|_{-s}\big).
\nonumber
\end{equation}
\end{lemma}

\begin{proof}
Letting $\varphi=1$ and $a = 0$ in Lemma \ref{Lemma6.5IabIneq1}, we have
\begin{equation}
\begin{split}
 \delta I_{(0,b)}(\frac{1}{r}w,w)   - M_b(v,v)
 \leq    C\big(\|v\|_{-s}^2 + \|w\|^2_{-s} + |(f,v)| + \|u\|_{-s}^2\big).
\end{split}
\nonumber
\end{equation}
Note that by Lemma \ref{psimlemma}, $V$ decays sufficiently rapidly.  We let $b \to \infty$ along a suitable sequence in the above inequality. By (\ref{CondMrvv=01}) and Lemma \ref{LemmaL2-1/2isBast0}, $M_b(v,v) \to 0$. Using Lemma  \ref{S4Apriori1}, we obtain
\begin{equation}
\|	w\|^2_{-1/2} \leq C\big(|(f,v)| + \|u\|^2_{-s}\big).
\label{|w|-1/2leqfvu-s}
\end{equation}  
Letting $\varphi=1$, $a= r$ and $b\to \infty$ along a suitable sequence in Lemma \ref{Lemma6.5IabIneq1}, we have
\begin{equation}
\begin{split}
  M_r(v,v) \leq  M_r(w,w)  +   C\big(\|v\|_{-s}^2 + \|w\|^2_{-s} + |(f,v)| + \|u\|_{-s}^2\big).
\end{split}
\nonumber
\end{equation}
By taking the integral mean in $r$,
\begin{equation}
\|v\|^2_{\mathcal B^{\ast}} \leq \|w\|^2_{\mathcal B^{\ast}} + C\big(\|v\|_{-s}^2 + \|w\|^2_{-s} + |(f,v)| + \|u\|_{-s}^2\big).
\nonumber
\end{equation}
Using Lemma \ref{S4Apriori1} (2), we then have
\begin{equation}
\|v\|^2_{\mathcal B^{\ast}} \leq C\Big(\|w\|_{\mathcal B^{\ast}}^2 + |(f,v)| + \|u\|_{-s}^2\Big).
\nonumber
\end{equation}
Since $\|w\|_{\mathcal B^{\ast}} \leq C\|w\|_{-1/2}$, using (\ref{|w|-1/2leqfvu-s}), we have
$$
\|v\|^2_{\mathcal B^{\ast}} \leq C\big(|(f,v)| + \|u\|_{-s}^2\big).
$$
Finally, using
$$
|(f,v)| \leq \epsilon \|v\|^2_{\mathcal B^{\ast}} + C_{\epsilon}\|f\|^2_{\mathcal B},
$$
we obtain the lemma.
\end{proof}

\medskip
Lemma \ref{Lemma6.3w-1/2leqfu} can be improved as follows.


\begin{lemma}\label{Lemma6.7vestimateimprove}
For  any $0 < t < {\rm min}\, (2\beta_0,\gamma_0)$, there exists a constant $C > 0$ such that if $\displaystyle{\liminf_{r\to\infty}M_r(r^tv,v)=0}$,
$$
\|v\|_{(t-1)/2} + \|w\|_{(t-1)/2} \leq C\big(\|f\|_{(t + 1)/2} + \|u\|_{(t-1-\gamma_0)/2}\big).
$$
\end{lemma}

\begin{proof}
Take $t < \delta < 2\beta_0$ and put $\varphi = r^{t}$, $a=0$ in the inequality in Lemma  \ref{Lemma6.5IabIneq1}. We drop the term $M_b(\varphi w,w)$ from the left-hand side. For large $b > 0$, the first two terms of the right-hand side are absorbed into the left-hand side. For the 3rd and 4th terms, we apply the inequalities
$$
|I_{(0,b)}(\varphi f,v)| \leq \epsilon I_{(0,b)}(r^{t-1}v,v) + C_{\epsilon}\|f\|_{(t+1)/2},
$$
$$
|I_{(0,b)}(\varphi Vu,v)| \leq \epsilon I_{(0,b)}(r^{t-1}v,v)  + C_{\epsilon}\|u\|^2_{-s}, \quad s = (t-1-\gamma_0)/2.
$$
Here, we have used Lemma \ref{psimlemma} (2). Then, the term $\epsilon I_{(0,b)}(r^{t-1}v,v)$ is absorbed into the left-hand side. Letting $b \to \infty$, we obtain the lemma.
 \end{proof}

\begin{remark}
If $2\beta_0 \geq \gamma_0 > 1$, one can take 
$1 < t < \gamma_0$ in Lemma \ref{Lemma6.7vestimateimprove}. This is important to study the behavior of the resolvent as $r \to \infty$ (Lemma \ref{LemmaDpm(k9uvto0}), and explains the appearance of border-lines $\beta_0 = 1/2$ and $\gamma_0 = 1$.
\end{remark}

\begin{lemma}\label{Lemmauestimate}
For any $s > 1/2$, there exists a constant $C > 0$ such that 
$$
\|u\|_{\mathcal B^{\ast}} \leq C(\|f\|_{\mathcal B} + 
\|u\|_{-s}).
$$
\end{lemma}

\begin{proof}
By Lemma \ref{IntIdentity2} 
\begin{equation}
\begin{split}
- {\rm Im}\,M(v, u)\Big|_a^b = &\ \,M(({\rm Re}\,\psi)  u,u)\Big|_a^b 
+ 2I_{(a,b)}(({\rm Re}\,\psi{\rm Im}\,\psi)u,u) \\
& \ 
- I_{(a,b)}(u,({\rm Re}\, \psi')u) + {\rm Im}\, I_{(a,b)}( f, u) 
- {\rm Im}\,I_{(a,b)}(Vu,u).
\end{split}
\nonumber
\end{equation}
For $z \in J_{+}$, ${\rm Re}\,\psi\, {\rm Im}\,\psi \geq 0$ by (\ref{S5Definepsimrxk}). Hence
\begin{equation}
\begin{split}
M(({\rm Re}\,\psi) u, u)\Big|_a^b \leq& \ - {\rm Im}\,M( v,u)\Big|_a^b 
+ I_{(a,b)}(u,({\rm Re}\, \psi')u) \\
& \ \ - {\rm Im}\, I_{(a,b)}( f, u) 
+ {\rm Im}\,I_{(a,b)}(Vu,u).
\end{split}
\label{Lemma6.7Equation1}
\end{equation}
We let $0= a < 2 < b = r$. By Lemma \ref{psimlemma}, there exist constants $C, r_0 > 0$ such that
$$
C \leq {\rm Re}\,\psi, \quad {\rm for}\quad r > r_0.
$$
Since $\psi' \in S^{-1-\alpha_0}$ by Lemma \ref{psimlemma} (1), we have by (\ref{Lemma6.7Equation1})
$$
M_r( u,u)  \leq C\big(|M_r( v,u)| +  \|f\|_{\mathcal B}\|u\|_{\mathcal B^{\ast}} + \|u\|^2_{-s}\big),
$$
which, together with $|M_r(v,u)| \leq \epsilon M_r(u,u) + C_{\epsilon}M_r(v,v)$,  yields
$$
M_r(u,u) \leq C\Big(M_r(v,v) + \|f\|_{\mathcal B}\|u\|_{\mathcal B^{\ast}} + \|u\|^2_{-s}\Big).
$$
Integrating with respect to $r$ over $(0,R)$ and dividing by $R$, we have
$$
\|u\|_{\mathcal B^{\ast}} \leq C(\|v\|_{\mathcal B^{\ast}} + \|f\|_{\mathcal B} + \|u\|_{-s}).
$$
The lemma then follows from Lemma \ref{Lemma6.3w-1/2leqfu}. 
\end{proof}


\subsection{Cusp ends}
We consider the case in which $\mathcal M = (0,\infty)\times M$ is a cusp. We assume as in the beginning of this section.
By (\ref{A-2CuspCondition}), there exists a constant $C_0 > 0$ such that 
\begin{equation}
\rho(r)^{-1} \geq C_0(1 + r)^{|\beta|}.
\label{rho'/rhaopm1}
\end{equation}
Let $\Lambda(r)$ and $\Lambda_0$ be the Laplace-Beltrami operators on $M$ associated with the metric $h(r,x,dx)$ and $h_M(x,dx)$, respectively. We put $B(r) = -\rho(r)^{-2}\Lambda(r)$. Let $P_0(r)$ be the projection associated with the $0$ eigenvalue of $\Lambda(r)$, and put $P_1(r) = 1 - P_0(r)$. 


\begin{lemma}
\label{S6P0(r)estimate}
There exists $r_0 > 0$ such that for $r > r_0$, $P_0(r)$ satisfies
\begin{equation}
\Big\|\big(\frac{d}{dr}\big)^m P_0(r)\Big\|  \leq 
C_m (1 + r)^{-m - \gamma_0}, \quad \forall m \geq 1,
\label{S6P0derivativedecay}
\end{equation}
\begin{equation}
\Big\|\Big[ P_0(r),\frac{g'}{g}\Big]\Big\| \leq 
C_m (1 + r)^{-1 - \gamma_0}. 
\label{S6P0commutatordecay}
\end{equation}
The same inequalities hold for $P_1(r)$.
\end{lemma}

\begin{proof}
Since $-\Lambda(r)$ and $-\Lambda_0$ have compact resolvents, and the 2nd eigenvalue of $-\Lambda_0$ is positive, there exists a constant $\delta_0 > 0$ and $r_0 > 0$ such that the second eigenvalue of $-\Lambda(r)$ is greater than $2\delta_0$ for all $r > r_0$. We take $\varphi \in C_0^{\infty}({\mathbb R})$ such that $\varphi(r) = 1$ for $|r| < \delta_0/2$, $\varphi(r) = 0$ for $|r| >\delta_0$. Then, we have
$P_0(r) = \varphi(-\Lambda(r))$ for $r > r_0$. The assertions for $P_0(r)$ then follows from 
 Lemma \ref{varphiLambdaderivative}. Since $P_1(r) = 1 - P_0(r)$, the lemma also holds for $P_1(r)$.
\end{proof}


\begin{lemma}\label{Cuspestimate}
Let $u$ be a solution to the equation
$$
\Big( - \partial_r^2 - \frac{g'}{2g}\partial_r + B(r) - z\Big)u =f, \quad 
z \in J_{\pm}
$$
satisfying $u = 0$ for $r < 2r_0$, 
and put $u_0 = P_0(r)u, u_1 = P_1(r)u$. 
\\
\noindent
(1) For any $s' > 1/2$, there exists a constant $C > 0$ such that for $s = |\beta| + 1 + \gamma_0 - s'$
\begin{equation}
\|u_1'\|_{s} +  \|\sqrt{B(r)}u_1\|_{s} + \|u_1\|_{s} \leq C\big(\|P_1(r)f\|_{s} + \|u\|_{-s'}\big),
\nonumber
\end{equation}
\begin{equation}
\|u_1''\|_{s} + \|B(r)u_1\|_{s} \leq C\big(\|P_1(r)f\|_{s} + \|u\|_{-s'}\big).
\nonumber
\end{equation}
(2) For $s' > 1/2$
$$
\|u_0\|_{\mathcal B^{\ast}} +\|D(k)u_0\|_{\mathcal B^{\ast}} \leq C\big(\|P_0(r)f\|_{\mathcal B} + \|u\|_{-s'}\big),
$$
(3) For $0 < t < \gamma_0$,  there exists a constant $C > 0$ such that
\begin{equation}
\|D(k)u_0\|_{(t-1)/2} \leq C\big(\|P_0(r)f\|_{(t+1)/2} + \|u\|_{(t-1 - \gamma_0)/2}\big). 
\nonumber
\end{equation}
\end{lemma}

\begin{proof}
Letting $B = B(r), P_0 = P_0(r), P_1 = P_1(r)$, we have
\begin{equation}
\begin{split}
& - u_0'' -\frac{g'}{2g}u_0'  - zu_0 = f_0, \\
& f_0  = P_0f 
- P_0''u - 2P_0'u' - \frac{g'}{2g}P_0'u + \Big[P_0,\frac{g'}{2g}\Big]u',
\end{split}
\label{CuspEqforu0}
\end{equation}
\begin{equation}
\begin{split}
& - u_1'' -\frac{g'}{2g}u_1'  + P_1Bu_1 - zu_1 = f_1, \\
& f_1 = P_1f 
- P_1''u - 2P_1'u' - \frac{g'}{2g}P_1'u + \Big[P_1,\frac{g'}{2g}\Big]u'.
\end{split}
\label{CuspEqforu1}
\end{equation}
Assume that $u \in L^{2,s}$ for some $s$. 
Integration by parts in the equation (\ref{CuspEqforu1}) gives 
\begin{equation}
\begin{split}
& \big((1+r)^{2s}u_1',u_1'\big) + {\rm Re}\,(2s(1+r)^{2s-1}u_1',u_1) + \big((1 + r)^{2s}Bu_1,u_1\big) \\
& = {\rm Re}\,z\, \|u_1\|_{s}^2 + {\rm Re}\,((1 +r)^{2s}f_1,u_1),
 \end{split}
\nonumber
\end{equation}
which yields
\begin{equation}
\|u_1'\|^2_{s} - 2s\|u_1'\|_{s}\|u_1\|_{s-1}+ \|\sqrt{B}u_1\|_{s}^2  \leq C\big(\|f_1\|^2_{s} + \|u_1\|^2_{s}\big).
\label{Lemma6.10ProofInequality1}
\end{equation}
Noting that $P_1(-\Lambda) \geq \delta P_1$,   by (\ref{rho'/rhaopm1}) we have
\begin{equation}
(Bu_1,u_1)_{{\bf h}(r)} \geq C(1+r)^{2|\beta|}\|u_1\|_{{\bf h}(r)}^2, \quad C > 0.
\nonumber
\end{equation}
This and  (\ref{Lemma6.10ProofInequality1}) imply
\begin{equation}
\|u_1'\|_{s} + \|u_1\|_{|\beta| +s} \leq C\big(\|f_1\|_{s} + \|u_1\|_{s}\big).
\nonumber
\end{equation}
By Lemma \ref{S6P0(r)estimate},
$$
\|f_1\|_s \leq C(\|f\|_s + \|u\|_{s-1-\gamma_0} + \|u'\|_{s - 1 - \gamma_0}).
$$
In view of Lemma \ref{S4Apriori1} (2), we then have as long as 
$s - 1 - \gamma_0 \leq - s'$,
\begin{equation}
\|u_1'\|_{|\beta|+s} + \|u_1\|_{|\beta| +s} \leq C\big(\|f\|_{|\beta|+s} + \|u\|_{-s'}\big).
\nonumber
\end{equation}
Starting from $s = -s'$, we obtain for 
$s = |\beta| + 1 + \gamma_0 - s'$,
\begin{equation}
\|u_1'\|_{s} + \|u_1\|_{s} \leq C\big(\|f\|_{s} + \|u\|_{-s'}\big).
\nonumber
\end{equation}
We  have also obtained
\begin{equation}
\|u_1'\|_{s} + \|\sqrt{B(r)}u_1\|_{s} \leq C_s\big(\|P_1f\|_{s} + \|u_1\|_{-s'}\big).
\nonumber
\end{equation}
Arguing in the same way as in the proof of Lemma \ref{Lemma2.10}, we can show
$$
\|u_1''\|_s + \|B(r)u_1\|_s \leq C(\|f_1\|_s + \|u_1\|_s) \leq
C(\|f_1\|_s + \|u_1\|_{-s'}).
$$
We  have thus proven the assertion (1).

As for $u_0$,  since $Bu_0 = 0$,  the problem is reduced to the one-dimensional case, and we can argue in the same way as in the previous subsection. This proves the  assertions (2) and (3). 
\end{proof}

Let us note that Lemma \ref{Lemmauestimate} also holds for the cusp end, since the proof for the regular end applies to $u_0$ as well and $u_1$ belongs to $L^{2, -1/2 + \epsilon}$ for some $\epsilon > 0$ by virtue of Lemma \ref{Cuspestimate} (1).


\begin{lemma}\label{LemmauestimateCusp}
For any $s > 1/2$, there exists a constant $C > 0$ such that 
$$
\|u\|_{\mathcal B^{\ast}} \leq C(\|f\|_{\mathcal B} + 
\|u\|_{-s}).
$$
\end{lemma}

\bigskip


\section{Spectrum of the Laplacian on $\mathcal M$}\label{SctionSpectreLaplacian}

We are now ready to study the spectral theory of $- \Delta_{\mathcal M}$.  Let $\mathcal M$ be a connected $n$-dimensional CMAG of the  form
\begin{equation}
 \mathcal M = \mathcal K\cup\mathcal M_1\cup\cdots\cup\mathcal M_{N+N'}
\label{M=kM1toMN+N'}
\end{equation}
satisfying (\ref{rho'rho-c0}), (\ref{h-hMgamma}), (\ref{S1growthCond}) and (\ref{A-2CuspCondition}). 
In this section, we assume
\begin{equation}
\alpha_{0,i} > 0, \quad \beta_{0,i} > 0,  \quad \gamma_{0,i} > 0
\label{betaigammaiassume}
\end{equation}
for all regular ends and (A-3) for cusp ends.


\subsection{Essential spectrum} \label{SubseqEssSpectrum}
Assuming that 
\begin{equation}
\mathcal K \cap \big((1/2,\infty)\times M_i\big) = \emptyset, \quad  i \neq 0,
\nonumber
\end{equation}
we take a partition of unity $\{\chi_i\}_{i=0}^{N+N'}$ on $\mathcal M$ such that
\begin{equation}
\left\{
\begin{split}
& \sum_{i=0}^{N+N'}\chi_i = 1, \quad {\rm on} \quad \mathcal M,\\
& \chi_0 \in C_0^{\infty}(\mathcal M), \quad \chi_0 = 1 \quad {\rm on} \quad \mathcal K, \\
& \chi_i = 0 \quad {\rm on} \quad (0,2)\times M_i, \quad i \neq 0,\\
&\chi_i = 1 \quad {\rm on} \quad (3,\infty)\times M_i, \quad i \neq 0.
\end{split}
\right.
\label{S8Partitionunity}
\end{equation} 
We also take $\widetilde\chi_i \in C^{\infty}((0,\infty)\times M_i)$ such that
\begin{equation}
\left\{
\begin{split}
& \widetilde\chi_i = 1, \quad {\rm  on} \quad (1,\infty)\times M_i, \\
& \widetilde\chi_i = 0, \quad {\rm  on} \quad (0,1/2)\times M_i.
\end{split}
\right.
\nonumber
\end{equation}

  For $i=1,\cdots,N$ and $s \in {\mathbb R}$, the Banach spaces $L^{2,s}(\mathcal M_i)$, $\mathcal B(\mathcal M_i)$, 
$\mathcal B^{\ast}(\mathcal M_i)$, $\mathcal B^{\ast}_0(\mathcal M_i)$ are defined in the same way as in Subsection \ref{spaceBBast}. For $i = N+1,\cdots,N+N'$, these space are defined similarly with 
${\bf h}(r)$ replaced by ${\mathbb C}$. We put for $f \in L^2_{loc}({\mathcal M})$ and $s \in {\mathbb R}$,
\begin{equation}
\|f\|_{L^{2,s}(\mathcal M)} = \|\chi_0 f\|_{L^2(\mathcal M)} + \sum_{i=1}^{N+N'}\|\chi_if\|_{L^{2,s}(\mathcal M_i)},
\nonumber
\end{equation}
\begin{equation}
\|f\|_{\mathcal B(\mathcal M)} =  \|\chi_0 f\|_{L^2(\mathcal M)} + \sum_{i=1}^{N+N'}\|\chi_if\|_{\mathcal B(\mathcal M_i)},
\nonumber
\end{equation}
\begin{equation}
\|f\|_{\mathcal B^{\ast}(\mathcal M)} =  \|\chi_0 f\|_{L^2(\mathcal M)} + \sum_{i=1}^{N+N'}\|\chi_if\|_{\mathcal B^{\ast}(\mathcal M_i)},
\nonumber
\end{equation}
and also define
\begin{equation}
f \in \mathcal B_0^{\ast}(\mathcal  M) \Longleftrightarrow 
\chi_i f \in \mathcal B_0^{\ast}(\mathcal M_i), \quad  1 \leq \forall i \leq N+N'.
\nonumber
\end{equation}
In the following, we often denote these norms by $\|f\|_s, \|f\|_{\mathcal B}, \|f\|_{\mathcal B^{\ast}}$ omitting the end $\mathcal M_i$ or $\mathcal M$, which will not confuse our arguments. Finally, we define $L^2_{comp}(\mathcal M)$ to be the set of all compactly supported $L^2$-functions on $\mathcal M$.

We define
\begin{equation}
H = - \Delta_{\mathcal M}
\nonumber
\end{equation}
to be the Friedrichs extension of $-\Delta_{\mathcal M}$. Therefore, $D(\sqrt{H}) = H^1(\mathcal M)$, the Sobolev space of order 1 on $\mathcal M$. 
We introduce the following two operators defined on $\widetilde{\mathcal M}_i = (2,\infty)\times M_i$ with Dirichlet boundary condition at $r=2$ :
\begin{equation}
H_i = - \Delta_{\mathcal M}, \quad 
 H_{free(i)} = - \partial_r^2 - \frac{(n-1)\rho_i'(r)}{\rho_i(r)}\partial_r - \rho_i(r)^{-2}\Lambda_{i},
\label{Hfreei}
\end{equation}
where $\Lambda_{i}$ is the Laplace-Beltrami operator on $M_i$ equipped with the metric $h_{M_i}(x,dx)$. 
Let for $z \not\in {\mathbb R}$
\begin{equation}
R(z) = (H -z)^{-1}, \quad 
R_i(z) = (H_i-z)^{-1}, \quad
R_{free(i)}(z) = (H_{free(i)}-z)^{-1}.
\nonumber
\end{equation}
Note that 
\begin{equation}
R(z) \in {\bf B}(L^2(\mathcal M);H^1(\mathcal M)), \quad 
R_i(z),  R_{free(i)}(z) \in {\bf B}(L^2(\widetilde{\mathcal M}_i);H^1(\widetilde{\mathcal M}_i)).
\nonumber
\end{equation}
The following formula holds:
\begin{equation}
R(z) = \sum_{i=1}^{N+N'}\chi_iR_i(z)\widetilde\chi_i + R(z)(\chi_0 - S(z)),
\label{R(z)Reprsent1}
\end{equation}
\begin{equation}
S(z) = \sum_{i=1}^{N+N'}S_i(z)\widetilde \chi_i,
\quad
S_i(z) = [H,\chi_i]R_i(z).
\nonumber
\end{equation}
In fact, we have
$$
(H-z)\sum_{i=1}^{N+N'}\chi_iR_i(z)\widetilde \chi_i = \sum_{i=1}^{N+N'}S_i(z)\widetilde \chi_i + 
\sum_{i=1}^{N+N'}\chi_i\widetilde\chi_i,
$$
from which (\ref{R(z)Reprsent1}) follows. Similary, one can show the following formula.
\begin{equation}
R(z) = \sum_{i=1}^{N+N'}\chi_iR_{free(i)}(z)\widetilde\chi_i + R(z)(\chi_0 - T(z)),
\label{S7Rz=Rfreei+chi0}
\end{equation}
\begin{equation}
T(z) = \sum_{i=1}^{N+N'}T_i(z)\widetilde \chi_i,
\nonumber
\end{equation}
\begin{equation}
T_i(z) = [H,\chi_i]R_{free(i)}(z) + \chi_i(H - H_{free(i)})\widetilde\chi_iR_{free(i)}(z).
\nonumber
\end{equation}

We put 
\begin{equation}
E_{0,i} = \Big(\frac{(n-1)c_{0,i}}{2}\Big)^2, \quad E_{0,tot} = \min_{1\leq i \leq N+N'}E_{0,i}.
\label{DefE0i}
\end{equation}
\index{$E_{0,i}$} \index{$E_{0,tot}$}


\begin{lemma}
\label{EssentialSpectrum}
(1) $\ \sigma_e(H_i) = [E_{0,i},\infty)$, $\ i = 1,\cdots, N + N'.$ \\
\noindent
(2) 
$\ 
\sigma_{e}(H) = [E_{0,tot},\infty).$
\end{lemma}

\begin{proof}
Let $0 = \lambda_0^{(i)} < \lambda_1^{(i)} \leq \cdots$ be the eigenvalues of $- \Lambda_{i}$. Then, $H_{free(i)}$ is unitarily equivalent to the direct sum $ \displaystyle{{\mathop\oplus_{\ell=0}^{\infty} }L_{\ell}^{(i)}}$, where
$$
L_{\ell}^{(i)} = - \partial_r^2 - \frac{(n-1)\rho_i'(r)}{\rho_i(r)}\partial_r + \frac{\lambda_{\ell}^{(i)}}{\rho_i(r)^{2}}.
$$
By the transformation $u \to \rho_i(r)^{(n-1)/2}u$, $L_{\ell}^{(i)}$ is unitarily equiavalent to $\widetilde L_{\ell}^{(i)}$, where
$$
\widetilde L_{\ell}^{(i)} = - \partial_r^2  + \frac{\lambda_{\ell}^{(i)}}{\rho_i(r)^2} + Q_i, \quad
Q_i = \Big( \frac{n-1}{2}\frac{\rho_i'}{\rho_i}\Big)^2 + 
\Big(\frac{n-1}{2}\frac{\rho_i'}{\rho_i}\Big)'.
$$
By (\ref{rho'rho-c0}), we have $Q_i = E_{0,i} + O(r^{-\alpha_{0,i}})$ as $r \to \infty$. Since $\rho_i(r)\to \infty$ for the case of regular end, we have  $\sigma_{e}(\widetilde L_{\ell}^{(i)})= [E_{0,i},\infty)$. For the case of cusp end, we have $\rho_i^{-1}(r) \to \infty$, hence
$\sigma_e(\widetilde L_{\ell}^{(i)}) = \emptyset$ if $\lambda_{\ell}^{(i)} > 0$, and 
$\sigma_e(\widetilde L_{\ell}^{(i)}) = [E_{0,i},\infty)$ if $\lambda_{\ell}^{(i)} = 0$. 
This proves that $\sigma_e(H_{free(i)}) = [E_{0,i},\infty)$. By Weyl's theorem for the perturbation of the essential spectrum, we have 
$\sigma_e(H_i) = \sigma_e(H_{free(i)})$.

Applying well-known Weyl's method of singular sequence, we can show 
$\sigma_e(H_i) \subset \sigma_e(H)$. Therefore
$[E_{0,tot},\infty) \subset \sigma_e(H)$. 
To prove the converse inclusion relation, take any compact interval $I \subset (-\infty,E_{0,tot})$, and $f(\lambda) \in C_0^{\infty}({\mathbb R})$ such that
 $f(\lambda)=1$ on $I$ and $f(\lambda) = 0$ on $[E_{0,tot},\infty)$. Using (\ref{R(z)Reprsent1}), one can show
 $$
 f(H) = \sum_{i=1}^{N+N'}\chi_if(H_i)\widetilde\chi_i + K,
 $$
where $K$ is a compact operator (see e.g. the proof of Chap. 3, Theorem 3.2 of \cite{IsKu10}). Since $I\cap\sigma_e(H_i) = \emptyset$, $f(H_i)$ is compact. Therefore, $f(H)$ is also compact. This proves that 
$I\cap\sigma_e(H) = \emptyset$, which implies $\sigma_e(H) \subset [E_{0,tot},\infty)$. 
\end{proof}

\subsection{Embedded eigenvalues}

We put
\index{$\mathcal T$}
\begin{equation}
\mathcal T = \{E_{0,1},\cdots,E_{0,N+N'}\}.
\label{E01cdotsE0N+N'}
\end{equation}


\begin{theorem}\label{noembededeigenval}
If there exists a regular end $\mathcal M_i$ with $\beta_i > 1/3$, we have
$$
 \sigma_p(H)\cap (E_{0,i},\infty) = \emptyset.
$$
\end{theorem}

\begin{proof}
By (\ref{S7CondC}) and (\ref{LemmatB(t)leqct-epsilon}), 
we can apply Theorem \ref{Growthproperty}. Therefore, if $u$ is an eigenfunction of $- \Delta_{\mathcal M}$ with eigenvalue in $(E_{0,i},\infty)$, $u$ vanishes near infinity of $\mathcal M_i$. By the unique continuation theorem, 
$u$ vanishes identically on $\mathcal M$.
\end{proof}

Recall that $\dfrac{\rho'_i}{\rho_i} \to c_{0,i}$ on each end.
Therefore, if $c_{0,i} > 0$, $\rho_i(r)$ is exponentially growing and $\beta_i > 1/3$ holds, hence there is no embedded eigenvalue in $(E_{0,i},\infty)$. If $c_{0,i} = 0$ for some regular end and $\beta_i > 1/3$, again there is no embedded eigenvalue in $(0,\infty)$. 
The remaining cases are the ones in which either $\beta_i \leq 1/3$ for all regular ends or all the ends are cusp. In these cases the essential spectrum is $[0,\infty)$, and the embedded eigenvalues are discrete in the following sense.

\begin{theorem}
\label{Theorem7.3E1/3case}
Assume that $\beta_i \leq 1/3$ for all regular ends. 
Then, the eigenvalues in $(0,\infty)\setminus\mathcal T$ are of finite multiplicities with possible accumulation points at $\mathcal T$ and $\infty$. 
\end{theorem}

In particular, this theorem holds when all ends are cusp. 

To prove this theorem, we first show that the eigenfunctions associated with embedded eigenvalues decay faster than $L^2$. 

\begin{lemma}
\label{S7Eigenfunctioncompact}
Assume that $\beta_i \leq 1/3$ for all regular ends. 
Let $I$ be a compact interval in $(0, \infty)\setminus\mathcal T$. Then, there exist constants $s, C > 0$ such that 
\begin{equation}
\|u\|_{H^{1,s}(\mathcal M)} \leq C \|u\|_{L^2(\mathcal M)}
\end{equation}
holds for any eigenfunction $u$ of $- \Delta_{\mathcal M}$ with eigenvalue in $I$.
\end{lemma}

\begin{proof}
We have proven that $u \in L^{2,s}$ for any $s > 0$ on the regular end $\mathcal M_i$, moreover, by Lemma \ref{S11polynomweighinequality},
$$
\|\chi u\|_{H^{1,s}(\mathcal M_i)} \leq C_s\|u\|_{L^2}
$$
where, $\chi = 1$ for $r > 2$, $\chi = 0$ for $r < 1$. On the cusp end, multiplying $\chi$ to $u$, which we denote by $u$ and  use Lemma \ref{Cuspestimate}. Then, $P_1(r)u$ obeys the desired estimate. Let $\widetilde u_0 = g^{1/4}P_0(r)u$. Then, it satisfies
\begin{equation}
- \widetilde u_0'' - k^2 \widetilde u_0 = P_0(r)f -
\Big(\big(\frac{g'}{4g}\big)' + \big(\frac{g'}{4g}\big)^2\Big)\widetilde u_0,
\nonumber
\end{equation}
where $f$ is  compactly supported. By virtue of (\ref{Subsection2.8Inequality2.5}), $P_0(r)u$ satisfies the desired inequality.
\end{proof}

Let us now prove Theorem \ref{Theorem7.3E1/3case}.

\begin{proof}
Take a compact interval $I \subset (0,\infty)\setminus \mathcal T$ and suppose that there exists an infinite number of eigenvalues 
$\{\lambda_{\ell}\}_{\ell=1}^{\infty} \subset I$. Let $\{\varphi_{\ell}\}_{\ell=1}^{\infty}$ be the orthonormal system of eigenvectors associated with $\{\lambda_{\ell}\}_{\ell=1}^{\infty}$. By Lemma \ref{S7Eigenfunctioncompact}, $\{\varphi_{\ell}\}_{\ell=1}^{\infty}$ is bounded in $H^
{1,s}(\mathcal M)$ for some $s > 0$. This implies that for any $\epsilon > 0$ there is a constant $r_0>0$ such that on each end
\begin{equation}
\|\varphi_{\ell}\|_{L^2((r_0,\infty)\times M_i)} < \epsilon, \quad \forall \ell \geq 1.
\nonumber
\end{equation}
By Rellich's selection theorem, $\{\varphi_{\ell}\}_{\ell=1}^{\infty}$ contains a subsequence which is convergent in $L^2(\mathcal M)$, which is a contradiction.
\end{proof}


\subsection{Radiation condition and uniqueness}
The radiation condition is the boundary condition at infinity to guarantee the uniqueness of solutions to  the reduced wave equation.  It is closely related to the Rellich-Vekua theorem. 

We define $D_j(k)$ in the same way as in (\ref{DefDpmk}) on each end $\mathcal M_j$ with $E_0$ replaced by $E_{0,j}$. Let $\psi_j(k) = \psi_j(r,x,k)$ be an approximate solution of the equation
$$
-i\psi' + \psi^2 + Q_j-k^2=0, \quad 
Q_j = \Big(\frac{g_j'}{4g_j}\Big)^2 + \Big(\frac{g_j'}{4g_j}\Big)' - E_{0,j}.
$$
By Lemma \ref{psimlemma}, it behaves like
$$
\psi_j(k) = \sqrt{z - E_{0,j}} + O(r^{-\alpha_{0,j}}), \quad {\rm as} \quad r \to \infty.
$$
For $z = re^{i\theta}$, $0 < \theta < 2\pi$, we defined  $\sqrt{z} = \sqrt{r}e^{i\theta/2}$. Therefore, for $\lambda > 0$,
$$
\sqrt{\lambda \pm i\epsilon} \to \pm \sqrt{\lambda}, \quad {\rm as} \quad \epsilon \to 0.
$$
Hence, for $\lambda > E_{0,j}$, we have two $\psi_j$'s, denoted by $\psi_j^{(\pm)}$, where
$$
\psi_j^{(\pm)}(\sqrt{\lambda - E_{0,j}}) = \pm \sqrt{\lambda - E_{0,j}} + O(r^{-\alpha_{0,j}}), \quad {\rm as} \quad r \to \infty.
$$
We put 
\begin{equation}
D_j^{(\pm)}(k) = \partial_r + \frac{g_j'}{4g_j} - i\psi_j^{(\pm)}(k),
\label{S7DjkpmDefine}
\end{equation}
\begin{equation}
\mathcal E = \mathcal T \cup \sigma_p(H). \index{$\mathcal E$}
\label{S9ExceptionalSet}
\end{equation}

\begin{definition}
\label{TotalRadCond}
Let $\lambda \in \sigma_e(H)\setminus\mathcal E$. 
A solution $u \in \mathcal B^{\ast}(\mathcal M)$ of the equation
\begin{equation}
\big(-\Delta_{\mathcal M} - \lambda\big)u = f \quad on \quad \mathcal M
\label{EquationonM}
\end{equation}
is said to satisfy the  outgoing  radiation condition  on the end $\mathcal M_j$, 
if it has the following properties.

\noindent
(1) For $\lambda < E_{0,j}$ and $1 \leq j \leq N + N'$,
$$
u, u' \in \mathcal B^{\ast}_0(\mathcal M_j).
$$ 
(2) For $\lambda > E_{0,j}$ and  $1 \leq j \leq N$,
\begin{equation}
\begin{split}
&  D_j^{(+)}(k)u \in \mathcal B^{\ast}_0(\mathcal M_j), \quad {\it if} \quad 1/3 < \beta_{0,j}, \\
&  D_j^{(+)}(k)u \in L^{2, -(1-\epsilon)/2}(\mathcal M_j), \quad {\it if}
\quad 0 < \beta_{0,j} \leq 1/3,
\end{split}
\nonumber
\end{equation}
for some $\epsilon > 0$.

\noindent
(3) For $\lambda > E_{0,j}$ and $N+1 \leq j \leq N+N'$,
\begin{equation}
D_j^{(+)}(k)u \in \mathcal B_0^{\ast}(\mathcal M_j).
\nonumber
\end{equation}

When $D_j^{(+)}(k)$ is replaced by $D_j^{(-)}(k)$, $u$ is said to satisfy the {\it incoming radiation condition}. 
We say that $u$ satisfies the outgoing (incoming) radiation condition on $\mathcal M$, if it satisfies the outgoing (incoming) radiation condition on all $\mathcal M_j$. 
\end{definition}
\index{radiation condition}

The purpose of this subsection is to prove the following theorem.


\begin{theorem}
\label{Theorem7.6}
Assume that $\lambda \in \sigma_{e}(H)\setminus\mathcal E$,  and let $u \in \mathcal B^{\ast}(\mathcal M)$ be a solution to the equation $(- \Delta_{\mathcal M}-\lambda)u = 0$ satisfying the radiation condition. Then, $u = 0$.
\end{theorem}

The starting point of the proof of this theorem is the following formula  (\ref{Lemma7.7RadCondFormula}) in Lemma \ref{RadConduniqueformula}.
 Let $\{\chi_j\}_{j=0}^{N+N'}$ be the partition of unity on $\mathcal M$ satisfying (\ref{S8Partitionunity}).
Take  $\varphi(r) \in C^{\infty}({\mathbb R})$ such that $\varphi(r) = 1$ for $r < 1$, $\varphi(r)=0$ for $r > 2$, and put 
$$
\phi(r) = \int_r^{\infty}\varphi(t)dt, \quad 
\phi_R(r) = \phi\big(\frac{r}{R}\big).
$$
Then, $\phi_R(r) = 0$ for  $r > 2R$ and $\phi'_R(r) = - \frac{1}{R}\varphi(\frac{r}{R})$. 


\begin{lemma}
\label{RadConduniqueformula}
Let $u$ be a solution to the equation $(-\Delta_{\mathcal M} - \lambda) u = 0$. Then, for any constant $\sigma \in {\mathbb R}$,
\begin{equation}
\sum_{j=1}^{N+N'}{\rm Im}\,\frac{1}{R}
\big(\varphi(\frac{r}{R})r^{\sigma}\chi_ju,u') = 
\frac{1}{2}{\rm Im}\,
\big([H,r^{\sigma}(1 - \chi_0)]\phi_Ru,u).
\label{Lemma7.7RadCondFormula}
\end{equation}
\end{lemma} 

\begin{proof}
Consider the equation
$$
{\rm Im}\,\big(\phi_Rr^{\sigma}\chi_ju,(\Delta_{\mathcal M} + \lambda)u\big)=0.
$$
Then by integration by parts, 
$$
{\rm Im}\, \big((\phi_Rr^{\sigma}\chi_j)'u,u'\big) = 0,
$$
which implies
\begin{equation}
{\rm Im}\,\frac{1}{R}\Big(\varphi\big(\frac{r}{R}\big)r^{\sigma}\chi_j u,u'\Big) = {\rm Im}\, \big(\phi_R(r^{\sigma}\chi_j)'u,u'\big), \quad \forall j \geq 1.
\nonumber
\end{equation}
Summing up with respect to $j$, we have
\begin{equation}
\sum_{j=1}^{N+N'}{\rm Im}\,\frac{1}{R}\Big(\varphi\big(\frac{r}{R}\big)r^{\sigma}\chi_ju,u'\Big) ={\rm Im}  \big(\phi_R(r^{\sigma}(1-\chi_0))'u,\partial_ru\big).
\label{S8Equation10}
\end{equation}
Noting that on each end for any smooth function $w(r)$
$$
2w'\partial_r = [w,H] - w'' - \frac{g'}{2g}w',
$$
we have 
$$
{\rm Im}\, \big(\phi_Rw'u,\partial_ru\big) 
= \frac{1}{2}{\rm Im}\, \big([H,w]\phi_Ru,u\big).
$$
Letting $w = r^{\sigma}(1 - \chi_0)$, and taking notice of (\ref{S8Equation10}), we obtain  the lemma. 
\end{proof}

 We consider the case of outgoing radiation condition, and  divide the proof of Theorem \ref{Theorem7.6} into three cases. 


\begin{lemma}\label{D+ktoB0ast}
 Let $u \in \mathcal B^{\ast}(\mathcal M)$ be a solution to  $(- \Delta_{\mathcal M} - \lambda)u = 0$ satisfying the radiation condition. Suppose there exists a regular end $\mathcal M_j$ such that $\lambda > E_{0,j}$ and  $\beta_{0,j} > 1/3$. Then, $u = 0$ on $\mathcal M$.
\end{lemma}

\begin{proof}
Take $\sigma = 0$ in (\ref{Lemma7.7RadCondFormula}).  Using the equation $(H - \lambda)u = 0$, we have
\begin{equation}
([H,- \chi_0']\phi_Ru,u) = (\chi_0'[H,\phi_R]u,u).
\nonumber
\end{equation}
The right-hand side tends to 0 as $R \to \infty$. 
Using the radiation condition for the left-hand side of (\ref{Lemma7.7RadCondFormula}) and noting that $u, u' \in \mathcal B_0^{\ast}(\mathcal M_j)$ if $\lambda < E_{0,j}$, we then have
$$
\lim_{R\to\infty}{\sum_{j}}'\frac{1}{R}(\varphi\big(\frac{r}{R}\big)\chi_ju,\sqrt{\lambda - E_{0,j}}\chi_ju) = 0,
$$
where the sum $\sum'_{j}$ ranges over all $j$ such that $\lambda > E_{0,j}$. 
This yields
\begin{equation}
\lim_{R\to\infty}\frac{1}{R}(\varphi\big(\frac{r}{R}\big)\chi_ju,\chi_ju) = 0, \quad \forall j \geq 1.
\label{S8Equatio11}
\end{equation}
Hence $u \in \mathcal B^{\ast}_0(\mathcal M_j),  \forall j \geq 1$. 
 Then, letting $S(r) = \{r\}\times M_j$, we have 
$$
\liminf_{r\to\infty}
\int_{S(r)}\big(|u'|^2 + |u|^2\big)dS(r) = 0.
$$
 By Theorem \ref{S2RellichTh1}, $u=0$ near infinity of $\mathcal M_j$ with $\beta_{0,j} > 1/3$.  Therefore, $u =0$ on $\mathcal M$ by the unique continuation theorem. 
\end{proof}

We next consider the case in  which all ends are cusp.


\begin{lemma}\label{LemmaRadCondCusp}
Suppose all ends are cusp, and $\lambda \in \sigma_e(H)\setminus {\mathcal E}$. 
If $u \in \mathcal B^{\ast}$ satisfies  $(-\Delta_{\mathcal M} - \lambda)u = 0$ and the radiation condition, then $u=0$.  
\end{lemma}

\begin{proof}
We pick up one end $\mathcal M_i$ and drop the subscript $i$. Take $\chi \in C^{\infty}((0,\infty))$ such that $\chi = 1$ for $r > r_0+1$ and $\chi=0$ for $r < r_0$, where $r_0$ is chosen large enough.  Then, $U = \chi u$ satisfies
$$
\big(-\partial_r^2 - \frac{g'}{2g}\partial_r + B(r) - \lambda\big)U = f,
$$
where $f = -\chi'' u - 2\chi'u' - \frac{g'}{2g}\chi'u$. Recalling that now $\mathcal M = \mathcal M_i$, we show that $U \in L^2(\mathcal M_i)$. If $\lambda < E_{0,i}$, this is true. Below, we consider the case in which $\lambda > E_{0,i}$. Again dropping the subscript $i$,  let  $P_0(r)$ be the projection to the 0-eigenvalue of $\Lambda(r)$ and $P_1(r) = 1- P_0(r)$.  Then, $u_0 = P_0(r)U$ and $u_1 = P_1(r)U$ satisfy
\begin{equation}
\big( - \partial_r^2 - \frac{g'}{2g}\partial_r - \lambda\big)u_0 = f_0,
\nonumber
\end{equation}
\begin{equation}
\big( - \partial_r^2 - \frac{g'}{2g}\partial_r + B(r) - \lambda\big)u_1 = f_1,
\nonumber
\end{equation}
\begin{equation}
f_i = P_i(r)f - P_i''(r)U - 2P_i'(r)U' - \frac{g'}{2g}P_i'(r)u + [P_i(r),\frac{g'}{2g}]U', \quad i = 0, 1.
\nonumber
\end{equation}
Letting $v_0 = g^{1/4}u_0$ and $k^2 = \lambda - E_0$,  we have
\begin{equation}
- v_0'' - k^2 v_0 = f_2, 
\nonumber
\end{equation}
\begin{equation}
f_2 = g^{1/4}f_0 - \Big(\big(\frac{g'}{4g}\big)' + \big(\frac{g'}{4g}\big)^2 - E_0\Big)v_0.
\nonumber
\end{equation}
Note that $f$ is compactly supported.
Lemma \ref{S6P0(r)estimate} implies 
\begin{equation}
\|P_i'(r)\| + \|P_i''(r)\| + \|[P_i(r),\frac{g'}{2g}]\| \leq C(1 + r)^{-1 - 2\epsilon},
\nonumber
\end{equation}
and the assumption (\ref{A-1alpha11/2}) yields 
$$
\Big|\big(\frac{g'}{4g}\big)' + \big(\frac{g'}{4g}\big)^2 - E_0\Big| \leq C(1 + r)^{-1 - 2\epsilon}
\nonumber
$$
for a small $\epsilon > 0$. Note also that
$$
U \in \mathcal B^{\ast} \subset L^{2, - (1 + \epsilon)/2}.
$$
Therefore, $f_1 \in L^{2,(1 + \epsilon)/2}$.
Lemma \ref{Cuspestimate} (1) then implies that $u_1 \in L^{2,(1+\epsilon)/2}$.

Similarly, $f_2 \in L^{2,(1+ \epsilon)/2}$. 
In the proof of Lemma \ref{D+ktoB0ast}, we have already proven that $u_0 \in \mathcal B^{\ast}_0(\mathcal M)$, which implies 
$\liminf_{r\to\infty}|v_0(r)| = 0$. 
Then, by Lemma \ref{Lemma2.181dimHelmholtz},  $v_0 \in L^{2,(\epsilon - 1)/2}$. This fact yields $f_0 \in L^{2,(1 + 2\epsilon)/2}$. Then, by the same argument, we have 
$u_0 \in L^{(2\epsilon - 1)/2}$. Repeating this procedure, we obtain $u_0 \in L^{2,(1 + \epsilon)/2}$. We have thus proven that $U \in L^2$. Therefore, $u \in L^2(\mathcal M_i)$ on all ends $\mathcal M_i$. Since $\lambda$ is not an eigenvalue of $- \Delta_{\mathcal M}$, we have $u = 0$.
\end{proof}

We consider the remaining case.

 \begin{lemma}
Assume that $0 < \beta_i \leq 1/3$ on all regular ends.  Let $u \in \mathcal B^{\ast}(\mathcal M)$ be a solution to  $(- \Delta_{\mathcal M} - \lambda)u = 0$ on $\mathcal M$ satisfying 
 the radiation condition and  $\lambda \in \sigma_e(H)\setminus \mathcal E$. Then, $u= 0$. 
\end{lemma}

\begin{proof} 
By  the arguments in the proof of 
Lemma \ref{LemmaRadCondCusp},   on all cusp ends $\mathcal M_i$, $u \in L^{2,(1 + \epsilon)/2}(\mathcal M_i)$, hence 
$D^{(+)}_i(k)u \in L^{2,(1 + \epsilon)/2}(\mathcal M_i)$.

Take $\sigma = \epsilon$ in (\ref{Lemma7.7RadCondFormula}), and let $w = r^{\epsilon}(1 - \chi_0)$. By the equation $- \Delta_{\mathcal M}u = \lambda u$, we have
$$
([H,w]\phi_Ru,u) = -(w[H,\phi_R]u,u)
$$
which tends to 0 as $R \to \infty$. In fact, since $\phi'_R(r) = - \frac{1}{R}\varphi(\frac{r}{R})$,  
$$
|(w[H,\phi_R]u,u)| \leq \frac{1}{R}\int_{r < 2R}(1 + r)^{\epsilon-1}|f(r,x)|^2 \rho^{n-1}(r)drdx
$$
for some $f \in \mathcal B^{\ast}$. Then, by (\ref{Lemma7.7RadCondFormula}), we have
$$
\sum_{j=1}^{N+N'}{\rm Im}\frac{1}{R}\Big(\varphi\big(\frac{r}{R}\big)r^{\epsilon}\chi_ju,u'\Big) = 0.
$$
By the radiation condition, we then have
$$
\sum_{j=1}^{N+N'}\frac{1}{R}\Big(\varphi\big(\frac{r}{R}\big)r^{\epsilon}\chi_ju,u\Big) = 0,
$$
which implies that 
$$
\liminf_{r\to\infty}\int_{S(r)}r^{\epsilon}\big(|u'|^2 + |u|^2\big)dS(r) =0.
$$
Theorem \ref{S2Rapiddecay}  then yields $u \in L^2$, hence $u = 0$.
\end{proof}


\section{Limiting absorption principle}
\label{SectionLAP}
For a self-adjoint operator $A$ on a Hilbert space $\mathcal H$, $(A - \lambda)^{-1}$ does not exist for $\lambda \in \sigma(A)$. However, when $\lambda \in \sigma_c(A)$, it often happens that  $\lim_{\epsilon\to 0}(A - \lambda \mp i\epsilon)^{-1}$ exists as an operator from $\mathcal X$ to $\mathcal Y$, where $\mathcal X$ and $\mathcal Y$ are Banach spaces rigging $\mathcal H$, i.e. $\mathcal X \subset \mathcal H \subset \mathcal Y$, with dense and continuous imbedding. This is called the
 {\it limiting absorption principle} and used as a fundamental tool 
in the study of continuous spectrum. The purpose of this section is to prove this limiting absorption principle for the Laplacian $H = - \Delta_{\mathcal M}$ on $\mathcal M$. In this section, we assume (\ref{betaigammaiassume}). 
We put $R(z) = (H - z)^{-1}$. Take a compact interval  $I \subset \sigma_e(H)\setminus \mathcal E$, and define
$J_{\pm} = \{z \in {\mathbb C}\, ; \, {\rm Re}\,z \in I, \ 0 < {\rm Im}\, z < 1\}$.

 
\begin{lemma} \label{ModelLAPinL2-s}
Let $s > 1/2$.\\
\noindent
(1) There exists a constant $C_s > 0 $ such that
$$
\|R(z)f\|_{-s} \leq C_s\|f\|_s, \quad z \in J_{\pm}.
$$
(2) For any $\lambda \in I$ and $f \in L^{2,s}$, there exists a strong limit ${\rm s-lim}\,_{\epsilon\to 0}
R(\lambda \pm i\epsilon)f$ in $L^{2,-s}$. Moreover, $R(\lambda \pm i0)f$ is an $L^{2,-s}$-valued strongly continuous function of $\lambda \in I$. \\
\noindent
(3) For $\lambda \in I$, $R(\lambda + i0)f$ satisfies the outgoing radiation condition, and
  $R(\lambda - i0)f$ satisfies the  incoming radiation condition.
\end{lemma}

\begin{proof}
Let us prove (1) for the case $z \in J_+$. If the assertion (1) does not hold, there exist sequences $z_j \in J_+$, $f_j \in L^{2,s}$ such that $\|f_j\|_s\to 0$, and $u_j=R(z_j)f_j$ satisfies $\|u_j\|_{-s}=1$. Then, there exists a subsequence, which is denoted by $\{z_j\}$ again, such that $z_j \to \lambda$. If $\lambda \not\in {\mathbb R}$, we easily arrive at a contradiction. Assume that  $\lambda \in I$. Fix any $i = 1,\cdots, N+N'$, and take $s>s'>1/2$. Note the inequality
\begin{equation}
\begin{split}
& \int_{R}^{\infty}(1 + r)^{-2s}\|\chi_iu_j\|^2_{{\bf h}_i(r)}\rho_i(r)^{n-1}dr \\
& \leq (1+R)^{-2(s-s')}\int_R^{\infty}(1 + r)^{-2s'}\|\chi_i u_j\|^2_{{\bf h}_i(r)}\rho_i(r)^{n-1}dr .
\end{split}
\nonumber
\end{equation}
Taking account of Lemmas \ref{Lemmauestimate} and \ref{Cuspestimate}, we have
\begin{equation}
\int_R^{\infty}(1+ r)^{-2s}\|\chi_iu_j\|^2_{{\bf h}_i(r)}\rho_i^{n-1}(r)dr \leq
C(1 + R)^{-2(s-s')}, \quad 1/2 < s' < s,
\label{S5Inequality2}
\end{equation}
where the constant $C$ does not depend on $j$. By the a-priori estimate (Lemma \ref{S4Apriori1}), 
 $\{u_j\}$ is bounded in $H^{1,-s}(\mathcal M)$. By Rellich's selection theorem and (\ref{S5Inequality2}), we can choose a subsequence, which is denoted by $\{u_j\}$ again,  and $u \in L^{2,-s}$ such that $u_j \to u$ in $L^{2,-s}$. Then $\|u\|_{-s}=1$, and $u$ satisfies
\begin{equation}
\big(-\Delta_{\mathcal M} - \lambda\big)u=0.
\nonumber
\end{equation}
By Lemmas \ref{Lemma6.3w-1/2leqfu}, \ref{Lemma6.7vestimateimprove}  and \ref{Cuspestimate}, $u$ satisfies the outgoing radiation condition.
Theorem \ref{Theorem7.6} then implies that $u=0$, which is a contradiction.

To prove (2), take a sequence $z_j \in J_+$ such that $z_j \to \lambda \in J_+$, and put $u_j = R(z_j)f$. Arguing as above, we see that $\{u_j\}$ contains a subsequence $\{u_{j'}\}$, which is convergent in $L^{2,-s}$, and the limit $u$ satisfies the outgoing radiation condition as well as the equation
$$
\big(- \Delta_{\mathcal M} - \lambda\big)u=f.
$$
Such a solution $u$ is unique. We thus see that any subsequence of $\{u_j\}$ contains a sub-sub sequence which converges to one and the same limit. This shows that $\{u_j\}$ itself converges in $L^{2,-s}$. The strong continuity is proven similarly. The assertion (3) is already proven. 
\end{proof}

We extend Lemma \ref{ModelLAPinL2-s} to  $\mathcal B, \mathcal B^{\ast}$ spaces.


\begin{theorem}\label{ModelLAPforbfh}
Let $I \subset \sigma_e(H)\setminus\mathcal E$ be a compact interval and $J_{\pm} = \{z \in {\mathbb C}\, ; \, {\rm Re}\,z \in I, \ 0 < {\rm Im}\, z < 1\}$. \\
\noindent
(1) There exists a constant $C > 0$ such that
$$
\|R(z)f\|_{\mathcal B^{\ast}} \leq C\|f\|_{\mathcal B}, \quad z \in J_{\pm}.
$$
(2) For any $f \in \mathcal B$, and $\lambda \in I$,  there exists a $\ast$-weak limit $\lim_{\epsilon \to 0}
R(\lambda \pm i\epsilon)f$, i.e. for any $f, g \in \mathcal B$,
$$
\lim_{\epsilon \to 0}(R(\lambda \pm i\epsilon)f,g) = (R(\lambda \pm i0)f,g).
$$
Moreover
$$
\|R(\lambda \pm i0)f\|_{\mathcal B^{\ast}} \leq C\|f\|_{\mathcal B}, \quad \lambda \in I,
$$
and $(R(\lambda \pm i0)f,g)$ is a continuous function of $\lambda \in I$ for any $f, g \in \mathcal B$. \\
\noindent
(3) For any $f \in \mathcal B$,  and $\lambda \in I$, $R(\lambda \pm i0)f$ satisfies $D^{(\pm)}(k)R(\lambda \pm i0)f \in \mathcal B^{\ast}_0(\mathcal M_i)$ on each end $\mathcal M_i$, $i = 1, \cdots, N + N'$.
  \end{theorem}

\begin{proof}
Lemma \ref{Lemmauestimate} and Lemma \ref{ModelLAPinL2-s} (1) prove $\|R(z)f\|_{\mathcal B^{\ast}} \leq C\|f\|_s$, $s > 1/2$. Taking the adjoint, we have $\|R(z)f\|_{-s} \leq C\|f\|_{\mathcal B}$, which, combined with Lemma \ref{Lemmauestimate}, yields the assertion (1).
 The assertion (2) is proved by approximating $f, g$ by $f_j, g_j \in L^{2,s}$ and using Lemma \ref{ModelLAPinL2-s}  (2).
Letting $u_j(z) = R(z)f_j$, $u(z) = R(z)f$, we then see that $(u_j(z),g_j)$ is uniformly continuous with respec to $z \in J_{\pm}$. Since $(u_j(z),g_j)$ converges uniformly to $(u(\lambda \pm i0),g)$, the continuity of $(u(\lambda \pm i0),g)$ follows. 
 Since $u^{(\pm)}_j \to u^{(\pm)}$ in $\mathcal B^{\ast}$ and $u^{(\pm)}_j \in \mathcal B^{\ast}_0$, and $\mathcal B^{\ast}_0$ is a closed subspace of $\mathcal B^{\ast}$, the assertion (3) follows. 
\end{proof}

\medskip
Let us remark that Lemma \ref{Lemma6.7vestimateimprove}  also holds for $R(\lambda \pm i0)$. In fact, by computing the commutator $[R(z), (1 + r^2)^{s/2}]$, it is easy to see that for $z \not\in {\mathbb R}$, $R(z) \in {\bf B}(L^{2,s};L^{2,s})$ for any $s \in {\mathbb R}$. Then, the inequality in Lemmas \ref{Lemma6.7vestimateimprove}  holds for $u = R(\lambda \pm i\epsilon)f$ with $f \in L^{2,(t+1)/2}$, $\epsilon > 0$, where the constant $C$ is independent of $\epsilon > 0$. Letting $\epsilon \to 0$, Lemma 
\ref{Lemma6.7vestimateimprove} is proven for $R(\lambda \pm i0)f$.

\medskip
Let $E_H(\lambda)$ be the resolution of the identity (spectral decomposition) for $H$. 
The absolutely continuous subspace for $H$, denoted by $\mathcal H_{ac}(H)$, is defined as the set of all $f \in L^2(\mathcal M)$ such that $(E_H(\lambda)f,f)$ is absolutely continuous with respect to the Lebesgue measure $d\lambda$. Thanks to the limiting absorption principle, for any $f \in \mathcal B$ and compact interval $I \subset (E_0,\infty)\setminus\mathcal E$, we have 
\begin{equation}
(E_H(I)f,f) = \frac{1}{2\pi i}\int_I\big([R(\lambda + i0)-R(\lambda - i0)]f,f\big)d\lambda. 
\label{S8EHI=Rff}
\end{equation}
This implies that $E_H(I)f \in \mathcal H_{ac}(H)$ for any $f \in L^2(\mathcal M)$ and compact interval $I \subset (E_0,\infty)\setminus \mathcal E$. One can then prove the  following theorem.

\begin{theorem}
Letting $P_{ac}(H) $ be the orthogonal projection onto $\mathcal H_{ac}(H)$, and  $\mathcal H_{pp}(H)$ is the closure of the linear hull of eigenvectors, we have
$$
P_{ac}(H) = E_H((E_0,\infty)\setminus\sigma_p(H)),
$$
\begin{equation}
\sigma_e(H) = \sigma_{ac}(H) = [E_0,\infty),
\nonumber
\end{equation}
\begin{equation}
L^2(\mathcal M) = \mathcal H_{ac}(H) \oplus \mathcal H_{pp}(H).
\nonumber
\end{equation}
\end{theorem}

\section{Resolvent asymptotics - Non perturbative approach}\label{ResolventAsymptotics1}


\subsection{Asymptotic expansion at infinity}
We study the asymptotic behavior of the resolvent 
$R(\lambda \pm i0)f$ as $r \to \infty$ on each end, where $\lambda \in (E_0,\infty)\setminus\mathcal E$. In this section, we consider the metric of long-range behavior, using the method developed for the Schr\"odinger operator $- \Delta + V(x)$ in ${\mathbb R}^n$ (see e.g. \cite{Saito77}, \cite{Saito}, \cite{Is80}).
We pick up a regular end $\mathcal M_j$ satisfying 
\begin{equation}
\alpha_{0,j} > 0, \quad \beta_{0,j} > 1/2, \quad \gamma_{0,j} > 1.
\label{betaj>1/2}
\end{equation}
The aim of this section is to prove the following theorem.


\begin{theorem}\label{TheoremexistlimitR(lambdapmi0)f}
Let $\mathcal M_j$ be a regular end satisfying (\ref{betaj>1/2}), and put
\begin{equation}
\Phi_j(r,\lambda) = \int_0^r\phi_j(t,\lambda)dt,
\label{DefPjij(r)}
\end{equation}
\begin{equation}
\phi_j(r,\lambda) = \sqrt{\lambda - \frac{(n-1)^2}{4}\Big(\frac{\rho'_j(r)}{\rho_j(r)}\Big)^2}.
\label{Th9.1phijrlambda}
\end{equation}
Then, for $\lambda > E_{0,j}$, there exists a strong limit
\begin{equation}
\lim_{r\to\infty}e^{\mp i\Phi_j(r,\lambda)}\rho_j(r)^{(n-1)/2}\big(R(\lambda \pm i0)f\big)(r,\cdot) \quad 
{\it in }  \quad L^2(M_j),
\label{S11LRlimitbetageq1}
\end{equation}
 for $f \in L^2_{comp}(\mathcal M)$.
\end{theorem}

\medskip
Note that 
\begin{equation}
\frac{(n-1)^2}{4}\Big(\frac{\rho'_j(r)}{\rho_j(r)}\Big)^2- E_{0,j} \in 
\left\{
\begin{split} 
& S^{-\alpha_j}, \quad {\rm if} \quad c_{0,j}  \neq 0, \\
& S^{-2\alpha_j}, \quad {\rm if} \quad c_{0,j} = 0,
\end{split}
\right.
\nonumber
\end{equation}
hence $\Phi_j(r,\lambda) \sim \sqrt{\lambda - E_{0,j}}\, r = kr$ as $r \to \infty$. 


\subsection{Reduction to each end}
Let $H_j$ be the Laplacian on 
$
\widetilde{\mathcal M}_j =(2,\infty)\times M_j
$
 associated with the metric $(dr)^2 + \rho_j(r)^2h_j(r,x,dx)$ with Dirichlet boundary condition at $r=2$. By the well-known perturbation theory, one can show that $\sigma_e(H_j) = [E_{0,j},\infty)$, and by Theorem \ref{S2RellichTh1},  
$H_j$ has no eigenvalues in $(E_{0,j},\infty)$.  The limiting absorption principle in \S 8 can also be proved for $H_j$, since our concern is only the estimates of the resolvent at infinity. Hence, Theorem \ref{ModelLAPforbfh} holds for $H_j$ on $\widetilde{\mathcal M}_j$ as well. We put $R_j(z) = (H_j - z)^{-1}$.

By using the partition of unity $\{\chi_j\}_{j=1}^{N+N'}$ in (\ref{S8Partitionunity}), we have for $\lambda \in \sigma_e(H)\setminus \mathcal E$
\begin{equation}
\chi_j R(\lambda \pm i0) = R_j(\lambda \pm i0)\chi_j + R_j(\lambda \pm i0)[H_j,\chi_j]R(\lambda \pm i0).
\label{S9chiR=Hjreolventchij}
\end{equation}
Hence, the behavior of $R(\lambda \pm i0)f$ on the end $\mathcal M_j$ is reduced to that of $R_j(\lambda \pm i0)f$. 
Let $u_j^{(\pm)} = R_j(\lambda \pm i0)f$, where $f \in C_0^{\infty}(\widetilde{\mathcal M}_j)$. 
Let $g_j= \det\big(\rho_j(r)^2h_j(r,x,dx)\big)$ and
\begin{equation}
Q_j= \Big(\frac{g'_j}{4g_j}\Big)^2 + \Big(\frac{g'_j}{4g_j}\Big)' - E_{0,j}.
\nonumber
\end{equation}
Note that 
\begin{equation}
D_j^{(\pm)}(k) = \partial_r + \frac{g'_j}{4g_j} - i\psi_j^{(\pm)},
\nonumber
\end{equation}
where $\psi_j^{(\pm)}$ is defined by $\psi_m$ in Lemma \ref{psimlemma}. Since $Q_j = \frac{(n-1)^2}{4}\big(\frac{\rho'_j}{\rho_j}\big)^2 - E_{0,j} + O(r^{-1-\epsilon})$ for some $\epsilon > 0$, we have by 
(\ref{Th9.1phijrlambda})
\begin{equation}
\psi_j ^{(\pm)}(r,x,\lambda)  =\pm \phi_j(r,\lambda) + O(r^{-1-\epsilon}).
\nonumber
\end{equation}
Letting
\begin{equation}
\Psi_j^{(\pm)}(r,x,\lambda) = \int_0^r\psi_j^{(\pm)}(t,x,\lambda)\, dt, 
\nonumber
\end{equation}
we have only to show Theorem \ref{TheoremexistlimitR(lambdapmi0)f} with $\pm \Phi_j$ replaced by $\Psi_j^{(\pm)}$.


\subsection{Utility of radiation condition}
Let us note that if $\beta_{0,j} > 1$, $\gamma_{0,j} > 2$, the proof of Theorem \ref{TheoremexistlimitR(lambdapmi0)f} is easy. In fact, in this case,  in Lemma \ref{Lemma6.7vestimateimprove}, 
one can choose $2 < t < \gamma_{0,j}$ so that $(t-1)/2 > 1/2$ and $(t-1-\gamma_0)/2 < - 1/2$. Note that
$$
\partial_r\Big(g_j^{1/4}e^{-i\Psi_j^{(\pm)}}u_j^{(\pm)}\Big) =
g_j^{1/4}e^{-i\Psi_j^{(\pm)}}D_j^{(\pm)}(k)u_j^{(\pm)},
$$
and that 
the right-hand side belongs to $L^{2,s}(\widetilde{\mathcal M}_j)$ for some $s > 1/2$. This yields
\begin{equation}
\begin{split}
& \int_2^{\infty}\|\partial_r\Big(g_j^{1/4}e^{-i\Psi_j^{(\pm)}}u_j^{(\pm)}\Big)\|_{L^2(M_j)}dr\\
&  \leq 
C
\left(\int_2^{\infty}r^{2s}\|\partial_r\Big(g_j^{1/4}e^{-i\Psi_j^{(\pm)}}u_j^{(\pm)}\Big)\|_{L^2(M_j)}^2dr\right)^{1/2},
\end{split}
\nonumber
\end{equation}
which proves Theorem \ref{TheoremexistlimitR(lambdapmi0)f}.


We turn to the general case. 

\begin{lemma}\label{LemmaDpm(k9uvto0}
Assume (\ref{betaj>1/2}), and  for $f \in L^2_{comp}(\mathcal M)$, let $u_j^{(\pm)} = R_j(\lambda \pm i0)f$. Let $v_j^{(\pm)} \in H^1_{loc}(\widetilde{\mathcal M}_j)\cap  \mathcal B^{\ast}(\widetilde{\mathcal M}_j)$ be such that 
\begin{equation}
D_j^{(\pm)}(k)v_j^{(\pm)}, \sqrt{B_j}v_j^{(\pm)} \in L^{2,-s}(\widetilde{\mathcal M}_j)
\label{Lemma9.2Assumption}
\end{equation}
 for some $0 < s < \min\frac{1}{2}\big(2\beta_j - 1, \gamma_j-1\big)$. We put
$$
a_j^{(\pm)}(r) = \int_{M_j}\big(D_j^{(\pm)}(k)u_j^{(\pm)}(r,x)\big)\overline{v_j^{(\pm)}(r,x)}\,\sqrt{g_j(r,x)}\,dx.
$$
Then, we have
\begin{equation}
\frac{d}{dr}a_j^{(\pm)}(r) = \mp 2i\phi_j(r,\lambda)a_j^{(\pm)}(r) + F_j^{(\pm)}(r), \quad
\int_1^{\infty}|F_j^{(\pm)}(r)|dr < \infty,
\label{Fjpm(r)inL1}
\end{equation}
and also
\begin{equation}
\lim_{r\to\infty}a_j^{(\pm)}(r) =0.
\label{S10phijpmrto0}
\end{equation}
\end{lemma}

\begin{proof}
 We prove the $+$ case.  Using (\ref{DefVpm}) and (\ref{S3Eq1}), we have
\begin{equation}
\begin{split}
\partial_r\big(D_j^{(+)}(k)u_j^{(+)}\big) & = -\frac{g_j'}{4g_j}D_j^{(+)}(k)u_j^{(+)} + B_ju_j^{(+)} \\
& - i\psi_jD_j^{(+)}(k)u_j^{(+)} - f + V_ju_j^{(+)},
\end{split}
\nonumber
\end{equation}
\begin{equation}
\partial_r\big(v_j^{(+)}\sqrt{g_j}\big) = \big(D_j^{(+)}(k)v_j^{(+)}\big)\sqrt{g_j} + 
\big(\frac{g_j'}{4g_j} + i\psi_j^{(+)}\big)v_j^{(+)}\sqrt{g_j},
\nonumber
\end{equation}
where we put
$$
V_j^{(+)} = - i(\psi_j^{(+)})' + (\psi_j^{(+)})^2 + Q_j - k^2,
$$
$$
Q_j = \big(\frac{g_j'}{4g_j}\big)^2 + \big(\frac{g_j'}{4g_j}\big)' - E_{0,j}.
$$
We then have, by a straightforward computation,
\begin{equation}
\begin{split}
& \partial_r\Big((D_j^{(+)}(k)u_j^{(+)})\overline{v_j^{(+)}}\sqrt{g_j}\Big)  \\
 =&   - i(\psi_j^{(+)}+ \overline{\psi_j^{(+)}})(D_j^{(+)}(k)u_j^{(+)})\overline{v_j^{(+)}}\sqrt{g_j} + (D_j^{(+)}(k)u_j^{(+)})\overline{( D_j^{(+)}(k)v_j^{(+)})}\sqrt{g_j} \\ 
&+ (B_ju_j^{(+)})\overline{v_j^{(+)}}\sqrt{g_j} + V_j^{(+)}u_j^{(+)}\overline{v_j^{(+)}}\sqrt{g_j} - f\overline{v_j^{(+)}}\sqrt{g_j}.
\end{split}
\nonumber
\end{equation}
This implies
$$
F_j^{(+)} (r) = \int_{M_j} G_j(r,x)\sqrt{g_j}dx, 
$$
\begin{equation}
\begin{split}
G_j(r,x) =& - i(\psi_j^{(+)} + \overline{\psi_j^{(+)}} - 2\phi_j(r,\lambda))(D_j^{(+)}(k)u_j^{(+)})\overline{v_j^{(+)}}\\
&+ (D_j^{(+)}(k)u_j^{(+)})\overline{( D_j^{(+)}(k)v_j^{(+)})} + (B_ju_j^{(+)})\overline{v_j^{(+)}} \\
& + V_j^{(+)}u_j^{(+)}\overline{v_j^{(+)}}- f\overline{v_j^{(+)}}.
\end{split}
\end{equation}
Since $\psi_j^{(+)} = \phi_j(r,\lambda) + O(r^{-1-\epsilon})$, we have
$$
\psi_j^{(+)} + \overline{\psi_j^{(+)}} - 2\phi_j(r,\lambda) = O(r^{-1-\epsilon}), \quad V_j^{(+)} =  O(r^{-1-\epsilon}).
$$
 By virtue of the assumption $\beta_j > 1/2$ and Lemmas \ref{Lemma6.7vestimateimprove}  and \ref{Cuspestimate}, we have 
 $$
 D_j^{(+)}(k)u_j^{(+)}, \ \sqrt{B_j}u_j^{(+)} \in L^{2,s}(\widetilde M_j), \quad {\rm for} \quad 0 < s < \min\frac{1}{2}\big(2\beta_j - 1, \gamma_j-1\big).
 $$
 This and the assumption (\ref{Lemma9.2Assumption}) imply that
 $F_j^{(+)}(r) \in L^1((2,\infty))$, which proves 
(\ref{Fjpm(r)inL1}). Letting
$$
b_j(r,\lambda) = 2\int_0^r\phi_j(t,\lambda)dt,
$$
we then have
$\frac{d}{dr}\big(e^{ib_j}a_j^{(+)}(r)\big)= e^{ib_j}F_j^{(+)}(r)$, which shows the existence of the limit
$\lim_{r\to\infty}e^{ib_j}a_j^{(+)}(r)$. However, this limit is equal to 0, since $r^{-\alpha}a_j^{(+)}(r) \in L^1((2,\infty))$ for some $0 < \alpha < 1$. This proves (\ref{S10phijpmrto0}). 
\end{proof}

In the course of the proof, we have obtained
$$
a_j^{(+)}(r) = - \int_r^{\infty}e^{i(b_j(t,\lambda)-b_j(r,\lambda))}F_j^{(+)}(t)dt,
$$
hence
\begin{equation}
\big|a_j^{(+)}(r)\big| \leq \int_r^{\infty}|F_j^{(+)}(t)|dt.
\label{Lemma9.3after}
\end{equation}


\begin{lemma}\label{S10Lemma10.3}
Let $u_j^{(\pm)}$ be as in Lemma \ref{LemmaDpm(k9uvto0}. Then, we have
\begin{equation}
 \big(R_j(\lambda + i0)f - R_j(\lambda-i0)f,f\big)
 = \lim_{r\to\infty}
2ik\int_{M_j}|u_j^{(\pm)}(r,x)|^2\, \sqrt{g_j(r,x)}\, dx 
\nonumber
\end{equation}
\end{lemma}

\begin{proof}
We prove the case for $u_j^{(+)}$. Integrating $u_j^{(+)}\overline{f} - f\overline{u_j^{(+)}}$ on the region
\begin{equation}
\widetilde{\mathcal M}_{j,\leq t} := (2,t)\times M_j,
\nonumber
\end{equation}
 we obtain
\begin{eqnarray}
 & & (u_j^{(+)},f)_{L^2(\widetilde{\mathcal M}_{j,\leq t})} - (f,u_j^{(+)})_{L^2(\widetilde{\mathcal M}_{j,\leq t})} 
\label{Lemma10.4Integbyparts} \\ 
& = &  
\int_{M_j}\left((\partial_tu_j^{(+)})\overline u_j^{(+)}- u_j^{(+)}\overline{\partial_tu_j^{(+)}}\right)\sqrt{g_j}\,dx
\nonumber
\\
&= &\int_{M_j}\left((D_j^{(+)}(k)u_j^{(+)})\overline{u_j^{(+)}} - u_j^{(+)}\overline{D_j^{(+)}(k)u_j^{(+)}}+ 2i\,{\rm Re}\,\psi_j^{(+)}|u_j^{(+)}|^2\right)\sqrt{g_j}\,dx.
\nonumber
\end{eqnarray}
We use Lemma \ref{LemmaDpm(k9uvto0} with $v_j^{(\pm)} = u_j^{(\pm)}$. 
Letting $t \to \infty$, we get the present lemma.   
\end{proof}

Let us note that  (\ref{Lemma9.3after}) and (\ref{Lemma10.4Integbyparts}) yield
\begin{equation}
\|u_j^{(+)}(r,\cdot)\|_{L^2(M_j)} \leq C\rho_j(r)^{-(n-1)}\|f\|_{L^{2,s}}, \quad {\rm if} \quad s = 2\beta_j.
\label{Lemma9.3afterafter}
\end{equation}

We show the existence of the limit of
\begin{equation}
\mathcal F_j^{(\pm)}(\lambda,r)f  = C_j(\lambda)^{-1} \rho_j(r)^{(n-1)/2}e^{\mp i\Phi_j(r,\lambda)}\big(R_j(\lambda \pm i0)f\big)(r,x),
\nonumber
\end{equation}
\begin{equation}
C_j(\lambda) = \Big(\frac{\pi}{\sqrt{\lambda - E_{0,j}}}\Big)^{1/2}, \quad 
E_{0,j} = \Big(\frac{(n-1)c_{0,j}}{2}\Big)^2
\label{S9Cj(lambda)Define}
\end{equation}
as $r\to\infty$, where $\Phi_j(r,\lambda)$ is defined by (\ref{DefPjij(r)}).


\begin{lemma}\label{S10weaklimitlemma}
For $f \in L^2_{comp}(\widetilde{\mathcal M}_j)$, there exists a weak limit
\begin{equation}
{\mathop{\rm w-lim}_{r\to\infty}}\, \mathcal F_j^{(\pm)}(\lambda,r)f =: \mathcal F_j^{(\pm)}(\lambda)f, \quad {\rm in} \quad L^2(M_j).
\nonumber
\end{equation}
\end{lemma}

\begin{proof}
 We prove the $+$ case.  
 Since $\sup_{r>2}\|\mathcal F_j^{(+)}(\lambda,r)f\|_{L^2(M_j)}$ is finite by  (\ref{Lemma9.3afterafter}), we have only to show the existence of
\begin{equation}
\lim_{r\to\infty}(\mathcal F_j^{(+)}(\lambda,r)f,\varphi)_{L^2(M_j)}
\label{Lemma9.4weaklimit}
\end{equation}
for $\varphi \in D(\Lambda_j)$, where $\Lambda_j$ is the Laplace-Beltrami operator on $M_j$. We put
$$
v = \rho_j(r)^{-(n-1)/2}e^{i\Phi_j}\varphi, \quad \varphi \in D(\Lambda_j).
$$
By a direct calculation, we have,
\begin{equation}
\begin{split}
\rho_j^{(n-1)/2}e^{-i\Phi_j}v' & = \Big(- \frac{n-1}{2}\frac{\rho_j'}{\rho_j} + i\Phi_j'\Big)\varphi,\\
\rho_j^{(n-1)/2}e^{-i\Phi_j}v'' & = \Big(\frac{(n-1)^2}{4}\Big(\frac{\rho_j'}{\rho_j}\Big)^2 - \frac{n-1}{2}\Big(\frac{\rho_j'}{\rho_j}\Big)' - i(n-1)\frac{\rho_j'}{\rho_j}\Phi_j' \\
&\ \ \ \ \  + i\Phi_j'' - \big(\Phi_j'\big)^2\Big)\varphi.
\end{split}
\nonumber
\end{equation}
By our assumption, 
\begin{equation}
\frac{\rho'_j}{\rho_j} = c_{0,j} + O(r^{-\alpha_{0,j}}), \quad
\frac{g_j'}{2g_j} = (n-1)c_{0,j} + O(r^{-\alpha_{0,j}}).
\nonumber
\end{equation}
Therefore, modulo a term of $O(r^{-1-\epsilon})$
\begin{equation}
 \rho_j^{(n-1)/2}e^{-i\Phi_j}\Big( v'' + \frac{g_j'}{2g_j}v'\Big)  \equiv - \lambda \varphi
\label{S10Lemma10.4Proofvj''etc}
\end{equation}
in the sense of $L^2(M_j)$-norm.
Since $2\beta_j > 1$, we have
\begin{equation}
\begin{split}
B_j\varphi & = O(\rho_j(r)^{-2}\rho_j(r)^{-(n-1)/2}) = O(r^{-2\beta_j}\rho_j(r)^{-(n-1)/2}) \\
& = O(r^{-1-\epsilon}\rho_j(r)^{-(n-1)/2}),
\end{split}
\nonumber
\end{equation}
hence
$$
\big(- \Delta_{\mathcal M_j} - \lambda\big)v = O(r^{ - 1 - \epsilon}\rho_j(r)^{-(n-1)/2}).
$$ 
We let 
$$
g = \big(- \Delta_{\mathcal M_j} - \lambda\big)v, 
$$
and integrate by parts as in (\ref{Lemma10.4Integbyparts}) to obtain
\begin{equation}
\begin{split}
& (u_j^{(+)},g)_{L^2(\widetilde{\mathcal M}_{j,\leq t})} - (f,v)_{L^2(\widetilde{\mathcal M}_{j,\leq t})} \\
&= \int_{M_j}\left((D_j^{(+)}(k)u_j^{(+)})\overline{v} - u_j^{(+)}\overline{D_j^{(+)}(k)v}+ 2i\,({\rm Re}\,\psi_j^{(+)})\,u_j^{(+)}\overline{v}\right)\sqrt{g_j}\,dx.
\end{split}
\label{Lemma10.5Integbyparts}
\end{equation}
Letting $t \to \infty$ and using Lemma \ref{LemmaDpm(k9uvto0}, we have proven the existence of the limit (\ref{Lemma9.4weaklimit}).  
\end{proof}


\begin{lemma}\label{Lemma10.5}
Let $f \in L^2_{comp}(\widetilde{\mathcal M}_j)$, and $u_j^{(+)} = R_j(\lambda + i0)f$. \\
\noindent
(1) There exist $\epsilon > 0$ and  a sequence $\{r_p\}_{p=1}^{\infty}$ tending to infinity such that
\begin{equation}
r_p^{1+\epsilon}\int_{M_j}|(D_j^{(+)}(k)u_j^{(+)})(r_p,x)|^2\sqrt{g_j(r_p,x)}\, dx \to 0.
\nonumber
\end{equation}
(2) Let $\{r_p\}_{p=1}^{\infty}$ be as above. Then, for any $r_q  > r_p$ and $\varphi \in H^1(M_j)$
\begin{equation}
\begin{split}
& \big|\big(\mathcal F_j^{(+)}(\lambda,r_p)f - \mathcal F_j^{(+)}(\lambda,r_q)f,\varphi\big)_{L^2(M_j)}\big| \\
\leq & \;  C(r_p)\Big(\|\varphi\|_{L^2(M_j)} + r_p^{- \epsilon}\|\sqrt{\Lambda_j(r_p)}\varphi\|_{L^2(M)}\Big),
\end{split}
\nonumber
\end{equation}
 where the constant $C(r_p)$ is independent of $\varphi$,  and $C(r_p) \to 0$ as $r_p \to \infty$.
\end{lemma}

\begin{proof}
By Lemmas \ref{Lemma6.7vestimateimprove} and \ref{Cuspestimate}, there exists an $\epsilon > 0$ such that
$$
\int_2^{\infty}\int_{M_j}r^{\epsilon}|D_j^{(+)}(k)u_j^{(+)}|^2\sqrt{g_j}\, drdx < \infty,
$$
which implies (1) immediately. 

Let $v$ and $g$ be as in the proof of Lemma \ref{S10weaklimitlemma}. Letting
$$
\mathcal M_{s,t} = (s,t)\times M, \quad
F\Big|_s^t = F(t) - F(s),
$$
we have, by integration by parts,
\begin{equation}
\begin{split}
& (v,f)_{L^2(\mathcal M_{r_p, r_q})} - (g,u_j^{(+)})_{L^2(\mathcal M_{r_p, r_q})} \\
& + \int_{M_j}(D_j^{(+)}(k)v)\overline{u_j^{(+)}}\sqrt{g_j}\Big|_{r_p}^{r_q}dx - 
\int_{M_j}v\overline{D_j^{(+)}(k)u_j^{(+)}}\sqrt{g_j}\Big|_{r_p}^{r_q}dx \\
&=  2i \int_{M_j}({\rm Re}\,\psi_j^{(+)})v\overline{u_j^{(+)}}\sqrt{g_j}\Big|_{r_p}^{r_q}dx.
\end{split}
\label{S10Lemma10.5DifferennceFormula}
\end{equation}
By (\ref{Lemma9.3afterafter}),  $\|\mathcal F_j^{(+)}(\lambda,r)f\|_{L^2(M_j)}$ is uniformly bounded for $r > 2$. Since $\psi = k + O(r^{-\epsilon})$, the right-hand side of (\ref{S10Lemma10.5DifferennceFormula}) is estimated from below as follows:
\begin{equation}
\begin{split}
&  \Big| 2i \int_{M_j}({\rm Re}\,\psi_j^{(+)})v\overline{u_j^{(+)}}\sqrt{g_j}\Big|_{r_p}^{r_q}dx\Big|\\
& \ \geq  C\Big|\big(\mathcal F_j^{(+)}(\lambda,r_p)f - \mathcal F_j^{(+)}(\lambda,r_q)f,\varphi\big)_{L^2(M_j)}\Big| \\
&\  - C'r_p^{-\epsilon}\big(\sup_{q\geq p}\|\mathcal F_j^{(+)}(\lambda,r_q)f\|_{L^2(M_j)}\big)\|\varphi\|_{L^2(M_j)}.
 \end{split}
 \label{Lemma10.5kefhandsidefrombelow}
\end{equation}
We estimate the left-hand side of (\ref{S10Lemma10.5DifferennceFormula}) from above. 
The 4th term is estimated as follows:
\begin{equation}
\begin{split}
& \Big|\int_{M_j}v\overline{D_j^{(+)}(k)u_j^{(+)}}\sqrt{g_j}\Big|_{r_p}^{r_q}dx\Big| \\
& \leq C\sup_{q\geq p}\Big(\int_{M_j}|D_j^{(+)}(k)u_j^{(+)}(r_q,x)|^2\sqrt{g_j(r_q,x)}dx\Big)^{1/2}
\|\varphi\|_{L^2(M_j)}.
\end{split}
\label{Lemma105.3rdtermbelow}
\end{equation}
Since
$$
D^{(+)}(k)v = \Big(\frac{g_j'}{4g_j} - \frac{n-1}{2}\frac{\rho_j'}{\rho_j} + i(\Phi_j' - \psi_j)\Big)v = O(r^{-1-\epsilon})v,
$$
the 3rd term is dominated by $Cr_p^{-1-\epsilon}\|\varphi\|_{L^2(M_j)}$. 
The 1st term is estimated as follows
\begin{equation}
\big|(v,f)_{L^2(\mathcal M_{r_p,r_q})}\big| \leq Cr_p^{-\epsilon/2}\|\varphi\|_{L^2(M_j)}
\|f\|_{L^{2,(1+\epsilon)/2}}.
\label{Lemma10.5Proof1sttermestimate}
\end{equation}
We estimate the 2nd term.
Recalling (\ref{S10Lemma10.4Proofvj''etc}), we have
$$
g = (-\Delta_{{\mathcal M}_j}-\lambda)v = \alpha(r)\varphi + \beta(r)B_j(r)\varphi,
$$
$$
\alpha(r) = O(r^{-1-\epsilon}\rho_j(r)^{-(n-1)/2}), \quad
\beta(r) = O(\rho_j(r)^{-(n-1)/2}).
$$
The term $(\alpha(r)\varphi,u_j^{(+)})_{L^2(\mathcal M_{r_p,r_q})}$ is  estimated as
\begin{equation}
\big|(\alpha(r)\varphi, u_j^{(+)})_{L^2(\mathcal M_{r_p,r_q})}\big| \leq Cr_p^{-\epsilon/2}\|\varphi\|_{L^2(M_j)}
\|u_j^{(+)}\|_{-(1+\epsilon)/2}.
\label{Lemma10.5Proof1sttermestimate2}
\end{equation}
The term $(\beta(r)B_j(r)\varphi,u_j^{(+)})_{L^2(\mathcal M_{r_p,r_q})}$ is rewritten as
\begin{equation}
\begin{split}
(\beta(r)B_j(r)\varphi,u_j^{(+)})_{L^2(\mathcal M_{r_p,r_q})}= 
(\beta(r)\sqrt{B_j(r)}\varphi,\sqrt{B_j(r)}u_j^{(+)})_{L^2(\mathcal M_{r_p,r_q})}.
\end{split}
\nonumber
\end{equation}
By virtue of Lemma \ref{Lemma6.7vestimateimprove}, letting $t = (s-1)/2$ where $0 < s < \min\{2\beta_{0,j},\gamma_{0,j}\}$, this is estimated from above as
\begin{equation}
\begin{split}
&\big|(\beta(r)\sqrt{B_j(r)}\varphi,\sqrt{B_j(r)}u_j^{(+)})_{L^2(\mathcal M_{r_p,r_q})}\big| \\
& \leq C\|\sqrt{\Lambda_j(r_p)}\varphi\|_{L^2(M_j)}\|r^{-t}\beta(r)\rho_j(r)^{-1}\|_{L^2((r_p,r_q);\rho_j^{n-1}dr)}\\
& \ \ \ \ \ \ \times \|r^{t}\sqrt{B_j(r_p)}u_j^{(+)}\|_{L^2(\mathcal M_{r_p,r_q})}\\
&\leq C\|\sqrt{\Lambda_j(r_p)}\varphi\|_{L^2(M_j)}\|r^{-t}\rho_j(r)^{-(n+1)/2}\|_{L^2((r_p,r_q);\rho_j^{n-1}dr)}\|f\|_{(s+1)/2} \\
& \leq C\|\sqrt{\Lambda_j(r_p)}\varphi\|_{L^2(M_j)}\|f\|_{(s+1)/2}\Big(\int_{r_p}^{r_q}r^{-2t-2\beta_{0,j}}dr\Big)^{1/2} \\
& \leq Cr_p^{-t}\|\sqrt{\Lambda_j(r_p)}\varphi\|_{L^2(M_j)}\|f\|_{(s+1)/2}.
\end{split}
\nonumber
\end{equation}
This and the inequalities (\ref{Lemma10.5kefhandsidefrombelow}),  (\ref{Lemma105.3rdtermbelow}), (\ref{Lemma10.5Proof1sttermestimate}),  (\ref{Lemma10.5Proof1sttermestimate2}), together with (1),  prove the lemma. 
\end{proof}

\begin{lemma}\label{Lemma10.6}
Let $f \in L^2_{comp}(\widetilde{\mathcal M_j})$. Then, there exists a sequence $\{r_p\}_{p=1}^{\infty}$ such that $\mathcal F_j^{(\pm)}(\lambda,r_p)f$ converges to $\mathcal F_j^{(\pm)}(\lambda)f$ strongly on $L^2(M_j)$.
\end{lemma}

\begin{proof}
 Let $\{r_p\}_{p=1}^{\infty}$ be the sequence in Lemma \ref{Lemma10.5} (1). We can also assume that
$$
r_p^{1 + \epsilon}\|\sqrt{B_j(r_p)}u_j^{(+)}\|^2_{L^2(M_j)} \to 0, \quad u_j^{(+)} = R_j(\lambda + i0)f.
$$
Using Lemma \ref{Lemma10.5} with $\varphi$ replaced by $\mathcal F_j^{(+)}(\lambda,r_p)f$, we have
\begin{equation}
\begin{split}
& \Big|\|\mathcal F_j^{(+)}(\lambda,r_p)f\|^2_{L^2(M_j)} - \big(\mathcal F_j^{(+)}(\lambda,r_q)f,\mathcal F_j^{(+)}(\lambda,r_p)f\big)_{L^2(M_j)}\Big| \\
& \leq C(r_p)\left(\|\mathcal F_j^{(+)}(\lambda,r_p)f\|_{L^2(M_j)}
+ r_p^{-\epsilon}\|\sqrt{B_j(r_p)}\mathcal F_j^{(+)}(\lambda,r_p)f\|_{L^2(M_j)}
\right).
\end{split}
\nonumber
\end{equation}
We first let $q$ tend to $\infty$, and then $p$ to $\infty$. Then, we see that the norm of  $\mathcal F_j^{(+)}(\lambda,r_p)f$ converges to that of $\mathcal F_j^{(+)}(\lambda)f$. This proves the strong convergence, since we already know its weak convergence in Lemma \ref{S10weaklimitlemma}. 
\end{proof}


\begin{lemma}\label{ParsevalFlambda}
Let $f \in L^2_{comp}(\widetilde{\mathcal M_j})$. Then,  $\mathcal F_j^{(\pm)}(\lambda,r)f$ converges to $\mathcal F_j^{(\pm)}(\lambda)f$ strongly on $L^2(M_j)$, and we have
\begin{equation}
\frac{1}{2\pi i}\big(R_j(\lambda + i0)f-R_j(\lambda-i0)f,f\big) = 
\|\mathcal F_j^{(\pm)}(\lambda)f\|^2_{L^2(M_j)}
\nonumber
\end{equation}
\end{lemma}

\begin{proof}
Sinec $\mathcal F_j^{(+)}(\lambda)f$ is the weak limit of $\mathcal F_j^{(+)}(\lambda,r)f$, we have
$$
\|\mathcal F_j^{(+)}(\lambda)f\|_{L^2(M_j)} \leq \liminf_{r\to\infty}\|\mathcal F_j^{(+)}(\lambda,r)f\|_{L^2(M_j)}.
$$
Using the sequence in Lemma \ref{Lemma10.6}, we have
$$
\liminf_{r\to\infty}\|\mathcal F_j^{(+)}(\lambda,r)f\|_{L^2(M_j)} \leq 
\lim_{r_p\to\infty}\|\mathcal F_j^{(+)}(\lambda,r_p)f\|_{L^2(M_j)} = \|\mathcal F_j^{(+)}(\lambda)f\|_{L^2(M_j)}.
$$
By Lemma \ref{S10Lemma10.3}, we already know the existence of the limit of $\|\mathcal F_j^{(+)}(\lambda,r)f\|_{L^2(M_j)}$. Therefore, we have $
\|\mathcal F_j^{(+)}(\lambda)f\|_{L^2(M_j)} = \lim_{r\to\infty}\|\mathcal F_j^{(+)}(\lambda,r)f\|_{L^2(M_j)}$. This proves the existence of the strong limit $\lim_{r\to\infty}\mathcal F_j^{(+)}(\lambda,r)f$ in $L^2(M_j)$. In view of Lemma \ref{S10Lemma10.3}, we obtain the lemma. 
\end{proof}
\subsection{Asymptotic expansion of the resolvent}
Lemma \ref{ParsevalFlambda} implies that
\begin{equation}
\|\mathcal F_j^{(\pm)}(\lambda)f\|_{L^2(M_j)} \leq C\|f\|_{\mathcal B}, \quad f \in L^2_{comp}(\mathcal M),
\label{FboudedBtoh}
\end{equation}
where the constant $C$ is independent of $\lambda$ varying over a compact set in $\sigma_e(H)\setminus \mathcal E$. Therefore, it is uniquely extended on $\mathcal B$ and 
(\ref{FboudedBtoh}) holds also for $f \in \mathcal B$. Making use of (\ref{S9chiR=Hjreolventchij}), we compute the asymptotics of $R(\lambda \pm i0)$ on each  end. 
For $f, g \in \mathcal B^{\ast}$, we use the folloing notation: \index{$f \simeq g$}
\begin{equation}
f \simeq g \Longleftrightarrow f-g \in \mathcal B^{\ast}_0.
\label{Definefsimeqg}
\end{equation}


\begin{theorem}\label{S10ResovAsympReg}
For any $f \in \mathcal B$ and $\lambda \in (E_{0,j},\infty)\setminus \mathcal E$, we have on the regular ends $\mathcal M_j$ satisfying (\ref{betaj>1/2})
\begin{equation}
R(\lambda \pm i0)f \simeq C_j(\lambda)\rho_j(r)^{-(n-1)/2}e^{\pm i\Phi_j(r,\lambda)}\mathcal F_j^{(\pm)}(\lambda)f,
\label{EndResolvAsym}
\end{equation}
where $C_j(\lambda)$ is defined by (\ref{S9Cj(lambda)Define}).
\end{theorem}

\begin{proof}
For $f \in L^2_{comp}(\mathcal M)$, (\ref{EndResolvAsym}) is already proven.  Since $R(\lambda \pm i0) \in {\bf B}(\mathcal B;\mathcal B^{\ast})$, in view of (\ref{FboudedBtoh}), we have only to approximate $f \in \mathcal B$ by compactly supported functions to prove (\ref{EndResolvAsym}). 
\end{proof}


\subsection{Cusp ends}
For the cusp ends, one can argue as above without any change. Thus, we have also proven the following theorem.


\begin{theorem}\label{S10ResovAsympCusp}
For any $f \in \mathcal B$ and $\lambda \in (E_{0,j},\infty)\setminus \mathcal E$, we have on  the cusp ends $\mathcal M_j$,
\begin{equation}
R(\lambda \pm i0)f \simeq C_j(\lambda)P_{0,j}\otimes\rho_j(r)^{-(n-1)/2}e^{\pm i\Phi_j(r,\lambda)}\mathcal F_j^{(\pm)}(\lambda)f, 
\label{CuspEndResolAsym}
\end{equation}
 $P_{0,j}$ being the projection associated with the $0$-eigenvalue for $-\Delta_{M_j}$. Here,  
$C_j(\lambda)$ is defined by (\ref{S9Cj(lambda)Define}).
\end{theorem}

\begin{proof}
 By Lemma \ref{Cuspestimate}, $(1-P_{0,j}\otimes 1)R(\lambda\pm i0)f \in L^2(\mathcal M_j)$ for the cusp end, from which (\ref{CuspEndResolAsym}) follows. 
\end{proof}


\section{Resolvent asymptotics - Perturbative approach}\label{ResolventAsympto2}
We study the remaining case in this section. 
We fix one regular end $\mathcal M_j = (0,\infty)\times M_j$, (we shifted the interval $(2,\infty)$ to $(0,\infty)$, which does not matter at all), and observe  the  asymptotic behavior of the resolvent on $\mathcal M_j$ under the assumption
\begin{equation}
c_{0,j}=0, \quad \alpha_{0,j} > 0,  \quad 1/2 \geq \beta_{0,j}  > 0, \quad  \gamma_{0,j} > 1.
\label{RegCuspSRAssumption}
\end{equation}
 In this section, we drop the subscript $j$ in  $\rho_j(r)$, $c_{0,j}, \alpha_{0,j}, \beta_{0,j}, \gamma_{0,j}$. 

We reduce the problem to the  one-dimensional case, and apply the perturbation technique. 
Hence, letting  $E\geq 0$ be an arbitrary constant, we start with the  equation
\begin{equation}
\Big(- {\partial_r}^2 - \frac{(n-1) \rho'(r)}{ \rho(r)}\partial_r + \frac{E}{ \rho(r)^{2}} - z\Big)u  = f, \quad {\rm on} \quad (0,\infty).
\label{S111-dimeeq}
\end{equation}

\subsection{WKB method}

\subsubsection{Asymptotic solutions}
We seek an asymptotic solution of  the equation
\begin{equation}
\left(- {\partial_r}^2 - \frac{(n-1) \rho'(r)}{ \rho(r)}\partial_r + \frac{E}{ \rho(r)^{2}}\right)u = \lambda u, \quad \lambda > 0
\label{S2Homogeneouseq}
\end{equation}
in the form 
$u = e^{\varphi}a$. 
A direct computation yields
\begin{equation}
\begin{split}
& e^{-\varphi}\left(- {\partial_r}^2 - \frac{(n-1) \rho'}{ \rho}\partial_r + \frac{E}{\rho^{2}} - \lambda\right)e^{\varphi}a \\
 =&  -\Big((\varphi')^2 + \frac{(n-1)\rho'}{\rho}\varphi' + \lambda - \frac{E}{\rho^2}\Big)a   - \left\{\Big(2\varphi'+ \frac{(n-1)\rho'}{\rho}\Big)a' + \varphi''a \right\} - a''.
\end{split}
\label{S2Equation10}
\end{equation}
By  the assumptions (A-1) and (A-2), we have 
$$
\frac{\rho'}{\rho} \in S^{-\alpha_0}, \quad 
\partial_r^m\big(\frac{1}{\rho}\big) \in S^{-\beta - m}, \quad \forall m \geq 0.
$$
Letting $\epsilon = \min(2\alpha_0,2\beta_0)$, we then have
\begin{equation}
\Big(\frac{(n-1)\rho'}{2\rho}\Big)^2 + \frac{E}{\rho^2} \leq
Cr^{-\epsilon}(1 +E).
\nonumber
\end{equation}
We define $r_0(\lambda,E)$ by 
\begin{equation}
r_0(\lambda,E) = \left(\frac{2C(1+E)}{\lambda}\right)^{1/\epsilon}
\label{Definer0lambdaE}
\end{equation}
so that $Cr^{-\epsilon}(1 +E) < \lambda/2$ for $r > r_0(\lambda,E)$. We put
\begin{equation}
\alpha(\lambda,E,r) = \sqrt{\lambda - \Big(\frac{(n-1)\rho'}{2\rho}\Big)^2 - \frac{E}{\rho^2}}, 
\label{S2Definealphalambda}
\end{equation}
\begin{equation}
\varphi_{\pm} = - \frac{n-1}{2}\log \rho(r) \pm i
\int_{r_0(E,\lambda)}^r\alpha(\lambda,E,s)ds.
\label{S2Varphi'define}
\end{equation}
\begin{equation}
a_{\pm} = \exp\left(\int_{r}^{\infty}\frac{\varphi_{\pm}''}{2\varphi_{\pm}' + (n-1)\rho'/\rho}ds\right).
\label{S10Defineapm}
\end{equation}
They satisfy the eikonal equation
\begin{equation}
\big(\varphi')^2 + \frac{(n-1)\rho'}{\rho}\varphi' + \lambda - \frac{E}{\rho^2} = 0,
\nonumber
\end{equation}
and the transport equation
\begin{equation}
\Big(2\varphi'+ \frac{(n-1)\rho'}{\rho}\Big)a' + \varphi''a = 0.
\nonumber
\end{equation}
Hence the equation (\ref{S2Equation10}) becomes
\begin{equation}
 e^{-\varphi}\left(- {\partial_r}^2 - \frac{(n-1) \rho'}{ \rho}\partial_r + \frac{E}{\rho^{2}} - \lambda\right)e^{\varphi}a 
 = - a''.
\nonumber
\end{equation}
Fix a compact interval $I \subset (0,\infty)$ arbitrarily.
In the following, $C$'s denote constants independent of $\lambda \in I$, $E \geq 0$ and $r > r_0(\lambda,E)$. Similarly, the various estimates are uniform with respect to $\lambda, E$ and $r$ satisfying these conditions.
The following lemma is proven by a direct computation.


\begin{lemma}\label{S10rhoinversinS-beta-m}
For $\alpha$, $\varphi_{\pm}$, $a_{\pm}$ defined by (\ref{S2Definealphalambda}), (\ref{S2Varphi'define}), (\ref{S10Defineapm}), we have
$$
\alpha \in S^0, \quad  \alpha' \in S^{-1-\epsilon}, \quad
\varphi_{\pm}' \in S^0,  \quad\varphi_{\pm}'' \in S^{-1-\epsilon}, \quad 
a_{\pm} - 1 \in S^{-\epsilon}.
$$
\end{lemma}
Noting that
\begin{equation}
e^{\varphi_{\pm}} = \rho(r)^{-(n-1)/2}\exp\Big(\pm i\int_{r_0(E,\lambda)}^r\alpha(\lambda,E,s)ds\Big),
\nonumber
\end{equation}
and summarizing the above computation, we have  proven the following lemma. 
 

\begin{lemma} \label{S2Lemma1}
There  exist asymptotic solutions $a_{\pm}e^{\varphi_{\pm}}$ of (\ref{S2Homogeneouseq}) satisfying
\begin{equation}
\Big|a_{\pm}e^{\varphi_{\pm}}  - \rho(r)^{-(n-1)/2}\exp\Big(\pm i\int_{r_0(E)}^r\alpha(\lambda,E,s)ds\Big)\Big|
\leq C\rho(r)^{-(n-1)/2}r^{-\epsilon}, 
\label{Lemma11.2Estimate}
\end{equation}
\begin{equation}
\Big(-\partial_r^2 - \frac{(n-1)\rho'}{\rho}\partial_r + \frac{E}{\rho^2}- \lambda\Big)a_{\pm}e^{\varphi_{\pm}} =
- a_{\pm}''e^{\varphi_{\pm}}, 
\label{Lemma11.2Equationforaevaphi}
\end{equation}
\begin{equation}
\partial_r^m(a_{\pm} - 1) = O(r^{-m-\epsilon}), \quad m \geq 0,
\nonumber
\end{equation}
uniformly for $ r > r_0(\lambda,E)$.
\end{lemma}


\subsubsection{Exact solutions}
Next let us construct the exact solutions to  (\ref{S2Homogeneouseq}) which behave like $a_{\pm}e^{\varphi_{\pm}}$ as $r \to \infty$. Putting $a = a_{\pm}, \varphi = \varphi_{\pm}$, $u = ae^{\varphi}(1 + v)$ and using (\ref{Lemma11.2Equationforaevaphi}), we have for $r > r_0(E,\lambda)$
\begin{equation}
v'' + \Big(2(\frac{a'}{a} + \varphi') + \frac{(n-1)\rho'}{\rho}\Big)v' + \frac{a''}{a} v =- \frac{a''}{a}.
\nonumber
\end{equation}
Putting 
${\bf v}= \left(\begin{array}{c}
v \\ v'\end{array}\right)$, ${\bf f}= \left(\begin{array}{c} 0 \\ - a''/a\end{array}\right)$, we transform it into the 1st order system: 
\begin{equation}
\frac{d{\bf v}}{dr} = K(r){\bf v}+ L(r){\bf v}+ {\bf f}(r),
\label{S2DiffEq}
\end{equation}
\begin{equation}
K(r) = \left(\begin{array}{cc}
0 & 1 \\ 0 & - \Big(2(\frac{a'}{a} + \varphi') + \frac{(n-1)\rho'}{\rho}\Big)
\end{array}
\right), \quad 
L(r) = \left(\begin{array}{cc}
0 & 0 \\ -a''/a& 0
\end{array}\right).
\nonumber
\end{equation}
A fundamental matrix for the equation $d{\bf v}/dr =K(r){\bf v}$ is
\begin{equation}
W(r,s) = V(r)V(s)^{-1}, \quad V(r) = \left(\begin{array}{cc}1& F \\  0 &a^{-2}e^{-2\varphi}\rho^{-(n-1)}\end{array}
\right),
\nonumber
\end{equation}
where
$$
F = \int_{r_0(E,\lambda)}^ra^{-2}e^{-2\varphi}\rho^{-(n-1)}ds. 
$$ 
Then by (\ref{Lemma11.2Estimate}),  
$$
F' = a^{-2}e^{-2\varphi}\rho^{-(n-1)} = \exp\Big(\mp 2i\int_{r_0}^r\alpha \, ds\Big)
\Big( 1 + O(r^{-\epsilon})\Big).
$$
Using
$$
\frac{1}{\mp 2i\alpha}\frac{d}{dr}e^{\mp 2i\int_{r_0}^r\alpha ds} = 
e^{\mp 2i\int_{r_0}^r\alpha ds}
$$
and  integrating by parts, we have $F = O(1)$.
  Then $W(r,s)$ is uniformly bounded for $r_0(\lambda,E) \leq r \leq s$. The equation (\ref{S2DiffEq}) is  rewritten as  the integral equation
\begin{equation}
{\bf v}(r) = - \int_r^{\infty}W(r,s)L(s){\bf v}(s)ds - \int_r^{\infty}W(r,s){\bf f}(s)ds,
\nonumber
\end{equation}
which is solved by iteration, since $L(r) = O(r^{-1-\epsilon})$, ${\bf f}(r) = O(r^{-1-\epsilon})$. 
We have thus proved the  following lemma.


\begin{lemma} \label{S2Lemma2}
There  exist exact solutions $\Psi^{(\pm)}(\lambda,r,E)$  to (\ref{S2Homogeneouseq}) on $[r_0(\lambda,E),\infty)$ such that 
$\Psi^{(\pm)} = a_{\pm}e^{\varphi_{\pm}}\big(1 + O(r^{-\epsilon})\big)$ as $r \to \infty$.
\end{lemma}

We extend the solutions  $\Psi^{(\pm)}(\lambda,r;E)$  to the whole interval $[0,\infty)$.
The following lemma is an immediate consequence.


\begin{lemma}\label{S2Lemma3}
If $u$ satisfies (\ref{S2Homogeneouseq}) and
$$
\frac{1}{R}\int_0^R|u(r)|^2 \rho(r)^{n-1}dr \to 0, \quad R \to \infty,
$$
then $u$ is identically equal to 0.
\end{lemma}

\begin{proof}
 Since $\Psi^{(\pm)}$ in Lemma \ref{S2Lemma2} are linearly independent, $u$ is written as $u = c_+\Psi^{(+)} + c_-\Psi^{(-)}$ for some constants $c_{\pm}$. The assumption of the lemma then implies
$$
\frac{1}{R}\int_0^R|c_+e^{i\int_0^{r}\alpha ds} + c_-e^{-i\int_0^r \alpha ds}|^2dr \to 0.
$$
For large $r$, we can make the change of variable $r \to t$ by $\int_0^r\alpha ds = t$. We then have
$$
\frac{1}{R}\int_0^R|c_+e^{it} + c_-e^{-it}|^2dt \to 0,
$$
which implies $c_+ = c_- = 0$. 
\end{proof}


\subsubsection{Green function}
The Green operator $G^{(\pm)}(\lambda,E)$ for (\ref{S2Homogeneouseq}) with Dirichlet condition at $r=0$ is defined by
\begin{eqnarray*}
\big(G^{(\pm)}(\lambda,E)f\big)(r) &= &
\int_0^{\infty}G^{(\pm)}(r,s,\lambda,E)f(s)ds, \\
G^{(\pm)}(r,s,\lambda,E) &=& \frac{1}{-W^{(\pm)}(\lambda,s,E)}
\left\{
\begin{split}
& \Psi_0(\lambda,r,E)\Psi^{(\pm)}(\lambda,s,E), \quad 0 < r < s, \\
& \Psi^{(\pm)}(\lambda,r,E)\Psi_{0}(\lambda,s,E), \quad 0 < s < r,
\end{split}
\right.
\label{S11G(rslambdaB)} 
\\
W^{(\pm)}(\lambda,r,E) &=& \Psi_0(\lambda,r,E)\Psi^{(\pm)}(\lambda,r,E)' - \Psi_{0}(\lambda,r,E)'\Psi^{(\pm)}(\lambda,r,E),
\label{S11Wronskian}
\end{eqnarray*}
where  $\Psi_0(r) = \Psi_0(\lambda,r,E)$ is the solution of (\ref{S2Homogeneouseq}) satisfying $\Psi_0(0) = 0,  \Psi_0'(0) = 1$.
Since ${W^{(\pm)}}' = - \frac{(n-1)\rho'}{\rho}W^{(\pm)}$, we have
\begin{equation}
W^{(\pm)}(\lambda,r,E) = - \Psi^{(\pm)}(\lambda,0,E)\Big(
\frac{\rho(0)}{\rho(r)}\Big)^{n-1}.
\nonumber
\end{equation}
Note that $\Psi^{(\pm)}(\lambda,0,E)\neq0$. In fact, if it vanishes, $\Psi^{(\pm)}(\lambda,r,E)$ is a solution to (\ref{S2Homogeneouseq}) satisfying the boundary condition and the radiation condition.  Arguing in the same  way as in the proof of Lemma \ref{D+ktoB0ast} using Lemma \ref{S2Lemma3} (actualy much simpler), we see that $\Psi^{(\pm)}(\lambda,r,E) = 0$, which is a contradiction.

The following Lemma can be proven easily by using the Green function.


\begin{lemma}\label{S2Resolventlimit}
For $f \in L^2_{comp}((0,\infty))$, we put
\begin{equation}
\begin{split}
\widetilde f^{(\pm)}(\lambda,E) & = - \int_0^{\infty}\frac{\Psi_0(\lambda,s,E)}{W^{(\pm)}(\lambda,s,E)}f(s)ds \\
 & = \frac{\rho(0)^{1-n}}{\Psi^{(\pm)}(\lambda,0,E)}\int_0^{\infty}\Psi_0(\lambda,s,E)f(s)\rho(s)^{n-1}\,ds.
\end{split}
\label{S2ftildelambdamudefine}
\end{equation}
Then  if $f(r) = 0$ for $r > r'$, 
\begin{equation}
\big(G^{(\pm)}(\lambda,E)f\big)(r) = \Psi^{(\pm)}(\lambda,r;E)\widetilde f^{(\pm)}(\lambda,E), 
\quad r > r'
\label{S2GpmlambdamufPsitildef}
\end{equation}
holds, and the following limit exists
\begin{equation}
\lim_{r\to\infty}\rho(r)^{(n-1)/2}e^{\mp i\int_0^r\alpha(\lambda,E,t)dt}\big(G^{(\pm)}(\lambda,E)f\big)(r)= \widetilde f^{(\pm)}(\lambda,E).
\label{GpmlambdaEflimit}
\end{equation}
\end{lemma}


\subsection{Fourier transform}
Let $E \geq 0$ be a constant,  $L(E)$ the differential operator
$$
L(E) = -\partial_r^2 - \dfrac{(n-1)\rho'}{\rho}\partial_r + 
\dfrac{E^2}{\rho^2}
$$
 with Dirichlet boundary condition at $r=0$, and
\begin{equation}
R_{free}(z,E) = (L(E)-z)^{-1}.
\nonumber
\end{equation}
All results in the previous sections  hold for $L(E)$. 
In particular, Theorem \ref{ModelLAPforbfh} holds for $L(E)$ with $E \geq 0$. 
Letting
\begin{equation}
\varphi(\lambda,E,r) = \int_{r_0(\lambda,E)}^r\alpha(\lambda,E,s)ds,
\label{Definevarphi(lambdaBr)}
\end{equation}
we define for $f \in L^2_{comp}((0,\infty))$
\begin{equation}
\begin{split}
 F_{free,E}^{(\pm)}(\lambda)f  = \Big(\frac{\sqrt{\lambda}}{\pi}\Big)^{1/2}
\lim_{r\to\infty}\rho(r)^{(n-1)/2}e^{\mp i  \varphi(\lambda,E,r)}\left(R_{free}(\lambda \pm i0,E)f\right)(r).
\end{split}
\label{S8Ffreelambdadef}
\end{equation}
The existence of the limit (\ref{S8Ffreelambdadef}) is guaranteed by Theorem \ref{S2Resolventlimit}.


\begin{lemma}\label{1dimFfreemuresolvent}
For $f \in L^2_{comp}((0,\infty))$, we have
$$
|F_{free,E}^{(\pm)}(\lambda)f|^2 = \frac{1}{2\pi i}\left(\big[R_{free}(\lambda+i0,E)- R_{free}(\lambda-i0,E)\big]f,f\right).
$$
\end{lemma}

\begin{proof}
Let $u_{\pm} = R_{free}(\lambda \pm i0,E)f$. Multiply the equation 
$$
- u_{\pm}'' - \frac{(n-1)\rho'}{\rho}u_{\pm}' + \Big(\frac{E}{\rho^2}-\lambda\Big)u_{\pm} = f
$$
by $\overline{u_{\pm}}\rho^{n-1}$, integrate by parts over $(0,r)$, and take the imaginary part. Then we have
\begin{equation}
{\rm Im}\,u_{\pm}'(r)\overline{u_{\pm}(r)}\rho^{n-1} = - {\rm Im}\,\int_0^rf\overline{u_{\pm}}\rho^{n-1}dt.
\label{S8u'ubar=intfubar}
\end{equation}
The left-hand side is equal to 
$$
{\rm Im}\, (D_{\pm}(k)u_{\pm})\overline{u_{\pm}}\rho^{n-1} + {\rm Re}\,\psi_j^{(\pm)}|u_{\pm}|^2\rho^{n-1}.
$$
By Lemma \ref{Cuspestimate}, $D_{\pm}(k)u \in L^2$. Therefore, the 1st term tends to 0 along a suitable sequence 
$r_1 < r_2 < \cdots \to \infty$. Taking the limit in (\ref{S8u'ubar=intfubar}) along $\{r_n\}$, we then have
$$
k\lim_{r\to\infty}|u_{\pm}|^2\rho^{n-1} = \mp {\rm Im}\, (f,u_{\pm}) = 
\frac{1}{2i}\left([R_{free}(\lambda +i0,E) - R_{free}(\lambda-i0,E)]f,f\right).
$$
Noting $k = \sqrt{\lambda}$ and (\ref{S8Ffreelambdadef}), we prove the lemma 
\end{proof}

Lemma \ref{1dimFfreemuresolvent} and Theorem \ref{ModelLAPforbfh} imply
\begin{equation}
|F^{(\pm)}_{free,E}(\lambda)f| \leq C\|f\|_{\mathcal B},
\nonumber
\end{equation}
where the constant $C$ is independent of $E$ and $\lambda$ when they vary over a compact set in $(0,\infty)$.
 In view of (\ref{S8Ffreelambdadef}),  (\ref{S2ftildelambdamudefine}) and (\ref{GpmlambdaEflimit}), we have
\begin{equation}
\begin{split}
F_{free,E}^{(\pm)}(\lambda)f & =  \Big(\frac{\sqrt{\lambda}}{\pi}\Big)^{1/2}\widetilde f^{(\pm)}(\lambda,E) \\
& =    \Big(\frac{\sqrt{\lambda}}{\pi}\Big)^{1/2}\frac{\rho(0)^{1-n}}{\Psi^{(\pm)}(\lambda,0,E)}\int_0^{\infty}\Psi_0(\lambda,s,E)f(s)\rho(s)^{n-1}ds.
\nonumber
\end{split}
\end{equation}
This implies that $F^{(\pm)}_{free,E}(\lambda)^{\ast} \in {\bf B}({\mathbb C};{\mathcal B}^{\ast})$ is the operator  of multiplication by the function
$$ \Big(\frac{\sqrt{\lambda}}{\pi}\Big)^{1/2}\overline{\left(\frac{\rho(0)^{1-n}}{\Psi^{(\pm)}(\lambda,0,E)}\right)}\Psi_0(\lambda,r,E).
$$
Since $\overline{\Psi^{(+)}(\lambda,r,E)} = \Psi^{(-)}(\lambda,r,E)$, there exist constants $c_{\pm}(\lambda,E)$ such that
$$
\Psi_0(\lambda,r,E) = c_+(\lambda,E)\Psi^{(+)}(\lambda,r,E) +
c_-(\lambda,E)\Psi^{(-)}(\lambda,r,E),
$$
$$
\overline{c_+(\lambda,E)} = c_-(\lambda,E).
$$
 Therefore, letting
\begin{equation}
\left\{
\begin{split}
&\omega_+(\lambda,E) =   \Big(\frac{\sqrt{\lambda}}{\pi}\Big)^{1/2}\overline{\left(\frac{\rho(0)^{1-n}}{\Psi^{(-)}(\lambda,0;E)}\right)}c_+(\lambda,E), \\
&\omega_-(\lambda,E) =   \Big(\frac{\sqrt{\lambda}}{\pi}\Big)^{1/2}\overline{\left(\frac{\rho(0)^{1-n}}{\Psi^{(-)}(\lambda,0;E)}\right)}c_-(\lambda,E),
\end{split}
\right.
\nonumber
\end{equation}
we have the following lemma. 

\begin{lemma}
 $F^{(-)}_{free,E}(\lambda)^{\ast}$ is the operator of multiplication by the function
\begin{equation}
\omega_+(\lambda,E)\Psi^{(+)}(\lambda,r,E) + 
\omega_-(\lambda,E)\Psi^{(-)}(\lambda,r,E).
\nonumber
\end{equation}
\end{lemma}


\subsection{Warped product metric}
We equip $(0,\infty)\times M$ with the metric
\begin{equation}
ds^2 = (dr)^2 + \rho(r)^2h_M(x,dx),
\nonumber
\end{equation}
and let
\begin{equation}
H_{free} = - \partial_r^2 - \frac{(n-1)\rho'}{\rho}\partial_r - \frac{\Lambda}{\rho^2},
\nonumber
\end{equation}
assuming the Dirichlet boundary condition at $r=0$, where $\Lambda = \Delta_M$ is the Laplace-Beltrami operator on $M$. We also let
\begin{equation}
R_{free}(z) = (H_{free} - z)^{-1}.
\nonumber
\end{equation}
 Let $0 = \lambda_0 \leq \lambda_1 \leq \cdots \to \infty$ be the eigenvalues of $- \Lambda$, $P_{\ell}$ the eigenprojection associated with $\lambda_{\ell}$,  and $\varphi_{\ell}(x)$ the associated normalized eigenvector. 
We then have
$$
R_{free}(\lambda \pm i0) = \sum_{\ell=0}^{\infty}P_{\ell} \otimes R_{free}(\lambda \pm i0,\lambda_{\ell}).
$$
Let ${\bf h}_{\infty} = L^2(M)$ and put 
\begin{equation}
\mathcal F_{free}^{(\pm)}(\lambda) = \sum_{\ell=0}^{\infty}P_{\ell}\otimes  F^{(\pm)}_{free,\lambda_{\ell}}(\lambda),
\nonumber
\end{equation}
Letting $E = \lambda_{\ell}$ in Lemma \ref{1dimFfreemuresolvent} and summing up with respect to $\ell$, we obtain the following lemma.


\begin{lemma}\label{infinitedimFfreemuresolvent}
For $f \in \mathcal B$ and $\lambda > 0$, we have
$$
\|\mathcal F_{free}^{(\pm)}(\lambda)f\|_{\bf h_{\infty}}^2 = \frac{1}{2\pi i}\left(\big[R_{free}(\lambda+i0)- R_{free}(\lambda-i0)\big]f,f\right).
$$
\end{lemma}
It then follows that
\begin{equation}
\|\mathcal F_{free}^{(\pm)}(\lambda)f\|^2_{L^2(M)} = 
\sum_{\ell=0}^{\infty}|F_{free,\lambda_{\ell}}^{(\pm)}(\lambda)f|^2 \leq C\|f\|_{\mathcal B}^2,
\nonumber
\end{equation}
where the constant $C$ is independent of $\lambda$ when $\lambda$ varies ove a compact set in $(0,\infty)$.
Take $\chi \in C^{\infty}({\mathbb R})$ such that $\chi(r) = 0$ for $r < 1$, and $\chi(r)=1$ for $r>2$, and put
\begin{equation}
c_{\ell}(\lambda,r) = \Big(\frac{\pi}{\sqrt{\lambda}}\Big)^{1/2}\chi\Big(\frac{r}{r_0(\lambda,\lambda_{\ell})}\Big),
\nonumber
\end{equation}
where $r_0(\lambda,E)$ is given in  (\ref{Definer0lambdaE}). 
Recall that $\varphi(\lambda,E,r)$ is defined by (\ref{Definevarphi(lambdaBr)}).


\begin{theorem}\label{RfeeAsympto}
For $f \in \mathcal B$, we have 
\begin{equation}
R_{free}(\lambda \pm i0)f  \simeq \sum_{\ell=0}^{\infty}
c_{\ell}(\lambda,r)\rho(r)^{-(n-1)/2}e^{\pm i \varphi(\lambda,\lambda_{\ell},r)}
P_{\ell}\otimes F^{(\pm)}_{free,\lambda_{\ell}}(\lambda)f.
\label{Thm6.4Rasympto}
\end{equation}
\end{theorem}

\begin{proof}
Since both sides in (\ref{Thm6.4Rasympto}) are bounded operators from $\mathcal B$ to $\mathcal B^{\ast}$, we have only to prove it for $f$ of the form $f = \sum_{\ell=0}^m \varphi_{\ell}(x)f_{\ell}(r)$, where $f_{\ell}(r) \in L^2_{comp}((0,\infty))$.
 Assume that $f_{\ell}(r) = 0$ for $r > a>0,\  0 \leq \ell \leq m$. Then by (\ref{S2GpmlambdamufPsitildef}), we have
$$
R_{free}(\lambda \pm i0)f = \sum_{\ell=0}^{m}\Psi^{(\pm)}(\lambda,r;\lambda_{\ell})\widetilde f^{(\pm)}_{\ell}(\lambda,\lambda_{\ell}), \quad r > a,
$$
from which (\ref{Thm6.4Rasympto}) follows immediately. 
\end{proof}

Note that if $\rho(r) = O(r^{\beta})$, $r_0(\lambda,\lambda_{\ell}) = O((\lambda_{\ell})^{1/2\beta})$. This shows the subtlety of the expansion (\ref{Thm6.4Rasympto}). 


\subsection{Perturbed metric}
We return to the perturbed metric
$$
ds^2 = (dr)^2 + \rho(r)^2h(r,x,dx),  
$$
$$
h(r,x,dx) - h_{M}(x,dx) \in S^{-\gamma}, \quad \gamma > 1.
$$
 Let $R_{free}(z) = (H_{free}-z)^{-1}$ and $R(z) = (H - z)^{-1}$.
Take $\chi \in C^{\infty}(0,\infty)$ such that $\chi(r) = 0$ for $r<1$, $\chi(r)=1$ for $r > 2$. Then, we have

\begin{equation}
\chi R(z) = R_{free}(z)V(z),
\label{chijRz=RfreejQz}
\end{equation}
\begin{equation}
V(z) = \chi + \Big([H_{free},\chi]-\chi\widetilde V\Big)R(z),
\quad 
\widetilde V = H - H_{free}.
\nonumber
\end{equation}
Here we use the assumption $\gamma > 1$ to have 
$$
V(\lambda \pm i0)  \in {\bf B}(\mathcal B^{\ast};\mathcal B),
$$
which implies that, by virue of (\ref{chijRz=RfreejQz}),  $R(\lambda \pm i0)$ has the same asymptotic expansion as in Theorem \ref{RfeeAsympto} where $f$ of the right-hand sides of (\ref{Thm6.4Rasympto})
is replaced by $V(\lambda \pm i0)f$. Therefore, the following theorem is proved.


\begin{theorem}\label{RAsympto}
Let $f \in \mathcal B$. Then, on the  regular end satisfying (\ref{RegCuspSRAssumption}) 
\begin{equation}
 R(\lambda \pm i0)f  \simeq \sum_{\ell=0}^{\infty}
c_{\ell}(\lambda,r)\rho(r)^{-(n-1)/2}e^{\pm i\varphi(\lambda,\lambda_{\ell},r)}
 P_{\ell}\otimes F^{(\pm)}_{free,\lambda_{\ell}}(\lambda)V(\lambda\pm i0)f.
\nonumber
\end{equation}
\end{theorem}


\section{Spectral representation}
\label{SectionSpectralrepresentation}
We return to our original manifold $\mathcal M = \mathcal K\cup\mathcal M_1\cup\cdots\cup\mathcal M_{N+N'}$. From here until the end of \S \ref{SectionPhysicalSmatrix}, we assume 
\begin{equation}
\left\{
\begin{split}
& \alpha_{0,j} > 0, \quad \beta_{0,j} > 0, \quad \gamma_{0,j} > 1 \quad 
on \ regular\  ends,\\
& (A-3) \quad  on \ cusp \ ends. 
\end{split}
\right.
\nonumber
\end{equation}
In this section, 
we construct a spectral representation for $H$ by observing the asymptotic expansion of the resolvent at infinity.


\subsection{Definition of $\mathcal F_j^{(\pm)}(\lambda)$}
Let ${\bf h}_j = L^2(M_j)$ and
\begin{equation}
{\bf h}_{\infty,j} = \left\{
\begin{split}
&{\bf h}_j \quad {\rm for} \quad 1 \leq j \leq N, \\
&{\mathbb C} \quad {\rm for} \quad N+1\leq j\leq N+N',
\end{split}
\right.
\nonumber
\end{equation}
\begin{equation}
	{\bf h}_{\infty} = {\mathop\oplus_{j=1}^{N+N'}}{\bf h}_{\infty,j},
\nonumber
\end{equation}
\begin{equation}
\widehat{\mathcal H} = {\mathop\oplus_{j=1}^{N+N'}}L^2((E_{0,j},\infty);{\bf h}_{\infty,j};d\lambda).
\label{S11DefinewidehatH}
\end{equation}
Let $\lambda_{\ell,j}$ and $P_{\ell,j}$ $(\ell = 0, 1, 2, \cdots)$ be the eigenvalues and the associated eigenprojections for $- \Delta_{M_j}$ with respect to the metric $h_{M_j}(x,dx)$.

For $f \in \mathcal B(\mathcal M)$ and $\lambda \in \sigma_e(H)\setminus\mathcal E$, we put
\begin{equation}
\mathcal F^{(\pm)}(\lambda)f = \big(\mathcal F^{(\pm)}_1(\lambda)f,\cdots,\mathcal F_{N+N'}^{(\pm)}(\lambda)f\big) \in {\bf h}_{\infty},
\nonumber
\end{equation}
where $\mathcal F^{(\pm)}_j(\lambda)$ is defined as follows.

\medskip
\noindent
(I) {\it Regular ends with $\beta_{0,j} > 1/2$} : For $1 \leq j \leq N$ and $\lambda > E_{0,j}$, we define it making use of the asymptotic expansion  (\ref{EndResolvAsym}) in
Theorem \ref{S10ResovAsympReg}:
\begin{equation}
R(\lambda \pm i0)f \simeq C_j(\lambda)\rho_j(r)^{-(n-1)/2}e^{\pm i\Phi_j(r,\lambda)}\mathcal F^{(\pm)}_j(\lambda)f \quad 
{\rm on} \quad \mathcal M_j,
\label{S12Fjpmregulardefine}
\end{equation}
where
\begin{equation}
C_j(\lambda) = \Big(\frac{\pi}{\sqrt{\lambda - E_{0,j}}}\Big)^{1/2},
\nonumber
\end{equation}
\begin{equation}
\Phi_j(r,\lambda) = \int_{0}^r\phi_j(t,\lambda)\, dt,
\nonumber
\end{equation}
\begin{equation}
\phi_j(r,\lambda) = \sqrt{\lambda - \frac{(n-1)^2}{4}\Big(\frac{\rho'_j(r)}{\rho_j(r)}\Big)^2}.
\nonumber
\end{equation}

\medskip
\noindent
(II) {\it Regular ends with $0 < \beta_{0,j} \leq 1/2$} :  For $1 \leq j \leq N$ and $\lambda > E_{0,j}$, we define
\begin{equation}
\mathcal F_j^{(\pm)}(\lambda) = \sum_{\ell=0}^{\infty}
\mathcal F_{\ell,j}^{(\pm)}(\lambda),
\nonumber
\end{equation}
where 
\begin{equation}
\mathcal F^{(\pm)}_{\ell,j}(\lambda) = P_{\ell,j}\otimes  F^{(\pm)}_{free,\lambda_{\ell,j}}(\lambda)V_j(\lambda \pm i0)
\nonumber
\end{equation}
 appearing in  Theorem \ref{RAsympto}:
\begin{equation}
R(\lambda \pm i0)f \simeq \sum_{\ell=0}^{\infty}
c_{\ell,j}(\lambda,r)\rho_j(r)^{-(n-1)/2}e^{\pm i\varphi_{j}(\lambda,\lambda_{\ell,j},r)}\mathcal F^{(\pm)}_{\ell,j}(\lambda)f, \quad {\rm on} \quad 
\mathcal M_j,
\label{DefineFjlSR}
\end{equation} 
\begin{equation}
c_{\ell,j}(\lambda,r) = \Big(\frac{\pi}{\sqrt{\lambda - E_{0,j}}}\Big)^{1/2}\chi\Big(\frac{r}{r_0(\lambda,\lambda_{\ell,j})}\Big),
\nonumber
\end{equation}
\begin{equation}
\varphi_{j}(\lambda,E,r) = \int_{r_0(\lambda,E)}^r\alpha_{j}(\lambda,E,s)ds,
\nonumber
\end{equation}
\begin{equation}
\alpha_{j}(\lambda,E,r) = \sqrt{\lambda - \Big(\frac{(n-1)\rho_j'}{2\rho_j}\Big)^2-\frac{E}{\rho_j^2}}.
\nonumber
\end{equation}

\medskip
\noindent
(III) {\it Cusp ends} : For $N+1 \leq j \leq N+N'$ and $\lambda > E_{0,j}$, we  use (\ref{CuspEndResolAsym}) in Theorem \ref{S10ResovAsympCusp}:
\begin{equation}
R(\lambda \pm i0)f \simeq C_j(\lambda)P_{0,j}\otimes \rho_j(r)^{-(n-1)/2}e^{\pm i\Phi_j(r,\lambda)}\mathcal F^{(\pm)}_j(\lambda)f \quad 
{\rm on} \quad \mathcal M_j.
\label{S12Fjpmregulardefinecusp}
\end{equation}

\medskip
Finally, we define
\begin{equation}
\mathcal F_j^{(\pm)}(\lambda) = 0, \quad {\rm if} \quad \lambda < E_{0,j}, \quad 1 \leq j \leq N + N'.
\nonumber
\end{equation}



%


\begin{lemma}\label{ParsevalDiff}
For $f, g \in \mathcal B$ and $\lambda \in \sigma_{e}(H)\setminus\mathcal E$
$$
\frac{1}{2\pi i}\left([R(\lambda + i0) - R(\lambda - i0)]f,g\right) 
= \big(\mathcal F^{(\pm)}(\lambda)f,\mathcal F^{(\pm)}(\lambda)g\big)_{\bf h_{\infty}}
$$
\end{lemma}

\begin{proof}
Take $\varphi \in C_0^{\infty}({\mathbb R})$ such that $\varphi(r) = 0$ for $r< 1/2$ and $r > 3$, and $\varphi(r) = 1$ for $1<r<2$. Put
$$
\widetilde \varphi(r) = \int_r^{\infty}\varphi(t)dt.
$$
Take $\varphi_R \in C_0^{\infty}(\mathcal M)$ such that $\varphi_R = \widetilde\varphi(0)$ on $\mathcal K$, and on each end $\mathcal M_j$
$$
\varphi_R(r) = \widetilde\varphi(r/R).
$$
Then, we have 
$\varphi_R' = - \varphi\big(r/R\big)/R$ on each end, hence
\begin{equation}
\Big[-\frac{\partial^2}{\partial r^2} - \frac{(n-1)\rho_j'}{\rho_j}\frac{\partial}{\partial r}, \varphi_R\Big] = \frac{2}{R}\varphi\big(\frac{r}{R}\big)
\Big(\frac{\partial}{\partial r} + \frac{(n-1)\rho_j'}{2\rho_j}\Big)+ \frac{1}{R^2}\varphi''\big(\frac{r}{R}\big).
\label{S12[HchiR]}
\end{equation}
 Let $u = R(\lambda + i0)f$, $v = R(\lambda + i0)g$. Then, we have
 \begin{equation}
 (\varphi_Ru,g) - (f,\varphi_Rv) = ([H,\varphi_R]u,v).
\label{ug-fvinnerproduct}
 \end{equation}
 As $R \to \infty$, the left-hand side tends to $\widetilde \varphi(0)\big([R(\lambda+i0)-R(\lambda-i0)]f,g\big)$. By (\ref{S7DjkpmDefine}), we have on each end, 
$$
\frac{\partial}{\partial r} + \frac{(n-1)\rho_j'}{2\rho_j} = D_{j}^{(+)}(k) 
+ ik + O(r^{-\epsilon}), \quad \epsilon>0.
$$
This, together with (\ref{S12[HchiR]}), implies that the right-hand side of (\ref{ug-fvinnerproduct}) is asymptotically equal to  
$$
\frac{2ik}{R}\big(\varphi\big(\frac{r}{R}\big)u,v\big). 
$$
By (\ref{S12Fjpmregulardefine}), (\ref{DefineFjlSR}) and (\ref{S12Fjpmregulardefinecusp}), this is asymptotically equal to
\begin{equation}
2ik\,\widetilde \varphi(0)\frac{\pi}{\sqrt{\lambda - E_{0,j}}}\big(\mathcal F_{j}^{(+)}(\lambda)f,\mathcal F_{j}^{(+)}(\lambda)g\big)_{{\bf h}_j}, 
\label{S12Parseval1}
\end{equation}
for the case of regular end with $\beta_j > 1/2$, 
\begin{equation}
2ik\,\widetilde \varphi(0)\sum_{\ell=0}^{\infty}\frac{\pi}{\sqrt{\lambda - E_{0,j}}}\big(\mathcal F_{\ell,j}^{(+)}(\lambda)f,\mathcal F_{\ell,j}^{(+)}(\lambda)g\big)_{{\bf h}_j},
\label{S12Parseval2}
\end{equation}
for the case of regular end with $0 < \beta_{j} \leq 1/2$, and
\begin{equation}
2ik\,\widetilde \varphi(0) \frac{\pi}{\sqrt{\lambda}}\,\mathcal F_{0,j}^{(+)}(\lambda)f\,\overline{ \mathcal F_{0,j}^{(+)}(\lambda)g},
\label{S12Parseval3}
\end{equation}
for the case of cusp end. 

Let us prove (\ref{S12Parseval2}). It is easy to show that if $w_1 \in \mathcal B^{\ast}_0$ and $w_2 \in \mathcal B$, then $\displaystyle{\frac{1}{R}\big(\varphi(\frac{r}{R})w_1,w_2) \to 0.}$
By (\ref{DefineFjlSR}), $\displaystyle{\frac{1}{R}\big(\varphi(\frac{r}{R})u,v)}$ is asymptotically equal to 
\begin{equation}
\begin{split}
& \frac{1}{R}\sum_{\ell}\big(\chi(\frac{r}{R})c_{\ell,j},c_{\ell,j}\big)a_{\ell,j} \\
& = \frac{\pi}{\sqrt{\lambda - E_{0,j}}}
\sum_{\ell}\frac{1}{R}\int_0^{\infty}\varphi\big(\frac{r}{R}\big)
\chi\big(\frac{r}{r_0(\lambda,\lambda_{\ell,j})}\big)^2dr \, a_{\ell,j},
\end{split}
\nonumber
\end{equation}
where $a_{\ell,j} = (\mathcal F^{(+)}_{\ell,j}(\lambda)f,\mathcal F^{(+)}_{\ell,j}(\lambda)g)_{{\bf h}_j}$.  Let $ r = tR$. Since $\sum_{\ell}|a_{\ell,j}| < \infty$ and $\chi\big(\frac{tR}{r_0(\lambda,\lambda_{\ell,j})}\big) \to 1$, the right-hand side converges to 
$$
\frac{\pi}{\sqrt{\lambda - E_{0,j}}}\sum_{\ell}\int_0^{\infty}\varphi(t)dt\, a_{\ell,j} = 
\frac{\pi}{\sqrt{\lambda - E_{0,j}}}\widetilde \varphi(0)\sum_{\ell} a_{\ell,j} .
$$
This proves (\ref{S12Parseval2}). The proof of (\ref{S12Parseval1}) and (\ref{S12Parseval3}) is similar and simpler.

Summing up (\ref{S12Parseval1}), (\ref{S12Parseval2}), (\ref{S12Parseval3}) with respect to $j$ and dividing by $\widetilde \varphi(0)$, we obtain the lemma for the case of $\mathcal F^{(+)}(\lambda)$. 
\end{proof}

As a corollary, we have
\begin{equation}
\|\mathcal F^{(\pm)}(\lambda)f\|_{{\bf h}_{\infty}}
\leq C\|f\|_{\mathcal B}, 
\label{S12F(lambda)bounded}
\end{equation}
where the constant $C$ does not depend on $\lambda$ when $\lambda$ varies over a compact set in $\sigma_e(H)\setminus \mathcal E$.


\subsection{Generalized Fourier transform}

We put 
$$
(\mathcal F^{(\pm)}f)(\lambda) = \mathcal F^{(\pm)}(\lambda)f
$$
 for $f \in \mathcal B$.  In view of (\ref{S8EHI=Rff}) and  Lemma \ref{ParsevalDiff}, we have, for $f, g \in \mathcal B$, 
\begin{equation}
(P_{ac}(H)f, g) = \int_{E_{0,tot}}^{\infty}(\mathcal F^{(\pm}(\lambda)f,\mathcal F^{(\pm)}(\lambda)g)_{{\bf h}_{\infty}}d\lambda = (\mathcal F^{(\pm)}f,\mathcal F^{(\pm)}g)_{\widehat{\mathcal H}}.
\nonumber
\end{equation}
Therefore, $\mathcal F^{(\pm)}$ is uniquely extended to a partial isometry on $L^2(\mathcal M)$ with initial set $\mathcal H_{ac}(H)$ and final set in $\widehat{\mathcal H}$, defined in (\ref{S11DefinewidehatH}), which is denoted by $\mathcal F^{(\pm)}$ again. We show the following lemma.


\begin{lemma}
\label{Lemma11.2}
(1) For any $f \in D(H)$ and  a.e. $\lambda \in (E_{0,tot},\infty)$, we have 
$$
(\mathcal F^{(\pm)}Hf)(\lambda) = \lambda(\mathcal F^{(\pm)}f)(\lambda).
$$
(2) For any bounded Borel function $\alpha(\lambda)$ on ${\mathbb R}$, any $f \in L^2(\mathcal M)$ and a.e. $\lambda \in (E_{0,tot},\infty)$, we have
\begin{equation}
\big({\mathcal F}^{(\pm)}\alpha(H)f\big)(\lambda) = \alpha(\lambda)\big({\mathcal F}^{(\pm)}f\big)(\lambda).
\label{S12FalphaH=alphalambdaF}
\end{equation}
\end{lemma}

\begin{proof}
For $f \in L^2_{comp}(\mathcal M) \cap D(H)$, let $u = R(\lambda +i0)f$ and $v = Hu$. Then, we have
$(H-\lambda)u = f, \quad v = R(\lambda+i0)Hf$. 
Observing the spatial asymptotics of $v = Hu = \lambda u + f$, we have $\mathcal F^{(+)}(\lambda)Hf = \lambda\mathcal F^{(+)}(\lambda)f$, which proves (1). It then follows that
$$
(\lambda - z)\big(\mathcal F^{(\pm)}(H-z)^{-1}f\big)(\lambda) = 
\big(\mathcal F^{(\pm)}f\big)(\lambda)
$$
for $z \not\in {\mathbb R}$, which shows the assertion (2) for $\alpha(\lambda) = (\lambda - z)^{-1}$. Then, by Stone's formula, 
(2) holds for any step function, hence for any bounded Borel function. 
\end{proof}


\begin{lemma}\label{Fpmisonto}
$\ {\rm Ran}\; \mathcal F^{(\pm)} = \widehat{\mathcal H}$.
\end{lemma}

\begin{proof}
 We have only to show that the range of $\mathcal F^{(\pm)}$ is dense in $\widehat{\mathcal H}$. For the sake of notational simplicity, we consider the case that $N=N'=1$, and assume that the volume of $M_2$ is equal to 1. Suppose 
$$
(\phi_1(\lambda),\phi_2(\lambda)) \in L^2((E_{0,1},\infty);L^2(M_1);d\lambda)\times L^2((E_{0,2},\infty);{\mathbb C};d\lambda)
$$
is orthogonal to ${\rm Ran}\,\mathcal F^{(+)}$.  Let $\lambda_{\ell,1}$, $\ell = 0, 1, 2, \cdots$, be the eigenvalues of $- \Delta_{M_1}$, and
 $e_{\ell,1}$ the associated complete orthonormal system of eigenvectors in $L^2(M_1)$. We put
\begin{equation}
\phi_{\ell,1}(\lambda) = (\phi_1(\lambda),e_{\ell,1})_{L^2(M_1)}.
\label{S12phiell1lambda}
\end{equation}
For $\psi \in L^1_{loc}((E_{0,tot},\infty))$, let $\mathcal L(\psi)$ be the set of Lebesgue points of $\psi$, i.e.
$$
\mathcal L(\psi) \ni \lambda \Longleftrightarrow 
\psi(\lambda) = \lim_{\epsilon\to0}\frac{1}{2\epsilon}\int_{\lambda-\epsilon}^{\lambda+\epsilon}
\psi(t)dt.
$$
It is well-known that $(E_{0,tot},\infty)\setminus\mathcal L(\psi)$ is a null set. Take an arbitrary point $\mu \in (E_{0,tot},\infty)$ satisfying
$$
\mu \in \Big({\mathop\cap_{\ell=0}^{\infty}}\mathcal L(\phi_{\ell,1})\Big)\cap\Big(\mathcal L(\|\phi_1\|^2_{L^2(M_1)})\Big)\cap \Big(\mathcal L(\phi_2)\Big)\cap
\Big(\mathcal L(|\phi_2|^2)\Big)\cap
\big((E_{0,tot},\infty)\setminus \mathcal E\big).
$$
Let $\{\chi_j\}_{j=0}^2$ be the partition of unity in (\ref{S8Partitionunity}). We fix $m$ arbitrarily, and put
$$
u_{\mu}(r) = \chi_1(r)\Psi^{(+)}_1(\mu,r,\lambda_{m,1})\alpha e_{m,1}(x)
+ \chi_2(r)\Psi^{(+)}_2(\mu,r,0)\beta,
$$
where $\Psi^{(+)}_j(\lambda,r,E)$ is the solution constructed in Lemma \ref{S2Lemma2} for the end $\mathcal M_j$, and $\alpha, \beta$ are arbitrary constants. We  put
$$
(H - \mu)u_{\mu}= g_{\mu}.
$$
Then, by virtue of Lemma \ref{S2Lemma1},  $g_{\mu} \in L^{2,1+\epsilon}(\mathcal M)$. Since $u_{\mu}$ is outging, we have $u_{\mu} = R(\mu+i0)g_{\mu}$. Moreover, letting $\mathcal F^{(+)}(\lambda)g_{\mu} = (C_1(\lambda),C_2(\lambda))$ and observing the behavior of $u_{\mu}$ at infinity, we see that $(C_1(\lambda),C_2(\lambda))$ is an $L^2(M_1)\times {\mathbb C}\, $-valued continuous function of $\lambda >0$ satisfying
\begin{equation}
(C_1(\mu),e_{\ell,1}) = \Big(\frac{\pi}{\sqrt{\mu-E_{0,1}}}\Big)^{-1/2}\delta_{\ell,m} \alpha, \quad 
C_2(\mu) = \Big(\frac{\pi}{\sqrt{\mu- E_{0,2}}}\Big)^{-1/2}\beta.
\label{S12Ctteminnerprod}
\end{equation}
By the assumption, $(\phi_1(\lambda),\phi_2(\lambda))$ is orthogonal to $\mathcal F^{(+)}E_H(I)g_{\mu}$, $I$ being any interval in $(E_{0,tot},\infty)\setminus\mathcal E$. Hence by Lemma \ref{Lemma11.2} (2)
$$
\int_I\left((\phi_1(\lambda),C_1(\lambda))_{L^2(M_1)} + \phi_2(\lambda)\overline{C_2(\lambda)}\right)d\lambda = 0
$$
for any interval $I \subset (E_{0,tot},\infty)\setminus\mathcal E$. Since $C_2(\lambda)$ is continuous, and $\mu$ is a Lebesgue point of $\phi_2(\lambda)$, we have
$$
\frac{1}{2\epsilon}\int_{\mu-\epsilon}^{\mu+\epsilon}\phi_2(\lambda)\overline{C_2(\lambda)}d\lambda \to \phi_2(\mu)\Big(\frac{\pi}{\sqrt{\mu- E_{0,2}}}\Big)^{-1/2}\overline{\beta}.
$$
The 1st term is computed as
\begin{equation}
\begin{split}
\frac{1}{2\epsilon}\int_{\mu-\epsilon}^{\mu+\epsilon}(\phi_1(\lambda),C_1(\lambda))_{L^2(M_1)}d\lambda &= 
\frac{1}{2\epsilon}\int_{\mu-\epsilon}^{\mu+\epsilon}(\phi_1(\lambda),C_1(\lambda)- C_1(\mu))_{L^2(M_1)}d\lambda \\
&+ \frac{1}{2\epsilon}\int_{\mu-\epsilon}^{\mu+\epsilon}(\phi_1(\lambda),C_1(\mu))_{L^2(M_1)}d\lambda.
\end{split}
\nonumber
\end{equation}
By (\ref{S12phiell1lambda}) and (\ref{S12Ctteminnerprod}), $(\phi_1(\lambda),C_1(\mu))_{L^2(M_1)} = \big(\frac{\pi}{\sqrt{\mu - E_{0,1}}}\big)^{-1/2} \phi_{m,1}(\lambda)\overline{\alpha}$, hence
$$
\frac{1}{2\epsilon}\int_{\mu-\epsilon}^{\mu+\epsilon}(\phi_1(\lambda),C_1(\mu))_{L^2(M_1)}d\lambda \to 
\big(\frac{\pi}{\sqrt{\mu-E_{0,1}}}\big)^{-1/2}\phi_{m,1}(\mu)\overline{\alpha}.
$$
We also have
\begin{equation}
\begin{split}
& \left|\frac{1}{2\epsilon}\int_{\mu-\epsilon}^{\mu+\epsilon}(\phi_1(\lambda),C_1(\lambda)- C_1(\mu))_{L^2(M_1)}d\lambda\right| \\
& \leq \left(\frac{1}{2\epsilon}\int_{\mu-\epsilon}^{\mu+\epsilon}\|\phi_1(\lambda)\|^2_{L^2(M_1)}d\lambda
\right)^{1/2} \times
\left(\frac{1}{2\epsilon}\int_{\mu-\epsilon}^{\mu+\epsilon}\|C_1(\lambda)- C_1(\mu)\|^2_{L^2(M_1)}d\lambda
\right)^{1/2}.
\end{split}
\nonumber
\end{equation}
The right-hand side tends to 0, since $\mu$ is a Lebesgue point of $\|\phi_1(\lambda)\|^2_{L^2(M_1)}$, and $C_1(\lambda)$ is an $L^2(M_1)$-valued continuous function of $\lambda >0$. Therefore, we have 
$$
\phi_{m,1}(\mu)\big(\frac{\pi}{\sqrt{\mu - E_{0,1}}}\big)^{-1/2}\overline{\alpha} + \phi_2(\mu)\Big(\frac{\pi}{\sqrt{\mu- E_{0,2}}}\Big)^{-1/2}\overline{\beta}=0.
$$
Since $\alpha, \beta$ can be chosen arbitrarily, we have proven that $\phi_1(\mu)=0$, $\phi_2(\mu)=0$.
\end{proof}

Now, we have arrived at the main theorem.


\begin{theorem}\label{GeneralizeFourierTh}
(1) The operator $\mathcal F^{(\pm)}$ is uniquely extended to a partial isometry with initial set $\mathcal H_{ac}(H)$ and final set $\widehat{\mathcal H}$.\\
\noindent
(2) For any $\lambda \in (E_{0,tot},\infty)\setminus \mathcal E$, $\mathcal F^{(\pm)}(\lambda)^{\ast} \in {\bf B}({\bf h}_{\infty};\mathcal B^{\ast})$. Moreover, for any $a \in {\bf h}_{\infty}$,
 $$
(- \Delta_{\mathcal M}-\lambda)\mathcal F^{(\pm)}(\lambda)^{\ast}a = 0.
$$
(3) Fore any $f \in L^2(\mathcal M)$, and a compact interval $I \subset (E_{0,j},\infty)\setminus\mathcal E$, the $\mathcal B^{\ast}$-valued integral
$\displaystyle{
\int_{I}
\mathcal F_j^{(\pm)}(\lambda)^{\ast}\Big(\mathcal  F_j^{(\pm)}f\Big)(\lambda)d\lambda}$
belongs to $L^2(\mathcal M)$. Letting $I_k \to (E_{0,j},\infty)\setminus\mathcal E$, the strong limit
$$
\lim_{k\to\infty}\int_{I_k}
\mathcal F_j^{(\pm)}(\lambda)^{\ast}\Big(\mathcal  F_j^{(\pm)}f\Big)(\lambda)d\lambda = \int_{E_{0,j}}^{\infty}
\mathcal F_j^{(\pm)}(\lambda)^{\ast}\Big(\mathcal  F_j^{(\pm)}f\Big)(\lambda)d\lambda
$$
exists in $L^2(\mathcal M)$.\\
\noindent
(4) For any $f \in \mathcal H_{ac}(H)$, the inversion formula holds :
$$
f = \Big(\mathcal F^{(\pm)}\Big)^{\ast}\mathcal F^{(\pm)}f = 
\sum_{j=1}^{N+N'}\int_{E_{0,j}}^{\infty}
\mathcal F_j^{(\pm)}(\lambda)^{\ast}\Big(\mathcal  F_j^{(\pm)}f\Big)(\lambda)d\lambda.
$$
\end{theorem}

\begin{proof}
 The assertion (1) is already proven in Lemma \ref{Fpmisonto}. 
By (\ref{S12F(lambda)bounded}), $\mathcal F^{(\pm)}(\lambda)^{\ast} \in {\bf B}({\bf h}_{\infty}\, ; \, \mathcal B^{\ast}_0)$ \footnote{One needs to be careful about the definition of $\mathcal F^{(\pm)}(\lambda)^{\ast}$. We discuss it in the next section.}. Taking the adjoint in Lemma \ref{Lemma11.2} (1), we obtain (2).  
We prove (3) and (4) at the same time. Take any compact inetrval $I \subset (E_0,\infty)\setminus\mathcal E$, and put
$$
u_I = \int_I\mathcal F^{(\pm)}(\lambda)^{\ast}\big(\mathcal F^{(\pm)}f\big)(\lambda)f\,d\lambda,
$$
which belongs to $\mathcal B^{\ast}$ by (2). Letting $c_I(\lambda)$ be the characteristic function of $I$, we have for any $g \in \mathcal B$,
$$
(u_I,g) = \int c_I(\lambda)(\mathcal F^{(\pm)}f)(\lambda),\mathcal F^{(\pm)}(\lambda)g)d\lambda 
= (E_H(I)f,g).
$$
Therefore $u_I = E_H(I)f \in L^2(\mathcal M)$. This implies that for any simple function $\alpha(\lambda)$
$$
\int \alpha(\lambda) \mathcal F^{(\pm)}(\lambda)^{\ast}\big(\mathcal F^{(\pm)}f\big)(\lambda)f\,d\lambda = \alpha(H)f.
$$
Finally, we approximate $(E_{0,tot},\infty)$ by a union of compact intervals in $(E_{0,tot},\infty)\setminus \mathcal E$ to complete the proof.
\end{proof}


\section{Helmholtz equation and S-matrix}\label{SectionPhysicalSmatrix}


\subsection{Eigenoperator} 
The adjoint of $\mathcal F^{(\pm)}(\lambda)$ is an eigenoperator of $H$ in the sense that $(- \Delta_{\mathcal M} - \lambda)\mathcal F^{(\pm)}(\lambda)^{\ast}a=0$ holds for $a \in {\bf h}_{\infty}$. We will derive the S-matrix by observing its asymptotic behavior at infinity. 

Letting $c_j(\lambda)$ be the characteristic function of the interval $(E_{0,j},\infty)$, we put
\begin{equation}
{\bf h}_{\infty}(\lambda) = {\mathop\oplus_{j=1}^{N+N'}}c_j(\lambda){\bf h}_{\infty,j}.
\nonumber
\end{equation}
 For $a = (a_1,\cdots,a_{N+N'}) \in {\bf h}_{\infty}(\lambda)$, we have
$$
\mathcal F^{(\pm)}(\lambda)^{\ast}a = \sum_{j=1}^{N+N'}\mathcal F^{(\pm)}_j(\lambda)^{\ast}a_j.
$$

Let $\{e_{\ell,j}(x)\}_{\ell=0}^{\infty}$ be a complete orthonormal system of eigenvectors of $- \Delta_{M_j}$ associated with eigenvalues $\{\lambda_{\ell,j}\}_{\ell=0}^{\infty}$. In particular, $e_{0,j}(x) = ({\rm vol}\,(M_j))^{-1/2}$.  Let 
$D_{finite}(\Lambda_j)$  be the set of $c \in L^2(M_j)$ such that $(c,e_{\ell,j}) = 0$ except for a finite numer of $\ell$. 
We define a subset ${\bf h}_{\infty,j}^{comp} \subset {\bf h}_{\infty,j}$ by
\begin{equation}
{\bf h}^{comp}_{\infty,j} = 
\left\{
\begin{split}
& D_{finite}(\Lambda_j) \quad {\rm for} \quad 1 \leq j \leq N, \\
& {\mathbb C} \quad {\rm for} \quad N+1 \leq j \leq N + N'.
\end{split}
\right.
\nonumber
\end{equation}
and ${\bf h}_{\infty}^{comp}(\lambda)$ by
\begin{equation}
{\bf h}_{\infty}^{comp}(\lambda) = {\mathop\oplus_{j=1}^{N+N'}}c_j(\lambda){\bf h}_{\infty,j}^{comp}.
\nonumber
\end{equation}

Taking $a = (a_1,\cdots,a_{N+N'}) \in {\bf h}_{\infty}^{comp}(\lambda)$, and using the partition of unity $\{\chi_j\}_{j=0}^{N+N'}$, we put $u_j^{(\pm)}$ as follows. 

\medskip
\noindent
(I) For regular ends with $\beta_{0,j} > 1/2$: 
\begin{equation}
u_j^{(\pm)} =  C_j(\lambda)\chi_j\rho_j(r)^{-(n-1)/2}e^{\pm i\Phi_j(r,\lambda)}a_j, \quad 1 \leq j \leq N, 
\nonumber
\end{equation}

\medskip
\noindent
(II) For regular ends with $0 < \beta_{0,j} \leq 1/2$:
\begin{equation}
u_j^{(\pm)} = \sum_{\ell=0}^{\infty} c_{\ell,j}(\lambda,r)\chi_j\rho_j(r)^{-(n-1)/2}e^{\pm i\varphi_j(\lambda,\lambda_{\ell,j},r)}a_{\ell,j}e_{\ell,j}(x), \quad 
1 \leq j \leq N, 
\nonumber
\end{equation} 
where $a_j = \sum_{\ell=0}^{\infty}a_{\ell,j}e_{\ell,j}(x)$, 

\medskip
\noindent
(III) For cusp ends: 
\begin{equation}
u_j^{(\pm)} = 
C_j(\lambda)\chi_j\rho_j(r)^{-(n-1)/2}e^{\pm i\Phi_j(r,\lambda)}a_{0,j}e_{0,j}(x), \quad N+1 \leq j \leq N+N'.
\nonumber
\end{equation}

\medskip
We put 
\begin{equation}
f^{(\pm)}_j = (- \Delta_{\mathcal M} - \lambda)u^{(\pm)}_j, \quad f^{(\pm)} = \sum_{j=1}^{N+N'}f_j^{(\pm)}, \quad u^{(\pm)} = \sum_{j=1}^{N+N}u_j^{(\pm)}.
\label{S13Definefpmetc}
\end{equation}


\begin{lemma}\label{FlambdaastExpand}
Let $a, u^{(\pm)}, f^{(\pm)}$ be as above. Then we have
$$
\mp 2\pi i\mathcal F^{(\pm)}(\lambda)^{\ast}a = u^{(\pm)} - R(\lambda \mp i0)f^{(\pm)}.
$$
\end{lemma}

\begin{proof}
We prove the $(+)$ case. Let $v = R(\lambda + i0)h$ for $h \in \mathcal B$. Take $\chi(r) \in C^{\infty}((0,\infty))$ such that $\chi(r)=1$ for $r < 1$, $\chi(r)=0$ for $r > 2$, and put $\chi_t(r) = \chi(r/t)$.    Then  by integration by parts, we have 
\begin{equation}
\begin{split}
 & (\chi_t u_j^{(+)},h) - (\chi_tf_j^{(+)},v) \\
 & =   \frac{1}{t}\int \chi' \Big(\frac{r}{t}\Big) \Big(u_j^{(+)}\overline{D_j^{(+)}(k)v} - (D_j^{(+)}(k)u_j^{(+)})\overline{v}\Big)\sqrt{g_j}\, drdx \\
& + 2i\,\, 
\frac{1}{t} \int \chi' \Big(\frac{r}{t}\Big)\big({\rm Re}\, \psi_j^{(+)}\big)u_j^{(+)}\overline{v}\sqrt{g_j}\, drdx,
\end{split}
\label{S13chirtint}
\end{equation}
where $\psi_j^{(+)} = \sqrt{\lambda - E_{0,j}} + O(r^{-\epsilon})$, $\epsilon > 0$.
The first  term of the right-hand side of (\ref{S13chirtint}) vanishes as $t \to \infty$. To compute the 2nd term, we use the asymptotic expansions (\ref{S12Fjpmregulardefine}), (\ref{DefineFjlSR}),  (\ref{S12Fjpmregulardefinecusp}). Since $\frac{1}{t}\int\chi'(r/t)dr = -1$, it tends to
$$
- 2\pi i(a_j,\mathcal F^{(+)}_j(\lambda)h)_{{\bf h}_{\infty,j}}.
$$
We thus obtain
$$
(u_j^{(+)},h) - (R(\lambda -i0)f_j^{(+)},h) = - 2\pi i
(\mathcal F_j^{(+)}(\lambda)^{\ast}a_j,h).
$$
Summing up with respect to $j$, we obtain the lemma. 
\end{proof}

In view of Lemma \ref{FlambdaastExpand} and the asymptotic expansion of the resolvent (Theorems \ref{S10ResovAsympReg}, \ref{S10ResovAsympCusp}, \ref{RAsympto}), we have the following lemma.


\begin{lemma}\label{Flambdaastexmandend}
For any  $a \in {\bf h}_{\infty}^{comp}(\lambda)$, there exists $b^{(\pm)} \in {\bf h}_{\infty}(\lambda)$ such that, 
letting $U^{(\pm)} \in {\mathcal B}^{\ast}$ be defined by (1), (2), (3) below, $2\pi i\mathcal F^{(\pm)}(\lambda)^{\ast}a$ satisfies
\begin{equation}
\mp 2\pi i\mathcal F^{(\pm)}(\lambda)^{\ast}a - U^{(\pm)}\in \mathcal B^{\ast}_0,
\label{Flambda-UpminB0ast}
\end{equation}
where

\noindent
(1) for regular ends with $\beta_{0,j} > 1/2$,

\begin{equation}
U^{(\pm)} = \Big(\frac{\pi}{\sqrt{\lambda - E_{0,j}}}\Big)^{1/2}\rho_j(r)^{-(n-1)/2}\left(e^{\pm i\Phi_j(r,\lambda)}a_j(x) - e^{\mp i\Phi_j(r,\lambda)}b_j^{(\pm)}(x)\right),
\nonumber
\end{equation}
(2) for regular ends with $0 < \beta_{0,j} \leq 1/2$, 
\begin{equation}
\begin{split}
 U^{(\pm)} = & \Big(\frac{\pi}{\sqrt{\lambda - E_{0,j}}}\Big)^{1/2}\rho_j(r)^{-(n-1)/2}\\
& \times \sum_{\ell=0}^{\infty}
\chi_{\ell,j}(r)\Big(e^{\pm i\varphi_j(\lambda,\lambda_{\ell,j},r)}a_{j,\ell} 
 - e^{\mp i\varphi_j(\lambda,\lambda_{\ell,j},r)}b^{(\pm)}_{j,\ell}\Big)e_{j,\ell}(x),
 \end{split}
\nonumber
\end{equation} 
where $\chi_{\ell,j}(r) = \chi\big(\frac{r}{r_0(\lambda,\lambda_{\ell,j})}\big)$ with $\chi(r)$ satisfying $\chi \in C^{\infty}({\mathbb R})$, $\chi(r)=0$ for $r < 1$, $\chi(r)=1$ for $r > 2$,

\noindent
(3) for cusp ends
\begin{equation}
U^{(\pm)} = \Big(\frac{\pi}{\sqrt{\lambda - E_{0,j}}}\Big)^{1/2}\rho_j(r)^{-(n-1)/2}\Big(e^{\pm i\Phi_j(r,\lambda)}a_{j,0} - e^{\mp i\Phi_j(r,\lambda)}b^{(\pm)}_{j,0}\Big)e_{j,0}(x).
\nonumber
\end{equation}
\end{lemma}


We equip ${\bf h}_{\infty}(\lambda)$ with the inner product
\begin{equation}
(a,b)_{{\bf h}_{\infty}(\lambda)} = \sum_{j=1}^N(a_j,b_j)_{L^2(M_j)} 
+ \sum_{j=N+1}^{N+N'}a_j\overline{b_j}\, {\rm vol}\,(M_j).
\label{innerproductnew}
\end{equation}


\begin{lemma}\label{Lemma13.3}
The operator 
$$
{\bf h}_{\infty}(\lambda) \supset {\bf h}^{comp}_{\infty}(\lambda) \ni a \to b^{(\pm)} \in {\bf h}_{\infty}(\lambda)
$$
 is isometric.
\end{lemma}

\begin{proof}
We prove the lemma for $\mathcal F^{(+)}(\lambda)$. Let $v = - 2\pi i\mathcal F^{(+)}(\lambda)^{\ast}a$. 
Take $\chi(r)  \in C^{\infty}({\mathbb R})$ such that $\chi(r)=1$ for $r<1$, and $\chi(r)=0$ for $r>2$. Let
\begin{equation}
\chi_t = \chi_0 + \sum_{j=1}^{N+N'}\chi_j\chi(r/t),
\label{Lemma13.3Definechit}
\end{equation}
where $\{\chi_j\}$ is the partition of unity on $\mathcal M$ as in (\ref{S8Partitionunity}).
 Since $(-\Delta_{\mathcal M}-\lambda)v=0$, we have by integration by parts
\begin{equation}
\begin{split}
0 &=- {\rm Im}\,(\chi_t\Delta_{\mathcal M}v,v) \\
&= 
\sum_{j=1}^{N+N'}\frac{1}{t}\int_{\mathcal M_j}\chi'\Big(\frac{r}{t}\Big)\, {\rm Im}\,\Big(\frac{\partial v}{\partial r}\overline{v} - v\overline{\frac{\partial v}{\partial r}}\Big)\sqrt{g_j}\, drdx.
\end{split}
\nonumber
\end{equation}
Replacing $v$ by the asymptotic expansion  in Lemma 
\ref{Flambdaastexmandend}, we compute the resulting integral. 
Note that 
\begin{equation}
\begin{split}
 & \Big(\frac{\partial}{\partial r}\big(e^{i\Phi_j}a_j - e^{-i\Phi_j}b_j\big)\Big)\overline{\Big(e^{i\Phi_j}a_j - e^{-i\Phi_j}b_j\Big)}\\
= & i\sqrt{\lambda - E_{0,j}}\Big(|a_j|^2 - |b_j|^2 - 2 {\rm Re}\,e^{2i\Phi_j}a_j\overline{b_j}\Big) + O(r^{-\epsilon}).
\end{split}
\nonumber
\end{equation}
Moreover, by integration by parts, we have for any $k > 0$,
$$
\frac{1}{t}\int_0^{\infty}\chi'\big(\frac{r}{t}\big)e^{ikr}dr \to 0, \quad 
t \to \infty.
$$
Therefore, 
for regular ends with $\beta_{0,j} > 1/2$, we have
\begin{equation}
\begin{split}
& \frac{1}{t}\int_{\mathcal M_j}\chi'\Big(\frac{r}{t}\Big)\, {\rm Im}\,\Big(\frac{\partial v}{\partial r}\overline{v} - v\overline{\frac{\partial v}{\partial r}}\Big)\sqrt{g_j}\, drdx \\
& =
\frac{2\pi}{t}
\int_0^{\infty}\chi'\Big(\frac{r}{t}\Big)
\left(\|a_j\|^2_{L^2(M_j)} - \|b_j\|^2_{L^2(M_j)}\right)dr + o(1).
\end{split}
\nonumber
\end{equation}
For regular ends with $0 < \beta_{0,j} \leq 1/2$, we have by a similar computation
\begin{equation}
\begin{split}
& \frac{1}{t}\int_{\mathcal M_j}\chi'\Big(\frac{r}{t}\Big)\, {\rm Im}\,\Big(\frac{\partial v}{\partial r}\overline{v} - v\overline{\frac{\partial v}{\partial r}}\Big)\sqrt{g_j}\, drdx \\
& =
\frac{2\pi}{t}
\int_0^{\infty}\chi'\Big(\frac{r}{t}\Big)
\left(\|a_j\|^2_{\ell^2} - \|b_j\|^2_{\ell^2}\right)dr + o(1).
\end{split}
\nonumber
\end{equation}
For cusp ends, 
\begin{equation}
\begin{split}
& \frac{1}{t}\int_{\mathcal M_j}\chi'\Big(\frac{r}{t}\Big)\, {\rm Im}\,\Big(\frac{\partial v}{\partial r}\overline{v} - v\overline{\frac{\partial v}{\partial r}}\Big)\sqrt{g_j}\, drdx \\
& =
\frac{2\pi}{t}
\int_0^{\infty}\chi'\Big(\frac{r}{t}\Big)
\left(|a_j|^2 - |b_j|^2\right)dr\, {\rm Vol}\,(M_j) + o(1).
\end{split}
\nonumber
\end{equation}
Adding these equalities and using $\frac{1}{t}\int_0^{\infty}\chi'(r/t)dt=-1$, we obtain the lemma. 
\end{proof}

\medskip
By Lemma \ref{Lemma13.3}, the expansion (\ref{Flambda-UpminB0ast}) in Lemma \ref{Flambdaastexmandend} is extended to all $a \in {\bf h}_{\infty}(\lambda)$.


\begin{lemma}\label{Lemma13.4}
For any $a \in {\bf h}_{\infty}(\lambda)$, we have
$$
\lim_{R\to \infty}\sum_{j=1}^{N+N'}\frac{1}{R}\int_0^R\|\mathcal F_j^{(\pm)}(\lambda)^{\ast}a\|^2_{{\bf h}_{\infty,j}(\lambda)}dr = \sum_{j=1}^{N+N'}\frac{1}{2\pi\sqrt{\lambda - E_{0,j}}}\|a_j\|^2_{{\bf h}_{\infty,j}(\lambda)}.
$$
\end{lemma}

\begin{proof}
Since $\mathcal F^{(\pm)}(\lambda)^{\ast} \in {\bf B}({\bf h}_{\infty};\mathcal B^{\ast})$, we have only to prove this lemma for $a \in {\bf h}_{\infty}^{comp}(\lambda)$. Letting $v = \mathcal F^{(\pm)}(\lambda)^{\ast}a$, we have only to compute
$$
\sum_{j=1}^{N+N'}\frac{1}{R}\int_{(0,R)\times M_j}|\chi_j v|^2\sqrt{g_j}\,drdx.
$$
We can then replace $v$ by the terms in the asymptotic expansion (\ref{Flambda-UpminB0ast}) in Lemma \ref{Flambdaastexmandend}. Arguing in the same way as above, we obtain the lemma. 
\end{proof}


\begin{lemma}\label{Flamdaastbothestimates}
There exists a constant $C = C(\lambda) > 0$ such that
$$
C^{-1}\|a\|_{{\bf h}_{\infty}(\lambda)} \leq \|\mathcal F^{(\pm)}(\lambda)^{\ast}a\|_{\mathcal B^{\ast}} \leq C\|a\|_{{\bf h}_{\infty}(\lambda)}, \quad \forall a \in {\bf h}_{\infty}(\lambda).
$$
\end{lemma}

\begin{proof}
 This lemma follows from the definition of $\mathcal B^{\ast}$ and Lemma \ref {Lemma13.4}. 
\end{proof}


\begin{lemma}\label{Lemma13.6}
Let $u \in \mathcal B^{\ast}$ satisfy the equation $(- \Delta_{\mathcal M} - \lambda)u = 0$. If $f \in \mathcal B$ satisfies $\mathcal F^{(+)}(\lambda)f=0$ or $\mathcal F^{(-)}(\lambda)f=0$, then 
$(u,f) =0$.
\end{lemma}

\begin{proof}
Suppose $\mathcal F^{(+)}(\lambda)f=0$, and let $v = R(\lambda + i0)f$. Let $\chi_t$ be as  in (\ref{Lemma13.3Definechit}). Then, by integration by parts
$$
(\chi_t u,(- \Delta_{\mathcal M} - \lambda)v) = \sum_{j=1}^{N+N'}\frac{1}{R}\int_{\mathcal M_j}\chi'\Big(\frac{r}{t}\Big)
\Big(u\overline{\partial_rv} - (\partial_ru)\overline v\Big)\sqrt{g_j}\, drdx + o(1).
$$
Since $\mathcal F^{(+)}(\lambda)f=0$, we have $v \in \mathcal B^{\ast}_0$ and also 
$\partial_r v \in \mathcal B^{\ast}_0$. Therefore, letting $t \to \infty$, the above integral vanishes, which proves the lemma. 
\end{proof}

Recall Banach's closed range theorem (see e.g. \cite{Yosida66}, p. 205). 


\begin{theorem}\label{cloedrangeth}
Let $X, Y$ be Banach spaces, and $T$ a densely defined closed operator from $X $ to $Y$. 
Then, the following 4 assertions are equivalent. \\
\noindent
(1) $R(T)$ is closed. \\
\noindent
(2) $R(T')$ is closed. \\
\noindent
(3) $R(T)  = N(T')^{\perp} = \{y \in Y\, ; \, \langle y,y^{\ast}\rangle=0, \ \forall y^{\ast}\in N(T')\}$. \\
\noindent
(4) $R(T') = N(T)^{\perp} = \{x^{\ast} \in X'\, ; \, \langle x,x^{\ast}\rangle=0, \ \forall x \in N(T)\}$.
\\
\noindent
Here, for an operator $T$ on a Banach space $X$, $R(T)$ and $N(T)$ are the range and nullspace for $T$,  $T'$ is the dual operator, and $X'$ is the dual space of $X$.
\end{theorem}


\begin{theorem}\label{S13SolspaceHelmholtzeq}
For $\lambda \in (E_{0,tot},\infty)\setminus \mathcal E$, let
\begin{equation}
\mathcal N(\lambda) = \{u \in \mathcal B^{\ast}\, ; \, (-\Delta_{\mathcal M}-\lambda)u=0\}.
\nonumber
\end{equation}
Then, we have
\begin{equation}
\mathcal F^{(\pm)}(\lambda)\mathcal B = {\bf h}_{\infty}(\lambda),
\label{Flambdaonto}
\end{equation}
\begin{equation}
\mathcal N(\lambda) = 
\mathcal F^{(\pm)}(\lambda)^{\ast}{\bf h}_{\infty}(\lambda).
\label{FlambdastSolspace}
\end{equation}
\end{theorem}

\begin{proof}
Take $X = \mathcal B$, $Y = {\bf h}_{\infty}(\lambda)$ and $T = \mathcal F^{(\pm)}(\lambda)$ in Theorem \ref{cloedrangeth}. Lemma \ref{Flamdaastbothestimates} shows that $R(T')$ is closed. For $a \in {\bf h}^{comp}_{\infty}(\lambda)$, define $f$ by (\ref{S13Definefpmetc}). 
 Then $\mathcal F^{(\pm)}(\lambda)f = a$, hence the range of $\mathcal F^{(\pm)}(\lambda)$ is dense in ${\bf h}_{\infty}(\lambda)$. However, Theorem \ref{cloedrangeth} (1) shows that it is closed, whence (\ref{Flambdaonto}) follows. Since $\mathcal F^{(\pm)}(\lambda)^{\ast}$ is an eigenoperator, $\mathcal N(\lambda) \supset 
\mathcal F^{(\pm)}(\lambda)^{\ast}{\bf h}_{\infty}(\lambda)$. If $u \in \mathcal N(\lambda)$, Lemma \ref{Lemma13.6} implies that $u \perp N(T)$. By Theorem \ref{cloedrangeth} (4), $u \in R(T')$, which proves (\ref{FlambdastSolspace}). 
\end{proof}
 
 
 \subsection{S-matrix}
 By (\ref{FlambdastSolspace}), any $u \in \mathcal N(\lambda)$ is written as $u = 2\pi i\mathcal F^{(-)}(\lambda)^{\ast}a$ for some $a = a^{(in)} \in {\bf h}_{\infty}(\lambda)$.
 Lemma \ref{Flambdaastexmandend} and the remark after Lemma \ref{Lemma13.3} imply that 
 there exists $a^{(out)} \in {\bf h}_{\infty}(\lambda)$ such that
 $u$ has the asymptotic expansion in Lemma \ref{Flambdaastexmandend} with $a = a^{(in)}$ and $b = a^{(out)}$. The S-matrix is the mapping between these asymptotic profiles $a^{(in)}$ and $a^{(out)}$. We make this fact more precise in the following theorem.
 
 
 \begin{theorem}\label{Smatrixtheorem}
Let $\lambda \in (E_{0,tot},\infty)\setminus \mathcal E$. Then, for any $a^{(in)} \in {\bf h}_{\infty}(\lambda)$, there exist unique $u \in \mathcal N(\lambda)$ and $a^{(out)} \in {\bf h}_{\infty}(\lambda)$ such that, letting $a = a^{(in)}$ and $b = a^{(out)}$, $u$ behaves  as follows on each end $\mathcal M_j$. 

\smallskip
\noindent
(1) For regular ends with $\beta_{0,j} > 1/2$,
\begin{equation}
u \simeq \Big(\frac{\pi}{\sqrt{\lambda - E_{0,j}}}\Big)^{1/2}\rho_j(r)^{-(n-1)/2}\left(e^{- i\Phi_j(r,\lambda)}a_j(x) - e^{ i\Phi_j(r,\lambda)}b_j(x)\right).
\nonumber
\end{equation}
(2) For regular ends with $0 < \beta_{0,j} \leq 1/2$,
\begin{equation}
\begin{split}
u\simeq & \Big(\frac{\pi}{\sqrt{\lambda - E_{0,j}}}\Big)^{1/2}\rho_j(r)^{-(n-1)/2}\\
& \times \sum_{\ell=0}^{\infty}
\chi_{\ell,j}(r)\Big(e^{- i\varphi_j(\lambda,\lambda_{\ell,j},r)}a_{j,\ell} - e^{i\varphi_j(\lambda,\lambda_{\ell,j},r)}b_{j,\ell}\Big)e_{\ell,j}(x),
 \end{split}
\nonumber
\end{equation} 
\noindent
where $\chi_{\ell,j}(r) = \chi(r/r_0(\lambda,\lambda_{\ell,j}))$ with $\chi \in C^{\infty}({\mathbb R})$, $\chi(r)=0$ for $r < 1$, $\chi(r)=1$ for $r > 2$.

\noindent
(3) For cusp ends
\begin{equation}
u \simeq \Big(\frac{\pi}{\sqrt{\lambda - E_{0,j}}}\Big)^{1/2}\rho_j(r)^{-(n-1)/2}\Big(e^{- i\Phi_j(r,\lambda)}a_{j,0} - e^{ i\Phi_j(r,\lambda)}b_j,{0}\Big)e_{j,0}(x).
\nonumber
\end{equation}

The operator
\begin{equation}
S(\lambda) : a^{(in)} \to a^{(out)}
\nonumber
\end{equation}
is unitary on ${\bf h}_{\infty}(\lambda)$.
 \end{theorem}

\begin{proof}
First we prove the existence of $u$ and $a^{(out)}$. Take $u = 2\pi i \mathcal F^{(-)}(\lambda)^{\ast}a^{(in)}$.  Then by Lemma \ref{Flamdaastbothestimates}, $\|u\|_{\mathcal B^{\ast}} \leq C\|a^{(in)}\|_{{\bf h}_{\infty}(\lambda)}$. For any $\epsilon > 0$, take $a_{\epsilon} \in {\bf h}^{comp}_{\infty}(\lambda)$ such that $\|a^{(in)} - a_{\epsilon}\|_{{\bf h}_{\infty}(\lambda)} < \epsilon$. Put $u_{\epsilon} = 2\pi i\mathcal F^{(-)}(\lambda)^{\ast}a_{\epsilon}$. Then, by Lemma \ref{Flambdaastexmandend}, there exists $b_{\epsilon} \in {\bf h}_{\infty}(\lambda)$ such that $u_{\epsilon}, a_{\epsilon}, b_{\epsilon}$ have the asymptotic expansion in the present lemma. By the isometric property in Lemma \ref{Lemma13.3}, $b_{\epsilon}$ tends to some $a^{(out)} \in {\bf h}_{\infty}(\lambda)$ as $\epsilon \to 0$. Thus, $a^{(in)},u, a^{(out)}$ have the desired properties. 

To show the uniqueness, suppose for a given $a^{(in)}$, there exit two such $u', u''$. Then, $w = u'- u''$ is a solution to the equation $(-\Delta_{\mathcal M} - \lambda)w=0$ satisfying the outgoing radiation condition. Therefore, $w=0$.

Lemma \ref{Lemma13.3} shows that $S(\lambda)$ is isometric. Arguing as above,  changing the roles of $a^{(in)}$ and $a^{(out)}$, one can prove that the range of $S(\lambda)$ is dense. This implies the unitarity. 
\end{proof}


\section{Radial solutions on cusp end}\label{sectAsympSolCusp}
In this section, we equip a cusp end with a warped product metric $ds^2 = (dr)^2 + \rho(r)^2h_M(x,dx)$, and 
construct (super) exponentially growing or decaying solutions to the equation
$(- \Delta_{\mathcal M} - \lambda)u=0$.  A  typical example of $\rho(r)$ is
\begin{equation}
\rho(r) = \left\{
\begin{split}
& e^{c_0r}, \quad c_0 < 0, \\
& e^{c_1r^{1-\alpha}}, \quad c_1 < 0, \quad 0 < \alpha < 1, \\
& r^{\beta}, \quad \beta < 0.
\end{split}
\right.
\label{S13Examplecusp}
\end{equation}
We formulate its perturbation in the form of asymptotic series.

For a real constant $\kappa$,  let $\widetilde A_{<\kappa}$ be the set of finite linear combinations of the following terms
$$
(\log r)^i r^{\alpha}, \quad i = 0, 1, 2, \cdots, \quad  \alpha < \kappa.
$$
Let $\widetilde{\mathcal A}^{\kappa}$ be the set of real functions $\widetilde f(r)$ satisfying
\begin{equation}
\widetilde f(r) - a_{0}r^{\kappa} \in \widetilde A_{<\kappa},
\label{S7f(r)=sumcirkapai}
\end{equation}
where $a_{0}$ is a constant.  
We define $\mathcal A^{\kappa}$ to be the set of real functions  $f \in C^{\infty}((0,\infty))$ such that for any $N > 0$, there exists $\widetilde f \in \widetilde{\mathcal A}^{\kappa}$ satisfying 
\begin{equation}
\partial_r^m\big(f(r) - \widetilde f(r)\big) = O(r^{-N- m}), \quad r \to \infty, \quad \forall m \geq 0.
\nonumber
\end{equation}
We put for $f \in \mathcal A^{\kappa}$
\begin{equation}
a_0(f) = a_{0},
\nonumber
\end{equation}
where $a_{0}$ is from (\ref{S7f(r)=sumcirkapai}). 
Finally we define $\mathcal A^{\kappa}_{+}$ by 
\begin{equation}
\mathcal A^{\kappa}_{+} = \{f \in \mathcal A^{\kappa}\, ; \,  f(r) > 0\}.
\nonumber
\end{equation}

The following lemma is obvious.


\begin{lemma} $\cup_{\kappa}\mathcal A^{\kappa}$ is an algebra in the following snese. 
\begin{equation}
\mathcal A^{\kappa}\mathcal A^{\kappa'} \subset \mathcal A^{\kappa+\kappa'},
\nonumber
\end{equation}
\begin{equation}
\partial_r^m\mathcal A^{\kappa} \subset  \mathcal A^{\kappa-m}, \quad \forall m \geq 0,
\nonumber
\end{equation}
\begin{equation}
\frac{1}{f} \in  \mathcal A^{-\kappa}_+ \quad {\it if} \quad f \in \mathcal A^{\kappa}_+.
\nonumber
\end{equation}
\end{lemma}

We solve the equation
\begin{equation}
\left(- {\partial_r}^2 - \frac{(n-1)\rho'(r)}{\rho(r)}\partial_r + \frac{B}{\rho(r)^{2}} - \lambda\right)u=0,
\label{S13Equation1}
\end{equation}
where $B \in [B_0,\infty)$ is a parameter, $B_0$ being a fixed positive constant. In our applications, $B_0$ is the least non-zero eigenvalue of $- \Delta_M$.


 \begin{lemma}
Let $u = \rho^{-(n-2)/2}w$ in (\ref{S13Equation1}) and introduce a new variable
\begin{equation}
t = \int_0^r\frac{\sqrt{B}}{\rho(\tau)}d\tau.
\nonumber
\end{equation}
Then, we have
\begin{equation}
\left\{
\begin{split}
& - \frac{d^2w}{dt^2} + \big(1 + V(t)\big)w = 0, \\
& V(t) = \frac{\rho^2}{B}\left(- \lambda + \frac{(n^2 -2n)}{4}\Big(\frac{\rho'}{\rho}\Big)^2 + 
\frac{(n-2)}{2}\Big(\frac{\rho'}{\rho}\Big)'\right),\\
&\rho' = d\rho/dr, \quad \rho'' = d^2\rho/dr^2.
\end{split}
\right.
\label{Lemma13.2d2wdt2}
\end{equation}
\end{lemma}
\begin{proof}
We put $u = \rho^{-(n-1)/2}v$ so that
\begin{equation}
- v'' + \left(\frac{B}{\rho^2} - \lambda + \Big(\frac{(n-1)\rho'}{2\rho}\Big)^2 + \Big(\frac{(n-1)\rho'}{2\rho}\Big)'\right)v=0.
\nonumber
\end{equation}
Passing to the variable $t$, and letting
\begin{equation}
q(t) = 1 + \frac{\rho^2}{B}\left(- \lambda  + \Big(\frac{(n-1)\rho'}{2\rho}\Big)^2+ \Big(\frac{(n-1)\rho'}{2\rho}\Big)'\right),
\nonumber
\end{equation}
where $' = d/dr$, 
we have
\begin{equation}
- \frac{d^2v}{dt^2} + \frac{\rho'}{\sqrt{B}}\frac{dv}{dt} + q(t)v=0.
\nonumber
\end{equation}
We put $v = e^{\varphi}w$, and obtain
\begin{equation} 
- \frac{d^2w}{dt^2} + \Big(-2\frac{d\varphi}{dt} + \frac{1}{\sqrt{B}}\frac{d\rho}{dr}\Big)\frac{dw}{dt} + \Big(q - \frac{d^2\varphi}{dt^2} - \big(\frac{d\varphi}{dt}\big)^2 + \frac{1}{\sqrt{B}}\frac{d\rho}{dr}\frac{d\varphi}{dt}\Big)w=0.
\nonumber
\end{equation}
We take $\varphi = \log\sqrt{\rho}$ so that the 2nd term vanishes. We thus arrive at
\begin{equation}
\left\{
\begin{split}
& - \frac{d^2w}{dt^2} + \Big(q + \frac{\rho'^2 - 2\rho\rho''}{4B}\Big)w = 0, \\
& u = \rho^{-(n-2)/2}w, \quad \rho' = d\rho/dr, \quad \rho'' = d^2\rho/dr^2.
\end{split}
\right.
\nonumber
\end{equation}
The term in the parenthesis is rewritten as in (\ref{Lemma13.2d2wdt2}).
\end{proof}

 Let us check the properties of the change of variable $[0,\infty) \ni  r \to t \in [0,\infty)$. We put
\begin{equation}
s = \frac{t}{\sqrt{B}} = \int_0^r\rho(\tau)^{-1}d\tau,
\label{S13s=intrho-1}
\end{equation}
and study the following three cases corresponding to (\ref{S13Examplecusp}): 
\begin{equation}
\rho = \left\{
\begin{split}
& e^{-\phi_1}f_1, \quad \phi_1 \in \mathcal A^{1}, \quad a_0(\phi_1) > 0,  \quad f_1 \in \mathcal A^{0}_+,  \\
&  e^{-\phi_2}f_2,  \quad \phi_2 \in \mathcal A^{1-\alpha}, \quad 0 < \alpha < 1, \quad a_0(\phi_2) > 0, \quad f_2 \in \mathcal A^0_+, \\
&  f_3, \quad f_3 \in \mathcal A_+^{\beta}, \quad \beta<0.
\end{split}
\right.
\label{S7rhoasymp1}
\end{equation}

We use the following notation. 
By $f(r) \sim g(r)$, we mean that
\begin{equation}
\frac{f(r)}{g(r)} \to 1, \quad {\rm as}\quad r \to \infty,
\label{Definef(r)simg(r)}
\end{equation}
and $f(r) \asymp g(r)$ means that the following inequalities hold
\begin{equation}
C^{-1}g(r) \leq f(r) \leq Cg(r), \quad r > r_0,
\label{S1AsypDefine}
\end{equation}
for some constants $C, r_0>0$.


\begin{lemma}
We  have as a function of $r$
\begin{equation}
s = s(r) = \left\{
\begin{split}
& e^{\phi_1}g_1,  \quad g_1 \in \mathcal A^{0}_+, \\
&  e^{\phi_2}g_2, \quad g_2 \in \mathcal A_+^{\alpha}, \\
& g_3 \in \mathcal A_+^{1-\beta}, 
\end{split}
\right.
\label{S7t(r)=}
\end{equation}
where $\phi_1, \phi_2$ are the same functions as in (\ref{S7rhoasymp1}). 
In particular, we have
\begin{equation}
\log s \asymp r, \quad \log s \asymp r^{1-\alpha}, \quad 
s \asymp r^{1-\beta},
\label{S13srasympresectively}
\end{equation}
respectively.
\end{lemma}
\begin{proof} Noting that $1/\rho$ has the form
$\dfrac{1}{\rho} = \dfrac{e^{\phi}}{f} =  \dfrac{1}{f\phi'}\dfrac{d}{d\tau}e^{\phi},$
we integrate by parts in (\ref{S13s=intrho-1}).  We then have
\begin{equation}
s \sim \frac{e^{\phi}}{f\phi'},
\label{Lemma13.3saympto}
\end{equation}
which proves the 1st two cases of (\ref{S7t(r)=}). The 3rd case is easy. 
For $0 < \kappa \leq 1$ and $\beta < 0$,  as $r \to \infty$,
$$
\int_1^re^{ct^{\kappa}}dt \asymp r^{1-\kappa}e^{cr^{\kappa}} \ \ (c  > 0), \quad 
\int_1^r t^{-\beta}dt \asymp r^{1-\beta}.
$$
Therefore, we have for the 1st two cases 
\begin{equation}
\log s \asymp (1-\kappa)\log r + r^{\kappa} \asymp r^{\kappa},
\nonumber
\end{equation}
where $\kappa = 1, 1- \alpha$ 
and for the 3rd case
\begin{equation}
s \asymp r^{1-\beta}, 
\label{S13baeta<0}
\nonumber
\end{equation}
which proves (\ref{S13srasympresectively}).   
\end{proof}


\begin{lemma}\label{Lemma13.4dsmh(r(s))}
 If $h(r) \in \mathcal A^{\kappa}$ and $m \geq1$, we have, as a function of $r>0$,
\begin{equation}
\Big(\big(\frac{d}{ds}\big)^mh(r(s))\Big)(r) 
\in 
\left\{
\begin{split}
& e^{-m\phi_1}\mathcal A ^{\kappa-1},  \\
& e^{-m\phi_2}\mathcal A^{\kappa-1  - (m-1)\alpha}, \\
&\mathcal A^{\kappa + m(\beta-1)}.
\end{split}
\right.
\label{d/dtmh(r)}
\nonumber
\end{equation}
\end{lemma}

\begin{proof}
Using $\dfrac{d}{ds}= \rho\dfrac{d}{dr}$, we have
$$
\frac{d}{ds}h(r) = \left\{
\begin{split}
& e^{-\phi_1}f_1h', \\
& e^{-\phi_2}f_2h', \\
& f_3h'.
\end{split}
\right.
$$
Differentiataing this equality, we get the lemma. 
\end{proof}

We return to the differential equation (\ref{Lemma13.2d2wdt2}). 


\begin{lemma} \label{Lemma13.5V(t)decay}
Fix $B_0 > 0,  C_0 > 1$ arbitrarily. Then, for $B \geq B_0$,  $\dfrac{t}{\sqrt{B}} \geq C_0$ and $m \geq 1$, 
\begin{equation}
\Big(\frac{d}{dt}\Big)^m V(t) = \left\{
\begin{split}
& O\Big(\big(\frac{\sqrt{B}}{t}\big)^mt^{-2}\Big), \\
& O\Big(\big(\frac{\sqrt{B}}{t}\big)^m\Big(\log\frac{t}{\sqrt{B}}\Big)^{\frac{2\alpha}{1-\alpha}}t^{-2}\Big), \\
& O\Big(\big(\frac{\sqrt{B}}{t})^{m-\frac{2}{1-\beta}}t^{-2}\Big).
\end{split}
\right.
\nonumber
\end{equation}
\end{lemma}

\begin{proof}
Recalling  $t = \sqrt{B}s$, we put
$$
V(t) = W(s) = \frac{\rho^2}{B}\left(- \lambda + \frac{(n^2 -2n)}{4}\Big(\frac{\rho'}{\rho}\Big)^2 + 
\frac{(n-2)}{2}\Big(\frac{\rho'}{\rho}\Big)'\right).
$$
Then, by Lemma \ref {Lemma13.4dsmh(r(s))}, we have
$$
\Big(\frac{d}{ds}\Big)^mW(s) = \frac{1}{B}
\left\{
\begin{split}
&e^{-(m+2)\phi_1}k_1(r), \quad k_1(r) \in \mathcal A^0, \\
&e^{-(m+2)\phi_2}k_2(r), \quad k_2(r) \in \mathcal A^{-m\alpha}, \\
&k_3(r), \quad k_3(r) \in \mathcal A^{2\beta + m(\beta-1)}.
\end{split}
\right.
$$
In view of (\ref{Lemma13.3saympto}), we have
$$
e^{-\phi_1} \sim \frac{1}{sf_1\phi_1'} \asymp \frac{1}{s}, \quad 
e^{-\phi_2} \sim \frac{1}{sf_2\phi_2'} \asymp \frac{r^{\alpha}}{s}, 
$$
By Lemma \ref{Lemma13.4dsmh(r(s))}, we then have for $m \geq 1$
$$
\Big(\frac{d}{ds}\Big)^mW(s) = \frac{1}{B}\left\{
\begin{split}
& O(s^{-(m+2)}), \\
& O(s^{-m-2}(\log s)^{\frac{2\alpha}{1-\alpha}}), \\
& O(s^{-m + \frac{2\beta}{1-\beta}}).
\end{split}
\right.
$$
Using $t = \sqrt{B}s$, we obtain the lemma. 
\end{proof}

In the equation (\ref{Lemma13.2d2wdt2}), we put $w = ae^{\varphi}$. Then
\begin{equation}
- w'' + (1 + V)w = e^{\varphi}\Big(a\left(1 + V - (\varphi')^2\right) 
- (2a'\varphi' + a\varphi'') -a''\Big).
\nonumber
\end{equation}
We take
\begin{equation}
\varphi(t) =  \int_{t_0}^t\sqrt{1 + V(\tau)}\, d\tau,
\nonumber
\end{equation}
\begin{equation}
a(t) = \varphi'^{-1/2} = \left(1 + V(t)\right)^{-1/4},
\nonumber
\end{equation}
where $t_0 = t_0(B)$ is chosen so that $|V(t)| \leq 1/2$ for $t \geq t_0$. Then, we have
\begin{equation}
1 +V - (\varphi')^2 = 0, \quad 2a'\varphi' + a\varphi'' = 0,
\label{S13varphiandaequation}
\end{equation}
 and by virtue of Lemma \ref{Lemma13.5V(t)decay}, there exists $p > 2$ such that 
\begin{equation}
a''(t) = O(t^{-p}).
\label{S13aprimeorime=t-2}
\end{equation}
if $\dfrac{t}{\sqrt{B}} \geq C_0 \max\{t^{\epsilon},1\}$ for some $0 < \epsilon < 1$.

We look for solutions of  (\ref{Lemma13.2d2wdt2}) in the form $ae^{\varphi}(1 + v)$. Then $v$ satisfies
\begin{equation}
v''+ 2\Big(\frac{a'}{a} + \varphi'\Big)v' + \frac{a''}{a}v =- \frac{a''}{a},
\end{equation}
where we have used (\ref{S13varphiandaequation}).
Putting $ f = \left(\begin{array} {c}v\\
v'\end{array}\right)$,  $g = \left(\begin{array} {c}0 \\
-a''/a\end{array}\right)$, we get the equation
\begin{equation}
\frac{df}{dt} = A(t)f + B(t)f + g(t),
\label{S7f(t)diffeq}
\end{equation}
\begin{equation}
A(t) = \left(
\begin{array}{cc}
0 & 1 \\
0 & - \dfrac{2a'}{a} - 2 \varphi'
\end{array}
\right), 
\quad 
B(t) =  \left(
\begin{array}{cc}
0 & 0 \\
- \dfrac{a''}{a} &0
\end{array}
\right).
\nonumber
\end{equation}
The fundamental matrix for the equation $dh/dt = A(t)h$ is
\begin{equation}
F(t,s) = C(t)C(s)^{-1},
\nonumber
\end{equation}
\begin{equation}
C(t) = \left(
\begin{array}{cc}
1 & v_0(t)\\
0 & v_0'(t)
\end{array}
\right), \quad 
v_0(t) = -\int_t^{\infty}\frac{e^{-2\varphi}}{a^2}d\tau.
\nonumber
\end{equation}
Then 
$$
F(t,s) = \frac{1}{v_0'(s)}\left(
\begin{array}{cc}
v_0'(s) & - v_0(s) + v_0(t) \\
0 & v_0'(t)
\end{array}
\right).
$$
Noting that 
$$
v_0(t) = e^{-2\varphi(t)}c_0(t), \quad c_0(t) \in \mathcal A_+^0, \quad \varphi(t) = t  + O(t^{1-\epsilon}), 
$$
we see that
 $F(t,s)$ is bounded for $t \geq s \geq t_0$. The differential equation 
 (\ref{S7f(t)diffeq}) is now transformed into the integral equation
\begin{equation}
f(t) = - \int_t^{\infty}F(t,s)B(s)f(s)ds - \int_t^{\infty}F(t,s)g(s)ds.
\nonumber
\end{equation}
By virtue of (\ref{S13aprimeorime=t-2}), one can solve it uniquely by iteration, and obtain $f(t) = O(t^{-1-\epsilon})$. Therefore  (\ref{Lemma13.2d2wdt2}) has a solution $w_+$ which behaves like 
$w_+ \sim ae^{\varphi}$. Another solution $w_-(t)$ is obtained by 
$$
w_-(t) = w_+ (t)\int_{t}^{\infty}w_+(\tau)^{-2}d\tau.
$$
We have thus constructed two solutions $w_{\pm}$ of  (\ref{Lemma13.2d2wdt2}) such that
\begin{equation}
w_{\pm}(t) \sim a_{\pm}e^{\pm \varphi}, \quad 
a_+ = a, \quad a_- \sim \frac{1}{a\varphi'}.
\nonumber
\end{equation}
We pass  to the variable $r$, and put $\psi(r) = \varphi(t)$. Then, we have
$$
\frac{d\psi}{dr} = \sqrt{\frac{B}{\rho^2}- \lambda + \frac{(n^2-2n)}{4}\Big(\frac{\rho'}{\rho}\Big)^2 + \frac{(n-2)}{2}\Big(\frac{\rho'}{\rho}\Big)'}, 
 $$
and proved the following theorem. 

\begin{theorem}\label{AsymptoticsolODE}
Assume (\ref{S7rhoasymp1}). Then, for any $B > 0$, there exist solutions $u_0^{(\pm)}$ of the equation
$$
- u'' - \frac{(n-1)\rho'}{\rho}u' + \Big(\frac{B}{\rho^2}-\lambda\Big)u = 0
$$
satisfying
$$
u_0^{(\pm)} \sim \rho(r)^{-(n-2)/2} e^{\pm \psi(r)},
$$
$$
\psi(r) = \int_{r_0}^r \sqrt{\frac{B}{\rho^2}- \lambda + \frac{(n^2-2n)}{4}\Big(\frac{\rho'}{\rho}\Big)^2 + \frac{(n-2)}{2}\Big(\frac{\rho'}{\rho}\Big)'}dr,
$$
$r_0 = r_0(B)$ being a sufficiently large constant. 
\end{theorem}

Note that by (\ref{Lemma13.3saympto}), we have the following asymptotics of $\psi$ as $r \to \infty$:
\begin{equation}
\frac{\psi(r)}{\sqrt{B}} \sim \left\{
\begin{split}
& \frac{e^{\phi_1(r)}}{f_1(r)\phi_1'(r)},\\
&  \frac{e^{\phi_2(r)}}{f_2(r)\phi_2'(r)}, \\
& \int f_3(r)^{-1}dr.
\end{split}
\right.
\nonumber
\end{equation}


\section{Generalized S-matrix}\label{S15GeneSmatrix}
We generalize the notion of  S-matrix by enlarging the solution spece of the equation $(-\Delta_{\mathcal M} - \lambda)u=0$ on the cusp end. To make the distinction clear, we call  the scattering data and the S-matrix constructed  in \S \ref{SectionPhysicalSmatrix} {\it physical}, and call the ones to be inroduced here {\it non-physical}.  To construct  these non-physical scattering data and S-matrix by the separation of variables, we assume that our end is a pure cusp. Namely, we impose the assumption (A-4-2).

Let $N+k \leq j \leq N+N'$, and $0 = \lambda_{0,j} \leq \lambda_{1,j} \leq \lambda_{2,j} \leq \cdots$ be the eigenvalues of $- \Delta_{M_j}$ with complete orthnormal system of eigenvectors $e_{\ell,j}(x), \ell = 0,1,2,\cdots$. 
We put
 \begin{equation}
 \Phi_j(r,B) = \int_{r_0}^r \sqrt{\frac{B}{\rho_j^2}- \lambda + \frac{(n^2-2n)}{4}\Big(\frac{\rho_j'}{\rho_j}\Big)^2 + \frac{(n-2)}{2}\Big(\frac{\rho_j'}{\rho_j}\Big)'}dr.
\nonumber
 \end{equation}
By Theorem \ref{AsymptoticsolODE}, there exist solutions $u_{\ell,j,\pm}$ to the equation
\begin{equation}
- u'' - \frac{(n-1)\rho_j'}{\rho_j}u' + \Big(\frac{\lambda_{\ell,j}}{\rho_j^2}-\lambda\Big)u = 0,
\label{S15radialequation} 
\end{equation}
which behave like
\begin{equation}
u_{\ell,j,\pm} \sim \rho_j(r)^{-(n-2)/2}e^{\pm \Phi_j(r,\lambda_{\ell,j})}, \quad r \to \infty.
\nonumber
\end{equation}
Take any solution $u$ of the equation
\begin{equation}
(- \Delta_{\mathcal M} - \lambda)u = 0, \quad {\rm on} \quad \mathcal M_j, \quad j = N+1,\cdots,N+N'.
\label{EquationincuspS15}
\end{equation}
Expanding it by $e_{\ell,j}$, we have
\begin{equation}
(u(r,\cdot),e_{\ell,j})_{L^2(M_j)} = a_{\ell,j}u_{\ell,j,+}(r) + 
b_{\ell,j}u_{\ell,j,-}(r).
\label{S1uexapnd}
\end{equation}

Here, we introduce two spaces of sequences ${\bf A}_{j,\pm}$ : 
\begin{equation}
{\bf A}_{j,\pm} \ni \{c_{\ell,\pm}\}_{\ell=0}^{\infty} 
\Longleftrightarrow \sum_{\ell=0}^\infty|c_{\ell,\pm}|^2|u_{\ell,j,\pm}(r)|^2 < \infty,\quad \forall r > 0.
\nonumber
\end{equation}


\begin{lemma}
For a solution $u$ of (\ref{EquationincuspS15}), let $a_{\ell,j}, b_{\ell,j}$ be defined by (\ref{S1uexapnd}). If $\{a_{\ell,j}\}_{\ell=0}^{\infty} \in {\bf A}_{j,+}$, then $\{b_{\ell,j}\}_{\ell=0}^{\infty} \in {\bf A}_{j,-}$.
\end{lemma}

\begin{proof}
Since $\sum_{\ell\geq 0}|(u(r,\cdot),e_{\ell,j})_{L^2(M_j)}|^2 < \infty$, the lemma follows from (\ref{S1uexapnd}).
\end{proof}

Any finite sequence belongs to ${\bf A}_{j,\pm}$. For the hyperbolic metric, one can find a more explicit subspace of ${\bf A}_{j,\pm}$ by using the asymptotic expansion of modified Bessel functions (see \cite{IKL11}).

\medskip
Using the partition of unity $\{\chi_j\}$ in (\ref{S8Partitionunity}),  we define the generalized incoming  solution on the cusp end $\mathcal M_j$ by
\begin{equation}
\Psi_j^{(in)} = \chi_j\sum_{\ell=0}^{\infty}a_{\ell,j}u_{\ell,j,+}(r)e_{\ell,j}(x),\quad \{a_{\ell,j}\}_{\ell=0}^{\infty}\in {\bf A}_{j,+},
\nonumber
\end{equation}
which is (super)-exponentially growing as $r \to \infty$, and the generealized outgoing solution by
\begin{equation}
\Psi_j^{(out)} = \chi_j\sum_{\ell=0}^{\infty}b_{\ell,j}u_{\ell,j,-}(r)e_{\ell,j}(x),\quad \{b_{\ell,j}\}_{\ell=0}^{\infty}\in {\bf A}_{j,-},
\label{S15Psiin}
\end{equation}
which is (super)-exponentially decaying as $r \to \infty$. We also define the spaces of generalized scattering data by
\begin{equation}
{\bf h}_{\infty}^{(in)}(\lambda) = \Big({\mathop\oplus_{j=1}^N}c_j(\lambda)L^2(M_j) \Big)\oplus
 \Big({\mathop\oplus_{j=N+1}^{N+k-1}}c_j(\lambda){\mathbb C}_j\Big)\oplus \Big({\mathop\oplus_{j=N+k}^{N+N'}}c_j(\lambda){\bf A}_{j,+}\Big),
\label{S15Psiout}
\end{equation}
\begin{equation}
{\bf h}_{\infty}^{(out)}(\lambda) = \Big({\mathop\oplus_{j=1}^N}c_j(\lambda)L^2(M_j) \Big)\oplus
 \Big({\mathop\oplus_{j=N+1}^{N+k-1}}c_j(\lambda){\mathbb C}_j\Big)\oplus \Big({\mathop\oplus_{j=N+k}^{N+N'}}c_j(\lambda){\bf A}_{j,-}\Big),
\nonumber
\end{equation}
where $c_j(\lambda)$ is the characteristic function of the interval $(E_{0,j},\infty)$.


\begin{theorem}\label{GeneralSmatrix}
For any generalized incoming data $a^{(in)} \in {\bf h}_{\infty}^{(in)}(\lambda)$, there exist a  unique solution $u$ of the equation $(-\Delta_{\mathcal M} - \lambda)u=0$, and the  outgoing data $a^{(out)} \in {\bf h}_{\infty}^{(out)}(\lambda)$ such that
\begin{equation}
u- \sum_{j=N+k}^{N+N'}\Psi_j^{(in)} \in \mathcal B^{\ast},
\nonumber
\end{equation}
\begin{equation}
u = \Psi_j^{(in)} - \Psi_j^{(out)}, \quad {\rm on} \quad \mathcal M_j, \quad j = N+k,\cdots,N+N',
\nonumber
\end{equation}
and on the ends $\mathcal M_j$, $1 \leq j \leq N+k-1$, $u$ has the asymptotic form in Theorem \ref{Smatrixtheorem}. Here, $\Psi_j^{(in)}$ and  $\Psi_j^{(out)}$ are written by (\ref{S15Psiin}), (\ref{S15Psiout}) with $a_{\ell,j}$, $b_{\ell,j}$ replaced by the associated components of $a_j^{(in)}$ and $a_j^{(out)}$.
\end{theorem}

\begin{proof}
Put $u^{(in)} = \sum_{j=N+k}^{N+N'}\Psi_j^{(in)}$ and $u = u^{(in)} - R(\lambda + i0)(- \Delta_{\mathcal M}-\lambda)u^{(in)}$. Then, $u$ has the desired properties. If $u_1$ and $u_2$ are two such solutions, $u_1 - u_2$ is an outgoing solution of the equation $(-\Delta_{\mathcal M} - \lambda)u=0$, hence vanishes identically. 
\end{proof}

We call the mapping
\begin{equation}
\begin{split}
\mathcal S(\lambda) : {\bf h}_{\infty}^{(in)}(\lambda) \ni a^{(in)} \to a^{(out)} \in {\bf h}_{\infty}^{(out)}(\lambda)
\end{split}
\nonumber
\end{equation}
the {\it generalized scatteing matrix}. \index{generalized scattering matrix}

\medskip
Let us remark here that in Theorem \ref{AsymptoticsolODE} 
the decaying solution $u_0^{(-)}$ is determined uniquely from its asymptotic behavior near infinity, while the growing solution $u^{(+)}_0$ is not unique, since $u_0^{(+)} + cu_0^{(-)}$, $c$ being any constant, is again an increasing solution. However, it gives no harm to the definition of the generalized S-matrix. In fact, given two incoming data $\Psi_j^{(in)}, {\Psi_j^{(in)}}'$, let $u, u'$ be the associated solutions to the Helmholtz equation. If $\Psi_j^{(in)} - {\Psi_j^{(in)}}' \in \mathcal B_0^{\ast}$, $u-u'$ is outgoing, hence vanishes identically. Therefore, $u$ and $u'$ give the same decaying solution in the end $\mathcal M_j$.

\chapter{Inverse scattering}



\section{From S-matrix to   source-to-solution map}
\label{SectionSmatrixtoDN1}
Let us  start the reconstruction of  the manifold from the (generalized) scattering matrix. We follow the arguments in \cite{IKL11} and \cite{IKL13(1)} with some modifications.  We  reduce the problem to the source-to-solution map in the {\it interior} domain,
see \cite{BKL2020,HLOS}.

 Let $\mathcal O\subset \mathcal M^{reg}$ be a relatively compact open set.
We  consider the problem with source $F$ supported in $\mathcal O$ and the radiation condition (see Definition \ref{TotalRadCond}):
\begin{equation}
\left\{
\begin{split}
& (- \Delta_{\mathcal M} - \lambda)u = F, \quad {\rm in} \quad {\mathcal M}, \\
&\hbox{$u$ satisfies the radiation condition.}
\end{split}
\right.
\label{Interior problem}
\end{equation}
We extend $F$ to be 0 outside $\mathcal O$. 
By Theorem \ref{ModelLAPforbfh}, for $\lambda \in \sigma_e(H)\setminus{\mathcal E}$, there exists a unique solution to this equation, denoted by
$$
u^F_+(\lambda) = (- \Delta_{\mathcal M} - \lambda - i0)^{-1}F, \quad u^F_-(\lambda) = (- \Delta_{\mathcal M} - \lambda + i0)^{-1}F.
$$
We define the {\it stationary source-to-solution operator}\index{stationary source-to-solution operator} by
\begin{equation}
U_{\mathcal O,{\pm}}(\lambda): L^2(\mathcal O) \ni F \to u^F_{\pm}(\lambda) \in L^2(\mathcal O).
\label{C2S1Sourcetosolop}
\end{equation}

Now, we enter into the first step of the inverse problem. Suppose we are given two manifolds $\mathcal M^{(1)}$ and  $\mathcal M^{(2)}$ satisfying the assumptions (A-1), (A-2), (A-3) and (A-4-1), (A-4-2).  
Then, $\mathcal M^{(i)}$ is written as 
\begin{equation}
\mathcal M^{(i)} = \mathcal K^{(i)}\cup{\mathcal M^{(i)}_1}\cup\cdots\cup{\mathcal M^{(i)}_{N_i+N_i'}},
\end{equation}
where $\mathcal K^{(i)}$ is bounded, and $\mathcal M^{(i)}_{j}$'s are non-compact. Note that the number of ends of $\mathcal M^{(i)}$ is not assumed to be the same for $i = 1, 2$ a-priori. 
Let $H_i$ be the Laplacian on $\mathcal M^{(i)}$, and $S^{(i)}(\lambda)$ the associated S-matrix, which is an $(N_i+N_i')\times (N_i+N_i')$ operator-valued matrix. Let $S^{(i)}_{pq}(\lambda)$ be its $(p,q)$-entry. 
Let $\mathcal E_i$ be the set of exceptional points $\mathcal E$ for $H_i$. 
First we consider the case for regular ends. Let $\mathcal M_1^{(i)}$ be the 1st regular end of $\mathcal M^{(i)}$.


\begin{theorem}\label{S15SmatrixdeterminesLW2}
Assume that $\mathcal M^{(1)}_{1}$ and $\mathcal M^{(2)}_{1}$ are isometric. 
If $S^{(1)}_{11}(\lambda) = S^{(2)}_{11}(\lambda)$ for some $\lambda \in \big(\sigma_e(H_1)\setminus\mathcal E_1\big)\cap\big(\sigma_e(H_2)\setminus \mathcal E_2\big)$, then  $U_{\mathcal O,\pm}^{(1)}(\lambda) = U_{\mathcal O,\pm}^{(2)}(\lambda)$.
\end{theorem}

\begin{proof}
We consider the outgoing case, and omit the subscript $+$. First we assume that $\beta_1 > 1/3$ on $\mathcal M_1^{(1)} = \mathcal M_1^{(2)}$. 
Take $a_1 \in L^2(M_1^{(1)}) = L^2(M_1^{(2)})$ and put $u^{(i)} = \mathcal F^{(i,+)}(\lambda)^{\ast}a^{(i)}$, where $a^{(i)} = (a_1,0,\cdots,0)$, and $\mathcal F^{(i,+)}(\lambda)$ is $\mathcal F^{(+)}(\lambda)$ for $\mathcal M^{(i)}$. 
Then 
\begin{equation}
(- \Delta_{\mathcal M}-\lambda)(u^{(1)}- u^{(2)}) = 0
\label{Th2.1.1Deltau1-u2=0}
\end{equation}
 on $\mathcal M^{(1)}_{1} = 
\mathcal M^{(2)}_{1}$.
In view of the asymptotic expansion in Lemma \ref{FlambdaastExpand} and the assumption 
$S^{(1)}_{11}(\lambda) = S^{(2)}_{11}(\lambda)$, we have $u^{(1)} - u^{(2)} \in \mathcal B^{\ast}_0$ on $\mathcal M^{(1)}_{1}$. Then $u^{(1)} - u^{(2)} = 0$ by Rellich-Vekua's theorem (Theorem \ref{S2RellichTh1}).
Take $F \in L^2(\mathcal M^{(1)}_{1})= L^2(\mathcal M^{(2)}_{1})$ with support in $\mathcal O\subset \mathcal M^{(1)}_{1,reg}=\mathcal M^{(2)}_{1,reg}$.
Let 
$$
w^{(i)} =R^{(i)}(\lambda +i0)F.
$$
Then for any $a_1 \in L^2(M_1^{(1)})$ and $a^{(i)} = (a_1,0,\dots,0)$,
\begin{eqnarray*}
(\mathcal F^{(1,+)}(\lambda)F,a_1)_{L^2(M_1)} &=&(F,\mathcal F^{(1,+)}(\lambda)^*a_1)_{L^2(\mathcal M^{(1)})}\\
&=&(F,u^{(1)})_{L^2(\mathcal O)} \\
&=&(F,u^{(2)})_{L^2(\mathcal O)} \\
&=&(F,\mathcal F^{(2,+)}(\lambda)^* a_1)_{L^2(\mathcal M^{(2)})} \\
&=&(\mathcal F^{(2,+)}(\lambda) F,a_1)_{L^2(M_1)} .
\end{eqnarray*}
As this holds for all $a_1 \in L^2(M_1^{(1)})$, this implies
$$
\mathcal F^{(1,+)}(\lambda) F =\mathcal F^{(2,+)}(\lambda) F.
$$
This implies that the far fields of $w^{(1)}$ and $w^{(2)}$ coincide.
%
Let $v=w^{(1)}-w^{(2)}$. Then
\begin{equation}
 (- \Delta_{\mathcal M} - \lambda)v = 0, \quad {\rm in} \quad {\mathcal M}_1^{(1)} = \mathcal M_1^{(2)}, \\
\label{C2S1Delta-lambdav=0}
\end{equation}
and the far field of $v$ in the end $M^{(1)}=M^{(2)}$ is zero. 
Then $v=0$ by Theorem \ref{S2RellichTh1}.
This implies that 
 $w^{(1)}=w^{(2)}$ in $\mathcal O$, and hence
 $U_{\mathcal O, +}^{(1)}(\lambda)F=U_{\mathcal O,+}^{(2)}(\lambda)F.$
 
 Next assume that $0 < \beta_1 \leq 1/3$. Then, by the assumption (A-4-1), on $\mathcal M_1^{(1)} = \mathcal M_1^{(2)}$, the metric has the warped product form (\ref{S1Intro1/3warpedproductmetric}). Then, the equation (\ref{Th2.1.1Deltau1-u2=0}) can be solved by the separtion of the variable: $u^{(1)} - u^{(2)} = \sum_jv_j(r)\varphi_j(x)$, where $\varphi_j(x)$ is the normalized eigenvector of the Laplace-Beltrami operator on $M_1^{(1)} = M_1^{(2)}$. Each $v_j(r)$ is a $\mathcal B_0^{\ast}$-solution of the radial equation. Hence $v_j(r) = 0$, and $u^{(1)} - u^{(2)} = 0$. By the same argument, $v = 0$ for (\ref{C2S1Delta-lambdav=0}). This completes the proof of the theorem.
\end{proof}

We next consider cusp ends. Assume that (A-4-2) is satisfied for cusp ends $\mathcal M_j^{(1)}, \mathcal M_j^{(2)}$. Assume further that $\mathcal M_j^{(1)}$ and $\mathcal M_j^{(2)}$ are isometric. Take a bounded open set $\mathcal O \subset \mathcal M_j^{(1)} = \mathcal M_j^{(2)}$ and 
define $U_{\mathcal O,\pm}^{(i)}(\lambda)$, $i = 1, 2$, as in (\ref{C2S1Sourcetosolop}). Let $\widetilde{\mathcal S^{(i)}}(\lambda) = \big(\widetilde{\mathcal S^{(i)}}_{pq}(\lambda)\big)$
be the associated generalized S-matrix. 

\begin{theorem}\label{S15SmatrixdeterminesLW2cusp}
Assume that cusp ends $\mathcal M^{(1)}_{j}$ and $\mathcal M^{(2)}_{j}$ are isometric, and that  $\widetilde{\mathcal S^{(1)}}_{jj}(\lambda) = \widetilde{\mathcal S^{(2)}}_{jj}(\lambda)$ for some $\lambda \in \big(\sigma_e(H_1)\setminus\mathcal E_1\big)\cap\big(\sigma_e(H_2)\setminus \mathcal E_2\big)$. Then,  $U_{\mathcal O,\pm}^{(1)}(\lambda) = U_{\mathcal O,\pm}^{(2)}(\lambda)$.
\end{theorem}

\begin{proof}
Take $F \in L^2(\mathcal O)$ and extend it to be 0 outside $\mathcal O$. Let $u^{(i)} = (H_i - \lambda - i0)^{-1}F$. By using the notation in \S \ref{S15GeneSmatrix}, one can expand it as 
\begin{equation}
u^{(i)} = \sum_{\ell \geq 0}\left(a_{\ell,j}^{(i)}u_{\ell,j,+}(r) + b^{(i)}_{\ell,j,-}(r)\right)e_{\ell,j}.
\nonumber
\end{equation}
Since $(- \Delta_g - \lambda)(u^{(1)} - u^{(2)}) = 0$ on $\mathcal M_j^{(1)} = \mathcal M_j^{(2)}$, and $u^{(1)} - u^{(2)} \in \mathcal B^{\ast}$, we then have
\begin{equation}
a_{\ell,j}^{(1)} = a_{\ell,j}^{(2)}, \quad \forall \ell.
\nonumber
\end{equation}
Since the incoming data and the generalized S-matrices  coincide, by virtue of Theorem \ref{GeneralSmatrix}, we have
$u^{(1)} = u^{(2)}$ on $\mathcal M_j^{(1)} = \mathcal M_j^{(2)}$ and $b_{\ell,j}^{(1)} = b_{\ell,j}^{(2)}$. This completes the proof.
\end{proof}

The physical as well as mathematical legitimacy
of the source-to-solution operator is easily seen in the following observation.
Let us consider the wave equation in ${\mathbb R}^3$:
\begin{equation}
v_{tt} - \Delta v = F(t,x), 
\label{Timedepedentequation}
\end{equation}
with the initial condition
$$
v(0) = v_t(0) = 0.
$$
The Duhamel principle gives the following solution 
\begin{equation}
v(t,x) = \frac{1}{4\pi}\int_{|y-x| < t}
\frac{F(t - |y-x|,y)}{|y- x|}dy.
\label{Duhamenlformula}
\end{equation}
If $F(t,x) = f(x)e^{-i\sqrt{\lambda} t}$ $(\lambda > 0)$, it is rewritten as 
\begin{equation}
u(t,x) = e^{-i\sqrt{\lambda} t}\frac{1}{4\pi}\int_{|y-x| < t}
\frac{e^{i\sqrt{\lambda}|y-x|}f(y)}{|y-x|}dy.
\nonumber
\end{equation}
Therefore, as $t \to \infty$, 
\begin{equation}
v(t,x) \sim e^{-i\sqrt{\lambda} t}\big( - \Delta - \lambda - i0\big)^{-1}f.
\label{Eq:u(t,x)asttoinfty}
\end{equation}
This means that if we apply the time-harmonic oscillation in a medium, the  wave motion is asymptotically equal to the solution of the Helmholtz equation with time-periodic factor. This is a well-known physical phenomenon called the {\it limiting amplitude principle}\index{limiting amplitude principle}. 

This is also a general fact for self-adjoint operators with absolutely continuous spectrum. Namely, the solution of the abstract wave equation
\begin{equation}
v_{tt}+ Hv = e^{-i\sqrt{\lambda} t}f, \quad u(0) = u_t(0) = 0
\nonumber
\end{equation}
behaves like 
$$
v(t) \sim e^{-i\sqrt{\lambda} t}\big(H - \lambda - i0)^{-1}f, \quad  t \to \infty,
$$
if the limiting absorption principle, i.e. the existence of the limit $\big(H - \lambda \mp i0)^{-1}$ and its H{\"o}lder continuity with respect to $\lambda$ are guaranteed. Thus, the  source-to-solution operator is the observation of the stationary wave by the time-periodic input. See e.g. 
\cite{Eidus}
and the references therein.

One can also consider the {\it time-dependent source-to-solution operator}\index{time-dependent source-to-solution operator}. For a bounded domain $\mathcal O \subset \mathcal M$, consider the wave equation
\begin{equation}
 v_{tt} +  Hv = F(t), \quad (t, x) \in {\mathbb R}\times \mathcal M,
\label{AbstsractWaveEqC2S1}
\end{equation}
where $F(t) = F(t,x)$ is assumed to be compactly supported in $(0,\infty)\times\mathcal M$. Then, the solution of this equation satisfying $v(t)=0$  for $t<0$
exists uniquely, which is denoted by $v_{\mathcal O,+}^F(t)$. We fix $I = (0, T)$ arbitrarily, and consider the operator
\begin{equation}
V_{\mathcal O,+}(T) : L^2(I\times\mathcal O) \ni F(t,x) \to v^F_{\mathcal O,+}(t,x) \in L^2(I\times\mathcal O).
\nonumber
\end{equation}
We call it an outgoing time-dependent source-to-solution operator. By the time reversal, one can also define an incoming time-dependent source-to-solution operator $V_{\mathcal O,-}(T)$.
Let $R(z) = (- \Delta_{\mathcal M} - z)^{-1}$ and $r_{\mathcal O}$ the operator of restriction to $\mathcal O$.

\begin{lemma}
\label{Lemma:rRlambda+i0rdeterminedrf(H)r}
Let $S$ be any set of positive measure in $\sigma_e(H)\setminus{\mathcal E}$. Then, for any $f(\lambda) \in C({\mathbb R})$ such that $f(\lambda) \to 0$ as $\lambda \to \infty$, the knowledge of $r_{\mathcal O}R(\lambda + i0)r_{\mathcal O}$ for all $\lambda \in S$ determines $r_{\mathcal O}f(H)r_{\mathcal O}$.
\end{lemma} 

\begin{proof}
Since $r_{\mathcal O}R(\lambda + i0)r_{\mathcal O}$ is the boundary value of an analytic function $r_{\mathcal O}R(z)r_{\mathcal O}$,  by Fatou's theorem, the knowledge of $r_{\mathcal O}R(\lambda + i0)r_{\mathcal O}$ for all $\lambda \in S$ determines $r_{\mathcal O}R(z)r_{\mathcal O}$ for all $z \in {\mathbb C}\setminus{\mathbb R}$. To prove the lemma, we have only to consider the case in which $f(\lambda) \in C_0^{\infty}({\mathbb R})$.  Then, letting $F(\lambda)$ be an almost analytic extension of $f(\lambda)$, we have the representaion formula (\ref{HelfferSjostrandFormula}). The lemma then readily follows. 
\end{proof}

\begin{lemma}
\label{Lemma:StationarySSopdeterminestimedependrntSS}
Take any subset $S$ of positive measure in $\sigma_e(H)$. 
Given a relatively compact open set $\mathcal O$ and a constant $T > 0$ arbitrarily, the knowledge of
$U_{\mathcal O,+}(\lambda)$ for all $\lambda \in S$ determines $V_{\mathcal O,+}(T)$ uniquely.
\end{lemma}

\begin{proof}
Let $\mathcal U(t) = H^{-1/2}\sin(tH^{1/2})$. Then, by Duhamel's principle, the solution of the equation (\ref{AbstsractWaveEqC2S1}) satisfying $v(0) = v_t(0) = 0$ is written as 
\begin{equation}
v(t) = \int_0^t\mathcal U(t-s)F(s)ds.
\label{C2S1DuhamelFormula}
\end{equation}
The Lemma then follows from Lemma \ref{Lemma:rRlambda+i0rdeterminedrf(H)r}.
\end{proof}

As is seen from the proof, Lemma \ref{Lemma:StationarySSopdeterminestimedependrntSS}  holds for 4 choices of the mapping
$$
U_{\mathcal O,\pm}(\lambda) \to V_{\mathcal O,\pm}(T), \quad U_{\mathcal O,\pm}(\lambda) \to V_{\mathcal O,\mp}(T).
$$

The converse of Lemma \ref{Lemma:StationarySSopdeterminestimedependrntSS} is also true.

\begin{lemma}
Given a relatively compact open set $\mathcal O$, the knowledge of $V_{\mathcal O,+}(T)$ for all $T > 0$ determines 
$U_{\mathcal O,+}(\lambda)$ for all $\lambda \in \sigma_e(H)\setminus{\mathcal E}$ uniquely.
\end{lemma}

\begin{proof}
In (\ref{C2S1DuhamelFormula}), take $F(t) = \varphi(t)f(x)$, where $\varphi(t) \in C_0^{\infty}((0,\infty))$ and $f \in L^2(\mathcal O)$ extended to be 0 outside $\mathcal O$. Let $v(t)$ be given in \eqref{C2S1DuhamelFormula}.
We put for $\epsilon > 0$
\begin{equation}
\widehat v_{\epsilon}(k) = \int_{-\infty}^{\infty}e^{i(k+i\epsilon)t}v(t)dt.
\label{C2S1Definewidehat vepsiklonk}
\end{equation}
Since $v(t) = 0$ for $t < 0$, the integral is convergent. Then, $\widehat v_{\epsilon}(k)$ satisfies
$$
\big(H - (k + i\epsilon)^2\big)\widehat v_{\epsilon}(k) = -
\widehat \varphi(k)f.
$$
One can show 
\begin{equation}
\|\widehat v_{\epsilon}(k)\|_{L^2(\mathcal M)} \leq \frac{C}{\epsilon^2}\|f\|.
\label{C2S1widehatvepsilonkestimate}
\end{equation}
In fact, noting that
$$
v(t) = \int_0^t\frac{\sin((t-s)\sqrt{H}}{(t-s)\sqrt{H}}(t-s)F(s)ds,
$$
we have
$$
\|v(t)\| \leq C\int_0^t(t-s)|\varphi(s)|ds\|f\| \leq C(1 + t)^2\|f\|.
$$
Using (\ref{C2S1Definewidehat vepsiklonk}), we obtain (\ref{C2S1widehatvepsilonkestimate}). 
Then, we have
$$
\widehat v_{\epsilon}(k) = -\widehat\varphi(k)R((k + i\epsilon)^2)f.
$$
Letting $\epsilon \to 0$, we have computed 
$\widehat\varphi(k)R(k^2 + i0)f$ from $v(t)$. Since $\widehat\varphi(k)$ is analytic, its zeros are discrete on ${\mathbb R}$. We can then compute $R(\lambda + i0)f$ for $\lambda \in \sigma_e(H)\setminus{\mathcal E}$ except for some discrete points. By using analytic continuation with respect to $\lambda \in {\mathbb C}_+$, we can find  $R(\lambda + i0)f$
 for all $\lambda \in \sigma_e(H)\setminus{\mathcal E}$. This proves the lemma.
\end{proof}

 \def\btext{}

\def \ba {\begin {eqnarray*} }
\def \ea {\end {eqnarray*} }
\def \beq {\begin {eqnarray}}    
\def \eeq {\end {eqnarray}}         
\def \supp {\hbox{supp}\,}
\def \diam {\hbox{diam }}
\def \rad {\hbox{rad}}
\def \ind {\hbox{Ind}\,}
\def \dist {\hbox{dist}}
\def \diag{\hbox{diag }}
\def \det {\hbox{det}}
\def \bra{\langle}
\def \cet{\rangle}
\def \e {\varepsilon}
\def \p {\partial}
\def \a {\alpha}
\def \la{{\lambda}}
\def \Im{{\rm Im}\,}
\def \Re{{\rm Re}\,}
\def \im{{\rm Im}\,}
\def \re{{\rm Re}\,}
\def\R{\mathbb R}
\def\N{\mathbb N}
\def\C{\mathbb C}
\def\Z{\mathbb Z}
\def\O{{\mathcal O}}
\def\M{{\mathcal M}}
\def\G{{\mathcal G}}
\def\U{{\mathcal U}}
\newcommand{\nntext}{}
\newcommand{\ntekst}{}
\newcommand{\nntekst}{}

\def \ba {\begin {eqnarray*} }
\def \ea {\end {eqnarray*} }
\def \beq {\begin {eqnarray}}    
\def \eeq {\end {eqnarray}}         
\def \supp {\hbox{supp}\,}
\def \diam {\hbox{diam }}
\def \rad {\hbox{rad}}
\def \ind {\hbox{Ind}\,}
\def \dist {\hbox{dist}}
\def \diag{\hbox{diag }}
\def \det {\hbox{det}}
\def \bra{\langle}
\def \cet{\rangle}
\def \e {\varepsilon}
\def \p {\partial}
\def \a {\alpha}
\def \b {\beta}
\def \la{{\lambda}}
\def \Im{{\rm Im}\,}
\def \Re{{\rm Re}\,}
\def \im{{\rm Im}\,}
\def \re{{\rm Re}\,}
\def\R{\mathbb R}
\def\N{\mathbb N}
\def\C{\mathbb C}
\def\Z{\mathbb Z}
\def\O{{\mathcal O}}
\def\M{{\mathcal M}}
\def\G{{\mathcal G}}
\def\U{{\mathcal U}}


\section{Definitions}

\subsection{Metric tensor}

Let ${\mathcal M}$  be a  CMGA, that can be considered as an orbifold with $C^\infty$-smooth coordinates and a non-smooth metric on it.
{\mtext Before we define the distance functions and consider the finite speed of wave propagation, we recall some properties of local coordinates on CMGA.}
%



Let us consider conic coordinates
$(\tilde U,\Gamma, \pi)$ and the set
$U= \pi( \tilde U)\subset \mathcal M$. Here  $U\subset \mathcal M$ and $\tilde U\subset \R^n$.
Recall that
 $\tilde \Phi:\Gamma\backslash \tilde U\to U$ is a homeomorphism. In $\tilde U$  we consider coordinates $(y,z)$,
where $y=Y(x)$, $Y:\tilde U\to \tilde W$  and
 $z=Z(x)$, $Z=\tilde U\to \tilde V$, that is, for $x\in U$
 $$
x=(Y(x),Z(x))\in \tilde W\times \tilde V\subset  {\mathbb R}^{k}\times  {\mathbb R}^{n-k}.
 $$ 
 The action of elements $\gamma\in \Gamma$
 is such that $Y(\gamma x)=Y(x)$, $x\in \tilde U$, that is, $\gamma$ keeps 
the $Y$-coordinate  coordinate invariant.
   In
 $B^{n-k}(0,R_0)\setminus \{0\}$ we also use spherical coordinates
 $Z(x)=r\omega$, $r=r(x)$, $\omega=\omega(x)$ such that
$r\in (0,R_0)$ and $\omega\in S^{n-k-1}$ is a unit vector. 
We assume that $r(\gamma x)=r(x)$ for all  $\gamma\in \Gamma$
and moreover, the operators $\omega_*\gamma:S^{n-k-1}\to S^{n-k-1}$,
where
$(\omega_*\gamma)(\omega(x))=\omega(\gamma(x))$, satisfy $\omega_*\gamma\in SO(n-k-1)$.

We assume that on $$\tilde U^{reg}=\tilde W\times (B^{n-k}(0,R_0)\setminus \{0\}),\quad 
\tilde W\subset 
\R^{k},\ \
 B^{n-k}(0,R_0)\subset 
{\mathbb R}^{n-k},
$$
we have a  $C^\infty$-smooth  metric tensor $\tilde g$, that  in the coordinates $(y,r,\omega)$ has the form
\beq\label{radial coordinates}
& &\tilde g=dr^2+\sum_{j,k=1}^k a_{jk}(y,r,\omega)dy^jdy^k+
\sum_{\a,\b=1}^{n-k-1} r^2 \,b_{\a \b}(y,r,\omega)d\omega^\a d\omega^\b+\\
& &+ \nonumber
\sum_{j=1}^k 
\sum_{\a=1}^{n-k-1} r\,h_{j \a}(y,r,\omega)dy^j d\omega^\a.
\eeq
We denote by $ g_{S^{n-k-1}}$ the standard metric of ${S^{n-k-1}}$.
We assume that 
\beq\label{conditions 1}
& &a_{jk}(y,r,\omega)\to \hat a_{jk}(y),\\ \nonumber
& &b_{\a \b}(y,r,\omega)d\omega^\a d\omega^\b\to \hat b_{\a \b}(y,\omega)d\omega^\a d\omega^\b,\\ \nonumber
& &h_{j\a}(y,r,\omega)dy^j d\omega^\a\to 0,
\eeq
uniformly in $(y,\omega)$, in compact subsets of  $\tilde W\times S^{n-k-1}$, as $r\to 0$. Moreover,
we assume that there are $c_0,c_1>0$ such that
\beq\label{conditions 2}
c_0  g_{S^{n-k-1}} \leq \hat b_{\a \b}(y,\omega)d\omega^\a d\omega^\b\leq   T(y)^2 g_{S^{n-k-1}}
\eeq
and $0<T(y)\leq c_1$. 
Also, we assume that
  the metric $\tilde g$ is $\Gamma$-invariant, that is, $\gamma_*\tilde g=\tilde g$
  on $\tilde U^{reg}$ for all $\gamma\in \Gamma$.

\subsection{Distance function}

We have that for any $R\in (0,R_0)$ there are $c_0(R)$, $c_1(R)>0$ such
 that for all $x\in \tilde W\times (B^{n-k}(0,R)\setminus \{0\})\subset {\mathbb R}^n$  the metric tensor 
 $\tilde g=\tilde g_{jk}(x)dx^jdx^k$ satisfies
 \beq\label{tilde g bounded A}
 c_0(R)\,g_E\leq \tilde g\leq c_1(R) g_E,
 \eeq
where $g_E$ is the Euclidean metric on ${\mathbb R}^k\times {\mathbb R}^{n-k}$.
However, in the 
Euclidean coordinates $x=(y,z)$ the map $x\mapsto \tilde g(x)$ is not Lipschitz
due to its behaviour near $z=0$.  This is due to the fact that the radial projection map $z\to P(z)$,
that is the matrix valued map that  has at $z$ the value
$$
P(z)=\bigg [\frac {z^jz^k}{|z|^2}\bigg]_{j,k=1}^{n-k},\quad P:{\mathbb R}^{n-k}\setminus \{0\}\to {\mathbb R}^{(n-k)\times (n-k)},
$$
is not Lipschitz on the set ${\mathbb R}^{n-k}\setminus \{0\}$ having the distance function induced from ${\mathbb R}^{n-k}$.

We assume  that there are conic coordinates $(\tilde U_{\ell},\Gamma^{(\ell)}, \pi^{(\ell)})$
such that  the sets $U_{\ell}= \pi^{(\ell)}( \tilde U_{\ell})\subset \mathcal M$,
 $\ell=1,2,\dots$ are an open covering of $\mathcal M$.
 We recall that $ \pi^{(\ell)}= \tilde \Phi^{(\ell)}\circ \tilde  \pi^{(\ell)}$ where
 $\tilde \Phi^{(\ell)}: \Gamma^{(\ell)}\backslash \tilde U_{\ell}\to U_{\ell}$ are homeomorphisms
 and $ \pi^{(\ell)}: \Gamma^{(\ell)}\backslash \tilde U_{\ell}^{reg}\to U_{\ell}^{reg}=U_{\ell}\cap \mathcal M^{reg}$
 are $C^\infty$-smooth. The sets $\pi^{(\ell)}(\tilde V)$, where $\tilde V \subset\tilde U_{\ell}$ is open,
 form a basis of the topology of ${\mathcal M}$.
 
We define metric tensors $\tilde g^{(\ell)}$ on $ \tilde U_{\ell}$ and
$$
g^{(\ell)}=({\pi}^{(\ell)})_*\tilde g^{(\ell)}\quad\hbox{ on }U_{\ell}^{reg}.
$$
We assume that for indexes $\ell$ and $\ell'$, we have that on the sets 
 $$
 U^{reg}_{(\ell)}\cap  U^{reg}_{\ell'}=
({\\pi}^{{(\ell)}}( \tilde U^{reg}_{\ell}))\cap ({\pi}^{{(\ell')}}( \tilde U^{reg}_{\ell'})) 
$$ 
the metric tensors  $g^{(\ell)}=({\pi}^{{(\ell)}})_*\tilde g^{(\ell)}$
 and $g^{(\ell')}=({\pi}^{{(\ell')}})_*\tilde g^{(\ell')}$
 coincide.

These
metric tensors  $g^{(\ell)}$
define a smooth metric on regular part  ${\mathcal M}^{reg}$ of ${\mathcal M}$ and we denote this metric by $g$.
We say that $\mu:[0,1]\to {\mathcal M}$  that
is piecewise $C^1$-smooth on the lifted local coordinates, if for all $[s_1,s_2]\subset [0,1]$
such that $\mu([s_1,s_2])\subset \tilde U_{\ell}$ there is a piecewise $C^1$-smooth  path
 $\tilde \mu:[s_1,s_2]\to \tilde U_{\ell}$ such that
\beq
\pi^{(\ell)}(\tilde \mu(s))= \mu(s),\quad\hbox{for }s\in [s_1,s_2].
\eeq

As ${\mathcal M}$ is the topological closure of  ${\mathcal M}^{reg}$, we define the distance on ${\mathcal M}$  by
\beq\label{extension of metric}
d_{\mathcal M}(x,y)=\inf_{\mu} \hbox{Length}_g(\mu([0,1])\cap {\mathcal M}^{reg}),
\eeq
where infimum is taken over paths $\mu:[0,1]\to {\mathcal M}$  that
are  piecewise $C^1$-smooth on the lifted local coordinates,
and satisfy $\mu([0,1)])\cap  {\mathcal M}^{sing}$ is a finite set, and 
$\mu(0)=x$, $\mu(1)=y$. Note that
$\mu([0,1])\cap {\mathcal M}^{reg,{(\ell)}}$ is rectifiable in all local coordinate charts $U^{reg,{(\ell)}}$. Also, 
when $\mu([s_1,s_2])\subset U_\ell^{reg}\subset {\mathcal M}^{reg}$,
  the length {is defined} using local coordinates as
\beq\label{path length}
\hbox{Length}_g(\mu([s_1,s_2]))=\int_{s_1}^{s_2}(g_{jk}(\mu(s))\dot\mu^j(s)\dot\mu^k(s))^{1/2}ds.
\eeq

Let $m_0\in \Z_+$ be such that $\frac 1{m_0}<R_0$, and let below $m\geq m_0$.
%
%
  
  \subsection{Approximation of $g$ by smooth metric tensors}
  
For constructions below, let us define an auxiliary metric that is smooth on lifted coordinates. Let $\tilde g^{(\ell),smooth}$  be a $C^\infty$-smooth metric tensor
defined on $\tilde U_{\ell}$ such that $\tilde g^{smooth,(\ell)}\geq 
\tilde g^{(\ell)}$ on $\tilde U_{\ell}^{reg}$, that means that the positive definite matrices satisfy $$(\tilde g^{smooth,(\ell)}_{jk})_{j,k=1}^n\geq 
(\tilde g^{(\ell)}_{jk})_{j,k=1}^n\hbox{ on $\tilde U_{\ell}^{reg}$,}$$
and
$\gamma_*\tilde g^{smooth,(\ell)}=\tilde g^{smooth,(\ell)}$  for all $\gamma\in \Gamma^{(\ell)}$. Moreover,
 \beq\label{tilde g bounded}
 c_0(R)\,g_E\leq \tilde g^{smooth,(\ell)}\leq 2c_1(R) g_E,
 \eeq
where $g_E$ is the Euclidean metric on $\R^k\times \R^{n-k}$.
The metric tensors $\tilde g^{smooth,(\ell)}$  define an ``orbifold metric''
on neighborhoods $U_{\ell}$ that we denote by
 $g^{smooth,(\ell)}=(\Phi^{{(\ell)}})_*\tilde g^{smooth,(\ell)}$.



%
%


We define a metric  $g^{smooth}$ by summing metric tensors $g^{smooth,\ell}$
together using a locally finite partition of unity $\phi_\ell:\mathcal M\to \R$,
such that $\supp(\phi_\ell)\subset U_{\ell}$,  $\sum_\ell \phi_\ell(x)=1$,
$\phi_\ell\geq 0$, 
and that there are  functions
$\tilde\phi_\ell\in C^\infty_0(\tilde U_{\ell})$ such that $\phi_\ell(\pi^{(\ell)}(x))=\tilde \phi_\ell(x)$.
Using such partition of unity we define
\ba
g^{smooth}(x)=\sum_{\ell=1}^L \phi_\ell(x)((\pi^{(\ell)})_*)g^{smooth,\ell}(x),\quad x\in {\mathcal M}^{reg}. 
\ea
 Strictly speaking, this metric tensor is defined only on ${\mathcal M}^{reg}$.
When
$\mu([s_1,s_2])\subset  U_{\ell}\subset {\mathcal M}$ and
 there is a piecewise $C^1$-smooth (on the lifted local coordinates)  path
 $\tilde \mu:[s_1,s_2]\to \tilde U_{\ell}$ such that
 $\pi^{(\ell)}(\tilde \mu(s))= \mu(s)$  for $s\in [s_1,s_2]$,
we define the length  of the path $\mu([s_1,s_2])$ with respect to the metric ${g^{smooth}}$ using local coordinates as
\beq\label{path length g smooth}
\hbox{Length}_{g^{smooth}}(\mu([s_1,s_2]))=\int_{s_1}^{s_2}(\tilde g^{(\ell)}_{jk}(\mu(s))\p_s \tilde 
\mu^j(s)\p_s \tilde \mu^k(s))^{1/2}ds.
\eeq
Decomposing a  path $\mu:[0,1]\to  {\mathcal M}$, that is  piecewise $C^1$-smooth on the lifted local coordinates, to a union of paths which
are all supported on some set $U_{\ell}$, we can define the length of 
arbitrary  path  $\mu:[0,1]\to  {\mathcal M}$ that is  piecewise $C^1$-smooth on the lifted local coordinates. Moreover,
we define
\beq\label{extension of metric}
d_{g^{smooth}}(x,y)=\inf_{\mu} \hbox{Length}_{g^{smooth}}(\mu([0,1]),
\eeq
where infimum is taken over paths $\mu:[0,1]\to {\mathcal M}$  that
are  piecewise $C^1$-smooth on the lifted local coordinates,
and satisfy 
$\mu(0)=x$, $\mu(1)=y$.
%
%
%
We say that a curve $\mu$   
is a geodesic  of the metric $g^{smooth}$ if it  is locally distance minimising.

Let 
$$K_{m}=\{x\in {\mathcal M};\ d_{\mathcal M}(x,{\mathcal M}^{sing})< \frac 1m\}.$$ 
Next
we modify the non-Lipschitz metric $g$  by defining a  metric $g^{(m)}$ that is smooth on
the lifted coordinates. We define that $g^{(m)}$ is in the
 set   ${\mathcal M}\setminus K_m$  equal to $g$ and in the set $K_m$,
 $$
 g_{jk}^{(m)}(x)=\psi_m(x) g^{smooth}_{jk}(x)+(1-\psi_m(x))g_{jk}(x),
 $$
 where $\psi_m\in C^\infty(K_m)$  is  equal to 1 in an open neighbourhood of
  ${\mathcal M}^{sing}$ and $0\le \psi_m(x)\le 1$. Then $g^{(m)}$ is $C^\infty$-smooth orbifold metric in ${\mathcal M}$.
We assume that  $\psi_{m+1}(x)\le \psi_m(x)$ so that
$g_{jk}^{(m+1)}(x)\leq g_{jk}^{(m)}(x)$.

Note that then $g^{(m)}\geq g$ on ${\mathcal M}^{reg}$ and using this  we will see that the
waves are slower when 
they propagate following the metric  $g^{(m)}$ than $g$.
Also,  $g^{(m)}\geq g^{(m+1)}$. In particular, the travel time of waves between points $x$ and $y$
with respect to $g^{(m)}$ is longer 
that the travel time of waves between points $x$ and $y$
with respect to $g^{(m+1)}$, i.e., $d_{g^{(m)}}(x,y) \ge d_{g^{(m+1)}}(x,y).$ 


 Below, we denote 
\ba
& &V_{m}={\mathcal M}\setminus K_m.
\ea

{\btext
 \begin{lemma}\label{deformation}
 Let $x,y\in {\mathcal M}$,  
$x\not=y$, $\e>0$ and let $\gamma:[0,1]\to {\mathcal M}$  be a
piecewise $C^1$-smooth path on the lifted local coordinates
 that connects $x=\gamma(0)$ and $y=\gamma(1)$.
 Then there is a path $\mu([0,1])$ that
 connects $x$ and $y$ so that 
   $\mu\cap {\mathcal M}^{sing}$  is a finite set, 
 and
$$
|\hbox{Length}_{g^{(m)}}(\mu([0,1])-
\hbox{Length}_{g^{(m)}}(\gamma ([0,1])|
\leq \e.
$$
 \end{lemma}
\noindent {\bf Proof.} 
 We assume that $\gamma$ is
parametrised so that it has constant speed 
with respect to $g^{(m)}$.


For let $s_j\in I$, $j=1,2,\dots,J$ be such that $s_1=0<s_2<s_3<\dots<s_J=1$ and  $[s_j,s_{j+1}]\subset I$ are intervals  such that  $\gamma(s_i)\in \mathcal M^{reg}$  for $i=2,3,\dots,J-1$,
and 
there are  projections from covering neighbourhoods ${\pi}_\ell :\tilde U_\ell \to U_\ell $, $\ell=\ell(j)$ such that
$\gamma([s_j,{s_{j+1}}])\subset U_{\ell(j)}$ for $j=1,2\dots,J-1$. Let
$\tilde \gamma_j:[s_j,{s_{j+1}}]\to \tilde U_{\ell(j)}$, $j=1,2\dots,J-1$ be
paths such that ${\pi}_{\ell (j)}(\tilde \gamma_j(s))=\gamma(s)$,
that is, $\tilde \gamma_j:[s_j,{s_{j+1}}]\to \tilde U_{\ell(j)}$  is the lift of  the path
$\gamma:[s_j,{s_{j+1}}]\to U_{\ell(j)}$.

Let 
\ba
& &\rho_2:\tilde U_{\ell(j)}=\tilde W_{\ell(j)}\times \tilde V_{\ell(j)}\to  \tilde V_{\ell(j)}\subset {\mathbb R}^{n-k},\\
& &\rho_2(y,z)=z
\ea
be the projection of $x$ to the $z$-coordinate.
Note that $n-k\geq 2$. Since $\tilde \gamma_j$ is a rectifiable path, we see that $\rho_2(\tilde \gamma_j([s_j,{s_{j+1}}]))$ has the Hausdorff dimension 1. Hence,
we see that {for any $i\in\Z_+$} there are vectors $v_{j,i}\in \R^{n-k},$    such that $|v_{j,i}|<1/i$  and that $-v_{j,i}\not \in \rho_2(\tilde \gamma_j([s_j,{s_{j+1}}]))$.
This implies that
$$
0\not \in v_{j,i}+\rho_2(\tilde \gamma_j([s_j,{s_{j+1}}]))
$$
and for $\tilde v_{j,i}=(0,v_{j,i})\in \R^n$
$$
(\tilde  v_{j,i}+\tilde \gamma_j([s_j,{s_{j+1}}]))\cap (\tilde W\times\{0\})=\emptyset.
$$
Denote 
$$
\eta_{j,i}:[s_j,{s_{j+1}}]\to {\mathcal M},\quad \eta_{j,i}(s)=\pi_{\ell (j)}(\tilde \eta_{j,i}(s)),\quad
\tilde \eta_{j,i}(s)=
\tilde  v_{j,i}+\tilde \gamma_j(s).
$$

For $j=1,2,\dots,J-1$, let $\tilde \alpha_{j,i}:[0,1]\to U_{\ell(j)}$  be the Euclidean line segment
$\tilde \alpha_{j,i}(t)=\tilde  v_{j,i}t+\tilde \gamma_j(s_j)$ that connects $\tilde \gamma_j(s_j)$ to 
$\tilde \eta_j(s_j).$ 
Similarly, let $\tilde \beta_{j,i}:[0,1]\to U_{\ell(j)}$  be the Euclidean line segment
that connects   
$\tilde \eta_j(s_{j+1})$ to $\tilde \gamma_j(s_{j+1})$.
Observe that when $i$ is large enough,
the line segments $\tilde \alpha_{j,i}([0,1])$  and  $\tilde \beta_{j,i}([0,1])$ are subsets of  $\tilde U_{\ell(j)}$ and intersect $\tilde W\times\{0\}$ {at most at the points $\tilde \alpha_{1,i}(0)$ or $\tilde \beta_{J-1,i}(1)$.}

Let  $\alpha_{j,i}(t)=\pi_{\ell(j)}(\tilde \alpha_{j,i}(t))$ and $\beta_{j,i}(t)=\pi_{\ell(j)}(\tilde \beta_{j,i}(t))$.
Then  $ \alpha_{j,i}([0,1])$  and  $\beta_{j,i}([0,1])$ intersect $\mathcal M^{sing}$ 
{at most at the points $\alpha_{1,i}(0)=\gamma(0)$ or $\beta_{J-1,i}(1)=\gamma(1)$.}

%
%
 
Let now $\mu_i:[0,1]\to {\mathcal M}$ be a path that is obtained by concatenating the paths
$\alpha_{1,i},\eta_{1,i},\beta_{1,i},\alpha_{2,i},\eta_{2,i},\beta_{2,i},\dots,\alpha_{J-1,i},\eta_{J-1,i},\beta_{J-1,i}$ 
that connects $\gamma(0)$ to $\gamma(1)$, that is,
\ba
\mu_i&=&
\bigcup_{j=1}^{J-1} (\alpha_{j,i}\cup \eta_{j,i}\cup\beta_{j,i}).
\ea
As the  metric tensor $g^{(m)}$  is a {smooth orbifold metric} in $\mathcal M$ 
we see that 
\ba
\lim_{i\to \infty}
\hbox{Length}_{g^{(m)}}(\mu_i(0,1]))=
\hbox{Length}_{g^{(m)}}(\gamma(0,1])).
\ea
By choosing $\mu=\mu_i$ with sufficiently large $i$, we prove the claim.
%


\hfill$\square$\medskip

%
%
%

 Let us consider the set ${\mathcal M}={\mathcal M}^{reg}\cup {\mathcal M}^{sing}$  as an orbifold with
a smooth metric  $g^{(m)}$  defined on the lifted coordinate neighbourhoods.
 
 \begin{lemma}\label{two distance functions}
 For $x,y\in {\mathcal M}$  we have
 $$
 d_{\mathcal M}(x,y)\leq d_{g^{(m)}}(x,y).
 $$ 
 \end{lemma}

\noindent {\bf Proof.}  Let $\gamma$ be a {$g^{(m)}-$length} minimizing curve
on $\mathcal M$ that connects $x$  to $y$.
By Lemma \ref{deformation} there is a path $\mu([0,1])$ that
 connects $x$ and $y$, has the property that
   $\mu\cap {\mathcal M}^{sing}$  is a finite set, 
 and
$$
|\hbox{Length}_{g^{(m)}}(\mu([0,1])-
d_{g^{(m)}}(x,y)|
\leq \e.
$$
Then
%
\ba
d_{g^{(m)}}(x,y)+\e&\ge &\hbox{Length}_{g^{(m)}}(\mu)\\
&\ge & \hbox{Length}_{g^{(m)}}(\mu\cap 
 {\mathcal M}^{reg})\\
 &\ge& \hbox{Length}_{g}(\mu\cap 
 {\mathcal M}^{reg})\\
 &\geq& d_{\mathcal M}(x,y).
\ea
As $\e>0$ is here arbitrary, the claim follows.
\hfill$\square$\medskip
}

\begin{lemma}\label{two distance functions opposite inequality}
 For $x,y\in {\mathcal M}$  we have
 $$
 d_{\mathcal M}(x,y)\geq \lim_{m\to \infty} d_{g^{(m)}}(x,y).
 $$ 
 \end{lemma}

\noindent {\bf Proof.} 
By above definitions,  
\beq\label{comparison of metrics}
d_{\mathcal M}(x,y)&=&\inf_{\gamma} \hbox{Length}_g(\gamma\cap {\mathcal M}^{reg}),
\eeq
where infimum is taken over {piecewise $C^1$-smooth paths $\gamma:[0,1]\to {\mathcal M}$ on
the lifted coordinates}
 such that
$\gamma(0)=x$, $\gamma(1)=y$, and $\gamma([0,1])\cap {\mathcal M}^{sing}$  is a finite set.

Let $\e>0$ and choose  $\gamma:[0,1]\to {\mathcal M}$  such that
$\gamma(0)=x$, $\gamma(1)=y$, and $\gamma([0,1])\cap {\mathcal M}^{sing}$  is a finite set and
\ba
\hbox{Length}_g(\gamma\cap {\mathcal M}^{reg})\leq d_{\mathcal M}(x,y)+\e.
\ea 
As the path $\gamma$ intersects ${\mathcal M}^{sing}$  finitely many times, 
 by definition there are $t_1,t_2,\dots,t_J$, $J\in \Z_+$, such that for $t\in [0,1]$ we have
 $$\gamma(t)\in {\mathcal M}^{sing}\hbox{ if and only if }t\in \{t_1,t_2,\dots,t_J\}.$$
Then for all $h>0$ there is $m_0(h,\gamma)$ such that for $m>m_0(h,\gamma)$ we have
\ba
\dist_{\R}(t,\{t_1,t_2,\dots,t_J\})\geq h\implies \gamma(t)\in {\mathcal M}\setminus K_m.
\ea
Indeed, if no such $m_0(h,\gamma)$  exists, then for all $m$ there is $t'_m\in [0,1]$
such that $d(t'_m,\{t_1,t_2,\dots,t_J\})\geq h$  and $\gamma(t'_m)\in K_m.$
By compactness, there is $t'$ such that $\lim_{k\to \infty} t'_{m_k}=t'$.
Then $\dist(t',\{t_1,t_2,\dots,t_J\})\geq h$. Also, for any $m_1>m_0(h,\gamma)$  we see that if $m>m_1$
then  $\gamma(t'_m)\in K_{m_1}.$ Hence $\gamma(t')\in K_{m_1}.$ As $m_1$ is arbitrary,
we see that   $\gamma(t')\in \bigcap_{m_1\geq 1}K_{m_1}={\mathcal M}^{sing}$ which is a contradiction as $\dist(t',\{t_1,t_2,\dots,t_J\})\geq h$. 

We see that there is $C_2=C_2(\gamma)$  such that
\ba
\sum_{j=1}^J \hbox{Length}_{g^{m}}(\gamma((t_j-h,t_j+h)\cap [0,1]))\leq C_2Jh
\ea
Hence, for $m>m_0(h,\gamma)$ we have
\beq\nonumber
d_{\mathcal M}(x,y)+\e&\geq&\hbox{Length}_g(\gamma\cap {\mathcal M}^{reg})\\ 
\label{E1}
&\geq&\hbox{Length}_g(\gamma([0,1]\setminus \cup_{j=1}^J (t_j-h,t_j+h)))\\ \nonumber
&\geq&\hbox{Length}_{g^{m}}(\gamma([0,1]\setminus \cup_{j=1}^J (t_j-h,t_j+h)))\\ \nonumber
& &+\bigg((\hbox{Length}_{g^{m}}(\gamma((t_j-h,t_j+h)\cap [0,1])))-C_2Jh\bigg)\\ \nonumber
&\geq&\hbox{Length}_{g^{m}}(\gamma([0,1]))-C_2Jh.
\eeq
Thus,
\ba
d_{\mathcal M}(x,y)+\e&\geq&\lim_{m\to \infty}\hbox{Length}_{g^{m}}(\gamma([0,1]))-C_2Jh.
\ea
As $h>0$ is here arbitrary, we have
\ba
d_{\mathcal M}(x,y)+\e&\geq&\lim_{m\to \infty}\hbox{Length}_{g^{m}}(\gamma([0,1]))\\
&\geq&\lim_{m\to \infty}d_{g^{m}}(x,y)
\ea
As $\e>0$ is arbitrary, we obtain the claim.
\hfill$\square$\medskip

For $x,y\in {\mathcal M}^{reg}$, define
 \beq
 d_{\mathcal M^{reg}}(x,y)=\inf_\gamma \hbox{Length}_{g}(\gamma)
\eeq
where the infimum is taken over rectifiable curves $\gamma : [0,1] \to {\mathcal M}^{reg}$
such that $\gamma(0)=x$ and $\gamma(1)=y$.

\begin{lemma}\label{two distance M and Mreg}
 For $x,y\in {\mathcal M}^{reg}$  we have
 $$
 d_{\mathcal M^{reg}}(x,y)= d_{\mathcal M}(x,y).
  $$ 
 \end{lemma}

\noindent {\bf Proof.} Assume that $x\not=y$, $\e>0$ and let $\gamma:[0,1]\to {\mathcal M}$  be a
piecewise $C^1$-smooth path 
 that connects $x=\gamma(0)$ and $y=\gamma(1)$, has the property
 that  $\gamma\cap {\mathcal M}^{sing}$  is a finite set, 
 and
$$
\hbox{Length}_g(\gamma([0,1])\cap {\mathcal M}^{reg})\leq d_{\mathcal M}(x,y)+\e.
$$
 We assume that $\gamma$ is
parametrised so that it has constant speed $1/\hbox{Length}_g(\gamma([0,1])\cap {\mathcal M}^{reg})$.

We start the proof by a modification of the proof of Lemma \ref{deformation}.
As the path $\gamma$ intersects ${\mathcal M}^{sing}$  finitely many times, 
 by definition there are $t_1,t_2,\dots,t_J$, $J\in \Z_+$, such that for $t\in [0,1]$ we have
 $$\gamma(t)\in {\mathcal M}^{sing}\hbox{ if and only if }t\in \{t_1,t_2,\dots,t_J\}.$$

For each $j=1,2,\dots,J$, let $[s_j,r_{j}]\subset I$ be intervals such that 
$s_j<t_j<r_{j}$
 and  ${\pi}_\ell :\tilde U_\ell \to U_\ell $, $\ell=\ell(j)$ be projections from covering neighbourhoods such that
$\gamma([s_j,r_j])\in U_\ell$. Let
$\tilde \gamma_j:[s_j,r_j]\to \tilde U_{\ell(j)}$
paths such that ${\pi}_{\ell (j)}(\tilde \gamma_j(s))=\gamma(s)$,
that is, $\tilde \gamma_j:[s_j,r_j]\to \tilde U_{\ell(j)}$  is the lift of  the path
$\gamma:[s_j,r_j]\to U_{\ell(j)}$.

%
Let 
\ba
& &\rho_2:\tilde U_{\ell(j)}=\tilde W_{\ell(j)}\times \tilde V_{\ell(j)}\to  \tilde V_{\ell(j)}\subset {\mathbb R}^{n-k},\\
& &\rho_2(y,z)=z
\ea
be the projection of $x$ to the $z$-coordinate.
Note that $n-k\geq 2$. Since $\tilde \gamma_j$ is a rectifiable path, we see that $\rho_2(\tilde \gamma_j([s_j,r_j]))$ has the Hausdorff dimension 1. Hence, for all $i\in\Z_+$ 
 that there are vectors $v_{j,i}\in \R^{n-k},$   such that $|v_{j,i}|<1/i$  and that $-v_{j,i}\not \in \rho_2(\tilde \gamma_j([s_j,r_j]))$.
This implies that
$$
0\not \in v_{j,i}+\rho_2(\tilde \gamma_j([s_j,r_j]))
$$
and for $\tilde v_{j,i}=(0,v_{j,i})\in \R^n$
$$
(\tilde  v_{j,i}+\tilde \gamma_j([s_j,r_j]))\cap (\tilde W\times\{0\})=\emptyset.
$$
Denote 
$$
\eta_{j,i}:[s_j,r_j]\to {\mathcal M},\quad \eta_{j,i}(s)=\pi_{\ell (j)}(\tilde \eta_{j,i}(s)),
\quad \tilde \eta_{j,i}(s)=
\tilde  v_{j,i}+\tilde \gamma_j(s).
$$

For $j=1,2,\dots,J$, let $\tilde \alpha_{j,i}:[0,1]\to U_{\ell(j)}$  be the Euclidean line segment
$\tilde \alpha_{j,i}(t)=\tilde  v_{j,i}t+\tilde \gamma_j(s_j)$ that connects $\tilde \gamma_j(s_j)$ to 
$\tilde \eta_j(s_j).$ 
Similarly, let $\tilde \beta_{j,i}:[0,1]\to U_{\ell(j)}$  be the Euclidean line segment
that connects   
$\tilde \eta_j(r_j)$ to $\tilde \gamma_j(r_j)$.
When $i$ is large enough,
the line segments $\tilde \alpha_{j,i}([0,1])$  and  $\tilde \beta_{j,i}([0,1])$ do not intersect ${\mathbb R}^k\times \{0\}$.

Let  $\alpha_{j,i}(t)=\pi_{\ell(j)}(\tilde \alpha_{j,i}(t))$ and $\beta_{j,i}(t)=\pi_{\ell(j)}(\tilde \beta_{j,i}(t))$.
Then  $ \alpha_{j,i}([0,1])$  and  $\beta_{j,i}([0,1])$ do not intersect $\mathcal M^{sing}$.

%
%
 
Let now $\mu_i:[0,1]\to {\mathcal M}$ be a path that is obtained by concatenating the paths
$\gamma([0,s_1]),\alpha_{1,i},\eta_{1,i},\beta_{1,i},\gamma([r_1,s_2]),\alpha_{2,i},\dots,\beta_{J,i}$ and $\gamma([r_J,1])$
that connects $\gamma(0)$ to $\gamma(1)$, that is,
\ba
\mu_i&=&\gamma([0,s_1])\cup \bigg(
\bigcup_{j=1}^{J-1} \alpha_{j,i}\cup \eta_{j,i}\cup\beta_{j,i}\cup \gamma([r_j,s_{j+1}])\bigg)\cup\\
& &\quad
\cup \,\alpha_{J,i}\cup \eta_{J,i}\cup\beta_{J,i}\cup \gamma([r_J,1])\subset  \mathcal M^{reg}.
\ea
The {\btext metric tensor $g$ satisfies \eqref{tilde g bounded A}
 and $g$ is smooth in $\mathcal M^{reg}$. Thus when
 we compute the lengths of curves using \eqref{path length}, we see using Lebesque dominated converge  theorem there are numbers $h_i>0$, $\lim_{i\to\infty} h_i=0$  such that
 \ba
 \sum_{j=1}^{J-1} |\hbox{Length}_{g^{}}(\eta_{j,i})-
\hbox{Length}_{g^{}}(\gamma([s_j,r_j])\cap \mathcal M^{reg} )
|<h_i
\ea
and there is $C_1$  such that
\ba
 \sum_{j=1}^J(\hbox{Length}_{g^{}}(\a_{j,i})+\hbox{Length}_{g^{}}(\beta_{j,i}))  \le C_1\frac Ji.
\ea
Hence,
\ba
|\hbox{Length}_{g^{}}(\mu_{i})
-  \hbox{Length}_{g^{}}(\gamma([0,1])\cap \mathcal M^{reg})|<h_i+C_1\frac Ji
\ea

Denote $I=[0,1]$.
Note that $\mu_i(I)\subset \mathcal M^{reg}$, $\mu_i(0)=\gamma(0)=x$ and $\mu_i(1)=\gamma(1)=y$. 
These  
 imply
\ba 
d_{\mathcal M}(x,y)+\e&\geq& \hbox{Length}_{g}(\gamma([0,1])\cap \mathcal M^{reg})\\
&\ge &\hbox{Length}_{g^{}}(\mu_i([0,1]))+h_i+C_1\frac Ji\\
& \ge & d_{ \mathcal M^{reg}}(x,y)+h_i+C_1\frac Ji.
\ea
By taking the limit $i\to \infty$,  we obtain
%
%
%
%
 $$
 d_{\mathcal M}(x,y)+\e\geq d_{\mathcal M^{reg}}(x,y).
  $$ 
  As $\e>0$ is arbitrary, we obtain}
   $$
 d_{\mathcal M}(x,y)\geq d_{\mathcal M^{reg}}(x,y).
  $$ 
The opposite inequality follows from the definitions.

\hfill$\square$\medskip

\section{Finite  speed of wave propagation}

{In this section we consider estimates on how the support of the solutions of the wave equation propagates in time.
The propagation of singularities has been analyzed extensively using microlocal analysis, see e.g. \cite{Baskin,deHoop-Uhlmann-Vasy,Ford1,Ford2,Lebeau,Melrose-Vasy-Wunsch2,Melrose-Vasy-Wunsch1,Vasy2}, and these techniques describe
the propagation of the singular support. However, to consider the
support of the wave we use more conventional methods with  $\Gamma$-convergence.}

\subsection{Auxiliary results based on $\Gamma$-convergence}

Let us now introduce {quadratic} forms
\begin{equation}
\mathcal Q^{(m)}(u,v) =\mathcal Q^{(m)}_1(u,v) +\mathcal Q^{(m)}_2(u,v)
\end{equation}
where
\ba
& &\mathcal Q^{(m)}_1(u,v) = 
\int_{\mathcal M^{reg}}  g^{(m),ij}\partial_iu\overline{\partial_j v}\,\det(g^{(m)})^{1/2}dx,\\
& &\mathcal Q^{(m)}_2(u,v) = 
\int_{\mathcal M^{reg}}  \lambda u\overline{v}\,\det(g^{(m)})^{1/2}dx,
\ea
and
\begin{equation}
\mathcal Q(u,v) =\mathcal Q_1(u,v) +\mathcal Q_2(u,v)
\end{equation}
where
\ba
& &\mathcal Q_1(u,v) = 
\int_{\mathcal M^{reg}}  g^{ij}\partial_iu\overline{\partial_j v}\,\det(g)^{1/2}dx,\\
& &\mathcal Q_2(u,v) = 
\int_{\mathcal M^{reg}}  \lambda u\overline{v}\,\det(g)^{1/2}dx.
\ea
All these quadratic forms are defined as (unbounded) closed quadratic forms in $L^2({\mathcal M})$
with the domain $\mathcal D(\mathcal Q^{(m)})=\mathcal D(\mathcal Q)=H^1({\mathcal M} ).$
Note that we next consider $L^2({\mathcal M})$ with the inner product given by metric $g$, that is,
\ba
(u,v)_{L^2}=  \int_{\mathcal M^{reg}} u\overline{v}\,\det(g)^{1/2}dx.
\ea
Next we consider quadratic forms as ${\mathbb R}\cup \{\infty\}$  valued functions $u\mapsto  \mathcal Q^{(m)}_1(u,u)$ in $L^2({\mathcal M})$
that have value $\infty$  in $L^2({\mathcal M})\setminus H^1({\mathcal M})$.

Let $C_0,C_1>0$ be such that for all $m\ge 1$,
\ba
& &C_0^{-1}\le |g^{(m)}|^{-1/2}|g|^{1/2}\le C_0,\\
& &C_1^{-1} I\le g^{(m)}\le C_1I,\\
& &C_1^{-1} I\le g\le C_1I.
\ea

As for almost all $x\in \mathcal M^{reg}$ we have the pointwise limit
\begin{equation}
 \lim_{m\to \infty} g^{(m),ij}(x)\det(g^{(m)}(x))^{1/2}=g^{ij}(x)\det(g(x))^{1/2},
\end{equation}
we see that the proof of \cite[Prop. 5.14]{DalMaso}, 
applied in local coordinates, implies that,
 $\mathcal Q^{(m)}_1$  $\Gamma-$converges to  $\mathcal Q_1$ in the weak topology of $H^1(\mathcal M, |g|^{1/2}dx)$ as $m\to \infty$.
As $\mathcal Q^{(m)}_1$ are quadratic forms, by  \cite[Prop. 13.12]{DalMaso} 
this implies that $\mathcal Q^{(m)}_1$  $\Gamma-$converges to  $\mathcal Q_1$ in the strong topology of $L^2(\mathcal M, |g|^{1/2}dx)$ as $m\to \infty$.
In addition, by Lebesque dominated convergence theorem, $\mathcal Q^{(m)}_1$  converges also pointwise to  $\mathcal Q_1$.


As
 \begin{equation}
\big(g^{(m)}_{ij}\big) \geq \big(g^{(m+1)}_{ij}\big) 
\geq \big(g_{ij}\big),
\label{gsubijdecrease B}
\end{equation}
we see that 
 \begin{equation}
\hbox{det}(g^{(m)}) \geq \hbox{det}(g^{(m+1)}) \geq  \hbox{det}( g),
\label{gsubijdecrease C}
\end{equation}
and hence the sequence $\mathcal Q^{(m)}_2(u,u)$ is decreasing.
By Lebesque dominated convergence theorem,  $\mathcal Q^{(m)}_2(u,u)$ converge 
pointwise to  $\mathcal Q_2(u,u)$ as $m\to \infty$, and  $\mathcal Q_2(u,u)$ is lower semi-continuous in the strong topology of $L^2(\mathcal M, |g|^{1/2}dx)$, and hence
 \cite[Prop.\ 5.7]{DalMaso} 
 implies that  $\mathcal Q^{(m)}_2$  $\Gamma-$converges to  $\mathcal Q_2$ in the strong topology of $L^2(\mathcal M, |g|^{1/2}dx)$ as $m\to \infty$.

As the sequences $\mathcal Q^{(m)}_i$, $i=1,2$ both $\Gamma-$converges to  $\mathcal Q_i$ in the strong topology of $L^2(\mathcal M, |g|^{1/2}dx)$ and pointwise converge  as $m\to \infty$,
  \cite[Prop. 6.25]{DalMaso}, 
implies that the sum of the sequences,  $\mathcal Q^{(m)}$,  $\Gamma-$converges to the sum $\mathcal Q=\mathcal Q_1+\mathcal Q_2$  in the strong topology of $L^2(\mathcal M, |g|^{1/2}dx)$ as $m\to \infty$.

  
 

%
%
%
 
In the Hilbert space $L^2(\mathcal M, |g|^{1/2}dx)$, the symmetric quadratic form $\mathcal Q$  is associated,
{in the sense of \cite[Sec.\ VI, Thm 2.6]{Kato66}, to the selfadoint} operator 
\ba
A^{(m)}=A^{(m)}_1+A^{(m)}_2
\ea
in the Hilbert space $L^2(\mathcal M, |g|^{1/2}dx)$ endowed with the volume form of the metric $g$.
{We denote the domain of this selfadjoing operator $A^{(m)}$  by $\mathcal D(A^{(m)})$. Then,
for $u\in \mathcal D(A^{(m)})$ we have 
\ba
A^{(m)}_1u(x)&=&- |g(x)|^{-1/2}\p_j( |g^{(m)}(x)|^{1/2}g^{(m),jk}(x)\p_k u(x)),\quad\hbox{for }x\in \mathcal M^{reg},\\
A^{(m)}_2(x)&=&\lambda  |g^{(m)}(x)|^{1/2}|g(x)|^{-1/2},\quad\hbox{for }x\in \mathcal M^{reg}.
\ea
}

In the Hilbert space $L^2(\mathcal M, |g|^{1/2}dx)$, the symmetric quadratic form $\mathcal Q$  is associated to the
 {selfadjoint operator
\ba
A=A_1+A_2,
\ea
having the domain $\mathcal D(A)$, and {for $u\in \mathcal D(A)$ we have}}
\ba
A_1u(x)&=&{- |g(x)|^{-1/2}\p_j( |g(x)|^{1/2}g^{jk}(x)\p_ku(x))=-\Delta_gu(x),\quad\hbox{for }x\in \mathcal M^{reg},}\\
A_2u(x)&=&\lambda u(x) .
\ea

By \cite[Prop. 13.12]{DalMaso} 
the operators 
$A^{(m)}$
$G-$converges in the Hilbert space $L^2(\mathcal M, |g|^{1/2}dx)$ to the operator $A$
 in the strong topology of $L^2(\mathcal M, |g|^{1/2}dx)$.
Moreover, by \cite[Prop. 13.12]{DalMaso} 
then $A^{(m)}$ converge to $A$ in the strong resolvent sense in the  strong  topology of $L^2(\mathcal M, |g|^{1/2}dx)$.

Note that the equation
\ba
(-\Delta_{g^{(m)}}+\lambda)u^{(m)}=f
\ea
is equivalent to 
\ba
A^{(m)}u^{(m)}=  |g^{(m)}|^{1/2}|g|^{-1/2}f
\ea
and if $\supp(f)$ does not intersect $\supp(g^{(m)}-g)\subset K_m$, we have that $  |g^{(m)}|^{1/2}|g|^{-1/2}f=f$.

Let $f\in L^2(\mathcal M, |g|^{1/2}dx)$ and $\lambda\in {\mathbb C}\setminus (-\infty,0]$.

Let $\e>0$ be arbitrary and choose $m_0$ such that
\beq
\bigg(\int_{K_{m_0}}|f(x)|^2\, |g|^{1/2}dx\bigg)^{1/2}<\frac \e{4C_0} \frac 1{\dist(\lambda, {\mathbb R}_-)}.
\eeq
Let 
\ba
f_1(x)=\chi_{K_{m_0}} (x)f(x),\quad f_2(x)=f(x)-f_1(x).
\ea
Then for all $m\ge m_0$
\beq
\||g^{(m)}|^{1/2}|g|^{-1/2}f_1\|_{L^2(\mathcal M, |g|^{1/2}dx)}\leq \frac \e{4} \frac 1{\dist(\lambda, {\mathbb R}_-)},
\eeq
so that
\ba
u^{(m)}_1=(-\Delta_{g^{(m)}}+\lambda)^{-1}f_1,\quad u_1=(-\Delta_{g}+\lambda)^{-1}f_1
\ea
satisfy
\ba
\|u^{(m)}_1\|_{L^2(\mathcal M, |g|^{1/2}dx)}&=& \|(-\Delta_{g^{(m)}}+\lambda)^{-1}f_1\|_{L^2(\mathcal M, |g|^{1/2}dx)}\\
&=&\|(A^{(m)})^{-1}( |g^{(m)}|^{1/2}|g|^{-1/2}f_1)\|_{L^2(\mathcal M, |g|^{1/2}dx)}\\
&\leq& \frac \e{4}
\ea
and 
\ba
\|u_1\|_{L^2(\mathcal M, |g|^{1/2}dx)}<\frac \e{4}.
\ea
Moreover, as $K_{n+1}\subset  K_{n}$, we see that for $m\ge m_0$
\beq
\supp(f_2)\cap K_{m}=\emptyset,
\eeq
so that 
\ba
|g^{(m)}|^{1/2}|g|^{-1/2}f_2=f_2.
\ea
Next we consider $m \ge m_0$. 
As $A^{(m)}$ converge to $A$ in the strong resolvent sense in topology of $L^2(\mathcal M, |g|^{1/2}dx)$, 
by \cite[Def. 13.3]{DalMaso} we have that 
\ba
u^{(m)}_2=(-\Delta_{g^{(m)}}+\lambda)^{-1}f_2,\quad u_2=(-\Delta_{g}+\lambda)^{-1}f_2
\ea
satifies
\ba
u^{(m)}_2=(-\Delta_{g^{(m)}}+\lambda)^{-1}f_2=(A^{(m)})^{-1}f_2\to A^{-1}f_2=(-\Delta_{g}+\lambda)^{-1}f_2=u_2
\ea
in the  strong  topology of $L^2(\mathcal M, |g|^{1/2}dx)$ as $m\to \infty$.
Thus there is $m_1>m_0$ such that for all $m>m_1$
\ba
\|u^{(m)}_2-u_2\|_{L^2(\mathcal M, |g|^{1/2}dx)}<\frac \e 2.
\ea

As above $\e>0$ is arbitrary and $u^{(m)}=u^{(m)}_1+u^{(m)}_2$ and $u=u_1+u_2$
satisfy $\|u^{(m)}-u\|_{L^2(\mathcal M, |g|^{1/2}dx)}<\e$ for $m>m_1$, we see that
for all $f\in L^2(\mathcal M, |g|^{1/2}dx)$ and $\lambda\in {\mathbb C} {\setminus (-\infty,0]}$
\ba
\lim_{m\to \infty}(-\Delta_{g^{(m)}}+\lambda)^{-1}f=(-\Delta_g+\lambda)^{-1}f,
\ea
 strong  topology of $L^2(\mathcal M, |g|^{1/2}dx)$.

In particular, when $\lambda=1$, we have that
 $\big(- \Delta_{g^{(m)}}+1\big)^{-1} $ converges to $\big(- \Delta_{\widetilde g}+1\big)^{-1}$ strongly in $L^2(\mathcal M, |g|^{1/2}dx)$.

%

\medskip

\noindent
{\bf Remark.}
In the n-dimensional case with $n\ge 3$ the above considerations can be simplified as follows. Let $\omega_1,\omega_2,\dots,\omega_n$ be the eigenvalues of $g$.
Then the eigenvalues of $\det(g)^{1/2}g^{-1}$ are 
\ba
& &s_1= \omega_1^{-1/2}\omega_2^{1/2}\omega_3^{1/2}\dots \omega_n^{1/2},\\
& &s_2= \omega_1^{1/2}\omega_2^{-1/2}\omega_3^{1/2}\dots \omega_n^{1/2},\\
& &\dots\\
& &s_n= \omega_1^{1/2}\omega_2^{1/2}\omega_3^{1/2}\dots \omega_n^{-1/2},\\
\ea
and we see that $|s_j|>C_0^{-n/2}$. 
Also, assume that $g^{smooth}=\hat g$ is constant metric $\hat g=\hat \omega I$ in the domain $K_1$ (Note that gluing different local 
coordinate charts has to be added, e.g. by using a partition of unity and summing up to local construction).
Then the eigenvalues of $\det(\hat g)^{1/2}\hat g^{-1}$ are constants
\ba
& &\hat s=\hat \omega^{(n-2)/2}.
\ea
If we change the definition of $g^{(m)}$ in local coordinates to be defined using powers of symmetric matrixes, that is,
 $$
 g_{jk}^{(m)}(x)=g^{smooth}_{jk}(x)^{\psi_m(x) }g_{jk}(x)^{1-\psi_m(x)},
 $$
 we see that then $\det(g^{(m)})^{1/2}(g^{(m)})^{-1}$ has the eigenvalues $$s_j^{(m)}=\hat \omega ^{\psi_m(x)/2 }\omega_j(x)^{(1-\psi_m(x))/2}.$$
 If $\hat \omega<C_0^{-n/(n-2)}$ so that $$\hat s=\hat \omega^{(n-2)/2}<C_0^{-n/2}\le s_j,$$
  we see that the sequences $s_j^{(m)}$  are decreasing, and thus the the sequence of the positive definite matrixes  
 $\det(g^{(m)})^{1/2}(g^{(m)})^{-1}$ is decreasing as $m\to \infty$. Then the quadratic forms 
$ \mathcal Q^{(m)}(u,v)$ and the corresponding selfadjoint operators are also decreasing,  and we can use the following monotone theorem on the quadratic forms. Recall that for for two qurdartic forms $q_1$ and $q_2$ on a Hilbert space $\mathcal H$ the inequality 
$q_1 \leq q_2$ means that their form domains satisfy 
$D(q_1) \supset D(q_2)$ and $q_1[u] \leq q_2[u]$.

\begin{theorem}
Let $\{q_n\}_{n=1}^{\infty}$ be a sequence of closed, positive definite quadratic forms satisfying $0 \leq q_1 \leq q_2 \leq \cdots$. Suppose that 
\begin{equation}
D(q_{\infty}) = \{u \in \mathcal H\, ; \, \sup_nq_n[u] < \infty\}
\end{equation}
is dense {in $\mathcal H$}. Then, the quadratic form $q_{\infty}$ defined by
\begin{equation}
q_{\infty}[u] = \lim_{n\to \infty} q_{n}[u] = \sup_{n\in \mathbb Z_+} q_n[u]
\end{equation}
with domain $D(q_{\infty})$ is closed. Moreover, if $A_n$ and $A_{\infty}$ are the self-adjoint operators associated with $q_{n}$ and $q_{\infty}$, then $A_n \to A_{\infty}$ in the strong resolvent sense.
\end{theorem}

For this result, see \cite{Reed-Simon}, Theorem S.14,  p.373. 

In our case, $\mathcal Q^{(m)}$ and ${\mathcal Q}$ have the common domain $H^1(\mathcal M)$. Since $g^{(m)}_{ij} \to  g_{ij}$, we have for $x\in \mathcal M^{reg}$
\begin{equation}
\big(g^{(m),ij}(x)\big) \leq \big(g^{(m+1),ij}(x)\big)
 \leq \big(g^{ij}(x)\big) {\quad\hbox{and $\quad \big(g^{(m),ij}(x)\big)
\to \big(g^{ij}(x)\big)$ as $m\to \infty$}.}
\end{equation}
Therefore, we have that
 $A_{g^{(m)}}^{-1}$ converges to $A^{-1}$ strongly in $L^2(\mathcal M, |g|^{1/2}dx)$.

\bigskip

\subsection{Strong resolvent convergence}
In this section, we study some results related to the strong resolvent convergence in a Hilbert space $\mathcal H$.

\begin{lemma}
\label{AnAGconvertgencelemma1}
Suppose $A_n, A \geq 1$ are self-adjoint in a Hilbert space $\mathcal H$, and $A_n^{-1} \to A^{-1}$ strongly. Then 
$$
(A_n - z)^{-1} \to (A-z)^{-1}, \quad n \to \infty
$$
strongly for $z \in {\mathbb C}\setminus [1,\infty)$.
\end{lemma}

\begin{proof}
Take $z_0 \not\in [1,\infty)$ and assume that $(A_n - z_0)^{-1}g \to (A-z_0)^{-1}g$ strongly for any $g \in \mathcal H$.  For an { arbitrary} $f \in \mathcal H$, put
$$
u_n = (A_n - z_0)^{-1}f, \quad u = (A - z_0)^{-1}f.
$$
Then, 
$$
(A_n - z_0)^{-2}f = (A_n - z_0)^{-1}u_n = (A_n - z_0)^{-1}(u_n - u) + 
(A_n-z_0)^{-1}u.
$$
Since $\|(A_n - z_0)^{-1}\| \leq C$ for a constant $C>0$ independent of $n$, we have $$\|(A_n - z_0)^{-1}(u_n - u)\| \to 0.$$ Therefore, by the above assumption, 
$$
(A_n - z_0)^{-2}f \to (A-z_0)^{-1}u = (A-z_0)^{-2}f \quad 
strongly.
$$
By induction
$$
(A_n - z_0)^{-k}f \to (A-z_0)^{-k}f \quad 
strongly \quad \forall k \geq 1.
$$
{If $|z - z_0|<\dist(z_0,[1,\infty))$,  we have $|z - z_0|<\|(A_n - z_0)^{-1}\|^{-1}$ and}
$$
(A_n - z)^{-1} = (1 - (z - z_0)(A_n - z_0)^{-1})^{-1}(A_n - z_0)^{-1},
$$
and the Neumann series
$$
 (1 - (z - z_0)(A_n - z_0)^{-1})^{-1} = \sum_{k=0}^{\infty}
(z-z_0)^k(A_n - z_0)^{-k}
$$
is norm convergent. This implies that $(A_n - z)^{-1}f$ converges strongly to $(A-z_0)^{-1}f$ when {$|z - z_0|
<\dist(z_0,[1,\infty))
$}. 
Starting from $z_0=0$, we obtain the desired strong convergence for any $z \not\in [1,\infty)$
{by iterating the above analysis finitely many times with $z_0$ that are on a path the connects
$0$ to $z$ in ${\mathbb C}\setminus [1,\infty)$.}
\end{proof}

\begin{lemma}
\label{f(An)tof(A)Lemma}
Let $A_n, A \geq 1$ be self-adjoint and $(A_n - z)^{-1} \to (A-z)^{-1}$ strongly for any $z \not\in [1,\infty)$. 
Let $f(\lambda) \in C({\mathbb R})$ satisfy $f(\lambda) \to 0$ as $\lambda \to \infty$. Then, we have
$$
f(A_n) \to f(A) \quad strongly \quad as \quad n \to \infty.
$$
\end{lemma}

\begin{proof}
We first prove this lemma for $f \in C_0^{\infty}({\mathbb R})$. Let $F(z) \in C_0^{\infty}({\mathbb C})$ be an almost analytic extension of $f$ constructed by Lemma \ref{C1S3AlmostAnalyticExt}.  Usimg Lemma \ref{C1S3HelfferSjostrandFormula}, 
we then have for any $u \in \mathcal H$
\begin{equation}
\begin{split}
\big(f(A_n) - f(A)\big)u = \frac{1}{2\pi i}\int_{{\mathbb C}}\overline{\partial_z}F(z)
\big((z - A_n)^{-1} - (z - A)^{-1}\big)u\, dzd\overline{z}.
\end{split}
\nonumber
\end{equation}
Taking $N = 1$ and $s = 3$, we have
$$
\|\overline{\partial_z}F(z)
\big((z - A_n)^{-1} - (z - A)^{-1}\big)u\| \leq C(1 + |z|)^{-3}\|u\|.
$$
Moreover, by the assumption, for any $z \not\in [1,\infty)$,
$$
\|\big((z - A_n)^{-1} - (z - A)^{-1}\big)u\| \to 0.
$$
By Lebesgue's dominated convergence theorem, we obtain
$$
\|(f(A_n) - f(A))u\| \to 0.
$$
To prove the general case, {let $\epsilon>0$ and  put $f_R(\lambda) = f(\lambda)$ for $\lambda < R$ and $f_R(\lambda) = 0$ for $\lambda > R$. 
For a sufficiently large $R> 0$, 
$|f(\lambda) - g_R(\lambda)|\leq \epsilon$ for all $\lambda\in \R$.
Moreover, we can approximate $f_R(\lambda)$ by a $C_0^{\infty}({\mathbb R})$-function $g_R(\lambda)$ so that 
$\|f_R - g_R\|_{L^{\infty}({\mathbb R})} \leq \epsilon$.  Then the above shows that
$g_R(A_n) - g_R(A) \to 0$ strongly. As $\epsilon>0$  is arbitrary. This proves the lemma.}
\end{proof}

\begin{lemma}
\label{AdeltaLemma}
Let $A_n, A \geq 1$ be self-adjoint and $(A_n - z)^{-1} \to (A-z)^{-1}$ strongly for any $z \not\in [1,\infty)$. Suppose  there exists $\delta > 0$ such that for {the function $F(\lambda)=\lambda^\delta$
the self-adjoint operators $F(A_n)=A_n^{\delta}$ satisfy} $D(A_n^{\delta}) = D(A^{\delta})$ for all $n$ and there exists  $u \in D(A^{\delta})$ such that
$$
\sup_n\|A_n^{\delta}u\| < \infty, \quad \forall n.
$$
Then,  we have 
$$
\|f(A_n)u - f(A)u\| \to 0, {\quad \hbox{as }n\to \infty}
$$
for any bounded continuous function $f$ on $[1,\infty)$. 
\end{lemma}

\begin{proof}
Take $\chi_{0,R}(\lambda) \in C({\mathbb R})$ such that $\chi_{0,R}(\lambda) = 1$ for $\lambda < R$ and $\chi_{0,R}(\lambda) = 0$ for $\lambda > R+1$. Put $\chi_{\infty,R}(\lambda) = 1- \chi_{0,R}(\lambda)$. Then, we have
$$
\|\chi_{\infty,R}(A_n)u\| = \|\chi_{\infty,R}(A_n)A_n^{-\delta} A_n^{\delta}u\| \leq C\sup_{\lambda}|\chi_{R,\infty}(\lambda)\lambda^{-\delta}| \leq CR^{-\delta}.
$$
The same inequality holds for $A$. 
We fix large $R > 0$ and put $g(\lambda) = f(\lambda)\chi_{0,R}(\lambda)$. Then, by Lemma \ref{f(An)tof(A)Lemma}, $g(A_n)u \to g(A)u$ strongly.
\end{proof}

We consider  an abstract wave equation.

\begin{definition}
Let $A$ be the unbounded self-adjoint operator associated to the quadratic from $Q$
in $L^2({\mathcal M},g)$ having the domain $H^1({\mathcal M})=\mathcal D(Q)$. We say that $u\in C^1([0,T];L^2({\mathcal M}))\cap C^0([0,T];H^1({\mathcal M}))$ {is a (finite energy) solution of the wave equation} 
\begin{equation}
\partial_t^2 u + Au = 0
\end{equation}
if there are $u_0\in H^1({\mathcal M})$  and  $u_1\in L^2({\mathcal M})$  and
$$
u(t)=\cos(tA^{1/2})u_0+\sin(tA^{1/2})A^{-1/2}u_1.
$$
{Here we define $\sin(t\lambda^{1/2})\lambda^{-1/2}=t$ for $\lambda=0$.}
\end{definition}

Next we consider corresponding abstract system.

\begin{lemma}
\label{Waveeqconverge1}
Let $A_n, A \geq 0$ be self-adjoint in a Hilbert spece $\mathcal H$ such that  $(A_n - z)^{-1} \to (A-z)^{-1}$ strongly for any $z \in {\mathbb C} \setminus [0,\infty)$.
Let $u(t)$ be the solution of the equation
\begin{equation}
\left\{
\begin{split}
& \partial_t^2u + Au = 0, \\
& u\big|_{t=0} = w_0 \in \mathcal H, \quad \partial_tu\big|_{t=0} = w_1 \in H^1({\mathcal M})=\mathcal D(A^{1/2}),
\end{split}
\right.
\end{equation}
and $u_n(t)$ the solution of the same equation with $A$ replaced by $A_n$. Assume that  $D(\sqrt{A_n}) = D(\sqrt{A})$ for all $n$ and $w_0 \in D(\sqrt{A})$ with 
$\sup_n \|\sqrt{A_n} w_0\| < \infty$. Then, we have. 
$$
\|u_n(t) - u(t)\| \to 0\quad{\hbox{as }n\to\infty}
$$
for any $t$. 
\end{lemma}

\begin{proof}
We have
$$
u(t) = \cos(t\sqrt{A})w_0 + \sin(t\sqrt{A})\sqrt{A}^{-1}w_1,
$$
and the similar formula for $u_n(t)$. Then, the lemma follows from Lemma \ref{AdeltaLemma}.
\end{proof}

\begin{lemma}
Let $A_n$ and $A$ be as in Lemma \ref{Waveeqconverge1}. Let $u(t)$ and $u_n(t)$ be the solutions of the equations
\begin{equation}
\left\{
\begin{split}
& \partial_t^2 u +Au = f(t), \\
& u\big|_{t=0} = \partial_tu\big|_{t=0} = 0,
\end{split}
\right.
\end{equation}
\begin{equation}
\left\{
\begin{split}
& \partial_t^2 u_n + A_nu = f_n(t), \\
& u_n\big|_{t=0} = \partial_tu_n\big|_{t=0} = 0.
\end{split}
\right.
\end{equation}
Assume that for any $T > 0$ there exists a constant $C > 0$ such that $\|f_n(t)\| \leq C$ for any $n$ and $0 \leq t \leq T$, and $\|f_n(t) - f(t)\| \to 0$ as $n \to \infty$ for any $t \geq 0$. 
Then $\|u_n(t) - u(t)\| \to 0$ as $n \to \infty$ for any $t \geq 0$.
\end{lemma}

\begin{proof}
We put
$$
F(t,\lambda) = \big(\sin(t\sqrt{\lambda})\big)\sqrt{\lambda}^{-1}.
$$
Then, 
$$
u(t) = \int_0^tF(t-s,A)f(s)ds,
$$
and the similar formula holds with $A$ and $f(s)$ replaced by $A_n$ and $f_n(s)$. The lemma then follows from Lemma \ref{Waveeqconverge1}. 
\end{proof}

\subsection{Domains of influence}
Let $u$ satisfy
\beq\label{IVP 1}
& &(\p_t^2-\Delta_{g})u(x,t)=0,\quad \hbox{on }{\mathcal M}\times {\mathbb R}_+\\
& &{u(x,0)=v_0(x)\in L^2({\mathcal M}),\quad \p_t u(x,0)=v_1(x)\in H^1({\mathcal M}).} \label{IVP 2}
\eeq

Next we prove the finite propagation of waves.
Below, for $W\subset  {\mathcal M}$, let
\beq
 {\mathcal M}(W,T)=\{x\in M;\ d_ {\mathcal M}(x,W)<T\}
 \eeq
  denote the (open) domain of influence.

\begin {proposition}
\label {prop:3.2}
{\it  Let  $W\subset {\mathcal M}$ be an open, relatively compact set and $u^f(t)$ be the solution of
initial boundary value problem (\ref{IVP 1})--(\ref{IVP 2}) with
$\supp(v_0)\cup\supp(v_1)\subset W.$
Then
\ba
\supp(u(\,\cdotp,T))\subset {\mathcal M}(W,T).
\ea}
\end {proposition}

\noindent{\bf Proof.}
Let $W\subset {\mathcal M}$ be an open, relatively compact set and 
$$\supp(v_0)\cup\supp(v_1)\subset W.$$
Let $u^{(m)}$ satisfy
\beq\label{IVP 1 with tilde}
& &(\p_t^2-\Delta_{g^{(m)}})u^{(m)}(x,t)=0,\quad \hbox{on }{\mathcal M}\times {\mathbb R}_+\\
& &u^{(m)}(x,0)=v_0(x),\quad \p_t u^{(m)}(x,0)=v_1(x), \label{IVP 2 with tilde}
\eeq
where the wave equation is defined in weak sense on all coordinate neighbourhoods.

As we can represent $v_0$ and $v_1$  as a sum of functions
supported in single coordinate neighborhood $U_{\ell}$,
without loss of generality we can below assume that
$$
W\subset U_{\ell}.
$$
Also, 
let
$$
0<T_0<\dist_{g^{smooth,(\ell)}}(W,U_{\ell}).
$$

The standard results on the finite speed of wave propagation is valid on the smooth lifted coordinate neighbourhood  $(\tilde U_{\ell},g^{(m)})$, {see \cite{HoIV}
(Here the metric $g^{(m)}$ is smooth. We note} that the considerations can may be simplified using Lipschitz-smooth metric
as  the generalized results for the finite speed of wave propagation
seem to be valid on  Lipschitz-smooth manifold, or with divergence form equations with
 log-Lipschitz coefficients 
see e.g. \cite[Thm 1.6 and Remark 1.8]{CM}.)
Then 
there is an open set $W_1\subset \overline{W_1}\subset W$  such that 
$\supp(v_0)\cup\supp(v_1)\subset \overline W_1$.
Also, 
there is an open set $W_2\subset \overline{W_2}\subset W_1$  such that 
$\supp(v_0)\cup\supp(v_1)\subset \overline W_2$.
Then
\ba& &
\{(x,t)\in {\mathcal M}\times [0,T_0);\ d^{(m)}(x,W_2)> t\} \\
&\subset& \{(x,t)\in {\mathcal M}\times [0,T_0);\ u^{(m)}(x,t)=0\},
\ea
so that
\ba
& &\supp(u^{(m)}(\,\cdotp,\cdotp ))\cap ({\mathcal M}\times [0,T_0))\\
&\subset& \{(x,t)\in {\mathcal M}\times [0,T_0);\ d^{(m)}(x,\overline W_2)\leq  t\}\\
&\subset & \{(x,t)\in {\mathcal M}\times [0,T_0);\ d^{(m)}(x,W_1)< t\}.
\ea

Let $u$ satisfy
\ba
& &(\p_t^2-\Delta_{g})u(x,t)=0,\quad \hbox{on }{\mathcal M}\times {\mathbb R}_+\\
& &u(x,0)=v_0(x),\quad \p_t u(x,0)=v_1(x).
\ea
As propagation of waves can be studied in separately on small time intervals $[jT_0,(j+1)T_0]$, that cover a longer interval $[0,T_1]$, without loss of generality we can consider the case when $0<T\leq T_0$
and initial data $(v_0,v_1)$  supported in $W$.

We see using
Lemma \ref{Waveeqconverge1}
that 
$$
\lim_{m\to \infty}u^{(m)}=u\quad\hbox{in }L^\infty([0,T];L^2({\mathcal M})).
$$
This implies 
\beq\label{L2 limit}
\lim_{m\to \infty}u^{(m)}=u\quad\hbox{in }L^2([0,T]\times {\mathcal M}).
\eeq

As  above   $g^{(m)}\geq g^{(m+1)}$, we have
\ba
& &\{(x,t)\in {\mathcal M}\times{[0,T]};\ d^{(m)}(x,W_1)< t\}\\
&\subset& \{(x,t)\in {\mathcal M}\times{[0,T]};\ d^{(m+1)}(x,W_1)< t\},
\ea
and thus equation \eqref{L2 limit} implies that 
\beq\nonumber
\supp(u(\,\cdotp,\cdotp )){\cap({[0,T]}\times {\mathcal M})}&\subset&\hbox{cl}\bigg(
\bigcap_{m'>m_0} \bigcup_{m>m'} \supp(u_{m'}(\,\cdotp,\cdotp )){\cap({[0,T]}\times {\mathcal M})}\bigg)
\\
 &\subset&\label{m supports}
\hbox{cl}\bigg( \bigcup_{m\in \Z_+}\{(x,t)\in {\mathcal M}\times{[0,T]};\ d^{(m)}(x,W_1)< t\}\bigg).\hspace{-1.5cm}
\eeq


Lemmas \ref{two distance functions} and \ref{two distance functions opposite inequality}
 imply that 
\ba
 & &\bigcup_{m\in \Z_+}\{(x,t)\in {\mathcal M}\times{[0,T]};\ d^{(m)}(x,W_1)< t\}
 \\
 &&=\{(x,t)\in {\mathcal M}\times{[0,T]};\ d_{\mathcal M}(x,W_1)< t\}.
\ea
Then (\ref{m supports}) implies that 
\beq\label{u supports}
\supp(u(\,\cdotp,\cdotp )){\cap( [0,T]\times {\mathcal M})}&\subset&
 \{(x,t)\in {\mathcal M}\times{[0,T]};\ d_{\mathcal M}(x,W_1)\leq t\} \hspace{-1cm}
 \\ \nonumber
 &\subset&
 \{(x,t)\in {\mathcal M}\times{[0,T]};\ d_{\mathcal M}(x,W)< t\}.\hspace{-1cm}
\eeq


The above yield that 
\ba
\supp(u(\,\cdotp,\cdotp )){\cap( [0,T]\times {\mathcal M})}\subset \{(x,t)\in {\mathcal M}\times{[0,T]};\ d_{\mathcal M}(x,W)< t\}.
\ea
This means  finite velocity of wave propagation (i.e., that the waves propagate with velocity one or slower) is valid on $({\mathcal M},g)$.
This proves Proposition \ref{prop:3.2}. \hfill$\square$\medskip


\section{Unique continuation}

Next we will show that Tataru's approximate controllability result, 
see e.g. \cite{Ta1},
is valid for {CMGAs}. As usual, we start with the observability result.

\begin{theorem}[Tataru's Unique Continuation Principle] \label{observability}
 Let  
{$u\in C^0([0,2T];H^1({\mathcal M}))\cap C^1([0,2T];L^2({\mathcal M}))$
satisfy} 
\ba
\left( \p_t^2 - \Delta_{g} \right)u(x, t)= 0, \quad (x, t) \in {\mathcal M}  \times (0, 2T).
\ea
  Assume, in addition, that 
 $u(x, t)=0$ in $ W  \times (0, 2T)$,
 where $W$ is an open subset of $ {\mathcal M}^{reg}$.
 Then, 
\ba
       u(x, t)=0 \quad \hbox{in}\,\, K(W, T),
\ea
where $K(W, T)$ is the double cone of influence,
\beq \label{double_cone}
K(W, T)=\{(x, t)\in {\mathcal M}  \times (0, 2T): d(x, W) < T-|t-T|\}.\hspace{-2cm}
\eeq
\end{theorem}

\noindent{\bf Proof.}
Let $V\subset {\mathcal M}^{reg}$ be open.
Assume that $u$  satisfies on ${\mathcal M}\times {\mathbb R}_+$ the wave equation
\ba
& &(\p_t^2-\Delta_{g})u(x,t)=0,\quad \hbox{on }{\mathcal {\mathcal M}}\times {\mathbb R}_+\\
& &u(x,0)=v_0(x),\quad \p_t u(x,0)=v_1(x).
\ea
Also, assume that $u(x,t)$ vanishes in the set $V\times (0,2T)$.


The  restriction of $u(x,t)$ on ${\mathcal M}^{reg}\times {\mathbb R}_+$ satisfies
\ba
& &(\p_t^2-\Delta_{g})u(x,t)=0,\quad \hbox{on }{\mathcal M}^{reg}\times {\mathbb R}_+,
\ea
and by {applying Tataru's theorem \cite{Ta1}  on the smooth manifold ${\mathcal M}^{reg}\times {\mathbb R}_+$ (see \cite{BKLlocal,BKLglobal} for the corresponding
stability results),} we see that
$u$ vanishes in
\ba
\Sigma_{V,T}^{reg}=\{(x,t)\in {\mathcal M}^{reg}\times {\mathbb R}_+;\ |d_{{\mathcal M}^{reg},g}(x,V\cap {\mathcal M}^{reg})-T|<T\}.
\ea

As $V$  is open, we see that $V\cap {\mathcal M}^{reg}\not =\emptyset$.
Next, let $y\in V\cap {\mathcal M}^{reg}$  and
let $x\in {\mathcal M}^{reg}$ be such that  $d_{\mathcal M}(x,y)<T$. Above in
Lemma \ref{two distance functions opposite inequality} we have shown
that there is $m$ such that  $d_{g^{(m)}}(x,y)<T$. Let
$\mu$ be a $g^{(m)}$-geodesic that connects $x$ and $y$ and has length
$L=d_{g^{(m)}}(x,y)<T$. Let $\e=(T-L)/2$. In the proof of 
Lemma \ref{two distance functions} we showed that there is a path
$\eta:I=[0,1]\to {\mathcal M}$  such that $\eta(0)=x$, $\eta(1)=y$, and
 $\eta(I^{int})\subset {\mathcal M}^{reg}$, and finally, $g^{(m)}$-length of $\eta$
 is at most $L+\e$. As $g^{(m)}\geq g$ in ${\mathcal M}^{reg}$,
 this shows that $g$-length of $\eta$
 is at most $L+\e$. Thus we see that 
\ba
x\in \{x'\in {\mathcal M};\ d_{{\mathcal M}^{reg},g}(x',V)<T\}
\ea
and we have {by Lemma \ref{two distance M and Mreg} that}
\ba
\{x'\in {\mathcal M}^{reg};\ d_{{\mathcal M}^{reg},g}(x,V)<T\}=
\{x'\in {\mathcal M}^{reg};\ d_{{\mathcal M}}(x,V)<T\}.
\ea
This also shows that the solution of the wave equation, $u(x,t)$, vanishes in the set
 \ba
({\mathcal M}^{reg}\times {\mathbb R})\cap \Sigma_{V,T}^g=\{(x,t)\in {\mathcal M}^{reg}\times {\mathbb R}_+;\ |d_{{\mathcal M},g}(x,V)-T|<T\}.
\ea
%
%
%
As ${\mathcal M}^{reg}$ is dense in ${\mathcal M}$ with  respect to the topology defined with metric $d_{\mathcal M}$, 
we have
\ba
\Sigma_{V,T}^{reg}=\Sigma_{V,T}\cap ({\mathcal M}^{reg}\times {\mathbb R}_+),
\ea
where
\ba
\Sigma_{V,T}=\{(x,t)\in {\mathcal M}\times {\mathbb R}_+;\ |d_{{\mathcal M}}(x,V)-T|<T\}.
\ea
Thus we see that $u$  vanishes a.e. in $\Sigma_{V,T}$,
as $u$ is in $L^2_{loc}({\mathcal M}\times{\mathbb R})$, this means that 
 $u$, considered  as a distribution,  vanishes in the set $\Sigma_{V,T}$.
This means that Tataru's theorem is valid on $({\mathcal M},g)$.
\hfill$\square$\medskip

\bigskip

\section{Controllability results}

%
%
%

%
%
%
%
%
%

Consider the initial value problem
\beq
\label{IBVP1}
& &\left( \p_t^2 - \Delta_{g} \right)u(x,t)= H(x,t), \quad \hbox{in } {\mathcal M}\times {\mathbb R}_+,
\\ \nonumber
& &  u|_{t=0}=0, \quad \partial_t u|_{t=0}=0,
\eeq
and denote its solution by $u^H(x,t)=u(x,t)$.

Next we prove Tataru's controllability theorem on CMGA.

\begin{theorem} \label{Tataru} Let $W\subset {\mathcal M}^{reg}$ be an open set. 
  Then the set $\{u^H(\cdot, T):\, H \in C^\infty_0(W \times (0,
T))\}$  is dense in $L^2({\mathcal M}(W, T))$.
\end{theorem}
\noindent{\bf Proof.} 
 Assume that $\eta \in  L^2({\mathcal M}(W, T))$ satisfies
\beq\label{eta orthogonal}
( \eta, \, u^H(T))_{L^2({\mathcal M}; \, dV_g)} =0
\eeq
 for all $H \in C^\infty_0(W \times (0, T))$.

We consider the approximate initial value problem
\beq
\label{IBVP1 m}
& &\left( \p_t^2 - \Delta_{g^{(m)}} \right)u^{(m)}(x,t)= H(x,t), \quad \hbox{in } {\mathcal M}\times {\mathbb R}_+,
\\ \nonumber
& &  u^{(m)}|_{t=0}=0, \quad \partial_t u^{(m)}|_{t=0}=0,
\eeq
and denote its solution also by $u^{H,(m)}(x,t)=u^{(m)}(x,t)$.

We consider also the dual problem
\beq \label{dual}
& &\left( \p_t^2 - \Delta_{g} \right) a(x, t)=0, \quad \hbox{in}\,\, {\mathcal M} \times {\mathbb R},
\\
\nonumber
& & a(x, T)=0,\, \p_t a(x, T)= \eta(x),
\eeq
and the approximate dual problems
\beq \label{dual m}
& &\left( \p_t^2 - \Delta_{g^{(m)}} \right) a^{(m)}(x, t)=0, \quad \hbox{in}\,\, {\mathcal M} \times {\mathbb R},
\\
\nonumber
& & a^{(m)}(x, T)=0,\, \p_t a^{(m)}(x, T)= \eta(x).
\eeq
Then,
\ba
a^{(m)}(x, t) \in C(\R,\, H^1({\mathcal M})) \cap C^1({\mathbb R}, \, L^2({\mathcal M})) \subset H^1_{loc}({\mathcal M} \times {\mathbb R}).
\ea
By energy conservation for the wave equation, we have
\ba
\|\p_t a^{(m)}(\,\cdotp,t)\|_{L^2(\mathcal M,|g^{(m)}|^{1/2}dx)}^2&\leq& 
\|\p_t a^{(m)}(\,\cdotp,0)\|_{L^2(\mathcal M,|g^{(m)}|^{1/2}dx)}^2+0\\
&\leq &\|\eta\|_{L^2(\mathcal M,|g^{(m)}|^{1/2}dx)}^2
\ea
and hence for all $t\in [0,T]$ we obtain by integrating in the time variable
\ba
\| a^{(m)}(\,\cdotp,t)\|_{L^2(\mathcal M,|g^{(m)}|^{1/2}dx)}&\leq&
\|\eta+\int_0^t a^{(m)}(\,\cdotp,t')dt'\|_{L^2(\mathcal M,|g^{(m)}|^{1/2}dx)}\\
&\leq & 
 (1+T)\|\eta\|_{L^2(\mathcal M,|g^{(m)}|^{1/2}dx)}.
\ea
{Letting $m\to \infty$, we obtain}
\beq\label{final energy estimate}
\| a^{(m)}(\,\cdotp,t)\|_{L^2(\mathcal M,|g|^{1/2}dx)}&\leq&
 C_0(1+T)\|\eta\|_{L^2(\mathcal M,|g|^{1/2}dx)}. 
\eeq

Then, by integrating by parts, 
\ba
& &\int_{{\mathcal M} \times (0, T)} H(x, t) a^{(m)}(x, t)\, \,dV_{g^{(m)}}(x) dt
\\ \nonumber
& &=
\int_{W \times (0, T)} H(x, t) a^{(m)}(x, t)\, \,dV_{g^{(m)}}(x) dt
\\ \nonumber
& &=
\int_{{\mathcal M^{reg}} \times (0, T)} \bigg[  \bigg( \p_t^2 - \Delta_{g^{(m)}}  \bigg) u^{H,(m)}(x, t) a^{(m)}(x, t)
\\ \nonumber & &\quad \quad\quad\quad\quad\quad-
u^{H,(m)}(x, t) \bigg( \p_t^2 - \Delta_{g^{(m)}}  \bigg) a^{(m)}(x, t) \bigg]\, \,dV_{g^{(m)}}(x) dt
\\ \nonumber
& &=
-\int_{\mathcal M^{reg}} \left[\p_t u^{H,(m)}(x, 0) a^{(m)}(x, 0)-u^{H,(m)}(x, 0) \p_t a^{(m)}(x, 0)  \right]\,\,dV_{g^{(m)}}(x)
\\ \nonumber
& &\quad
+\int_{\mathcal M^{reg}} \left[\p_t u^{H,(m)}(x, T) a^{(m)}(x, T)-u^{H,(m)}(x, T) \p_t a^{(m)}(x, T)  \right]\,\,dV_{g^{(m)}}(x)
\\ \nonumber
& &=-\bra u^{H,(m)}(\cdot, T), \eta(\cdot) \cet_{L^2({\mathcal M}); \, dV_{g^{(m)})}},
\ea
where we use equations (\ref{IBVP1 m}) and (\ref{dual m}). Thus,
\beq\label{endpoint identity}
& &\int_{{\mathcal M} \times (0, T)} H(x, t) a^{(m)}(x, t)\, |{g^{(m)}}(x)|^{1/2}dx dt=\\
\nonumber
& &=-\int_{\mathcal M}  u^{H,(m)}(x, T)\eta(x)\, |{g^{(m)}}(x)|^{1/2}dx.
\eeq
By the strong resolvent convergence that for all $t\in [0,T]$
we have $a^{(m)}(\,\cdotp, t)\to a(\,\cdotp, t)$ in $L^2({\mathcal M}); \, dV_g)$ 
and $u^{H,(m)}(\cdot, T)\to u^{H}(\cdot, T)$ in $L^2({\mathcal M}); \, dV_g)$
as $m\to \infty$.

By applying estimate (\ref{final energy estimate}) and Lebesque dominated convergence theorem
we can take limit of both sides of (\ref{endpoint identity}) as $m\to \infty$, and obtain 
\ba
\int_{{\mathcal M} \times (0, T)} H(x, t) a(x, t)\, \,dV_{g}(x) dt=-\bra u^{H}(\cdot, T), \eta(\cdot) \cet_{L^2({\mathcal M}); \, dV_g)}=0,
\ea
where we use (\ref{eta orthogonal}) in the last identity.

Therefore, 
\ba
\int_{W \times (0, T)} H(x, t) a(x, t)\, \,dV_{g}(x) dt=0,
\ea for all
$H \in C^\infty_0(W \times (0, T))$, and hence $a(x, t)=0$ in $W \times (0, T)$.
Since $a(x, t)$ satisfies  $a(x, t)=a(x,2T- t)$ we have
\ba
a(x, t)=0 \quad \hbox{for}\,\ (x, t) \in W \times (0, 2T).
\ea
By Theorem \ref{observability},
\ba
a(x, t)=0 \quad \hbox{for}\,\, (x, t) \in K(W, T),
\ea
in particular,
\ba
\eta(x)= \p_t a(x, T)=0.
\ea
This yields that the set $\{u^H(\cdot, T):\, H \in C^\infty_0(W \times (0,
T))\}$  is dense in $L^2({\mathcal M}(W, T))$.
\hfill$\square$\medskip

\section{Uniqueness of inverse scattering}
\label{Uniquenessofinversescattering}

{Next we introduce an additional assumption for the inverse problem we use to show its unique solvability. To this end,
we} consider the volume factor
\beq\label{volume factor}
\Lambda(x)=\lim_{r\to 0+}\frac{\hbox{vol}_{\mathcal M}(B_{\mathcal M}(x,r)\cap {\mathcal M}^{reg})}{\hbox{vol}_{\R^n}(B(0,r))}
\eeq
Note that $\Lambda(x)=1$  for all $x\in {\mathcal M}^{reg}$.
\medskip

{\bf Assumption (L). }We assume below that
\beq\label{volume factor2}
\Lambda(x)\not=1,\quad\hbox{for all $x\in {\mathcal M}^{sing}$.}
\eeq
Note that 
\beq
\Lambda(x)\leq \frac{T(x)}{\#\Gamma_x},
\eeq
see \eqref{conditions 2}.

We aim to prove the following 

\begin{theorem}\label{thm: main result for IP}
Suppose we are given two conic manifolds with group action $\mathcal M^{(1)}$ and $\mathcal M^{(2)}$ satifying the assumptions {\bf (A-1)} $\sim$ {\bf (A-4)},
 {\bf (C-1)} $\sim$ {\bf (C-4)}, 
 {\bf (D)},  {\bf (VG)}, 
and  {\bf (L)}. Let the (1,1) component of the (generalized) scattering matrix coincide: 
$$
{\bf S}_{11}^{(1)}(k) = {\bf S}_{11}^{(2)}(k), \quad \forall k > 0, \quad  k^2\not\in \sigma_p(- \Delta^{(1)})\cup\sigma_p(- \Delta^{(2)}),
$$ 
and $r_1^{(1)}=r_1^{(2)}$.
Then there is an isometry between $\mathcal M^{(1)}$ and $\mathcal M^{(2)}$ in the following sense.

\noindent
(1) There is a homeomorphism $\Phi : \mathcal M^{(1)} \to \mathcal M^{(2)}$. \\
\noindent
(2) $\ \Phi(\mathcal M_{sing}^{(1)}) = \mathcal M_{sing}^{(2)}$.\\
\noindent
(3) $\ \Phi : \mathcal M^{(1)}\setminus\mathcal M_{sing}^{(1)} \to 
\mathcal M^{(2)}\setminus\mathcal M_{sing}^{(2)}$ is a Riemannian isometry. \\
%
\end{theorem}


\subsection{Blagovestchenskii's identity}

To prove the uniqueness of the inverse scattering problem we start with
some auxiliary results.
{Let ${\mathcal M}$ be a (compact or complete) conical Riemannian manifold with group action.
Let $\mathcal O\subset {\mathcal M}$ be open. 
Consider the solution $u^f(x,t)=u$ of the initial boundary value problem
\begin{equation}
\left\{
\begin{split}
& \partial_t^2u - \Delta_g u =f, \quad {\rm in} \quad {\mathcal M} \times {\mathbb R}_+, \\
& u\big|_{t=0} = \partial_tu\big|_{t=0} = 0, \quad {\rm in } \quad {\mathcal M}.
\end{split}
\right.
\label{S5:IBVPWave A}
\end{equation}

Also, we define the source-to-solution map
$V_{\mathcal O,+} :C^\infty_0({\mathcal O} \times \R_+)\to C^\infty(\overline {\mathcal O} \times {\mathbb R}_+),$ given by
\ba
V_{\mathcal O,+} (f)=u^f|_{{\mathcal O} \times {\mathbb R}_+}
\ea
where we denote $V_{\mathcal O,\pm}(T)=V_{\mathcal O,\pm}^T$ and $V_{\mathcal O,\pm}=V_{\mathcal O,\pm}^\infty.$ 
We denote below also $U_{\mathcal O,+}(\lambda)=U^\lambda_{\mathcal O,+}.$

 For considerations below, we observe that when a subset $\mathcal O\subset \mathcal M^{reg}$ and the metric $g|_{\mathcal O}$ on it are given, the hyperbolic source-to-solution operators $ V_{\mathcal O,+} $  and $U^\lambda_{\mathcal O,+}$ determine their Schwartz kernels $G(x,t;x_0,t_0)=
 G(x,x_0,t-t_0)$  and $G(\lambda,x,x_0)$ 
 that satisfy
 \beq
& & V_{\mathcal O,+} f(x,t)=\int_{\mathcal O\times {\mathbb R}_+}G(x,t;x_0,t_0)f(x_0,t_0)dV_g(x_0)dt_0,\\
& &U^\lambda_{\mathcal O,+} F(x)=\int_{\mathcal O}G(\lambda,x,x_0)F(x_0)dV_g(x_0)dt.
 \eeq

}
%


{To construct the manifold  $({\mathcal M},g)$ from local measurements, we use
a version of the boundary control method. The method originates from results of Belsihev and Kurylev 
\cite{Be87,BeKu92} and it has been further developed for different linear equations in \cite{BKLS,HLOS,KKL01,KOP,Lassas-Oksanen,Oksanen1,Oksanen2}, see also \cite{deHoop-Uhlmann} on the related  scattering control method.
The numerical implementation of the method has been recently developed in \cite{deHoop1,deHoop2,deHoop3}.
The present version of the method is based on focusing of the waves so that at a given time moment $t$  the value
$u(\,\cdotp,t)$ of the wave is concentrated in a neighborhood of a point, \cite{BKLS,DKL1}, {\mtext and the detection of singular points using local source to solution map (see also \cite{Kyria} for related techniques).} This focusing technique has recently been
used to study non-linear wave equation \cite{KLU,LUW} and the relation of reconstruction methods for the linear and the non-linear equations have been discussed in \cite{Lassas}.

}

\begin{theorem}
\label{th_Blago}
Let $({\mathcal M},g)$ be a compact or complete Riemannian manifold. Let $T>0$, ${\mathcal O}\subset {\mathcal M}$ be open and bounded.  

(i) {\btext Let $f,h \in C^\infty_0({\mathcal O}\times {\mathbb R}_+)$, then 
\begin{equation}
\label{Blagovestchenskii identity}
\langle u^f( \cdot,T),u^h( \cdot,T)\rangle_{L^2({\mathcal O})}=\langle f,J V^{2T}_{\mathcal O,+}h\rangle_{L^2({\mathcal O}\times 
(0,T))}-
\langle V^{2T}_{\mathcal O,+} f,  J h\rangle_{L^2( {\mathcal O})\times 
(0,T)}
\end{equation}
where  the operator $J:L^2({\mathcal O}\times 
(0,2T))\rightarrow L^2({\mathcal O}\times 
(0,T))$ is defined as
$$
J\phi(x,t)=\frac{1}{2}\ \int_{t}^{2T-t}\phi(x,s) \;ds.
$$}

(ii)
 Let $f \in C^\infty_0({\mathcal O}\times {\mathbb R}_+)$, then
\ba
\bra u^f(t),1\cet_{L^2({\mathcal M})}=\int_{{\mathcal M}} u^f(x,t)\,dV(x)
\ea
is given by
\ba
\langle u^f(t),1\rangle_{L^2({\mathcal M})}=\int_0^t\int_0^{t'}
\langle f(t''),1\rangle_{L^2({\mathcal O})}dt''dt'.
\ea
\end{theorem}
\begin{proof}
(i) Let $f,h \in C^\infty_0({\mathcal O}\times {\mathbb R}_+)$ and consider the mapping $V:[0,2T]\times [0,2T]\rightarrow {\mathbb R}$,
$$
V(t,s)=\langle u^f(t),u^h(s)\rangle_{L^2({\mathcal M})}.
$$
Then using Green's formula {we obtain}
\ba
& &(\p^2_t-\p^2_s)V(t,s)=(\p^2_t-\p^2_s)\langle u^f(t),u^h(s)\rangle_{L^2({\mathcal M})}\\
& &
=\langle f(t), V_{\mathcal O,+} h(s)\rangle_{L^2(N)}-\langle V_{\mathcal O,+} f(t),h(s)\rangle_{L^2(N)}:=F(t,s).
\ea
 The function  $(t,s)\mapsto F(t,s)$ can be computed, if the
 {the source-to-solution map
$V_{\mathcal O,+}$} is given. Note that
$$
V(0,s)=0,\quad \p_t V(t,s)|_{t=0}=0.
$$
Thus $V$ is the solution of the following $(1+1)$-dimensional initial value problem:
\begin{equation}
\label{1D wave equ}
\left\lbrace \begin{array}{l}
(\p^2_t-\p_s^2)V(t,s)=F(t,s), \: \textrm{ in } (0,2T)\times {\mathbb R}
\\
V|_{t=0}=\p_tV|_{t=0}=0.
\end{array} \right.
\end{equation}
Recall that the following formula
\begin{equation}
\label{Solution formula for 1d zero initial inner source wave equ} 
V(t,s)=\frac{1}{2}\int_0^t \int_{s-\tau}^{s+\tau}F(t-\tau,y) \;dyd\tau, \: s\in \R, \: t\in [0,2T],
\end{equation}
solves \eqref{1D wave equ}. 
By the change of variables $T-s = \tau$, we conclude
$$
V(T,T)=\frac{1}{2}\int_0^T \int_{\tau}^{2T-\tau}F(\tau,y) \;dyd\tau.
$$
$$
=\langle f,J V^{2T}_{\mathcal O,+} h\rangle_{L^2({\mathcal O} \times (0,T))}-\langle  V^{2T}_{\mathcal O,+} f\bigg|_{{\mathcal O}\times(0,T)},J h\rangle_{L^2({\mathcal O}\times (0,T))}.
$$
This proves (i).

(ii)
Let $f \in C^\infty_0({\mathcal O}\times {\mathbb R}_+)$ and consider the mapping $I:[0,T]\rightarrow {\mathbb R}$,
$$
I(t)=\langle u^f(t),1\rangle_{L^2({\mathcal M})}.
$$
Then using Green's formula 
\ba
& &\p^2_tI(t)=\langle \p^2_t u^f(t),1\rangle_{L^2({\mathcal M})}\\
& &=\langle \Delta u^f(t),1\rangle_{L^2({\mathcal M})}+\langle f(t),1\rangle_{L^2({\mathcal O})}
=\langle f(t),1\rangle_{L^2({\mathcal O})}.
\ea
Also, we have $I(0)=\p_tI(t)|_{t=0}=0$.
By solving the ordinary differential equation for $I(t)$ with initial conditions, we obtain the claim.
\end{proof}

\medskip
Above result is a generalization of Blagovestchenskii identity (see 
\cite[Theorem 3.7]{KKL01}) for Riemannian manifolds with conic singularities.

Next we will apply these formulas to compute the volume of the (open) domain of influence
\begin{equation}
{\mathcal M}(\tilde {{\mathcal O}},T)=\{x \in {\mathcal M}:\, { d_{\mathcal M}}(x, \tilde {{\mathcal O}}) < T   \},
\quad \tilde {{\mathcal O}} \subset {\mathcal O},
\label{S5Domainofinfluence}
\end{equation}
where ${d_{\mathcal M}}$ denotes the distance in ${\mathcal M}$ with respect to $g$. We denote the volume {of ${\mathcal M}(\tilde {{\mathcal O}},T)\cap \mathcal M^{reg}$ by ${\hbox{Vol}_g}({\mathcal M}(\tilde {{\mathcal O}},T))$
and define that ${\hbox{Vol}_g}({\mathcal M}^{sing})=0$.}


\begin{lemma} \label{lem: volumes}
Assume that we are given ${\mathcal O}$, the metric $g|_{\mathcal O}$
and the map $V^{2T}_{\mathcal O,+}$. 
Then, for any given open set $\tilde {{\mathcal O}} \subset  {\mathcal O}$ and $T>0$, these data 
uniquely 
determine the volume ${\hbox{Vol}_g}({\mathcal M}(\tilde {{\mathcal O}},T))$ of ${\mathcal M}(\tilde {{\mathcal O}},T)$.
\end{lemma}

\noindent {\bf Proof.} Let $w\in L^2({\mathcal M})$ be a function such that $w=1$ in 
${\mathcal M}(\tilde {{\mathcal O}},T)$.
For $f\in  C^\infty_0(\tilde {{\mathcal O}}\times (0,T))$, real-valued, we define the quadratic functional 
\ba
I_T(f)&=& \|u^f(\cdot, T)-w\|^2_{L^2({\mathcal M})}-\|w\|^2_{L^2({\mathcal M})}.
\ea
Since
 $
 \supp(u^f(\cdot, T))\subset {\mathcal M}(\tilde {{\mathcal O}},T), 
$
 we have
\beq \label{5.5a}
I_T(f)
&=& \|u^f(\cdot, T)\|^2_{L^2({\mathcal M})}-2 \bra u^f(\cdot, T),1\cet _{L^2({\mathcal M})}.
\eeq
Hence, by Theorem \ref{th_Blago},   we can compute $I_T(f)$ for any 
$f\in C^\infty_0(\tilde {{\mathcal O}}\times (0,T))$ uniquely by using 
$({\mathcal O},g|_{\mathcal O})$ and $V^{2T}_{\mathcal O,+}$. In the sequel, this is phrased as {\it we can compute}  $I_T(f)$.

Now we use again the fact that, for $f\in  C^\infty_0(\tilde {{\mathcal O}}\times (0,T))$,
 $\supp(u^f(\cdot, T))\subset {\mathcal M}(\tilde {{\mathcal O}},T)$ so that (\ref{5.5a}) yields that
\beq\label{def IT}
I_T(f)&=& \|u^f(\cdot, T)-\chi_{{\mathcal M}(\tilde {{\mathcal O}},T)}\|^2_{L^2( {\mathcal M})}- 
\|\chi_{{\mathcal M}(\tilde {{\mathcal O}},T)}\|^2_{L^2( {\mathcal M})},
\eeq
where $\chi_{{\mathcal M}(\tilde {{\mathcal O}},T)}$ is the characteristic function of ${\mathcal M}(\tilde {{\mathcal O}},T)$.
Thus, 
\beq\label{minimimum}
I_T(f)\geq - {\hbox{Vol}_g}({\mathcal M}(\tilde {{\mathcal O}},T)), \quad\hbox{for all }f\in  C^\infty_0(\tilde {{\mathcal O}}\times (0,T)).\hspace{-2cm}
\eeq  
By Tataru's controllability theorem, Theorem \re{Ta1},  there is a sequence $h_j\in  C^\infty_0(\tilde {{\mathcal O}}\times (0,T))$, such that
\ba
\lim_{j\to \infty} u^{h_j}(\cdot, T)=\chi_{{\mathcal M}(\tilde {{\mathcal O}},T)}\quad \hbox{in }L^2({\mathcal M}).
\ea
For this sequence,
\beq\label{minimizing 0}
\lim_{j \to \infty}I_T(h_j) = - {\hbox{Vol}_g}({\mathcal M}(\tilde {{\mathcal O}},T)).
\eeq
On the other hand, if 
$f_j \in C^\infty_0(\tilde {{\mathcal O}}\times (0,T))$ is a minimizing sequence for $I_T$, i.e.,
\beq\label{minimizing I}
\lim_{j\to \infty}I_T(f_j)=m_0:=\inf \{I_T(f);\,  f\in C^\infty_0(\tilde {{\mathcal O}}\times (0,T))\},
\eeq
then, by using  (\ref{minimimum}) and (\ref{minimizing 0}) {and the definition \eqref{def IT} of $I_T(f)$
and the property that
 $\supp(u^f(\cdot, T))\subset {\mathcal M}(\tilde {{\mathcal O}},T)$ for all
  $f\in  C^\infty_0(\tilde {{\mathcal O}}\times (0,T))$ we see that}
\ba
\lim_{j\to \infty} u^{f_j}(\cdot, T)=\chi_{{\mathcal M}(\tilde {{\mathcal O}},T)}\quad \hbox{in }L^2({\mathcal M}).
\ea
Thus, using any sequence $(f_j)$ satisfying (\ref{minimizing I}), we can compute
\ba
{\hbox{Vol}_g}({\mathcal M}(\tilde {{\mathcal O}},T))=\lim_{j\to \infty} \bra u^{f_j}(\cdot, T),u^{f_j}(\cdot, T)\cet _{L^2({\mathcal M})}.
\ea
\qed


%
%
%
%
\subsection{Reconstruction near ${\mathcal M}  $} \label{Green}


To prove Theorem \ref{thm: main result for IP} our first aim is 
to show that $ {\mathcal M}_{reg}^{(1)}$ and $ {\mathcal M}_{reg}^{(2)}$
are isometric. The proof is based on the procedure 
of the continuation of 
Green's functions%

By the above considerations, the scattering operator $\mathcal S_{11}$ determines
in an open set $W\subset {\mathcal M}^{reg}\cap  {\mathcal M}_1$. the source-to-solution operator $V_{W,+}$,
and thus we next we assume that we are given $(W,g|_W)$  and the operator $V_{W,+}$.

We are going to prove the uniqueness for the inverse problem
 step by step by  constructing relatively open subsets 
 ${\mathcal M}^{(1),rec} \subset {\mathcal M}  ^{(1)}$ 
 and  ${\mathcal M}^{(2),rec} \subset {\mathcal M}  ^{(2)}$, which are isometric
 and enlarge these sets at each step.
 In the following, when  
 ${\mathcal M}^{(1),rec} \subset {\mathcal M}_{reg}^{(1)}\cap {\mathcal M}  ^{(1)}$ 
 and  ${\mathcal M}^{(2),rec} \subset {\mathcal M}_{reg}^{(2)}\cap {\mathcal M}  ^{(2)}$  
 are relatively open connected sets and 
 \ba
\Phi^{rec}:\, {\mathcal M}^{(1),rec} \rightarrow {\mathcal M}^{(2),rec},
\ea
is a diffeomorphism, we say that the triple 
  $({\mathcal M}^{(1),rec},{\mathcal M}^{(2),rec},\Phi^{rec})$ is {\it admissible}
  if 
  {\begin{itemize}
  
  \item [(i)] $\Phi^{rec}:\, {\mathcal M}^{(1),rec} \to {\mathcal M}^{(2),rec}$
 is {a diffeomorphism and} an isometry, that is, $(\Phi^{rec})_*g^{(1)}=g^{(2)} $,

  \item [(ii)] the source-to-solution maps  $V_{\mathcal O^{(1)},+}$
and $V_{\mathcal O^{(2)},+}$ are $\Phi^{rec}$-related
on ${\mathcal O}^{(1)}$ and  ${\mathcal O}^{(2)}$, that is,
\beq\label{eq: source to solution operators  related T is infty}
 V_{\mathcal O^{(1)},+}(f\circ \Phi^{rec})=
 (V_{\mathcal O^{(2)},+}(f))\circ \Phi^{rec}
\eeq
for all $f\in C^\infty_0(\mathcal O^{(2)}\times {\mathbb R}_+)$.
\end{itemize}}

By {the above, (\ref{eq: source to solution operators  related T is infty})
is equivalent to   that the source-to-solution operators for the spectral problems
$U^{(1)}_{\mathcal O^{(1)},+}(\lambda)$
and $U^{(2)}_{\mathcal O^{(2)},+}(\lambda)$ are $\Phi^{rec}$-related
for all 
 $\lambda \in \sigma_e(H^{(i)})\setminus{\mathcal E}^{(i)}$,
 $i=1,2$, that is,
 \beq\label{eq: source to solution operators  related spectral}
U^{(1)}_{\mathcal O^{(1)},+}(\lambda)\circ (\Phi^{rec})^*=
 (\Phi^{rec})^*\circ U^{(2)}_{\mathcal O^{(2)},+}(\lambda).
 \eeq
  Furthermore, 
  (\ref{eq: source to solution operators  related T is infty})
 implies that Schwartz kernels of the source-to-solution operators, that
 is, the} time-domain Green's functions
 $G^{(i)}(z, x, t)$ on
${\mathcal M}^{(i),rec} $ satisfy
  the relation
\beq\label{eq: time domain Greens functions related}
G^{(2)}(\Phi^{rec}(x), \Phi^{rec}(y),t)=G^{(1)}(x, y,t)
\eeq
 for $x, y \in {\mathcal M}^{(1),rec}$ and $t\in {\mathbb R}$.
%
%
%

 Note that then the values of Green's functions  $G^{(i)}(z, x, y)$ on
${\mathcal M}^{(i),rec} $ satisfy 
the relation
\beq\label{eq: Greens functions related}
G^{(2)}(z; \Phi^{rec}(x), \Phi^{rec}(y))=G^{(1)}(z; x, y),\quad\hbox{for }\
x, y \in {\mathcal M}^{(1),rec},\ 
\ z\in {\mathbb C}\setminus {\mathbb R}.\hspace{-2cm}
\eeq

First we consider  Green's functions in  the set  
$W$.


Our earlier considerations, Theorem  \ref{S15SmatrixdeterminesLW2} 
and Lemma \ref{Lemma:StationarySSopdeterminestimedependrntSS}
show the following lemma:

\begin{lemma}\label{lem. N1 admissible} When $W$ is considered both as  a subset
${\mathcal M}  ^{(1)}$ and ${\mathcal M}  ^{(2)}$ and $I:W\to W$
is the identity map, then the triple $(W,W,I)$ is admissible.
\end{lemma}
%
%


\subsection{Continuation by Green's functions} \label{Green2}

To reconstruct subsets of manifolds   ${\mathcal M}^{(i)}$, $i=1,2$,
by continuing Green's function, we need the the following result telling
 that the values  of Green's functions
identify the points of the manifold. 

	
\begin{lemma} \label{identification} Let  $x_1,x_2\in {\mathcal M}^{(i)}$ be such 
that
\beq \label{3.3.4}
G^{(i)}(z, x_1, y)=G^{(i)}(z, x_2, y)\,  
\eeq
for all $\, y \in W$ and some $z \in {\mathbb C} \setminus {\mathbb R}$.
Then $x_1=x_2.$
\end{lemma}
\noindent
{\bf Proof.} 
Using the unique continuation principle for the solutions
of elliptic equations, we see that (\ref{3.3.4}) implies that 
$G^{(i)}(z, x_1, y)=G^{(i)}(z, x_2, y)$, for all $y \in {\mathcal M}^{(i)}\setminus
\{x_1,x_2\}$.
As the map $y\mapsto G^{(i)}(z, x, y)$ is bounded in the compact
subsets of ${\mathcal M}^{(i)}\setminus\{x\}$ and tends
to infinity as $y$ approaches $x$, this proves that $x_1=x_2$. 	
\qed

\begin{remark} \label{Remark:3}
Lemma \ref{identification} has the following important
consequence: If the triples $(N_1^{(1)},N_1^{(2)},\Phi_1)$
and $(N_2^{(1)},N_2^{(2)},\Phi_2)$ are
admissible and $N_1^{(1)}\cap N_2^{(1)}\not = \emptyset,
$ then, by  Lemma \ref{identification}, the maps $\Phi_1(x)$ and 
$\Phi_2(x)$ have to coincide  in $N_1^{(1)}\cap N_2^{(1)}$.
Moreover, if $N_3^{(i)}=N_1^{(i)}\cup N_2^{(i)}$, $i=1,2,$ 
and
\beq\label{eq: glue}
\Phi_3(x)=\left\{\begin{array}{cl}
\Phi_1(x), & \hbox{for } x\in N_1^{(1)}, \\
\Phi_2(x), & \hbox{for } x\in N_2^{(1)},
\end{array}
\right.
\eeq
then, by  Lemma \ref{identification}, the map 
$\Phi_3:N_3^{(1)}\to N_3^{(2)}$ is bijective
and hence a diffeomorphims. This implies that the 
triple $(N_3^{(1)},N_3^{(2)},\Phi_3)$ is
admissible
\end{remark}


The procedure of constructing the isometry between $\mathcal M^{(1)}  $ and $\mathcal M^{(2)}  $ consists of extending the admissible triple $(\mathcal M^{(1),rec},\mathcal M^{(2),rec},\Phi^{rec})$. 
In the first step, we apply Lemma \ref{lem. N1 admissible} 
to the triple $(W,W,I)$. {\btext In the subsequent steps
we extend the sets $\mathcal M^{(1),rec}$ and $\mathcal M^{(2),rec}$
and use sets  ${\mathcal O}^{(i)}\subset {\mathcal M}^{(i),rec}$,
defined below, to have the role of the set $W$ above.}

{\ntekst Let 
$q_i\in  {\mathcal M}^{(i),rec}$, $i=1,2$, be such that
\beq \label{eq: basic assumptions}
\hspace{1cm}\Phi^{rec}(q_1)=q_2,
\eeq
and let $d^{(i)}=d_{g^{(i)}}{=d_{\mathcal M^{(i)}}}$ denote the distance on $ {\mathcal M}^{(i)}$.
{\btext Let $R>0$ be sufficiently small so that
\beq \label{eq: basic assumptions2}
& &B^{(i)}(q_i, 4R)\subset {\mathcal M}^{(i),rec},\quad i=1,2 
\eeq
and the set
$${\mathcal O}^{(i)}=B^{(i)}(q_i,R) $$  is geodesically convex set,
has smooth boundary and that all points in the closure of $B^{(i)}(q_i, 2R)$
can be joined by a unique length minimizing curve. Note that
when the set $ {\mathcal M}^{(i),rec}$ and the metric on it are known, one can verify for a given $q_i\in {\mathcal M}^{(i),rec}$ if a given value $R$ satisfies these assumptions.}
Note that then normal coordinates, centered at $q_i$, are well defined in $B^{(i)}(q_i, 2R)$.
Then  $\Phi^{rec}({\mathcal O}^{(1)})={\mathcal O}^{(2)}$ and
${\mathcal M}^{(i),rec}\setminus \overline {\mathcal O}^{(i)}$ are connected.}

{\btext Below, 
we 
 say 
that the source-to-solution maps  $V^{2T}_{\mathcal O^{(1)},+}$
and $V^{2T}_{\mathcal O^{(2)},+}$ are $\Phi^{rec}$-related
on ${\mathcal O}^{(1)}$ and  ${\mathcal O}^{(2)}$ if
\beq\label{eq: source to solution operators  related}
 V^{2T}_{\mathcal O^{(2)},+}(H\circ \Phi^{rec})=
 (V^{2T}_{\mathcal O^{(1)},+}(H))\circ \Phi^{rec}
\eeq
for all $H\in C^\infty_0(\mathcal O^{(2)}\times {\mathbb R}_+)$.}

\begin{lemma} \label{lem. determination of Green} 
Let $({\mathcal M}^{(1),rec},{\mathcal M}^{(2),rec},\Phi^{rec})$
be an admissible triple and ${\mathcal O}^{(i)}$, $i=1,2$ be relatively compact subsets
of $ {\mathcal M}^{(i),rec}$ such that  ${\mathcal O}^{(2)}= \Phi^{rec}({\mathcal O}^{(1)})$.
Then {\btext for all $T>0$ the source-to-solution maps  $V^{2T}_{\mathcal O^{(1)},+}$
and $V^{2T}_{\mathcal O^{(2)},+}$ are $\Phi^{rec}$-related
on ${\mathcal O}^{(1)}$ and  ${\mathcal O}^{(2)}$.}
\end{lemma}
\noindent
{\bf Proof.} 
%
%
%
 Our earlier considerations, Theorem  \ref{S15SmatrixdeterminesLW2} 
and Lemma \ref{Lemma:StationarySSopdeterminestimedependrntSS}, yield the claim.

\qed
\medskip



{\btext
\subsection{Source-to-solution maps for subdomains of $\mathcal M$ and recognition of singular points}
\label{subsec: Relation}
For
 $y\in  \p {\mathcal M}^{(i),reg}$ 
 we define the {singular set and cut locus distances}
 \ba
\tau^{(i)}(y)=\min(\tau^{(i)}_{cut}(y),\tau^{(i)}_{sing}(y))
\ea
where
\ba
\tau^{(i)}_{sing}(y)=\inf_{\xi\in S_y {\mathcal M}^{(i)}} \tau^{(i)}_{sing}(y,\xi),\quad
 \tau^{(i)}_{sing}(y,\xi)=
\inf\left\{t>0;\ \
\gamma_{y,\xi}^{(i)}(t)\in  {\mathcal M}^{(i)}_{sing}\right\},
\ea
where 
 $\gamma^{(i)}_{y,\xi}(t)$ is  the geodesic on ${\mathcal M}^{(i)}$,
 {and
  \ba
& &\tau^{(i)}_{cut}(y)=\inf_{\xi\in S_y {\mathcal M}^{(i)}} \tau^{(i)}_{cut}(y,\xi),\\
& &\tau^{(i)}_{cut}(y,\xi)=
\inf\bigg(\{t\in  (0,  \tau^{(i)}_{sing}(y,\xi)) ;\ \
 d^{(i)}(\gamma_{y,\xi}^{(i)}(t),y)<t\}\cup\{ \tau^{(i)}_{sing}(y,\xi)\}\bigg)
\ea}

When ${W}\subset  \O^{(i)}$ are open and $s>0$, we denote the (open) domain of influence by
$$
 {\mathcal M}^{(i)}({W},s)=\{x\in {\mathcal M}^{(i)};\
{d}^{(i)}(x,{W})<s\}.
$$
We also denote
$$
 {\mathcal M}^{(i)}_\O({W},s)=\{x\in {\mathcal M}^{(i)}\setminus
 {\mathcal O}^{(i)};\
{d}^{(i)}(x,{W})<s\}.
$$


Next we consider the points $q_i\in {\mathcal M}^{(i),reg}$ satisfying \eqref{eq: basic assumptions} and $R$ satisfying assumption
in the formula (\ref{eq: basic assumptions2}) and below it.
Moreover,
assume that 
\beq\label{condition for T}
T<\min_{i=1,2} \tau^{(i)}(q_i)-R.
\eeq
Below, we consider the sets 
\beq
{\mathcal N}^{(i)}=\{p\in {\mathcal M^{(i)}}\setminus
{\mathcal O}^{(i)}:\ d_{g^{(i)}}(p,\p {\mathcal O}^{(i)})<T\}
\eeq
and
 the families of the interior distance functions corresponding
 to interior points in $\mathcal N^{(i)}$,
\begin{equation}
\label{distance functions}
R_{\mathcal O^{(i)}}(\mathcal N^{(i)}):=\{d_{g^{(i)}}(x,\cdot)|_{{\mathcal O}^{(i)}}:\ x\in \mathcal N^{(i)}\}\subset C(\mathcal O^{(i)}).
\end{equation}

\def\Src{\mathcal X}
\def\P{\mathbb P}
\def\W{\mathbb W}

%
%
%

\begin{theorem}
\label{th_geom} 
Let $({\mathcal M}^{(1),rec},{\mathcal M}^{(2),rec},\Phi^{rec})$
be an admissible triple and ${\mathcal O}^{(i)}$, $i=1,2$ be relatively compact subsets
of $ {\mathcal M}^{(i),rec}$ such that  ${\mathcal O}^{(2)}= \Phi^{rec}({\mathcal O}^{(1)})$.
%
Then
\begin{equation}
\label{distance functions related}
\{(d_{g^{(2)}}(x,\cdot)|_{\mathcal O^{(2)}})\circ \Phi^{rec}:\ x\in N^{(2)}\}=
\{d_{g^{(2)}}(x',\cdot)|_{\mathcal O^{(1)}}:\ x'\in N^{(1)}\}.
\end{equation}

\end{theorem}

This is to be proved in several steps.

Next we drop the superindex $(i)$ for a while and consider manifold $({\mathcal M},g)$.

Let $T, \epsilon >0$. For each $r>\epsilon$ and  $x \in \mathcal N$ we define a set 
\beq\label{set S}
S_\epsilon(x,r):=(T-(r-\epsilon),T)\times B(x,\epsilon),\quad   { B(x,\epsilon)\subset \mathcal M.}
\eeq

We denote for any measurable $A\subset \mathcal N$ the function space 
$$
L^2(A):=\{u\in L^2(\mathcal M): \hbox{supp}(u) \subset \overline{A}\}.
$$


We use the following lemma, developed in \cite{HLOS}, to study interior distance functions
on smooth manifolds, and give its proof for convenience of the reader.

\begin{lemma}
\label{crossing balls and functions}
Let $p,y,z\in \mathcal O$, $\epsilon>0$ and $\ell_p,\ell_y,\ell_y \in (\epsilon,T)$. Then the following are equivalent: 
\begin{enumerate}[(i)]
\item We have
\begin{equation}
\label{cond 3}
B(p,\ell_p)\subset \overline{B(y,\ell_y) \cup B(z,\ell_z)}.
\end{equation}
\item Suppose that
\begin{equation}
\label{cond 4}
\begin{array}{c}
\forall f \in C^\infty_0(S_\epsilon(p,\ell_p))\, \exists (f_j)_{j=1}^\infty \subset C^\infty_0(S_\epsilon(y,\ell_y)\cup S_\epsilon(z,\ell_z)) 
\\
\textrm{ such that }\lim_{j\to \infty} \|u^f(\cdot,T)-u^{f_j}( \cdot,T)\|_{L^2(N)}=0 .
\end{array}
\end{equation}
Here $u^f, u^{f_j}$ are the solutions of \eqref{S5:IBVPWave A}
with metric $g$ and
with sources $f$ and $f_j$, respectively.
\end{enumerate}
\end{lemma}
\begin{proof}

Suppose that \eqref{cond 3} is valid. Let $ f \in C^\infty_0(S_\epsilon(p,\ell_p))$. Then by the finite speed of wave propagation it holds that 
$$
\supp (u^f(T) )\subset B(p,\ell_p) \subset  \overline{B(y,\ell_y) \cup B(z,\ell_z)}
$$

Let $\chi(x)$ 
be the characteristic function of the ball $B(y,\ell_y)$ and set \ba
u^f_{(in)}(x,T):=\chi(x) u^f(x,T),\quad u^f_{(ext)}(x,T):=u^f(x,T)-u^f_{(in)}(x,T).
\ea 
Then  $u^f_{(in)}(\cdot,T)\in L^2(B(y,\epsilon))$ and $u^f_{(ext)}(\cdot,T)\in L^2(B(z,\epsilon))$. 
By approximate controllability there exist sequences $(f_y^j)_{j=1}^\infty \subset C^\infty_0(S_\epsilon(y,\ell_y))$ and $(f_z^j)_{j=1}^\infty \subset C^\infty_0(S_\epsilon(z,\ell_z))$ such that sequences $(u^{f_y^j}(\cdot,T))_{j=1}^\infty$ and $(u^{f_z^j}(\cdot,T))_{j=1}^\infty$ converge to $u^f_{(in)}(\cdot,T)$ and $u^f_{(ext)}(\cdot,T)$, respectively, in $L^2(\mathcal M)$.
Therefore the sequence 
$$
f_j=f^j_y+f_z^j \in C^\infty_0(S_\epsilon(y,\ell_y)\cup S_\epsilon(z,\ell_z)), \: j=1,2,\ldots$$ 
satisfies \eqref{cond 4}. 

Suppose that \eqref{cond 3} is not valid. Then the open set
$$
U:=B(p,\ell_p)\setminus \overline{(B(y,\ell_y)\cup B(z,\ell_z))}
$$
is not empty. By approximate controllabilty, we can choose $f \in  C^\infty_0(S_\epsilon(p,\ell_p))$ such that $\|u^f(\cdot,T)\|_{L^2(U)}>0$. By finite speed of wave propagation it holds that 
\ba
& &\inf \{\|u^{f}(\cdot,T)-u^h(\cdot,T)\|_{L^2(N)}: h\in C^\infty_0(S_\epsilon(y,\ell_y)\cup S_\epsilon(z,\ell_z))\}
\\
& &{\ge \|u^f(\cdot,T)|_U\|_{L^2(U)}}>0.
\ea
Therefore \eqref{cond 4} is not true. This proves the claim.
%
\end{proof}

Next we consider the interior distance functions related to a point $p$ (in the unknown part of the manifold) that gives distances to points $z$ in the set $\mathcal O^{(i)}$.

\begin{proposition}
\label{finding the distance function}
Let $({\mathcal M}^{(1),rec},{\mathcal M}^{(2),rec},\Phi^{rec})$
be an admissible triple, and ${\mathcal O}^{(i)}=B^{(i)}(q_i,R)$, $i=1,2$ be 
the balls of radius $R$,
centered at $q_i$ in ${\mathcal M}^{(i),rec}$,
satisfying $\Phi^{rec}(q_1)=q_2$ and $T$  satisfy (\ref{condition for T}).
Assume that the source-to-solution maps  $V^{2T}_{\mathcal O^{(1)},+}$
and $V^{2T}_{\mathcal O^{(2)},+}$ are $\Phi^{rec}$-related.
Let $z_i \in \mathcal O^{(i)}$, $i=1,2$, 
be such that
$z_2=\Phi^{rec}(z_1)$.  
Let $\xi_i \in T_{q_i}\mathcal M^{(i)}$
be the unit vector, 
 $R\le \widetilde r <R+ T$ and $p_i = \gamma_{q_i, \xi_i}(\widetilde r)$. Then 
 \beq\label{distances are the same}
 d_{g^{(1)}}(p_1,z_1)=d_{g^{(2)}}(p_2,z_2).
 \eeq
\end{proposition}
\begin{proof} Note that we have $R+T< \tau^{(i)}_{cut}(q_i, \xi)$ for $i=1,2$  and hence $\widetilde r <\tau^{(i)}_{cut}(q_i, \xi)$.

Let $s\in (0,R)$ be such that 
$\gamma_{{{q}}_i,\xi_i}([0,s])\subset \mathcal O^{(i)}$. We denote $x_i=\gamma_{{{q}}_i,\xi_i}(s)$. 

Let $r:=\widetilde{r}-s$, ${{t}}>0$, and
 $0<\epsilon<\e_0=\min(R+T-\tilde r,s)$. Below we consider the balls
$B(x_i,r+\epsilon)\subset \mathcal M^{(i)}$  
and
the sets
$S_\epsilon(x,\ell)\subset \mathcal M^{(i)}$  defined in \eqref{set S}.

 By Lemma \ref{crossing balls and functions}
the balls $B(x_i,r+\epsilon)\subset \mathcal M^{(i)}$ satisfy  the inclusion 
\begin{equation}
\label{crossing balls with closure}
B(x_i,r+\epsilon)\subset \overline{B({{q}}_i,r+s) \cup B(z_i,{{t}})}
\end{equation}
 if and only if the equation \eqref{cond 4} is valid with $p=x_i,y={{q}}_i,z=z_i\in \mathcal M^{(i)},$ $\ell_p=r+\epsilon$, $\ell_y=r+s$ and $\ell_z={{t}}$. 
  
  Let us consider functions $f^i\in C^\infty_0(S_\epsilon(x_i,\ell_p))$, $i=1,2$
  and {the sequences} $ (f^i_j)_{j=1}^\infty \subset C^\infty_0(S_\epsilon({{q}}_i,\ell_y)\cup S_\epsilon(z_i,\ell_z))$, $i=1,2,$ 
      that satisfy 
   \beq\label{pushes}
    f^1=(\Phi^{rec})^*f^2,\quad f^1_j=(\Phi^{rec})^*f^2_j.
   \eeq

  Using the Blagovestchenskii identity, {Theorem \ref{th_Blago},} we can
compute  for functions $f^i$ and $f^i_j$ 
the norm $ \|u^{f^i}(\cdot,T)-u^{f^i_j}(\cdot,T)\|_{L^2(\mathcal M^{(i)})}$ using
the local source-to-solution data  $(\mathcal O^{(i)},{g|_{\mathcal O^{(i)}},}V^{2T}_{\mathcal O^{(i)},+})$. Hence,
as the source-to-solution maps  $V^{2T}_{\mathcal O^{(1)},+}$
and $V^{2T}_{\mathcal O^{(2)},+}$ are $\Phi^{rec}$-related,
 \eqref{pushes} implies
   \beq\label{pushes2}
\|u^{f^1}(\cdot,T)-u^{f^1_j}(\cdot,T)\|_{L^2(\mathcal M^{(1)})}
=\|u^{f^2}(\cdot,T)-u^{f^2_j}(\cdot,T)\|_{L^2(\mathcal M^{(2)})}.
\eeq
Hence for functions
$f^i\in C^\infty_0(S_\epsilon(x_i,\ell_p))$, $i=1,2$
  and the sequences $ (f^i_j)_{j=1}^\infty \subset C^\infty_0(S_\epsilon({{q}}_i,\ell_y)\cup S_\epsilon(z_i,\ell_z))$, $i=1,2,$ satisfying 
  \eqref{pushes}  we have
   \beq\label{pushes3}
& &\lim_{j\to \infty} \|u^{f^1}(\cdot,T)-u^{f^1_j}(\cdot,T)\|_{L^2(\mathcal M^{(1)})}=0
\quad\hbox{if and only if}\\ \nonumber
& &\lim_{j\to \infty} \|u^{f^2}(\cdot,T)-u^{f^2_j}(\cdot,T)\|_{L^2(\mathcal M^{(2)})}=0.
\eeq
Thus by  Lemma \ref{crossing balls and functions},
we have
\beq
\label{pushes crossing balls with closure}
& &B(x_1,r+\epsilon)\subset \overline{B({{q}}_1,r+s) \cup B(z_1,{{t}})}
\quad\hbox{if and only if}\\ \nonumber
& &B(x_2,r+\epsilon)\subset \overline{B({{q}}_2,r+s) \cup B(z_2,{{t}})},
\eeq
see \eqref{crossing balls with closure}.


 Let us define, for $i=1,2$,
 $$
 {{t}}_i^\ast:=\inf\{{{t}}>0: \textrm{Formula }\eqref{crossing balls with closure} \textrm{ is valid  with index $i$ for some } \epsilon\in (0,\e_0)\}.
$$
By \eqref{pushes crossing balls with closure}, we have $ {{t}}_1^\ast= {{t}}_2^\ast$. Due to this,  we denote below
$ {{t}}_1^\ast= {{t}}_2^\ast= {{t}}^\ast$.
 
Next we will show that
$$
d_{g^{(i)}}(p_i,z_i)={{t}}^\ast.
$$

Suppose $t$ is such that \eqref{crossing balls with closure} is valid
for some $\e>0$, {that we next denote by $\e_0$. Then, we have that \eqref{crossing balls with closure}  is
valid for all $\e\in (0,\e_0)$}.
For any $\e\in (0,\e_0)$,  let $$y_i^\e=\gamma_{x_i,\eta_i}(r+\e)=
\gamma_{{{q}}_i,\xi_i}(s+r+\e),\quad \eta_i=
\dot \gamma_{{{q}}_i,\xi_i}(s).$$ 
 Since
the right hand side of \eqref{crossing balls with closure} is a closed 
set, we have that $y^\e_i\in  \overline{B({{q}}_i,r+s) \cup B(z_i,{{t}})}$.
As
$x_i=\gamma_{{{q}}_i,\xi_i}(s)$ 
and we assumed that $$s+r+\e=\tilde r+\e< R+T\leq \tau^{(i)}_{cut}(q_i,\xi)$$ for $i=1,2$,
 it holds that $$d_{g^{(i)}}(y^\e_i,{{q}}_i)=r+s+\e>r+s,$$ so that 
  $y^\e_i \not \in \overline{B({{q}}_i,r+s)}$. As  \eqref{crossing balls with closure}
holds, we need to have
$y^\e_i  \in \overline{B(z_i,{{t}})}$. Thus $t\ge d_{g^{(i)}}(y^\e_i,z_i)$.
As this holds 
for all sufficiently small $\e>0$ and $y^\e_i \to p_i$ as $\e\to 0$,
we obtain $t\ge d_{g^{(i)}}(p_i,z_i)$.
This yields that 
${{t}}^\ast\ge d_{g^{(i)}}(p_i,z_i)$.

Next, suppose that there exists ${{t}} \in (d_{g^{(i)}}(p_i,z_i), {{t}}^\ast)$. Then for any sufficiently small $\epsilon>0$ the formula \eqref{crossing balls with closure} is not valid. Choose for every $k\in \N$  a point  
\beq\label{point pk}
p^k_i \in B(x_i,r+1/k) \setminus \overline{B({{q}}_i,r+s) \cup B(z_i,{{t}})}.
\eeq
By compactness of $\overline{B(x_i,r+1)}$ we may assume that $p^k_i \rightarrow \widetilde p_i \in \p B(x_i,r)$ as $k\rightarrow \infty$.

Next we will show that $\widetilde p_i=p_i$.
As $\widetilde p_i \in \p B(x_i,r)$ and $d_{g^{(i)}}(x_i,{{q}}_i)\leq s$, we have by triangle inequality $d_{g^{(i)}}(\widetilde p_i,{{q}}_i)\leq s+r$. Let $\alpha_i$ be a minimizing geodesic from $x_i$ to $\widetilde p_i$. Suppose first that $\alpha_i$ is not the geodesic continuation of the geodesic segment $\gamma_{{{q}}_i,\xi_i}([0,s])$. Since a curve $\gamma_{{{q}}_i,\xi_i}([0,s])\cup \alpha_i$ has a length $s+r$ and it is not smooth at $x_i$, it must hold that $d_{g^{(i)}}(\widetilde p_i,{{q}}_i)<s+r$. Then for sufficiently large $k$, we have $d_{g^{(i)}}( p^k_i,{{q}}_i)<s+r$ that is not possible because of 
\eqref{point pk}. This show that 
  $\alpha_i$ has to the geodesic continuation of segment $\gamma_{{{q}}_i,\xi_i}((0,s))$. This yields that  $\widetilde p_i=\gamma_{{{q}}_i,\xi_i}(s+r)=p_i$.


Since $p_i \in B(z_i,{{t}})$ we get a contradiction with 
the assumptions that
$p^k_i \rightarrow \widetilde p_i =p_i$ as $k\rightarrow \infty$
and $p^k_i\not \in \overline{B(z_i,{{t}})}.$ 
Therefore interval $(d_{g^{(i)}}(p_i,z_i),  {{t}}^\ast)=\emptyset$. This shows that ${{t}}^\ast =d_{g^{(i)}}(p_i,z_i)$ for both $i=1,2$.
Hence, $d_{g^{(1)}}(p_1,z_1)=d_{g^{(2)}}(p_2,z_2)$.
\end{proof}

For $i=1,2$, let $p_i\in \mathcal N^{(i)}$ and $z_i \in \mathcal O^{(i)}$
be such that $p_2=\Phi^{rec}(p_1)$  and $z_2=\Phi^{rec}(z_1)$.

 By definition of the set $ \mathcal N^{(i)}$,
we see that there exist 
 vectors $\xi_i\in S_{{{q}}_i}\mathcal M^{(i)}$ such that $p_i=\gamma_{{{q}}_i,\xi_i}(\widetilde{r})$, for some $\widetilde{r}\in [R,R+T)$. By Proposition  \ref{finding the distance function} 
 we see that \eqref{distances are the same} is valid.
 Therefore Theorem \ref{th_geom} is proved.\hfill $\square$
\bigskip

Consider {points}  $p_1\in \mathcal N^{(1)}$ and $p_2=\Phi^{rec}(p_1)$.
 Let $A(p_i)$ be the set of the points 
$y_i\in \p{\mathcal O^{(i)}}$ such that there are 
no points $z_i\in {\mathcal O}^{(i)}$ such that 
$d^{(i)}(p_i,z_i)=d^{(i)}(p_i,y_i)+d^{(i)}(y_i,z_i)$.  Note that if $y_i\in A(p_i)$,
then {any shortest geodesics in $\mathcal M^{(i)}$ from $p_i$ to $y_i$ does} not 
intersect ${\mathcal O}^{(i)}$.  Let us define the distance function
$d_{\mathcal N^{(i)}}(x_i,p_i)$ analogously to \eqref{extension of metric}, that is, as the infimum of length of paths connecting $x_i$  to $p_i$
in $\mathcal N^{(i)}$. As the ball 
$\mathcal O^{(i)}$  is convex, we see that for any $x_i\in \p\mathcal O^{(i)}$  and $p_i\in {\mathcal N}^{(i)}$
we have
\beq\label{boundary distance functions}
d_{\mathcal N^{(i)}}(x_i,p_i)=\inf_{y_i\in A(p_i)}
d_{\p\mathcal O^{(i)}}(x_i,y_i)+d_{\mathcal M^{(i)}}(y_i,p_i).
\eeq
{Here, $d_{\p\mathcal O^{(i)}}(x_i,y_i)$  is the intrinsic distance of points $x_i$ and $y_i$
along the boundary $\p\mathcal O^{(i)}$.}
By Theorem \ref{th_geom}, this and
$d_{\p\mathcal O^{(1)}}(x_1,y_1)=d_{\p\mathcal O^{(2)}}(\Phi^{rec}(x_i),
\Phi^{rec}(y_1))$ for $x_1,y_1\in \p\mathcal O^{(1)}$
yield that the following:

\begin{corollary}\label{coro on bdf}
Under assumptions of 
Theorem \ref{th_geom}, we have for all $p_1\in \mathcal N^{(1)}$ and $x_1\in  \p \mathcal O^{(1)}$ that
\beq\label{boundary distance functions in N}
d_{\mathcal N^{(2)}}(\Phi^{rec}(x_1),\Phi^{rec}(p_1))=d_{\mathcal N^{(1)}}(x_1,p_1).
\eeq
\end{corollary}}

Next we consider how the source-to-solution operator causes bounds for the cut locus functions.

\begin{lemma}
\label{Finding cut dist}
Let  $0<s<R$. {Assume that $r>0$  is such that  
 \beq\label{tau cut small}
0<r+s<\tau^{(i)}_{sing}(q_i).
\eeq
Let $\xi_i\in S_{q_i}\mathcal M^{(i)}$, and
$x_i=\gamma_{q_i,\xi_i}(s)$.

(i) If $\tau^{(i)}_{cut}(q_i,\xi_i)<s+r$,} then 
\begin{eqnarray}
\label{Ball condition for cut distance}
&&\textrm{there exists }\epsilon >0 \textrm{ such that } B(x_i,r+\epsilon)\subset \overline{B(q_i,s+r)}.
\end{eqnarray}

{(ii)} If \eqref{Ball condition for cut distance} is valid then $\tau^{(i)}_{cut}(q_i,\xi_i)\leq s+r.$

\end{lemma}
\begin{proof}
Denote $p_i=\gamma_{q_i,\xi_i}(s+r)$. 

\noindent
Suppose that \eqref{Ball condition for cut distance} is valid. Let $\delta\in (0,\epsilon)$ and consider a point 
$$
z_i=\gamma_{q_i,\xi_i}(s+r+\delta)\in B(x_i,r+\epsilon).
$$
By \eqref{Ball condition for cut distance}, $d^{(i)}(z_i,q_i)\leq s+r$. Thus $\tau^{(i)}_{cut}(q_i,\xi_i)< s+r+\delta$. 
Since $\delta$ was arbitrary we have $\tau^{(i)}_{cut}(q_i,\xi_i)\leq s+r$.

Suppose that {$\tau^{(i)}_{cut}(q_i,\xi_i)< s+r$}. We show first that 
\begin{equation}
\label{eq: closure inclusion}
\overline{B(x_i,{r})}\subset B(q_i,s+r).
\end{equation}
By triangle inequality it suffices to show that  $\p B(x_i,{ r})\subset B(y_i,s+r)$. Let $z_i \in \p B(x_i,r)$. By triangle inequality $d^{(i)}(z_i,q_i)\leq s+r$. Let $\alpha$ be a minimizing geodesic from $x_i$ to $z_i$. Suppose first that $\alpha$ is not the geodesic continuation of the geodesic segment $\gamma_{q_i,\xi_i}([0,s])$. Since a curve $\gamma_{q_i,\xi_i}([0,s])\cup \alpha$ has a length $s+r$ and it is not smooth at $x_i$, it must hold that $d^{(i)}(z_i,q_i)<s+r$. { Second, suppose} $\alpha$ is the geodesic continuation of segment $\gamma_{q_i,\xi}((0,s))$, then $z_i=\gamma_{q_i,\xi}(s+r)=p_i$. Since {$\tau^{(i)}_{cut}(q_i,\xi_i)< s+r$}, it holds that $d^{(i)}(q_i,p_i)<s+r$. {As $z_i \in \p B(x_i,r)$ above is arbitrary, the formula} \eqref{eq: closure inclusion} follows. Therefore dist$(\p B(x_i,r), \p B(q_i,s+r))>0$ and \eqref{Ball condition for cut distance} is valid.

\end{proof}

\begin{proposition}
\label{reconstruction of cut distance from data}
Let $({\mathcal M}^{(1),rec},{\mathcal M}^{(2),rec},\Phi^{rec})$
be an admissible triple, and ${\mathcal O}^{(i)}=B^{(i)}(q_i,R){ \subset {\mathcal M}^{(i),rec}}$, $i=1,2$ be 
the balls of radius $R$,
centered at $q_i$ in ${\mathcal M}^{(i),rec}$,
satisfying $\Phi^{rec}(q_1)=q_2$.
Assume that the source-to-solution maps  $V_{\mathcal O^{(1)},+}$
and $V_{\mathcal O^{(2)},+}$ are $\Phi^{rec}$-related.
Assume that  
\beq\label{tau cut small2}
& &\tau^{(1)}_{cut}(q_1)<\min(\tau^{(1)}_{sing}(q_1),\tau^{(2)}_{sing}(q_2)).
\eeq
Then
\beq
\tau^{(2)}_{cut}(q_2)\ge \tau^{(1)}_{cut}(q_1).
\eeq
\end{proposition}
\begin{proof}
Let $\xi_1\in S_{q_1}{\mathcal M}^{(1),rec}$ and 
$\xi_2=(\Phi^{rec})_*\xi_1$ and let $r<\tau^{(1)}_{cut}(q_1)$. 
Let $0<s<R$  and denote $x_i=\gamma_{q_i,\xi_i}({s})$. 
Also, let
\ba
& &a_i(\xi_i)=
\\
& &\inf\{s+r>0: \  \hbox{Formula \eqref{Ball condition for cut distance} with $i$ and $x_i$ holds for  some $r>0$  and $s\in (0,R)$}\}.
\ea
Then by
{Lemma \ref{Finding cut dist},} $a_1(\xi_1)\geq  \tau^{(1)}_{cut}(q_1).$
Note that here $\xi_1$ is an arbitrary unit vector and thus
$\inf_{\xi_1}a_1(\xi_1)\geq  \tau^{(1)}_{cut}(q_1).$

Choose $\epsilon\in (0,{R-s})$. Then, 
$$
B(q_i,\epsilon)\cup B(x_i,\epsilon)\subset \mathcal O^{(i)},\quad i=1,2.
$$

Now we can proceed as in the proof of Proposition \ref{finding the distance function}.

By applying the Blagovestchenskii identity \eqref{Blagovestchenskii identity}
and the fact that the source-to-solution maps  $V^{2T}_{\mathcal O^{(1)},+}$
and $V^{2T}_{\mathcal O^{(2)},+}$ are $\Phi^{rec}$-related with any $T>0$, we see that
formula \eqref{Ball condition for cut distance} holds for $r>0$ and $s\in (0,R)$ with the index $i=1$ if and only it holds with the index $i=2$.
Thus by taking $z=y=q_i$, $x=x_i$, $\ell_y=r+s=\ell_z$, $\ell_x=r+\epsilon$  and applying Lemma \ref{crossing balls and functions}, 
with $i=1$  and  $i=2$ we see that $a_1(\xi_1)=a_2(\xi_2)$.
 By
{Lemma \ref{Finding cut dist}}, $\tau^{(2)}_{cut}(q_2)\ge \inf_{\xi_2}a_2(\xi_2)=
\inf_{\xi_1}a_1(\xi_1).$
This yields the claim.
%

\end{proof}



\begin{theorem} \label{recognition}
Let $({\mathcal M}^{(1),rec},{\mathcal M}^{(2),rec},\Phi^{rec})$
be an admissible triple and 
${\mathcal O}^{(i)}=B^{(i)}(q_i,R) \subset {\mathcal M}^{(i),rec}$, $i=1,2,$
be a ball centered at $q_i$ and radius $R$ satisfying (\ref{eq: basic assumptions})
and (\ref{eq: basic assumptions2}). 
Then $\tau^{(1)}(q_1)=\tau^{(2)}(q_2)$.
Using the notation $\tau=\tau^{(1)}(q_1)$, then, for 
\ba
   {\widetilde{\mathcal M}}^{(i),rec}=
 {\mathcal M}^{(i),rec} \cup {\mathcal M}^{(i)}_{\mathcal O}( \p \O^{(i)},\tau-R),
 \quad i=1,2,
  \ea
  there is a map
  $  {\widetilde \Phi}^{rec}:{\widetilde{\mathcal M}}^{(1),rec}
 \to
 {\widetilde{\mathcal M}}^{(2),rec}
$ 
which is an extension
of  $\Phi^{rec}$. Moreover, the triple 
$( {\widetilde{\mathcal M}}^{(1),rec}, 
{\widetilde{\mathcal M}}^{(2),rec},  {\widetilde \Phi}^{rec})$ is  admissible.
  \end{theorem}
\noindent {\ntekst Note that $B^{(i)}(q_i,\tau)={\mathcal M}^{(i)}_{\mathcal O}( \p \O^{(i)},\tau-R)\cup B^{(i)}(q_i,R)$.}

\noindent
{\bf Proof.} Assume opposite to the claim  that we would have
{$\tau^{(1)}(q_1) > \tau^{(2)}(q_2)$.} 
Let 
\beq \label{5.21a}
a=\tau^{(1)}(q_1)-R,\quad b=\tau^{(2)}(q_2)-R, \quad 0<T<b<a.
\eeq

{\btext 
Recall that ${\mathcal N}^{(i)}=B^{(i)}(q^{(i)},T+R)\setminus O^{(i)}$.

Let us extend the manifolds ${\mathcal M}^{(i)}_{\mathcal O}( \p \O^{(i)},T)\subset
B^{(i)}(q^{(i)},T+R)$. First, observe that as $T+R<\tau^{(1)}(q_i)$,
$B^{(i)}(q^{(i)},T+R)$  is a diffeomorphic to a ball of an Euclidean space.

 On the surfaces
$\Gamma^{(i)}=\p 
B^{(i)}(q^{(i)},R)$ 
the boundary distance function 
$r^{(i)}_{p_i}\in C(\Gamma^{(i)})$ corresponding to the point
 $p_i\in {\mathcal N}^{(i)}$ is
$$
r^{(i)}_{p_i}(x)=d_{\mathcal N^{(i)}}(x,p_i),\quad \hbox{for }x\in \Gamma^{(i)}.
$$
By \eqref{boundary distance functions in N},
we have for all $p_1\in {\mathcal N}^{(1)}$
and $p_2=\Phi^{rec}(p_1)$
\beq\label{bdf coincide}
r^{(2)}_{p_2}(\Phi^{rec}(x))=r^{(1)}_{p_1}(x),\quad \hbox{for }x\in \Gamma^{(1)}.
\eeq
%
{For $x\in  \Gamma^{(i)}$ and $p\in \mathcal N^{(i)}$, denote
$$r^{(i),T}_{p}(x)=\max(r^{(i)}_{p}(x),T)$$
and let $R_{\Gamma^{(i)}}^T:\mathcal N^{(i)}\to C(\Gamma^{(i)})$ be the 
map $$R_{\Gamma^{(i)}}^T(p)=r^{(i),T}_{p}(\,\cdotp).$$
These and related boundary distance functions have been considered in \cite{deHoop-Holman1,deHoop-Holman2,deHoop-Saksala,Ivanov2018,
KatsudaKurylevLassas}.
Naturally,
the functions $r^{(i)}_{p}\in C(\Gamma^{(i)})$, $p\in {\mathcal N^{(i)}}$  determine the range of $R_{\Gamma^{(i)}}^T$, that is, the family functions 
$${R_{\Gamma^{(i)}}^T(\mathcal N^{(i)})}=\{r^{(i),T}_{p}\in C(\Gamma^{(i)}):\ p\in {\mathcal N}^{(i)}\}.$$
By \cite{KKL01}, Subsection 4.2.9, the family ${ R_{\Gamma^{(i)}}^T(\mathcal N^{(i)})}$ of  functions} determine
the topological and differentiable {type of the manifold} $ {\mathcal N}^{(i)}$
and the isometry type of  {the Riemannian manifold} $( {\mathcal N^{(i)}},g^{(i)}|_ {{\mathcal N}^{(i)}})$.
Moreover, when we identify the sets $\Gamma^{(1)}$  and
$\Gamma^{(2)}$ using the map $\Phi^{rec}$ and denote by
$\nu_i$  the unit exterior normal vector of $\Gamma^{(i)}$ and 
by $E_i:\Gamma^{(i)}\times [0,T)\to {\mathcal N^{(i)}}$
the normal exponential map,
\ba
E_i(x,t)=\gamma_{x,\nu_i(x)}(s),
\ea
 we see that when $J(x,s)=(\Phi^{rec}(x),s)$ {maps $J:\Gamma^{(1)}\times [0,T)\to \Gamma^{(2)}\times [0,T)$, then  the ``collar'' map}
 \beq\label{collar map}
  {\widetilde \Phi}_c=E_2\circ J\circ E_1^{-1}:{\mathcal N^{(1)}}\to {\mathcal N^{(2)}}
 \eeq
is a diffeomorphism and an isometry.

For the convenience of the reader, let us sketch the idea of the above cited construction
in  \cite{KKL01}, Subsection 4.2.9. There, the map   $ {\widetilde \Phi}_c:{\mathcal N^{(1)}}\to {\mathcal N^{(2)}}$  
gives
a diffeomorphism that can be used to identify the sets
${\mathcal N^{(1)}}$ and ${\mathcal N^{(2)}}$   and their
differentiable structures. Moreover, the gradients of the distance functions,
$$\nabla_p (r^{(i),T}_{p}(x_i))=\nabla_p (d^{(i)}(p,x_i))\in T_p\mathcal M^{(i)}$$   are unit length vectors.
When $x_i$ moves on the boundary near the closest point to $p$,
these unit vectors run over an open subset of the unit sphere in $T_p\mathcal M^{(i)}$.
As the differentials of the distance functions satisfy
(\ref{boundary distance functions in N}), we see that
\beq
(\Phi^{rec})_*(\nabla_pr^{(1),T}_{p}(x_1)){=}
\nabla_q r^{(2),T}_{q}(\Phi^{rec}(x_1))\bigg|_{q=
\Phi^{rec}(p)}.
\eeq
This implies that the linear map $(\Phi^{rec})_*$ in $T_p\mathcal M^{(1)}$ maps an open subset of unit
vectors {in $T_p\mathcal M^{(1)}$ to  unit vectors in $T_q\mathcal M^{(2)}$,} implying that 
 $ ({\widetilde \Phi}_c)^*g^{(2)}= g^{(1)}$.
 {This implies that}  $ {\widetilde \Phi}_c:{\mathcal N^{(1)}}\to {\mathcal N^{(2)}}$  
 is an isometry.   
 
 Now we return to the proof of the claim.
 As above  $T<b$ is arbitrary, 
{the fact that the map \eqref{collar map} is an isometry} implies that $ {\mathcal M}^{(1)}_{\mathcal O}(\p{\mathcal O}^{(1)},b)$ and $
 {\mathcal M}^{(2)}_{\mathcal O}(\p{\mathcal O}^{(2)},b)$  are 
isometric. We can extend this isometry to  the sets ${\mathcal O}^{(1)}$ and ${\mathcal O}^{(2)}$, and see
 there is an isometry 
 \beq\label{eq: new isometry}
 {\widetilde \Phi}:  B^{(1)}(q^{(1)},b+R) \to
 B(q^{(2)},b+R), 
 \eeq
 where we recall that $B^{(i)}(q,r)\subset {\mathcal M}^{(i)}$  denote the balls.
 {Let $b' <b$ be so
   small that $B^{(i)}(q^{(i)},b'+R) \subset
  {\mathcal M}^{(i),rec} $.}
{By the assumptions we made in the claim, }
 \beq \label{5.25b} 
 {\widetilde \Phi}(x)=\Phi^{rec}(x), \quad x \in  B^{(1)}(q^{(1)},b'+R).
 \eeq

By (\ref{eq: time domain Greens functions related}), the time domain Green's functions  $G^{(i)}(x, y,t)$, $i=1,2,$ satisfy
 the relation 
  \beq\label{eq: time domain Greens functions related2}
& &G^{(2)}({\widetilde \Phi}(x), {\widetilde \Phi} (y),t)=G^{(1)}(x, y,t),
\eeq 
in $\{(x,y,t)\in B^{(1)}(q^{(1)},b'+R)^2\times \R;\ x\not =y\}$.
Similarly, to the proof of Lemma \ref{lem. N1 admissible}, 
we use Tataru's unique continuation first in the $x$ and $t$ variables
to see that (\ref{eq: time domain Greens functions related2}) is valid in
$\{(x,y,t)\in B^{(1)}(q^{(1)},b+R)\times B^{(1)}(q^{(1)},b'+R)\times {\mathbb R};\ x\not =y\}$. As $G^{(1)}(x, y,t)=G^{(1)}(y, x,t)$, we can then 
use Tataru's unique continuation  in the $y$ and $t$ variables
to see that (\ref{eq: time domain Greens functions related2}) is valid in
$\{(x,y,t)\in B^{(1)}(q^{(1)},b+R)^2\times {\mathbb R};\ x\not =y\}$.
%
Thus {the triple} $(B^{(1)}(q^{(1)},b+R),
B^{(2)}(q^{(2)},b+R), {\widetilde \Phi})$
 is admissible. 

 {\ntekst Using (\ref{eq: new isometry}) and (\ref{eq: time domain Greens functions related2}), 
it follows from Remark \ref{Remark:3} that $\Phi^{rec}$ can be extended by ${\widetilde \Phi}$
as ${\widetilde \Phi}^{rec}$,
\beq \label{5.23a}
& &{\widetilde \Phi}^{rec}: {\widetilde \M}^{(1), rec} \rightarrow
{\widetilde \M}^{(2), rec};
\\\nonumber
& & {\widetilde \M}^{(i), rec}= \M^{(i), rec} 
\cup B^{(i)}(q^{(i)},b+R)
\eeq

}

Now, consider the case when 
\beq\label{cuts smaller}
\tau^{(i)}_{cut}(q_i)\le \tau^{(i)}_{sing}(q_i)
\eeq
for both $i=1,2$. {Then} we see using Proposition \ref{reconstruction of cut distance from data} that $\tau^{(1)}_{cut}(q_1)\leq \tau^{(2)}_{cut}(q_2)$, and using 
Proposition \ref{reconstruction of cut distance from data} with roles of $i=1$ and $i=2$  exchanged that $\tau^{(1)}_{cut}(q_1)\geq \tau^{(2)}_{cut}(q_2)$. Hence, we have $\tau^{(1)}(q_1)=\tau^{(2)}(q_2)$.
Then, applying Corollary \ref{coro on bdf} with all $T<\tau^{(1)}(q_1)=\tau^{(2)}(q_2)$ and applying results of 
 \cite{KKL01}, Subsection 4.2.9 {as described above}, yield the claim in the case \eqref{cuts smaller}. }

Recall that, by our assumption, {$\tau^{(1)}(q_1) > \tau^{(2)}(q_2)$}, that is, $a>b$.
{\nntext Due to  (\ref{eq: basic assumptions}),  this implies that
$B^{({\mtext 1})}_{sing}\cap \p B^{({\mtext 1})}(q^{({\mtext 1})},b+R)\not =\emptyset$.
Thus it remains to consider the case when
\beq
\tau^{({2})}(q_{2})=d^{({2})}(q_{2}, {\mathcal M}^{({ 2})}_{sing})<\tau^{({1})}(q_{1})\leq d^{({1})}(q_{1}, {\mathcal M}^{({1})}_{sing}).
\eeq
Next we show that this is not possible.
 

%

As the mapping
(\ref{eq: new isometry}) is an isometry between
$ B^{(1)}(q^{(1)},b+R)$
and $ B^{(2)}(q^{(2)},b+R)$,
we see that,
\ba
\min (\tau^{(1)}(q^{(1)}),b+R)=\min (\tau^{(2)}(q^{(2)}),b+R).
\ea

Next, any point $p^{(i)}\in \overline B^{(i)}(q^{(i)},b+R)$ 
 can be written in the form
 ${p^{(i)}=}\gamma^{(i)}_{q^{(i)}, \xi}(t)$ where $\xi$ is a unit vector and 
 $t  \leq 
\min (\tau^{(1)}(q^{(1)}),b+R)$.
 
  
Let 
\beq \label{5.??}
p^{(2)} \in {\mathcal M}^{(2)}_{sing}\cap \p B^{(2)}(q^{(2)},b+R).
\eeq
By the above, there is 
a point $\xi^{(2)}\in  S_{q^{(2)}}{\mathcal M}^{(2)}$ such that $p^{(2)}= \gamma^{(2)}_{q^{(2)}, \xi^{(2)}}(b+R)$.
Let $\xi^{(1)}=d(\Phi^{rec})^{-1}(\xi ^{(2)})$ and consider $p^{(1)}=\gamma^{(1)}_{q^{(1)}, \xi^{(1)}}(b+R)$.
Since {$\tau^{(1)}(q_1) > \tau^{(2)}(q_2)$}, that is, $a>b$, we have}
$p^{(1)} \notin \M^{(1)}_{sing}.$

Let 
\ba
p^{(i)}_{\e}:= \gamma_{y^{(i)}, \xi^{(i)}}(b+R-2\e),\quad \e >0,\ i=1,2.
\ea

We denote by ${\widetilde O}^{(i)}_\e= B^{(i)}(p^{(i)}_{\e}, \e)$  the 
metric ball in  ${\mathcal M}^{(i)}$ of radius $\e$.
By using (\ref{5.23a}) and choosing
$\e>0$ to be  small,  ${\widetilde \O}^{(i)}_\e$ satisfy the conditions of Lemma
\ref{lem. determination of Green} with ${\widetilde \O}^{(i)}_\e$ instead of $\O^{(i)}$.


Then, Lemma \ref{lem. determination of Green} implies that
\beq \label{7.3.2}\quad \quad
& &\hbox{Source-to-solution operators {$V^{2T}_{\tilde {\mathcal O}^{(i)},+}$ for }}  - \Delta^{(i)} \,\, \hbox{on} \,\,{\widetilde \O}_\e^{(i)},\,\, i=1,2,\,\,\\
\nonumber
& &
\hbox{are}\,\, 
{{} \Phi}^{reg}\hbox{-related}.
\eeq
Here $- \Delta^{(i)}$, s the Laplace operator
associated with ${{} {\mathcal M}}^{(i)}$. 
Equation (\ref{7.3.2}) together with Lemma \ref{lem: volumes} imply that
\ba
\hbox{Vol}^{(1)}({{} {\mathcal M}}^{(1)}({\widetilde \O}_\e^{(1)},r-\e))=
\hbox{Vol}^{(2)}({{} {\mathcal M}}^{(2)}({\widetilde \O}_\e^{(2)},r-\e))
\ea
when $r>0$. 
Therefore,
\beq \label{7.3.3}
\hbox{Vol}^{(1)}(B^{(1)}(p^{(1)}_\e, r))=\hbox{Vol}^{(2)}(B^{(2)}(p^{(2)}_\e, r)).
\eeq
Next, we observe that as
$d^{(i)}(p^{(i)}_{\e}, p^{(i)}) \leq 2\e,$ we have
$$
B^{(i)}(p^{(i)},r-2\e) \subset B^{(i)}(p^{(i)}_{\e},r) \subset 
B^{(i)}(p^{(i)},r+2\e),\quad \hbox{for }r>2\e.
$$
Thus,
by the continuity of the volume,
\ba
\hbox{Vol}^{(i)}(B^{(i)}(p^{(i)},r))
 = \lim_{\e \to 0} \hbox{Vol}^{(i)}(B^{(i)}(p^{(i)}_{\e},r)).
\ea
Together with (\ref{7.3.3}), this implies that, for $r>0$,
\beq \label{E2}
\hbox{Vol}^{(1)}(B^{(1)}(p^{(1)},r))=\hbox{Vol}^{(2)}(B^{(2)}(p^{(2)},r)).
\eeq

Let us now consider the conic  coordinates of $\mathcal M^{(i)}$ and the volume factor,
see (\ref{volume factor}) and (\ref{volume factor2}).
It then follows from (\ref{E2}), that
$$
\Lambda(p^{(i)})= 
 \lim_{r \to 0} \frac{1}{\hbox{Vol}_{\R^n}(B_{\R^n}(0,r))} \left[\hbox{Vol}^{(i)}(B^{(i)}(p^{(i)},r) )    \right] 
 $$
and thus we have
\beq \label{7.3.5}
\Lambda(p^{(1)})=\Lambda(p^{(2)}).
\eeq
Denote $\Lambda^{(1)}=\Lambda(p^{(i)})$.

Note that, 
if $p^{(i)}\in {\mathcal M}_{sing}^{(i)}$ we have $\Lambda^{(i)}\neq 1$
and if
  $p^{(i)}\in {\mathcal M}_{reg}^{(i)}$ then $\Lambda^{(i)}=1$.
  As we {assumed} that $a=\tau^{(1)}(q_1)-R>b=\tau^{(2)}(q_2)-R$, we have 
   $p^{(1)}\in {\mathcal M}_{reg}^{(1)}$ 
   and thus $\Lambda^{(1)}=1$. Hence, we also have   $\Lambda^{(2)}=1$, 
 and thus $p^{(2)}\in  {\mathcal M}_{reg}^{(2)}$, contradicting (\ref{5.??}).
%
{Thus the above assumption $\tau^{(1)}(q_1) > \tau^{(2)}(q_2)$ led to a contradiction. By changing
roles of indexes $1$ and $2$, the above considerations show} that we must have
\ba
\tau^{(1)}(q_1)=\tau^{(2)}(q_2).
\ea

Thus using Corollary \ref{coro on bdf} with all $T<\tau^{(1)}(q_1)=\tau^{(2)}(q_2)$ and applying results of 
 \cite{KKL01}, Subsection 4.2.9
 we prove the claim of the theorem.
\qed

\medskip
%

Let $\mathcal A$ be the collection of admissible triples
$({\mathcal W}^{(1)},\, {\mathcal W}^{(2)}, \Phi  )$ such that 
$W \subset {\mathcal W}^{(i)}$ {and
${\mathcal W}^{(i)}{\subset  \M^{(i)}_{reg}}$ are connected open sets for both indexes} $ i=1, 2$.
We define a partial order on $\mathcal A$ by setting
$({\mathcal W}^{(1)},\, {\mathcal W}^{(2)}, \Phi  )\leq
(\widetilde {\mathcal W}^{(1)},\, \widetilde{\mathcal W}^{(2)},\widetilde \Phi  )$
if ${\mathcal W}^{(1)}\subset \widetilde {\mathcal W}^{(1)}$ 
{\ntekst and
$
\Phi= {\widetilde \Phi}|_{{\mathcal W}^{(1)}}.
$

Note that, by Remark \ref{Remark:3}, if 
$({\mathcal W}^{(1)},\, {\mathcal W}^{(2)}, \Phi  )$ and 
$(\widetilde {\mathcal W}^{(1)},\, \widetilde{\mathcal W}^{(2)},\widetilde \Phi  )$
are admissible triples, {then the extended triple,} $({\mathcal W}^{(1), ex},\, {\mathcal W}^{(2), ex}, \Phi^{ex}  )$,
where
\ba
& &{\mathcal W}^{(i), ex}={\mathcal W}^{(i)} \cup {\widetilde {\mathcal W}}^{(i)},
\\ \nonumber
& & \Phi^{ex}|_{{\mathcal W}^{(1)}}=\Phi, \quad
\Phi^{ex}|_{{\widetilde {\mathcal W}}^{(1)}}={\widetilde \Phi},
\ea
is also an admissible triple. }
Therefore,  by Zorn's lemma, there exists a maximal element 
$({\mathcal W}_m^{(1)},\, {\mathcal W}_m^{(2)}, \Phi_m  ) \in \mathcal A$.


\begin{lemma} \label{maximal}
The maximal element $({\mathcal W}_m^{(1)},\, {\mathcal W}_m^{(2)}, \Phi_m  ) $ of $\mathcal A$ satisfies
 \beq \label{10.3.6}
 {\mathcal W}_m^{(1)}= \M^{(1),reg}.
 \eeq 
 \end{lemma}
 \noindent
{\bf Proof.} If the claim is not true, there exists
    $x^{(1)}_0 \in\mathcal M ^{(1),reg} \cap \p{\mathcal W}_m^{(1)}$.
    Let $\mu([0,1])$ be a smooth path from 
    $\mu(0)=z\in {W}$ to $\mu(1)=x^{(1)}_0$, such that 
    \ba
    \mu([0,1)) \subset \M^{(1),reg}.
    \ea
    Then $d_0=d^{(1)}(\mu, \M^{(1)}_{sing})>0$. Let 
     $c= \frac{d_0}{2}$.  
    We can cover $\mu([0, 1])$ by a finite number of balls $B^{(1)}_j=B^{(1)}(x^{(1)}_j, c/2) \subset 
    \M ^{(1),reg}$ 
    so that  
    \beq \label{10.3.5}
   { \overline B}_j^{(1)} \subset {\mathcal W}_m^{(1)}, \,\,
    x^{(1)}_{j+1} \in B_j^{(1)},
   \eeq 
   where we order them so that $x^{(1)}_0 \in B^{(1)}_1$.
   Let  ${\mathcal O}^{(1)}_1=B^{(1)}(x^{(1)}_1, R)$  be a small ball such 
   that $0<R<c/2$ satisfies  (\ref{eq: basic assumptions}), (\ref{eq: basic assumptions2}),
   and 
   ${\mathcal O}^{(1)}_1   \subset {\mathcal W}_m^{(1)}$. 
   As $d^{(1)}( x^{(1)}_1,{\mathcal M}^{(1)}_{sing})>\frac {d_0}2$,
    Theorem \ref{recognition} yields that
 {\ntekst we can extend 
 the admissible triple $({\mathcal W}_m^{(1)},\, {\mathcal W}_m^{(2)}, \Phi_m  ) $ onto
 \ba
 {\widetilde {\mathcal W}}^{(i)}={\mathcal W}_m^{(i)} \cup B^{(i)}(x^{(i)}_1, c),
 \quad x^{(2)}_1= \Phi_m(x^{(1)}_1).
 \ea}
 As $ x^{(1)}_0\in B( x^{(1)}_1, c)$ {and $x^{(1)}_0 \in \p{\mathcal W}_m^{(1)}$}, this contradicts 
the fact that $({\mathcal W}_m^{(1)},\, {\mathcal W}_m^{(2)}, \Phi_m  ) $ is a maximal element
 of $\mathcal A$, which completes the proof of (\ref{10.3.6}).
    \qed\medskip

   Lemma \ref{maximal} proves that there is a diffeomorphism
   \ba
\Phi_m: {\mathcal M} ^{(1),reg} \to {\mathcal W}_m^{(2)}, \quad 
{\mathcal W}_m^{(2)}= \Phi_m({\mathcal M} ^{(1),reg})
\subset  {\mathcal M} ^{(2),reg},
\ea
which is  a Riemannian isometry.
   Changing the role of indexes 1 and 2, we see that  there is also a  diffeomorphism
\ba
{\widetilde \Phi}_m: {\mathcal M} ^{(2),reg} \to {\widetilde {\mathcal W}}_m^{(1)}, 
\quad {\widetilde {\mathcal W}}_m^{(1)}
\subset  {\mathcal M} ^{(1),reg}
\ea
which is a Riemannian isometry. Moreover,
using Lemma \ref{lem. N1 admissible} we see that ${\widetilde \Phi}_m$
and $\Phi_m$ coincide with the identity map on $W$.

Using (\ref{eq: Greens functions related}) we see that for all  $z\in {\mathbb C}\setminus {\mathbb R}$,
$x\in {\mathcal M} ^{(2),reg} $ and $y\in {W}$.
\ba
G^{(1)}(z; \Phi_m({\widetilde \Phi}_m(x)), y)=G^{(2)}(z; x, y).
\ea
 By Lemma \ref{identification}, this implies that $\Phi_m({\widetilde \Phi}_m(x))=x$,
 that is, $\Phi_m\circ{\widetilde \Phi}_m=I$ on ${\mathcal M} ^{(1),reg}$.
 Similarly, we see that   $\widetilde \Phi_m\circ{ \Phi}_m=I$
 on ${\mathcal M} ^{(2),reg}$ and hence
 \ba
 {\mathcal W}_m^{(2)}={\mathcal M} ^{(2),reg},\,\, {\mathcal W}_m^{(1)}={\mathcal M} ^{(1),reg},
 \quad \hbox{and}\,\, {\widetilde \Phi}_m= \Phi_m^{-1}.
\ea 
Summarizing, we have shown that 
\ba
\Phi_m:({\mathcal M}_{reg}^{(1)}, g^{(1)}) \to ({\mathcal M}_{reg}^{(2)}, g^{(2)}),
\ea 
 {is a diffeomorphism and an isometry.}

{\nntext 
Skipping again the superscript $^{(i)}$, recall that by Lemma \ref{two distance M and Mreg},
\beq\label{distances are same}
 d_{\mathcal M}(x, y)=  d_{\mathcal M^{reg}}(x, y),\quad \hbox{for any $x, y \in {\mathcal M}^{reg}$},
\eeq
where $ d_{\mathcal M^{reg}}$ is the distance on $({\mathcal M}^{reg}, g)$,
 defined
as the infimum of the length of  rectifiable paths connecting $x$ to $y$.

The identity  (\ref{distances are same}) implies that $
({\mathcal M}  ,  d_{\mathcal M^{reg}})$, considered as a metric space, is isometric to
 the completion of the metric space 
$({\mathcal M}_{reg},\, d_{reg}).$}
Thus, we can uniquely extend $\Phi_m$ to a metric isometry
\begin{equation} \label{E3B}
 \Phi:({\mathcal M}^{(1)}  , d^{(1)}) \to ({\mathcal M}^{(2)}  , d^{(2)}).
\end{equation} 
Again, taking into account that ${\mathcal M}^{(1),reg}$ is mapped to  ${\mathcal M}^{(2),reg}$ 
we see that $ \Phi$ maps also singular points to singular points.  

These prove  conditions (1)--(3) of Theorem \ref{thm: main result for IP}.\hfill $\square$

\printindex
\end{document}